\documentclass[reqno]{amsart}
\usepackage{amsthm}
\usepackage{epsfig}
\usepackage{pst-grad} 
\usepackage{pst-plot} 
\usepackage[english]{babel}
\usepackage[latin1]{inputenc}
\usepackage{indentfirst}
\usepackage{graphicx,subfigure}
\usepackage{tabularx}
\usepackage{xspace} 
\usepackage{amsmath,amssymb}
\usepackage{color}
\definecolor{oneblue}{rgb}{0,0.0,0.75}

\usepackage{fancyhdr}
 \usepackage{setspace}

\makeatletter
\newcommand{\sech}{\mathop{\operator@font sech}}
\newcommand{\sign}{\mathop{\operator@font sign}}
\makeatother

\newtheorem{remark}{Remark}[section]
\numberwithin{equation}{section}

\begin{document}

\title[Numerical approximation to Benjamin type equations...]{Numerical approximation to Benjamin type equations. Generation and stability of solitary waves}

\author[V. A. Dougalis]{Vassilios A. Dougalis}
\address{Mathematics Department, University of Athens, 15784
Zographou, Greece \and Institute of Applied \& Computational
Mathematics, FO.R.T.H., 71110 Heraklion, Greece}
\email{doug@math.uoa.gr}

\author[A. Duran]{Angel Duran}
\address{ Applied Mathematics Department,  University of
Valladolid, 47011 Valladolid, Spain}
\email{angel@mac.uva.es}

\author[D. Mitsotakis]{Dimitrios Mitsotakis}
\address{School of Mathematics and Statistics, Victoria University of Wellington, PO Box 600
Wellington 6140
New Zealand}
\email{dmitsot@gmail.com}
\urladdr{http://dmitsot.googlepages.com/}



\begin{abstract}
This paper is concerned with the study, by computational means, of the generation and stability of solitary-wave solutions of generalized versions of the Benjamin equation. The numerical generation of the solitary-wave profiles is accurately performed with a modified Petviashvili method which includes extrapolation to accelerate the convergence. In order to study the dynamics of the solitary waves the equations are discretized in space with a Fourier pseudospectral collocation method and a fourth-order, diagonally implicit Runge-Kutta method of composition type as time-stepping integrator. The stability of the waves is numerically studied by performing experiments with small and large perturbations of the solitary pulses as well as interactions of solitary waves.
\end{abstract}

\maketitle

\section{Introduction}
\label{sec1}
The purpose of this paper is to study numerically some aspects of solutions of generalized versions of the Benjamin equation, mainly focused on the generation and dynamics of solitary wave solutions.

The Benjamin type equations under study are of the form
\begin{eqnarray}
u_{t}-\mathcal{L}u_{x}+F(u)_{x}=0,\label{gben1}
\end{eqnarray}
where $u=u(x,t), x\in\mathbb{R}, t\geq 0$,  $\mathcal{L}=\mathcal{L}(r,m)$ is the Fourier multiplier operator with symbol
\begin{eqnarray}
\widehat{\mathcal{L}u}(\xi)=(|\xi|^{2m}-\gamma
|\xi|^{2r})\widehat{u}(\xi),\quad \xi\in\mathbb{R},\label{gben2}
\end{eqnarray}
and
\begin{eqnarray}
F(u)=\frac{u^{p-1}}{p-1},\label{gben3}
\end{eqnarray}
In (\ref{gben2}), (\ref{gben3}) $m$ and $p$ are positive integers, $p>2$ and $r\in\mathbb{R}, 0\leq r<m$, with $\gamma$ a positive constant. Equations (\ref{gben1}) were introduced by Chen \& Bona, \cite{BonaC1998}, as generalizations of the Benjamin equation, \cite{Benjamin1992,Benjamin1996,AlbertBR1999}, for the propagation of unidirectional, one-dimensional internal waves along the interface of a two-layer system of incompressible, inviscid fluids under a rigid lid assumption for the upper layer while the lower layer is bounded from below by an impermeable, horizontal bottom. The model assumes that gravity and surface tension are not negligible and compete with nonlinear and dispersive effects. Equations of the type (\ref{gben1}) are essentially considered also in \cite{LinaresS2005} but written with
\begin{eqnarray}
&&\mathcal{L}=\delta D^{2m}+\gamma H^{2r},\nonumber\\
&&\widehat{D^{2m}u}(\xi)=|\xi|^{2m}\widehat{u}(\xi),\quad \widehat{H^{2r}u}(\xi)=-|\xi|^{2r}\widehat{u}(\xi)\quad \xi\in\mathbb{R},
\label{gben2b}
\end{eqnarray}
and
\begin{eqnarray}
F(u)=\frac{u^{q+1}}{q+1},\quad q\in \mathbb{N},\label{gben3b}
\end{eqnarray}
where now $q\geq 1$ (say $p=q+2$ in (\ref{gben2})) and $r, m$ as above, and $\delta, \gamma$ positive constants. The particular case $m=1, r=1/2$ is called the generalized Benjamin equation (gBenjamin) by Linares \& Scialom (reducing to the Benjamin equation when $q=1$) and so will be here, written in the form
\begin{eqnarray}
u_{t}+u^{q}u_{x}+\gamma\mathcal{H}u_{xx}+\delta u_{xxx}=0,\label{gben4}
\end{eqnarray}
where $\mathcal{H}$ is the Hilbert transform so that $H=-\mathcal{H}\partial_{x}$, including the case $\gamma=0$, that is the generalized KdV equation (gKdV) with positive dispersion coefficient. Linares and Scialom call the case $r=1/2, \delta=0$ the generalized Benjamin-Ono (gBO) equation, \cite{KenigPV1994,BonaK2004}. (In \cite{LinaresS2005}  $\gamma$ is assumed to be positive, while in \cite{BonaK2004} $\gamma=-1$.)

\begin{figure}[htbp]
\centering
\subfigure[]
{\includegraphics[width=6cm]{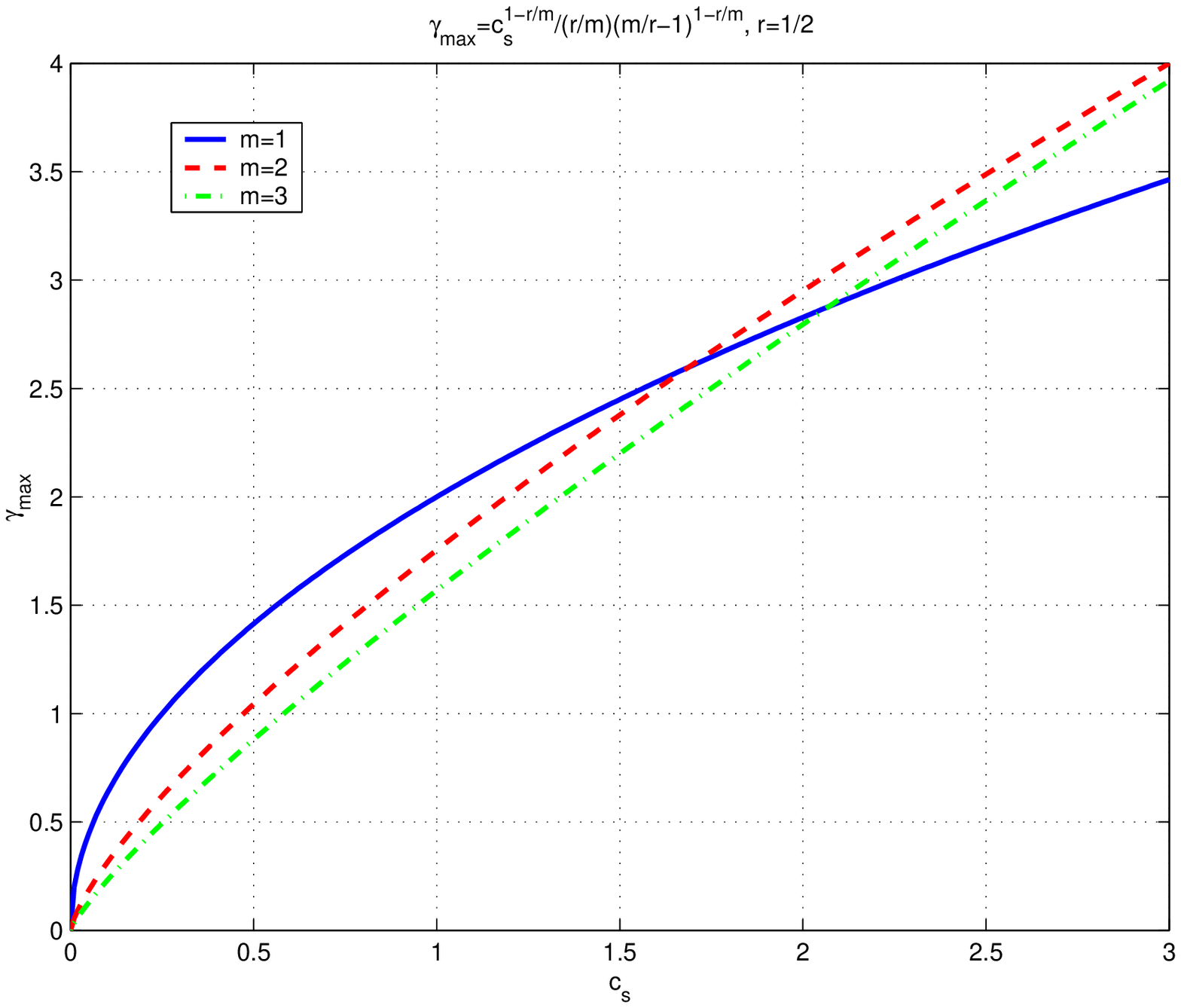}}
\subfigure[]
{\includegraphics[width=6cm]{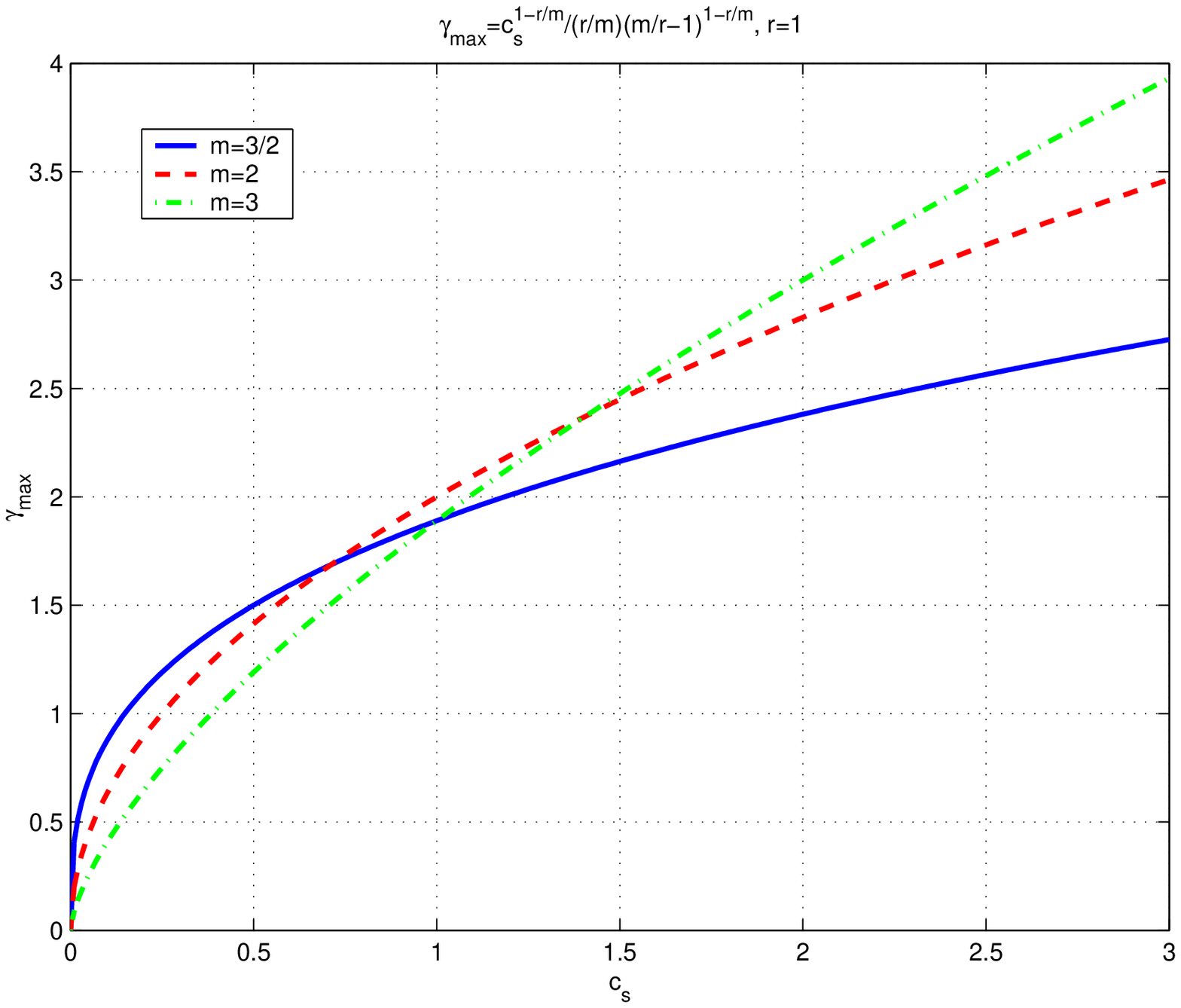}}
\subfigure[]
{\includegraphics[width=6cm]{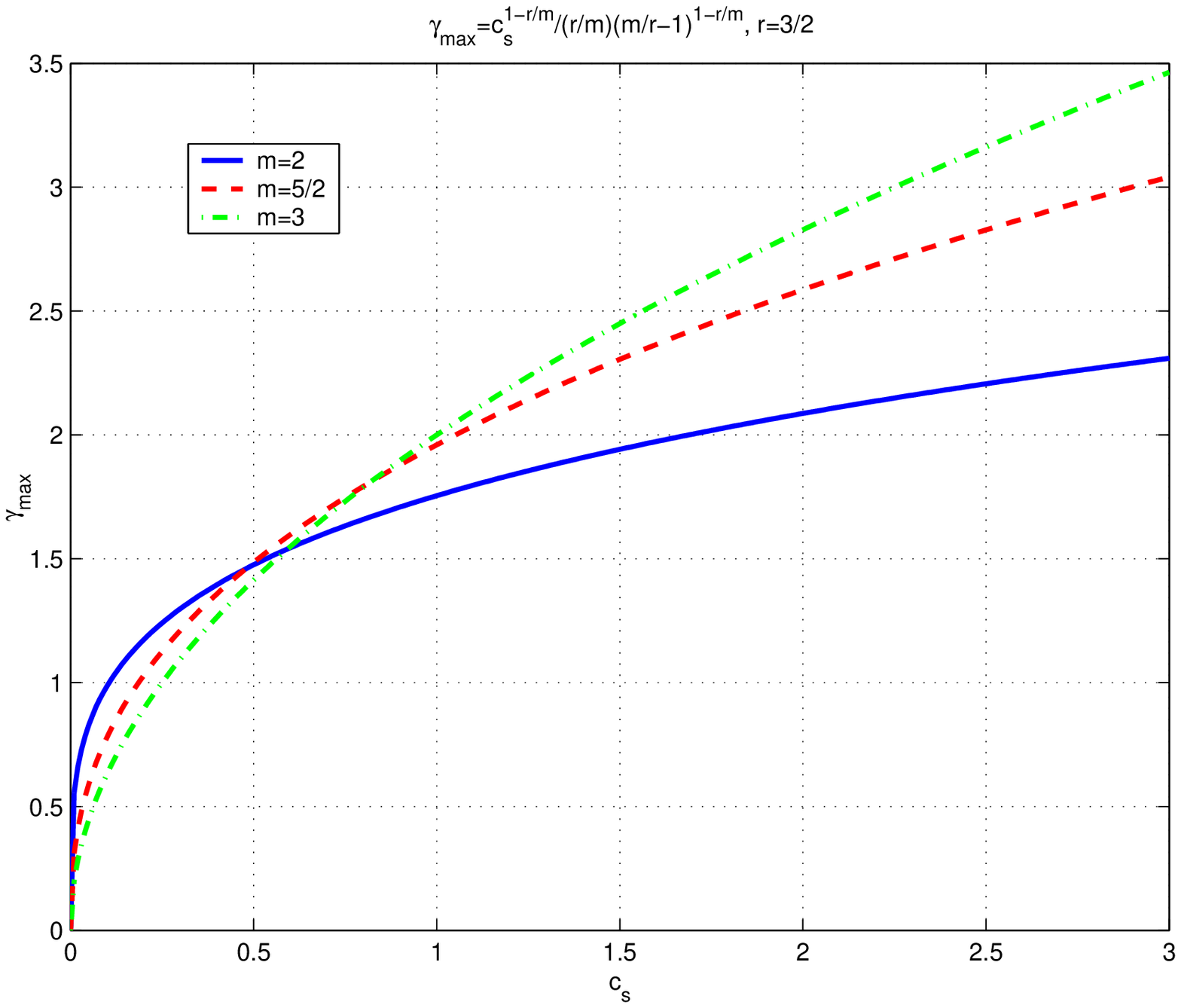}}
\caption{Parameter $\gamma_{max}(c_{s})$ as function of $c_{s}$. (a) $r=1/2$; (b) $r=1$; (c) $r=3/2$.}
\label{gbf1}
\end{figure}

Based on the theory developed for the gKdV and gBO equations in \cite{KenigPV1991,KenigPV1994}, Linares and Scialom, \cite{LinaresS2005}, establish local and global well-posedness results for the ivp for (\ref{gben1}), (\ref{gben2b}), (\ref{gben3b}). 
More specifically, the generalized Benjamin equation (\ref{gben4}) is studied first. In this case, for $u_{0}\in H^{s}(\mathbb{R}), s\geq 1$, the existence of $T=T(||u_{0}||_{H^{s}})$ and of a unique solution $u\in C([0,T],H^{s}(\mathbb{R})$ is established. Furthermore it is shown that the ivp is globally well-posed in $H^{1}(\mathbb{R})$ under the following conditions:
\begin{itemize}
\item For $q=2,3$, with no restriction on the initial data.
\item For $q\geq 4$ and the initial data small enough in $L^{2}(\mathbb{R})$.
\end{itemize}
In the general case, for $u_{0}\in H^{m}(\mathbb{R}), m>r>0, m>1$, the ivp for (\ref{gben1}), (\ref{gben2b}), (\ref{gben3b}) is globally well-posed in the energy space, i.~e.  the weakest space where the energy 
\begin{eqnarray}
E(u(t))=\int_{-\infty}^{\infty} \delta uD^{2m}u dx-\int_{-\infty}^{\infty} \gamma uH^{2r}u dx-\frac{2}{(q+1)(q+2)}\int_{-\infty}^{\infty} u^{q+2}dx.\label{gben7}
\end{eqnarray}
makes sense.
One has global solutions under the conditions:
\begin{itemize}
\item For $q<4m$, with no restriction on the initial data.
\item For $q\geq 4m$ the initial data must be small enough in $ H^{m}(\mathbb{R})$.
\end{itemize}
In addition to (\ref{gben7}) there are two conserved quantities for decaying, smooth enough solutions:
\begin{eqnarray*}
I(u)=\int_{-\infty}^{\infty} u^{2} dx,\quad C(u)=\int_{-\infty}^{\infty} u dx.
\end{eqnarray*}
This paper is mainly focused on the behaviour of solitary-wave solutions of (\ref{gben1}), (\ref{gben2b}), (\ref{gben3b}), i.~e. solutions of the form 
$$ u=\varphi(x-c_{s}t),\quad c_{s}>0,$$ where $\varphi$ and derivatives tend to zero as $X=x-c_{s}t\rightarrow \pm \infty$. Therefore, we have
\begin{eqnarray}
-c_{s}\varphi +F(\varphi)-\mathcal{L}\varphi=0.\label{gben8}
\end{eqnarray} 
{
We give a brief review on the literature about solitary wave solutions of (\ref{gben1}) focusing on \cite{BonaC1998} and \cite{Angulo2009} . In \cite{BonaC1998}, existence and asymptotic properties of solutions of (\ref{gben8}) are established (with $\delta=1$). Specifically, for $c_{s}>0$, $m$ and $q$ positive integers with $q>2$ and $0\leq r<m$, existence of solitary waves in $H^{\infty}$ is proved for $\gamma<\gamma_{max}(c_{s})$ where
\begin{eqnarray}
\gamma_{max}(c_{s})=\frac{c_{s}^{\frac{m-r}{m}}}{s(r,m)},\quad s(r,m)=\frac{r}{m}\left(\frac{m}{r}-1\right)^{\frac{m-r}{m}}.\label{gben9}
\end{eqnarray}
(See Figure \ref{gbf1}.) 
When $\delta\neq 1$, (\ref{gben9}) would have the form
\begin{eqnarray}
\gamma_{max}(c_{s})=\frac{\delta^{\frac{r}{m}}c_{s}^{\frac{m-r}{m}}}{s(r,m)},\quad s(r,m)=\frac{r}{m}\left(\frac{m}{r}-1\right)^{\frac{m-r}{m}}.\label{gben9b}
\end{eqnarray}
The main results concerning the asymptotic decay of the solitary waves, established in \cite{BonaC1998}, are:
\begin{itemize}
\item[(i)] (Exponential decay.) If $r$ is a positive integer (the $H$ term becomes local) or if $\gamma=0$ then there is a $\sigma_{0}>0$ such that for any $\sigma<\sigma_{0}$
$$
e^{\sigma |x|}\varphi(x)\rightarrow 0,\quad |x|\rightarrow\infty.
$$ 
\item[(ii)] (Algebraic decay.) Otherwise, there is a constant $\mu$ such that
$$
|x|^{2r+1}\varphi(x)\rightarrow \mu,\quad |x|\rightarrow\infty.
$$ 
\end{itemize}
\begin{remark}
These results are extended in \cite{BonaC1998} to a more general dispersion operator $\mathcal{L}$ of Fourier symbol
$$
\widehat{\mathcal{L}u}(\xi)=\sum_{j=1}^{K}\alpha_{j}|\xi|^{2r_{j}}\widehat{u}(\xi),
$$ where $\alpha_{j}\in\mathbb{R}, \alpha_{K}>0, 0<r_{1}<\cdots <r_{K}$, but $r_{K}$ is not an integer (i.~e.  the highest-order term is not local). 
\end{remark}
Another reference of interest is \cite{Angulo2009} (see also the references therein). In the case of the generalized Benjamin equation, the above described well-posedness results as well as the results on existence of solitary wave solutions are mentioned. On the other hand, a sort of instability of the solitary waves (in some orbital sense)  for $q>3$ and $\gamma$ small enough is established (see also \cite{Angulo2003}). These results (both existence and instability) are finally extended to the general case (\ref{gben1}).
}

We also mention that as in the case of the Benjamin equation, \cite{AlbertBR1999,DougalisDM2012}, a normalized form of (\ref{gben8}), (\ref{gben2b}), (\ref{gben3b}) may be derived. By using the change of variables $\varphi(X)=A(q)\psi(Z)$ where 
$$A(q)=((q+1)c_{s})^{1/q}, Z=BX, |B|=\left(\frac{c_{s}}{\delta}\right)^{\frac{1}{2m}},$$ we obtain the following ordinary differential equation (ode) for  $\psi$:
\begin{eqnarray}
\psi-\psi^{q+1}-D^{2m}\psi-\frac{1}{s(r,m)}\tilde{\gamma}H^{2r}\psi=0,\label{gben8b}
\end{eqnarray}
where $\tilde{\gamma}=\gamma/\gamma_{max}(c_{s})$ with $\gamma_{max}(c_{s})$ given by (\ref{gben9b}).

The numerical study of (\ref{gben1}), (\ref{gben2b}), (\ref{gben3b}) presented in this paper is concerned with the generation and dynamics of its solitary waves. As far as the generation is concerned, the nonexistence of exact formulas forces  to implement a numerical technique of approximation to the solitary profiles. The literature in this sense for (\ref{gben1}) is focused on the case of the Benjamin equation, where the main procedures to this end (see \cite{DougalisDM2012} and references therein) are based on iterative methods applied to the corresponding discretization of the equation for the profiles combined with numerical continuation on the parameter $\gamma$ that governs the nonlocal term and from the limiting case $\gamma=0$, that is the KdV equation, for which the profiles are known explicitly.

In a similar way, the references to the numerical integration of (\ref{gben1}) to explore numerically the dynamics of the solitary waves are concerned with schemes to approximate the Benjamin equation. In this sense, the corresponding periodic initial-value problem is discretized in \cite{DougalisDM2015} with a hybrid finite element-spectral method in space  and a two-stage, fourth-order, Gauss-Legendre implicit Runge-Kutta method as time integrator, while in \cite{DougalisDM2012} the accuracy of the solitary-wave profiles were verified by using a scheme consisting of a Fourier pseudospectral method in space and a third-order, simply diagonally implicit Runge-Kutta time-stepping integrator. 

In this paper we propose new schemes in order to perform a computational study of (\ref{gben1}), (\ref{gben2b}), (\ref{gben3b}) for the generation of the solitary waves and the analysis of their dynamics. Both are described in Section \ref{sec2}. As for the numerical generation, the system for the approximate profiles, obtained by discretizing (\ref{gben8}) on a long enough interval with periodic boundary conditions by using a Fourier pseudospectral collocation method, is numerically solved in Fourier space with the iterative method of Petviashvili, \cite{Petv1976}, along with an acceleration technique based on extrapolation. The purpose of the inclusion of acceleration is two-fold: avoiding the numerical continuation procedure (which is sometimes inefficient) and computing highly oscillatory profiles  that the simple Petviashvili method is not able to approximate, \cite{DougalisDM2015}. The numerical experiments aim to suggest some properties of the solitary waves. They are mainly concerned with the behaviour of the amplitude-speed relation and the asymptotic decay with respect to the parameters $r, m, q$ and $\gamma$. Once the accuracy of the method to generate numerically the solitary-wave profiles is verified, the next step of the study is to explore numerically their dynamics. To this end, a numerical scheme to integrate the periodic initial-value problem for these equations is proposed. The method is based on Fourier pseudospectral discretization in space and a fourth-order diagonally implicit Runge-Kutta method for time integration. The time-stepping scheme is a composition integrator with the implicit midpoint rule as basis method. It was derived as the first member of a family of high-order, symplectic, symmetric composition methods by Yoshida, \cite{Yoshida1990}, and its accuracy and favorable properties  for the integration of some nonlinear wave problems were observed in \cite{FrutosS1992}. With this scheme, in some numerical  experiments in Section \ref{sec2} we check the accuracy of the previously computed wave profiles. 

The experiments also serve to validate the proposed discretization as a stable and accurate way to approximate the dynamics of the solitary waves.
The corresponding computational study is carried out in Section \ref{sec3}, in
experiments with small and large perturbations as well as solitary wave interactions. Concluding remarks are summarized in Section \ref{sec4}.

The following notation will be used throughout the paper. We denote by $H^{s}(\mathbb{R}), 0\leq s\leq \infty$, the Sobolev space of order $s$, with $H^{0}(\mathbb{R})=L^{2}(\mathbb{R})$. For $s>0$, the norm in $H^{s}(\mathbb{R})$ is given by $||\cdot||_{s}$. The inner product for real or complex-valued functions in $L^{2}(\mathbb{R})$ is denoted by $\langle\cdot,\cdot\rangle$ and the associated norm by $||\cdot||$. Finally, for $T>0$, $C([0,T],H^{s})$ will stand for the space of continuous functions $u:[0,T]\rightarrow H^{s}(\mathbb{R})$.

\section{Numerical generation and propagation of solitary waves}
\label{sec2}
\subsection{Numerical generation. Accelerated Petviashvili method}
\label{sec21}
This section concerns the numerical generation of solitary wave solutions of (\ref{gben1}), (\ref{gben2b}), (\ref{gben3b}) by an iterative Petviashvili type method. We first review the basic iteration procedure for (\ref{gben8}) or (\ref{gben8b}) following \cite{DougalisDM2012}. The equation (\ref{gben8}) is discretized in some way on a sufficiently long interval $(-l,l)$ with periodic boundary conditions. For a uniform grid $x_{j}=-l+jh, j=0,\ldots,N, h=2l/N$ on $(-l,l)$ let $\varphi_{h}=(\varphi_{h,0},\ldots,\varphi_{h,N})$ be the vector approximating the values of $\phi$ at the grid points, with $\varphi_{h,j}\approx \varphi(x_{j}), j=0,\ldots, N$. Then $\varphi_{h}$ satisfies a system of the form
\begin{eqnarray}
-c_{s}\varphi_{h}+F_{h}(\varphi_{h})-\mathcal{L}_{h}\varphi_{h}=0,\label{gben10}
\end{eqnarray}
for some approximations $F_{h}, \mathcal{L}_{h}$ to the operators $F,\mathcal{L}$ respectively. As in \cite{DougalisDM2012} Fourier spectral approximation will be used to implement (\ref{gben10}) in the numerical experiments below.

In \cite{DougalisDM2012} the numerical generation of solitary wave profiles for the Benjamin equation was carried out by numerical continuation on $\gamma$ plus an iterative procedure for each step of the homotopic path. Here an alternative scheme will be implemented, replacing the numerical continuation procedure by an acceleration technique. This has been shown to be more efficient in computing profiles for values of $\gamma$ for which the first mentioned procedure does not seem to work, see \cite{AlvarezD2015a}. The underlying iterative technique will be the Petviashvili method, \cite{Petv1976}. Write (\ref{gben10}) in the fixed-point form
$$
(c_{s}+\mathcal{L}_{h})\varphi_{h}=F_{h}(\varphi_{h}),
$$ and introduce an iterative formulation putting
\begin{eqnarray}
m(\varphi_{h}^{[\nu]})&=&\frac{\langle (c_{s}+\mathcal{L}_{h})\varphi_{h}^{[\nu]}, \varphi_{h}^{[\nu]}\rangle}{\langle F_{h}(\varphi_{h}), \varphi_{h}^{[\nu]}\rangle},\label{gben11}\\
(c_{s}+\mathcal{L}_{h})\varphi_{h}^{[\nu+1]}&=&m(\varphi_{h}^{[\nu]})^{\epsilon}F_{h}(\varphi_{h}^{[\nu]}),\label{gben12}
\end{eqnarray}
In (\ref{gben11}), (\ref{gben12}) $\epsilon$ is a parameter that depends on the homogeneity degree $q+1$ of $F_{h}$ (assuming that the approximation to $F$ inherits this property) and $\varphi_{h}^{[0]}\neq 0$. For the continuous problem, the Fourier symbol of the operator $(c_{s}+\mathcal{L})$ is $(c_{s}+\delta |\xi|^{2m}-\gamma |\xi|^{2r})$ and the restriction $0<\gamma<\gamma_{max}(c_{s})$ for $c_{s}>0$ ensures the existence of an inverse. This property has to be preserved when choosing the discretization $\mathcal{L}_{h}$. (For instance, it is satisfied by the Fourier spectral approximation.) Formulas (\ref{gben11}), (\ref{gben12}) for this approximation and with $F_{h}(\phi)=\phi .^{q+1}/(q+1)$ (the dot stands for the Hadamard product) have the following form for the discrete Fourier coefficients:
\begin{eqnarray*}
m(\varphi_{h}^{[\nu]})&=&\frac{\sum_{k}(c_{s}+\delta |k|^{2m}-\gamma |k|^{2r})|\widehat{\varphi_{h}^{[\nu]}}(k)|^{2}}{\sum_{k}\widehat{F_{h}(\varphi_{h}^{[\nu]})}(k)\overline{\widehat{\varphi_{h}^{[\nu]}}(k)}},\nonumber\\
\widehat{\varphi_{h}^{[\nu+1]}}(k)&=&\frac{m(\varphi_{h}^{[\nu]})^{\epsilon}}{c_{s}+\delta |k|^{2m}-\gamma |k|^{2r}}\widehat{F_{h}(\varphi_{h}^{[\nu]})}(k),\quad k\in\mathbb{Z},\label{gben14}
\end{eqnarray*}
where $1<\epsilon<(q+2)/q$ (with $\epsilon=(q+1)/q$ as optimal value).
The iteration is controlled by using several strategies, namely by specifying:
\begin{itemize}
\item The maximum number of iterations.
\item The maximum tolerance $TOL$ for:
\begin{itemize}
\item The discrepancy
\begin{eqnarray*}
SFE_{\nu}=|1-m(\varphi_{h}^{[\nu]})|,\quad \nu=0,1,\ldots,\label{gben15}
\end{eqnarray*}
\item The residual error in some norm (Euclidean or maximum)
\begin{eqnarray}
RES_{\nu}=||(c_{s}+\mathcal{L}_{h})\varphi_{h}^{[\nu]}-F_{h}(\varphi_{h}^{[\nu])}||,
\quad \nu=0,1,\ldots,\label{gben15b}
\end{eqnarray}
or the relative version
\begin{eqnarray*}
RES_{\nu}=\frac{||(c_{s}+\mathcal{L}_{h})\varphi_{h}^{[\nu]}-F_{h}(\varphi_{h}^{[\nu])}||}{||\varphi_{h}^{[\nu]||}},\quad \nu=0,1,\ldots,\label{gben16}
\end{eqnarray*}
\end{itemize}
\end{itemize}
{As we mentioned previously, the difficulties that the Petviashvili method presents in generating solitary waves because of high oscillations, \cite{DougalisDM2012}, can be overcome by incorporating acceleration techniques in the iterative procedure. The most widely used methods in the literature for this purpose appear to be the so-called Vector Extrapolation Methods (VEM). They introduce an extrapolation procedure to transform  the original sequence of the iterative process by using different strategies. For a detailed analysis and implementation of the methods see e.~g.  \cite{smithfs,jbilous} and references therein. 
A description of the application of acceleration techniques for traveling wave computations may be found in \cite{AlvarezD2015a}. The study made there for the Benjamin equation leads to the choice of the so-called minimal polynomial extrapolation (MPE), \cite{smithfs,jbilous,sidifs,sidi}, as the VEM acceleration technique in the experiments to be described below.
}
\subsection{Numerical results}
\label{sec22}
In this section the accuracy of the iterative scheme (\ref{gben11}), (\ref{gben12}) accelerated with MPE is checked by generating several solitary wave profiles of (\ref{gben1}). From now on we fix $\delta=1$.
\subsubsection{Generalized Benjamin equation}
\begin{figure}[htbp]
\centering
{\includegraphics[width=0.8\textwidth]{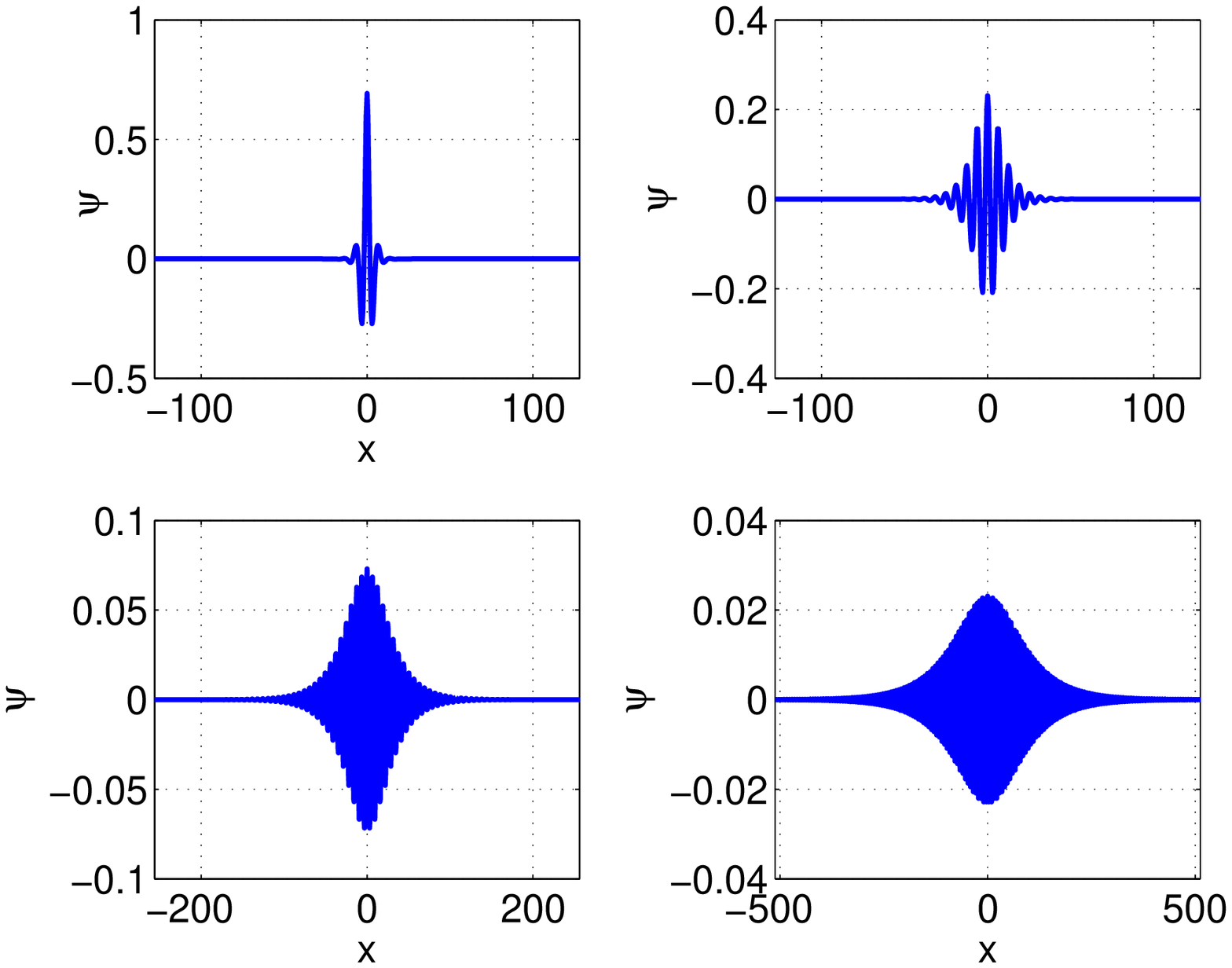}}
\caption{Computed solitary wave profiles for (\ref{gben8b}): $q=2, \widetilde{\gamma}=0.9,0.99,0.999,0.9999$.}
\label{gbenfig3}
\end{figure}

\begin{figure}[htbp]
\centering
{\includegraphics[width=0.8\textwidth]{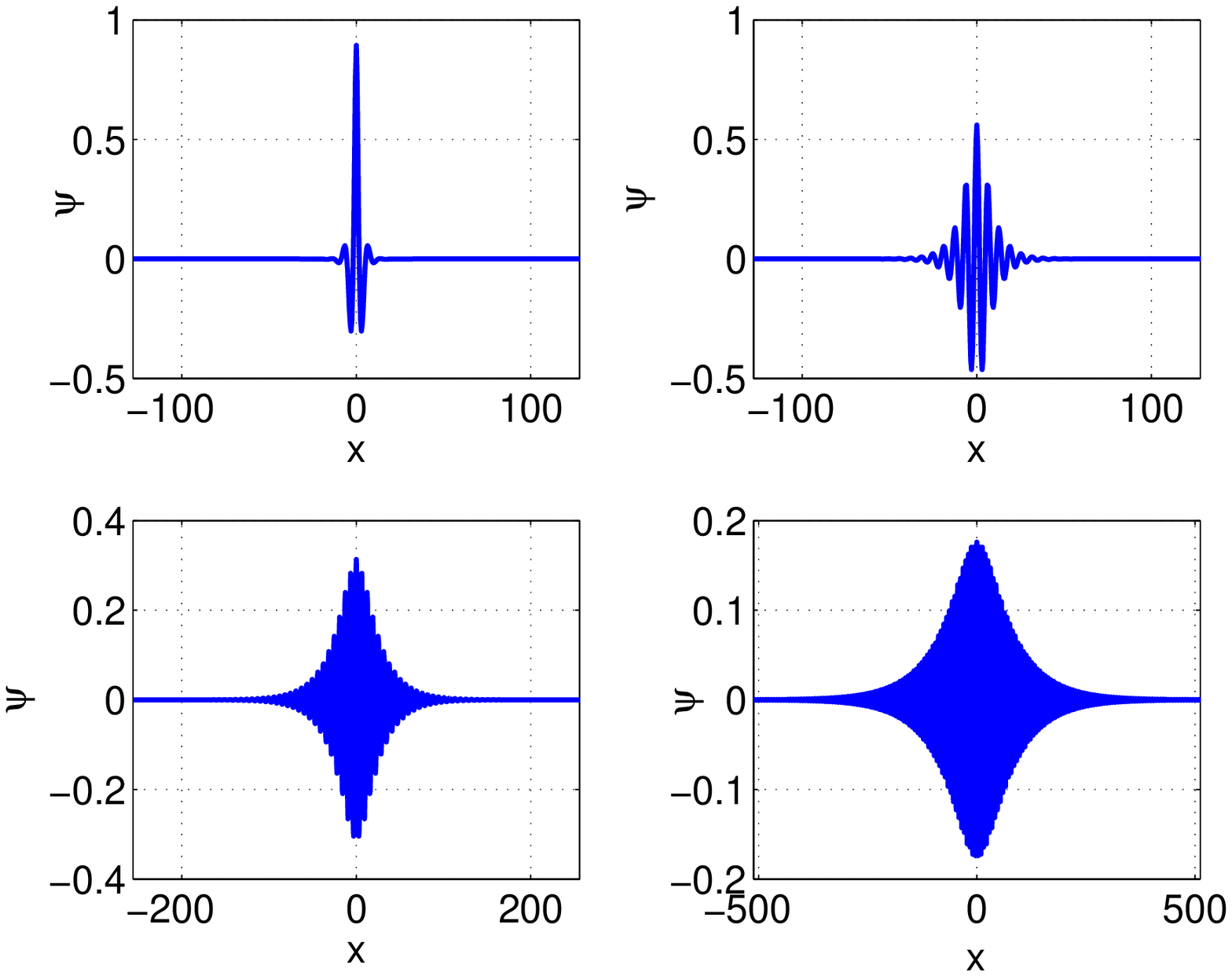}}
\caption{Computed solitary wave profiles for (\ref{gben8b}): $q=4, \widetilde{\gamma}=0.9,0.99,0.999,0.9999$.}
\label{gbenfig4}
\end{figure}

\begin{figure}[htbp]
\centering
{\includegraphics[width=0.8\textwidth]{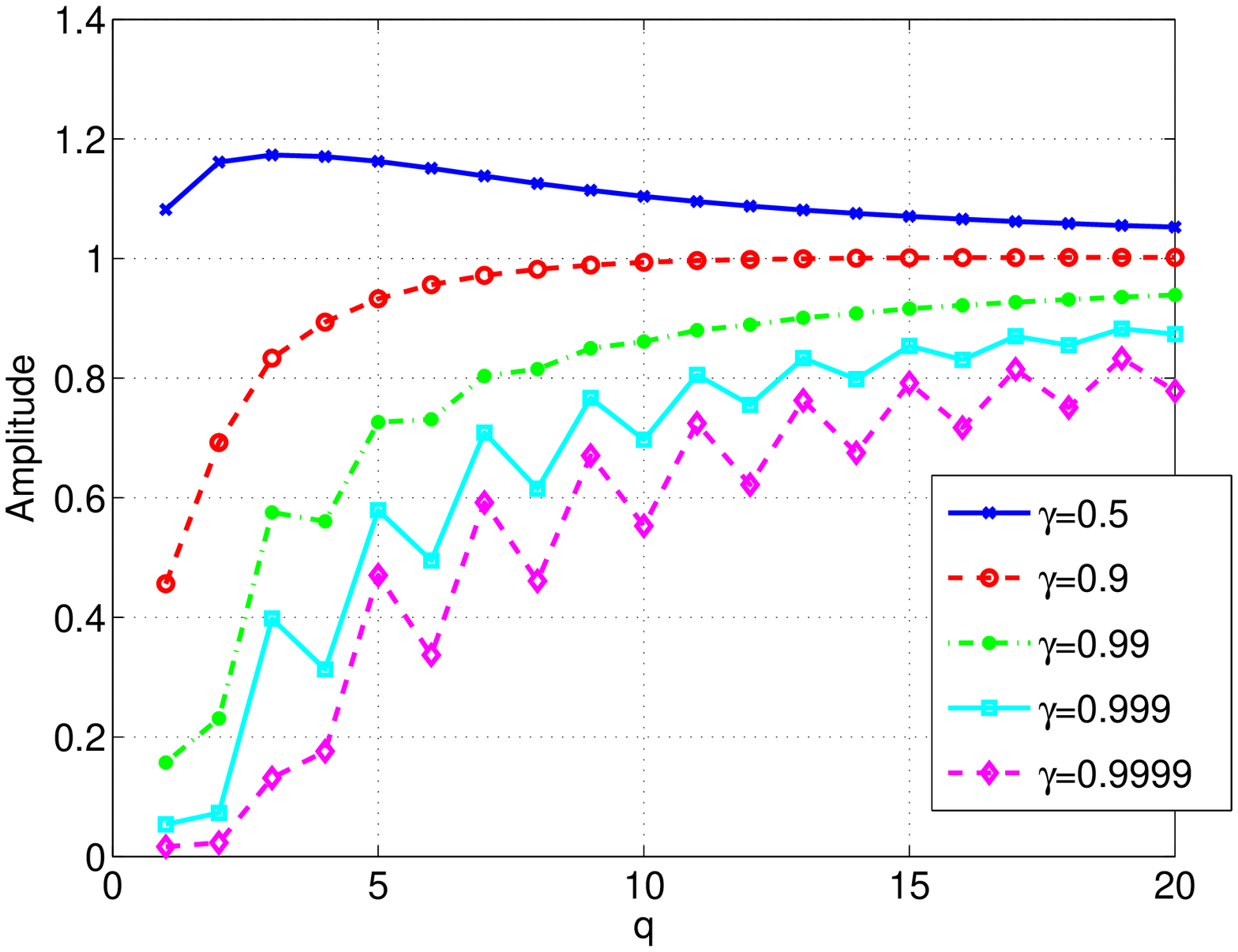}}
\caption{Amplitude vs. $q$ for approximate solitary wave profiles of (\ref{gben8b}).}
\label{gbenfig6}
\end{figure}

We start with the generalized Benjamin equation ($r=1/2, m=1$) in the normalized formulation (\ref{gben8b}). Our first experiment concerns the form of the solitary wave profiles for several values of the parameter $\widetilde{\gamma}\in [0,1)$ and $q\geq 1$. Two examples are shown in Figures \ref{gbenfig3} and \ref{gbenfig4}, corresponding respectively to $q=2$ and $4$ and for each $q$, to $\widetilde{\gamma}=0.9,0.99,0.999,0.9999$. The method is able to compute profiles for $\widetilde{\gamma}$ much closer to one as compared to other strategies, \cite{DougalisDM2012}. These results and also Figure \ref{gbenfig6} suggest that for those solitary waves with amplitude below one and fixed $\widetilde{\gamma}$, the amplitude grows with $q$, while for waves of amplitudes larger than one, the amplitude decreases with $q$. In other words, there seems to be a limiting value around one for the amplitude of the computed wave as $q$ grows. Figure \ref{gbenfig6} also shows that for some values of $\widetilde{\gamma}$ the amplitudes for odd $q$ are larger than those for the neighboring even $q$.

\begin{figure}[htbp]
\centering
\subfigure[]
{\includegraphics[width=6cm]{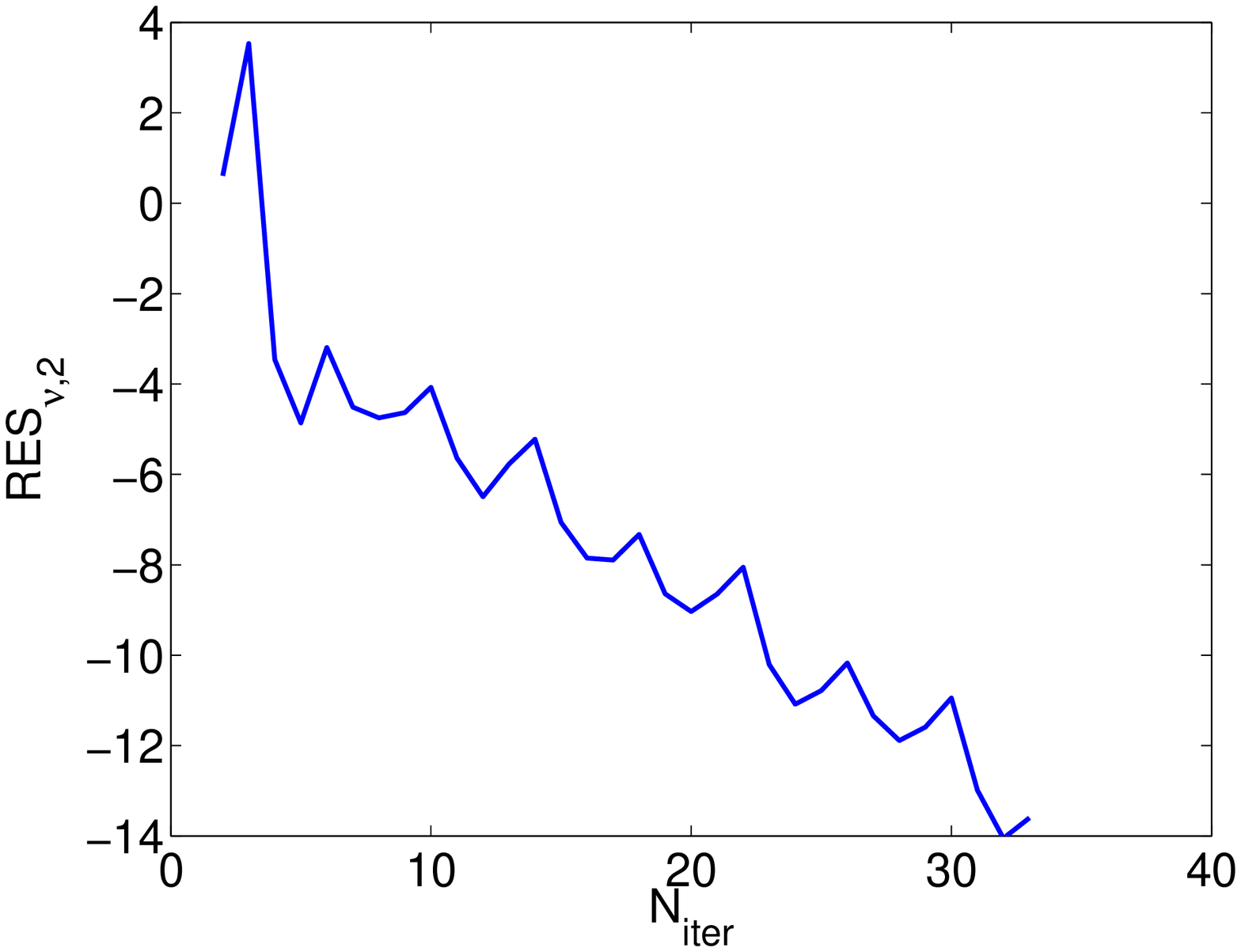}}
\subfigure[]
{\includegraphics[width=6cm]{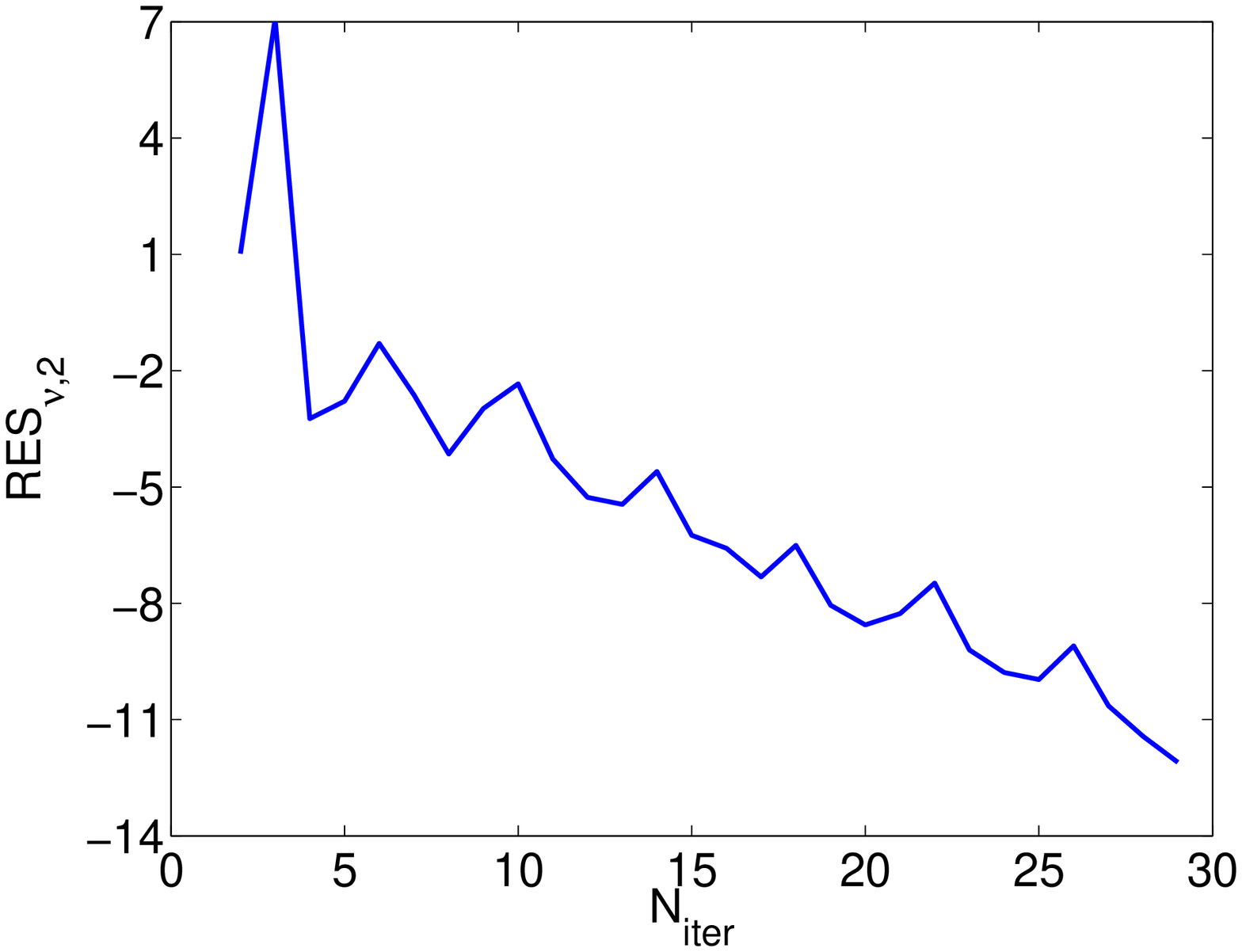}}
\caption{Residual error (\ref{gben15b}) in Euclidean norm as function of the number of iterations, $\widetilde{\gamma}=0.9999$. (a) $q=2$, (b) $q=4$.}
\label{gbenfig7}
\end{figure}

The performance of the method is checked in Figures \ref{gbenfig7}(a)-(b), which show the behaviour of the residual error (\ref{gben15b}) (in Euclidean norm) as a function of the number of iterations for $q=2,4$ and $\widetilde{\gamma}=0.9999$. Since the amplitude of the profiles is increasing with $q$, the method needs fewer iterations to achieve a tolerance of magnitude $TOL=10^{-12}$, see Table \ref{gbentab1}.



\begin{table}
\smallskip
\begin{center}
\begin{tabular}{|c|c|}
\hline
 $q$&$N_{iter}$\\
\hline
1&$49(1.911912E-14)$\\
2&$17(5.783437E-13)$\\
4&$16(4.838645E-13)$\\
6&$16(1.425915E-13)$\\
\hline
\end{tabular}
\end{center}
\caption{Number of iterations required to achieve a tolerance in the residual error of $TOL=1E-12$ (with the final residual error in parenthesis).}\label{gbentab1}
\end{table}



\begin{figure}[htbp]
\centering
{\includegraphics[width=0.8\textwidth]{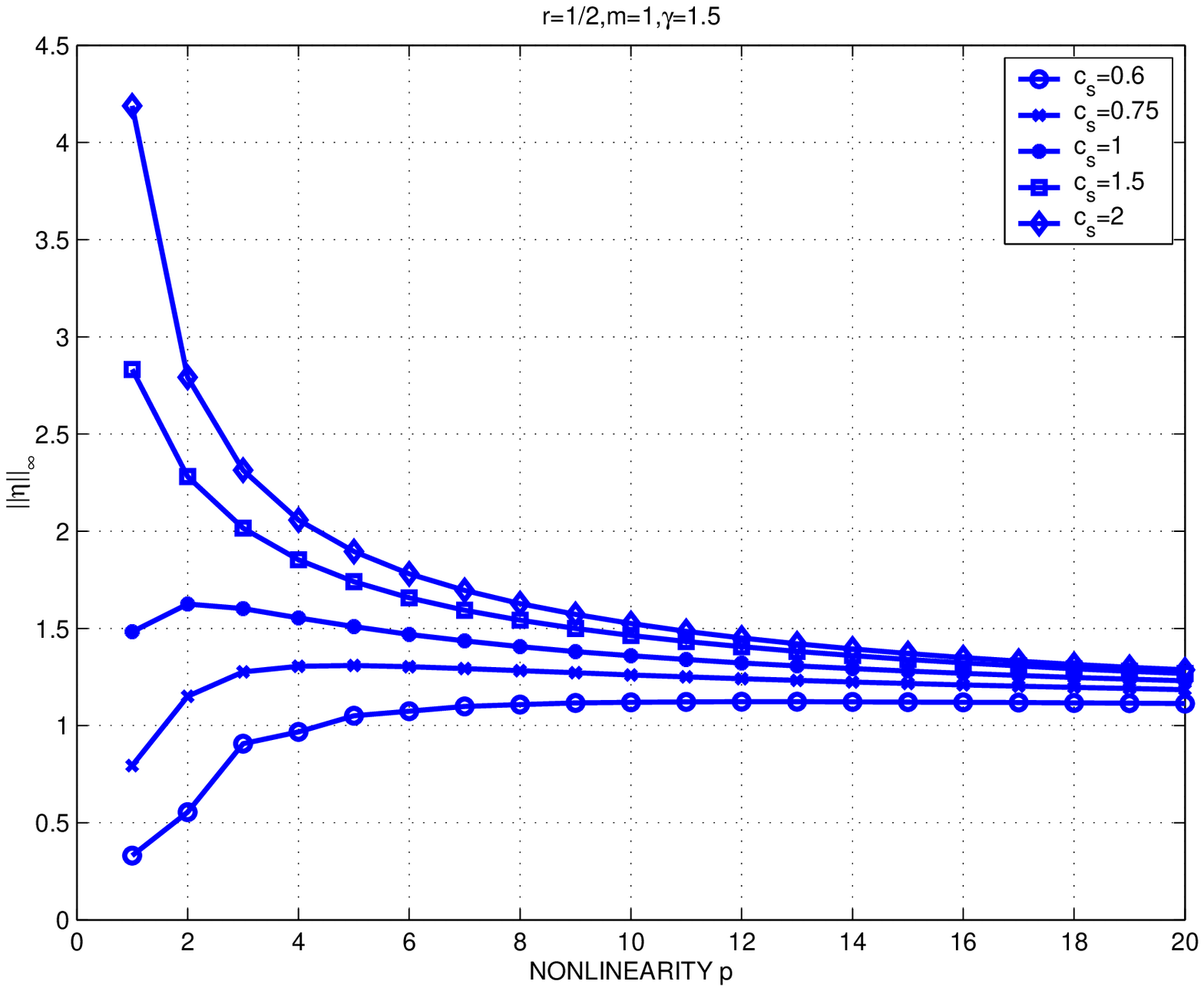}}
\caption{Amplitude vs. $q$ for approximate solitary wave profiles of (\ref{gben8b}).}
\label{gbenfig8}
\end{figure}

This behaviour of the amplitudes for the case of (\ref{gben8b}) is also observed in Figure \ref{gbenfig8}, which displays the amplitude of the computed wave as a function of $q$ when the generalized Benjamin equation without normalization (\ref{gben8}) is considered, for several values of $c_{s}$ and $\gamma=1.5$. 

\begin{figure}[htbp]
\centering
{\includegraphics[width=0.8\textwidth]{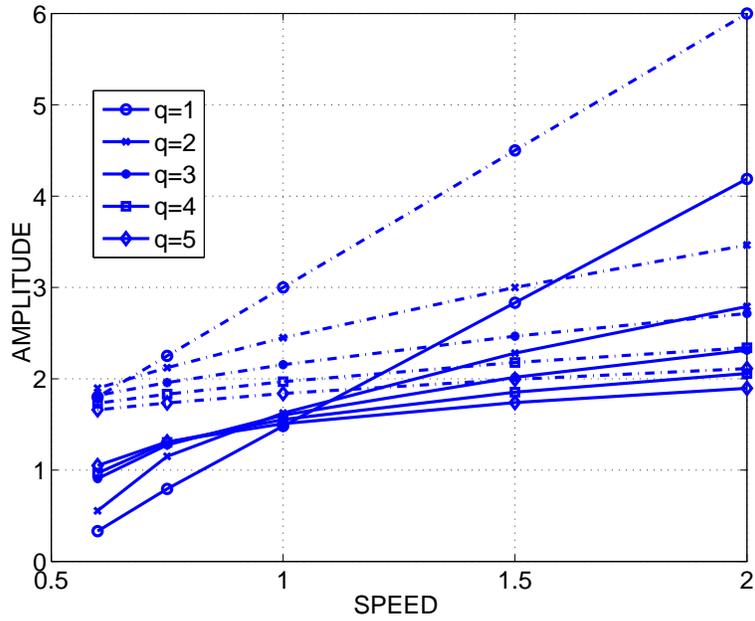}}
\caption{Amplitude-speed relation for computed solitary wave profiles of (\ref{gben8}) with $\gamma=1.5$ (solid lines) and for the gKdV equation (dashed lines).}
\label{gbenfig10}
\end{figure}

With the same value of $\gamma$, Figure \ref{gbenfig10} shows in solid lines the amplitude-speed relation for the first five integer values of $q$. Comparison with the corresponding relation for the generalized KdV (gKdV) equation (with dashed lines) reveals a very similar behaviour in all cases and suggests a similar amplitude-speed relation, \cite{BonaDKM1995}.
The corresponding speed-amplitude data (in logarithmic scale) have been fitted to a linear polynomial. (For the experiment, $c_{s}$ was varied from $1.1$ to $10$ with a stepsize of $\Delta c_{s}=0.01$.) The corresponding coefficients along with the goodness of fit parameters are given in Table \ref{gbentab2}. The coefficients suggest a relation $A(c_{s})\approx c_{s}^{\alpha}$ with $\alpha$ close to $1/q$ of the gKdV case.
\begin{table}
\smallskip
\begin{center}
\begin{tabular}{|c|c|c|c|c|}
\hline
 $q$&$K$&$\alpha$&$SSE$&$R^{2}$\\
\hline
1&$1.3042$&$1.177$&$1.7072E-01$&$9.9775E-01$\\
&$(1.3007,1.3076)$&$(1.173,1.18)$&&\\
2&$1.3243$&$0.557$&$1.5257E-01$&$9.9106E-01$\\
&$(1.3210,1.3276)$&$(0.5535,0.5605)$&&\\
4&$1.2659$&$0.2381$&$3.8525E-02$&$9.8769E-01$\\
&$(1.2643,1.2674)$&$(0.2363,0.2398)$&&\\
7&$1.1902$&$0.1241$&$4.4959E-02$&$9.9468E-01$\\
&$(1.1897,1.1908)$&$(0.1235,0.1247)$&&\\
\hline
\end{tabular}
\end{center}
\caption{Speed-amplitude relation, gBenjamin equation. Fitting curve $A(c_{s})\approx Kc_{s}^{\alpha}$ and goodness of fit parameters (with 95\% confidence bounds).}\label{gbentab2}
\end{table}
\subsubsection{General case}
\begin{figure}[htbp]
\centering
{\includegraphics[width=0.8\textwidth]{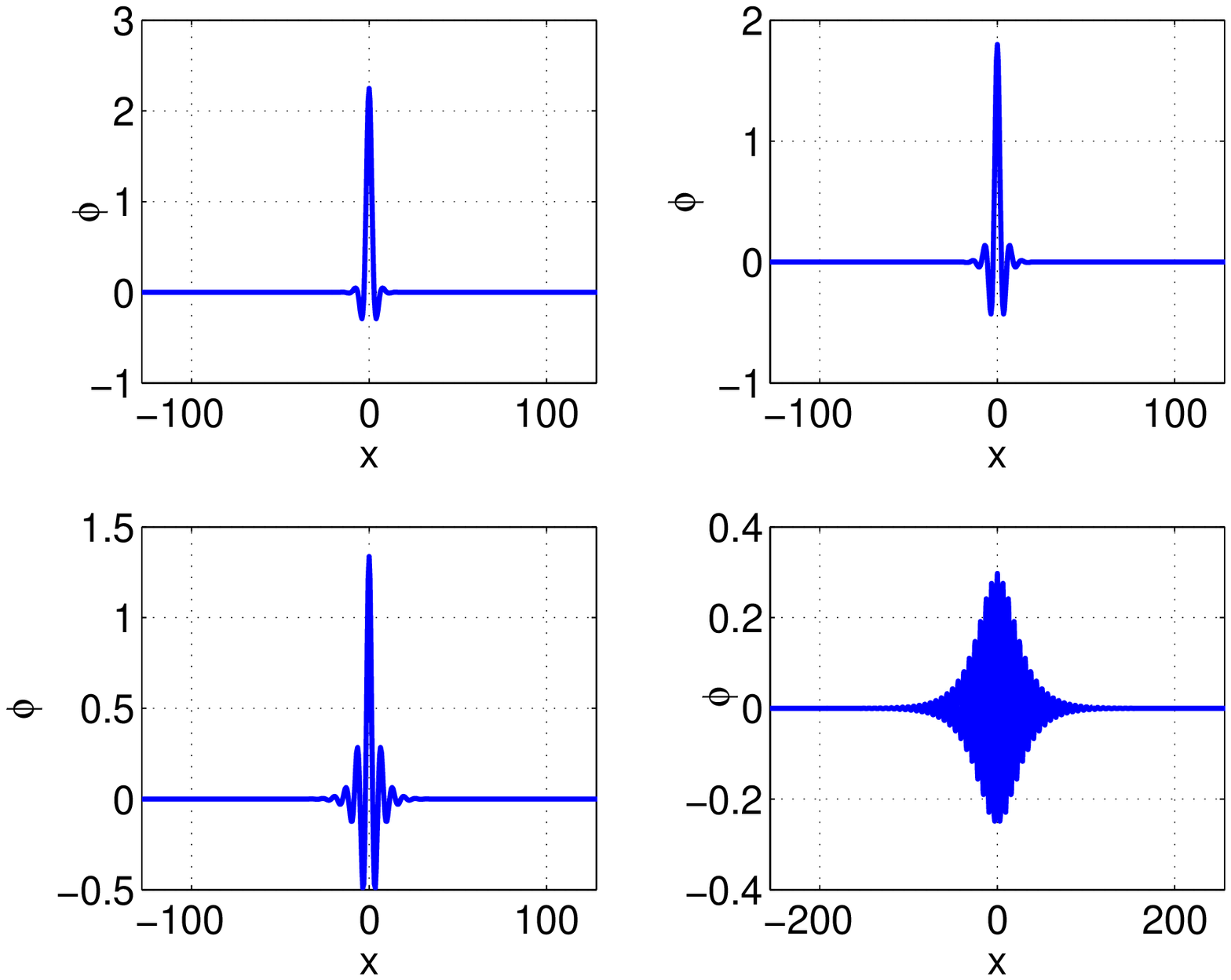}}
\caption{Approximate solitary wave profiles for (\ref{gben8}): $q=1, m=2, r=1, \gamma=1, 1.5, 1.8, 2$.}
\label{gbenfig11}
\end{figure}

\begin{figure}[htbp]
\centering
{\includegraphics[width=0.8\textwidth]{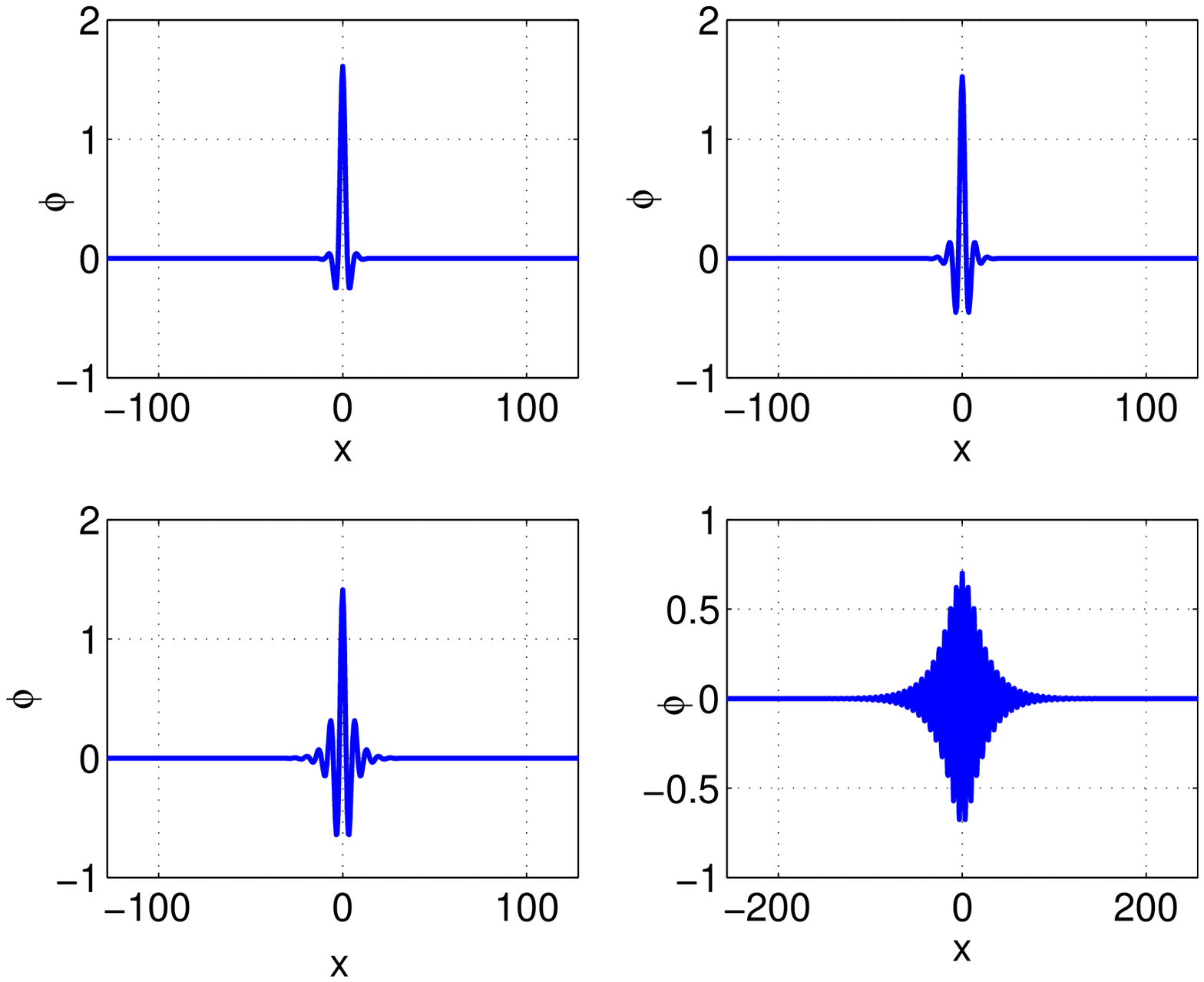}}
\caption{Approximate solitary wave profiles for (\ref{gben8}): $q=4, m=2, r=1, \gamma=1, 1.5, 1.8, 2$.}
\label{gbenfig13}
\end{figure}

Regarding the numerical generation of solitary-wave profiles, the case with general $m$ and $r$ shows a similar behaviour to that of the generalized Benjamin equation. By way of illustration, two examples have been considered. The first corresponds to the values $m=2, r=1, q=1$ with $c_{s}=1.01$. Figure \ref{gbenfig11} shows the form of the solitary-wave profiles generated with $\gamma=1, 1.5, 1.8, 2$; in all these cases $\gamma<\gamma_{max}(c_{s})\approx 2.009975$. As in previous examples, as $\gamma\rightarrow\gamma_{max}(c_{s})$ the profile oscillates more and more and the amplitude diminishes. This behaviour does not seem to depend on the nonlinearity parameter $q$, as exemplified in Figure \ref{gbenfig13}, corresponding to $q=4$. 

%
%

\begin{figure}[htbp]
\centering
{\includegraphics[width=0.8\textwidth]{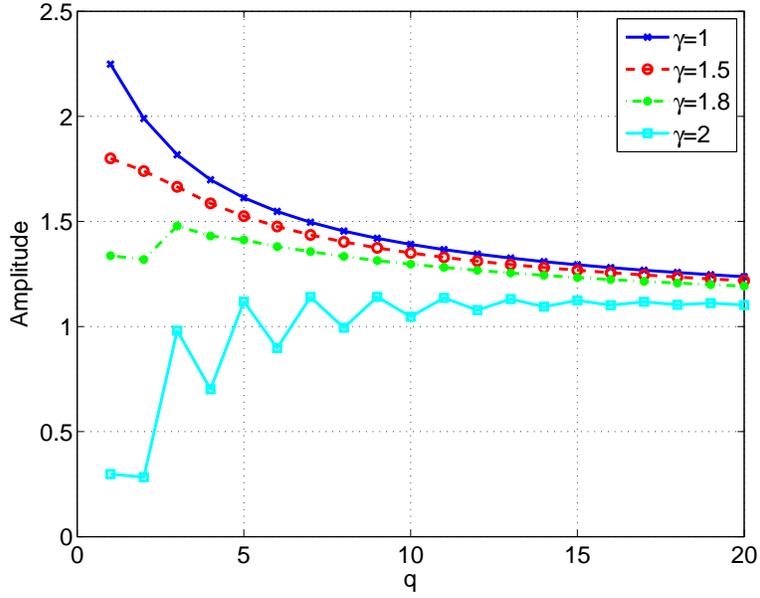}}
\caption{Amplitude vs. $q$ for (\ref{gben8}) with $m=2, r=1$ and several values of $\gamma$.}
\label{gbenfig15}
\end{figure}

\begin{figure}[htbp]
\centering
{\includegraphics[width=0.8\textwidth]{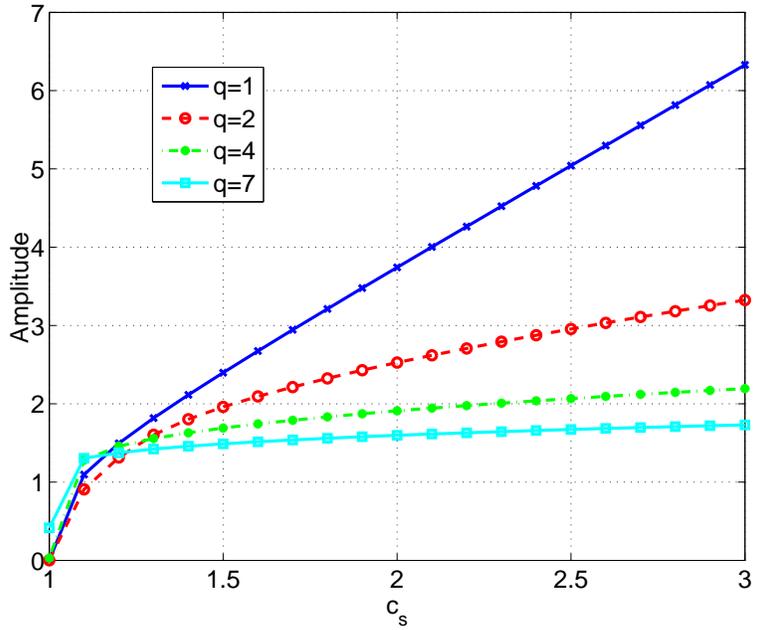}}
\caption{Amplitude vs. speed for (\ref{gben8}) with $m=2, r=1$ and several values of $q$.}
\label{gbenfig16}
\end{figure}

Figure \ref{gbenfig15} depicts the amplitude of the computed profile as function of $q$ and suggests again a similar asymptotic behaviour. 
The corresponding amplitude-speed relation is shown in Figure \ref{gbenfig16}. (Recall Figure \ref{gbf1}(b); the range of values of $c_{s}$ determines $\gamma_{max}(c_{s})$ and consequently the value of $\gamma$ to be taken.) A fit that may give an idea about an approximate amplitude-speed relation is presented in Table \ref{gbentab3}. It differs somewhat from the relation in the case of the gKdV solitary waves.
\begin{table}
\smallskip
\begin{center}
\begin{tabular}{|c|c|c|c|c|}
\hline
 $q$&$K$&$\alpha$&$SSE$&$R^{2}$\\
\hline
1&$1.2129$&$1.218$&$6.4179E-01$&$9.9212E-01$\\
&$(1.2066,1.2192)$&$(1.211,1.225)$&&\\
2&$1.2072$&$0.6501$&$6.4898E-01$&$9.7262E-01$\\
&$(1.2009,1.2135)$&$(0.6429,0.6573)$&&\\
4&$1.1960$&$0.3206$&$8.1200E-02$&$9.8571E-01$\\
&$(1.1938,1.1982)$&$(0.3181,0.3232)$&&\\
7&$1.1571$&$0.1825$&$2.3849E-02$&$9.8704E-01$\\
&$(1.1559,1.1582)$&$(0.1811,0.1839)$&&\\
\hline
\end{tabular}
\end{center}
\caption{Speed-amplitude relation, $m=2, r=1$. Fitting curve $A(c_{s})\approx Kc_{s}^{\alpha}$ and goodness of fit parameters (with 95\% confidence bounds).}\label{gbentab3}
\end{table}
\subsubsection{Asymptotic decay}
\begin{figure}[htbp]
\centering
{\includegraphics[width=0.8\textwidth]{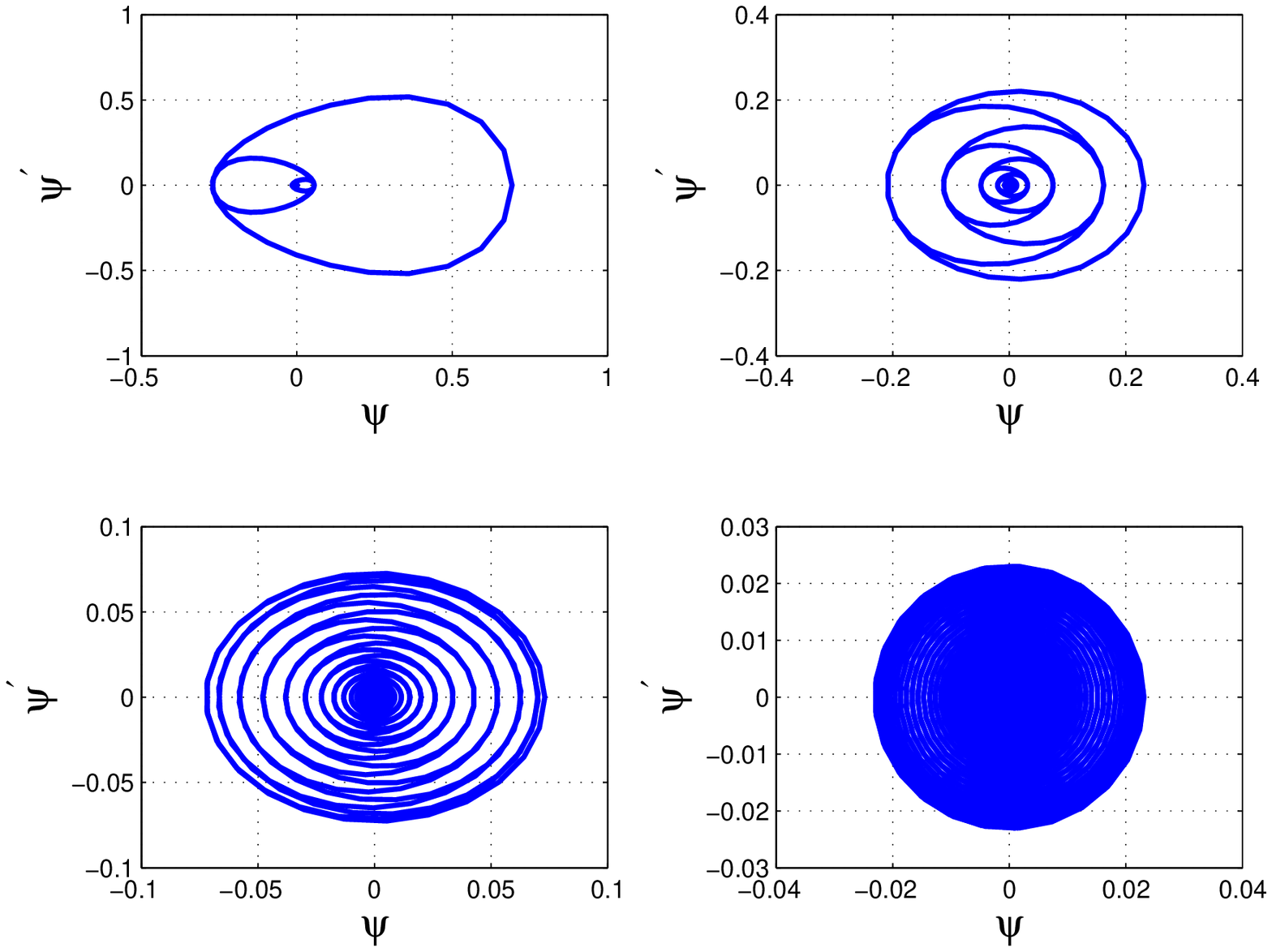}}
\caption{Phase plots of the profiles of Figure \ref{gbenfig3}, $q=2, \widetilde{\gamma}=0.9,0.99,0.999,0.9999$.}
\label{gbenfig3pp}
\end{figure}

\begin{figure}[htbp]
\centering
{\includegraphics[width=0.8\textwidth]{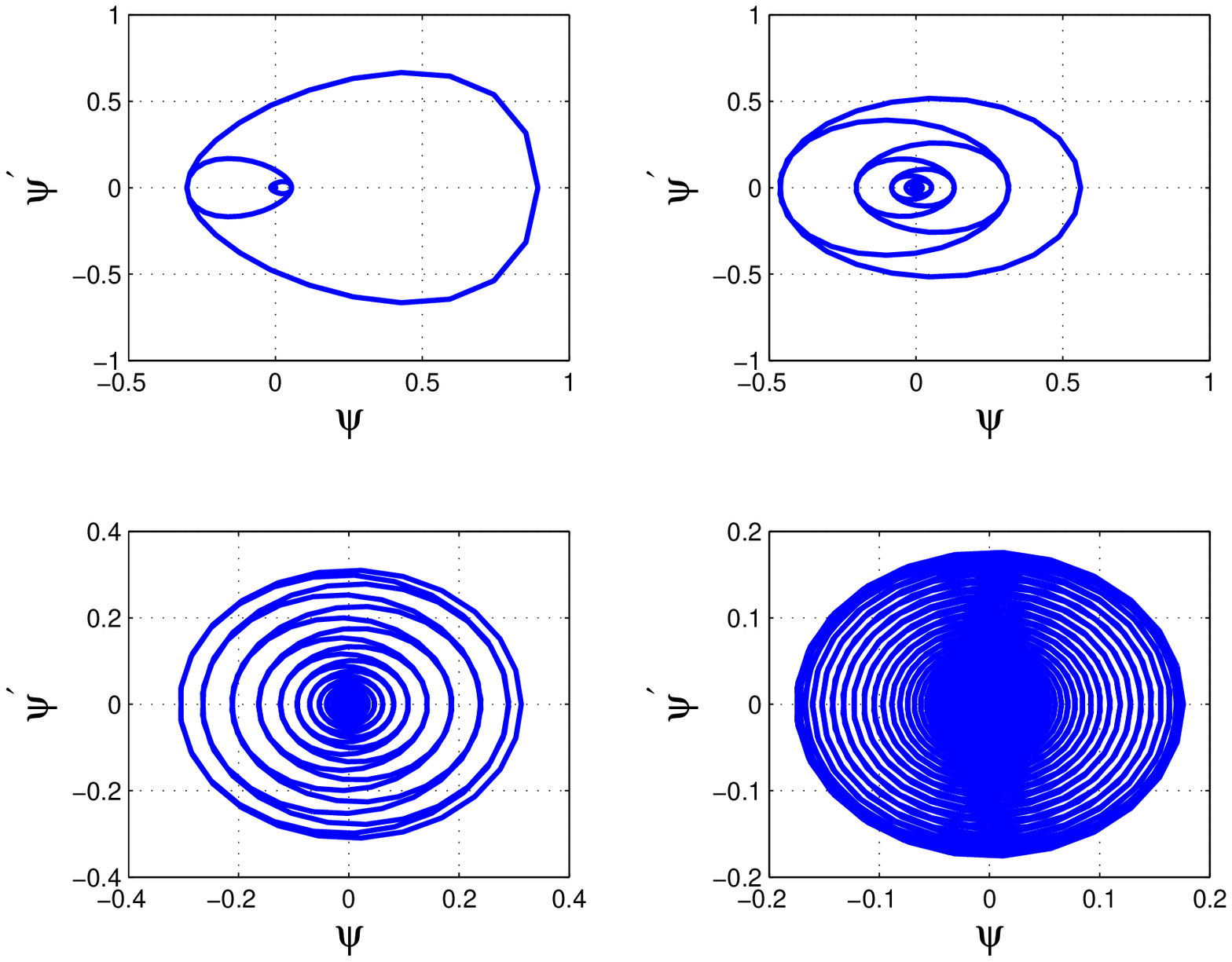}}
\caption{Phase plots of the profiles of Figure \ref{gbenfig4}, $q=4, \widetilde{\gamma}=0.9,0.99,0.999,0.9999$.}
\label{gbenfig4pp}
\end{figure}

The homoclinic to zero character of the computed waves as well as their asymptotic decay properties are illustrated in the following figures and computations.
In the case of the generalized Benjamin equation (normalized version) the phase plots of the profiles displayed in Figures \ref{gbenfig3}, \ref{gbenfig4} are shown in Figures \ref{gbenfig3pp}, \ref{gbenfig4pp}, respectively.

%

\begin{figure}[htbp]
\centering
{\includegraphics[width=0.7\textwidth]{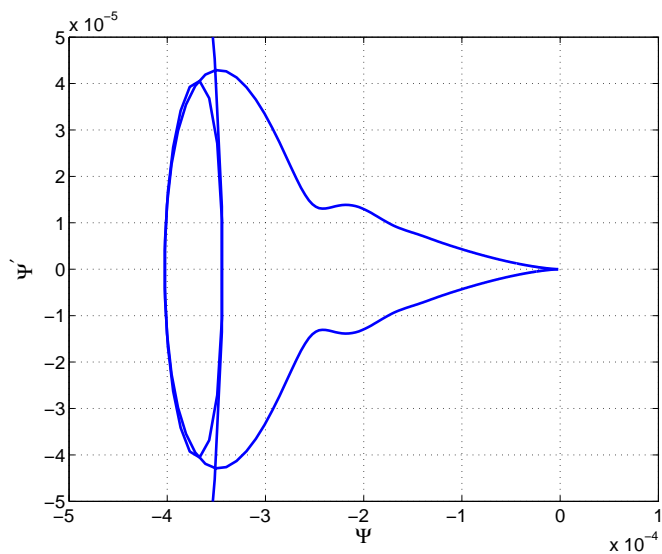}}
\caption{Phase plots of solitary-wave profiles. Magnified detail of Figure \ref{gbenfig3pp} near $(0,0)$ with $q=2, \widetilde{\gamma}=0.9$.}
\label{gbenfig3ppz}
\end{figure}

\begin{figure}
\centering
\includegraphics[width=0.9\textwidth]{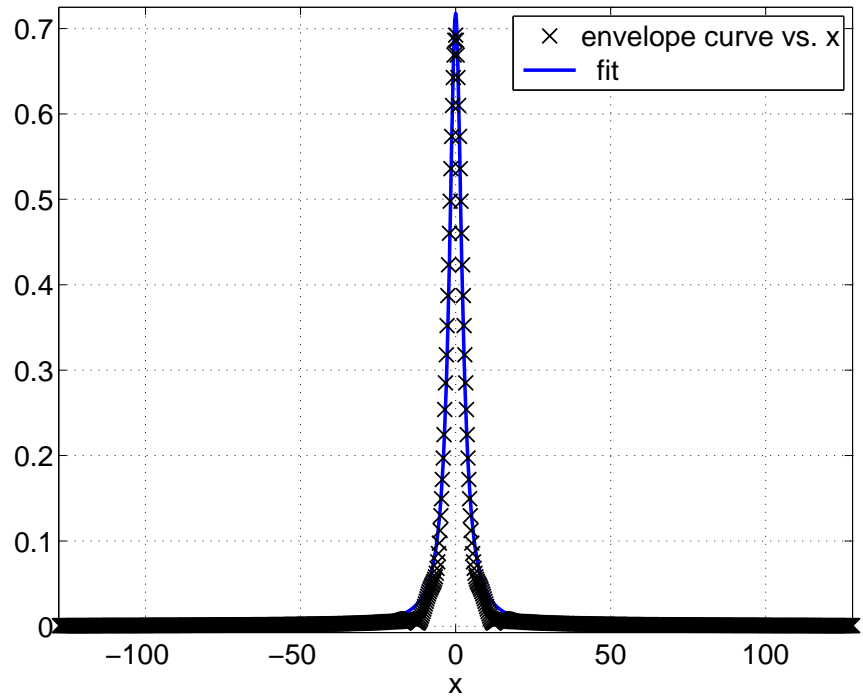}
\caption{Envelope of the absolute value of the computed profile and fitting curve: $q=2, \widetilde{\gamma}=0.9$, from Table \ref{gbtab3p1}.}
\label{gbenfig3ppzfit}
\end{figure}

In order to observe the character of asymptotic decay, a magnified detail is required. This is shown for example in Figure \ref{gbenfig3ppz}, corresponding to $q=2, \widetilde{\gamma}=0.9$, which suggests an algebraic decay, as the cusp behaviour at the origin implies. In fact, for the data 
%
of Figure \ref{gbenfig3ppz}, the envelope of the absolute value of the profile has been computed and fitted to a rational function of the form $f(x)=p_{1}/(x^{2}+q_{1}x+q_{2})$. The resulting coefficients and goodness of fit parameters are displayed in Table \ref{gbtab3p1}. The envelope data and the fitting curve are shown in Figure \ref{gbenfig3ppzfit}.

\begin{table}
\smallskip
\begin{center}
\begin{tabular}{|c|c|c|c|c|}
\hline
 $p_{1}$&$q_{1}$&$q_{2}$&$SSE$&$R^{2}$\\
\hline
$4.218$&$4.589E-05$&$5.876$&$6.019E-02$&$9.923E-01$\\
$(4.182,4.254)$&$(-1.835E-02,1.845E-02)$&$(5.813,5.939)$&&\\
\hline
\end{tabular}
\end{center}
\caption{Envelope of the absolute value of the computed profile for $q=2, \widetilde{\gamma}=0.9$, Figure \ref{gbenfig3ppz}. Fitting curve $f(x)=p_{1}/(x^{2}+q_{1}x+q_{2})$ and goodness of fit parameters (with 95\% confidence bounds).}
\label{gbtab3p1}
\end{table}



\begin{figure}[htbp]
\centering
\includegraphics[width=0.7\textwidth]{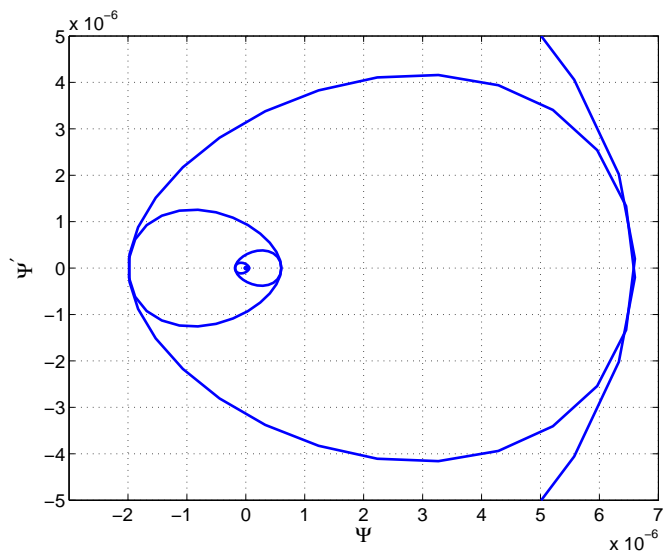}
\caption{Phase plots of the solitary-wave profile of (\ref{gben8b} with $r=1, m=2, q=2, c_{s}=1.01$ and $\gamma=1.5$). Magnified detail near $(0,0)$.}
\label{gbenfig3ppzz}
\end{figure}

\begin{figure}[htbp]
\centering
{\includegraphics[width=0.8\textwidth]{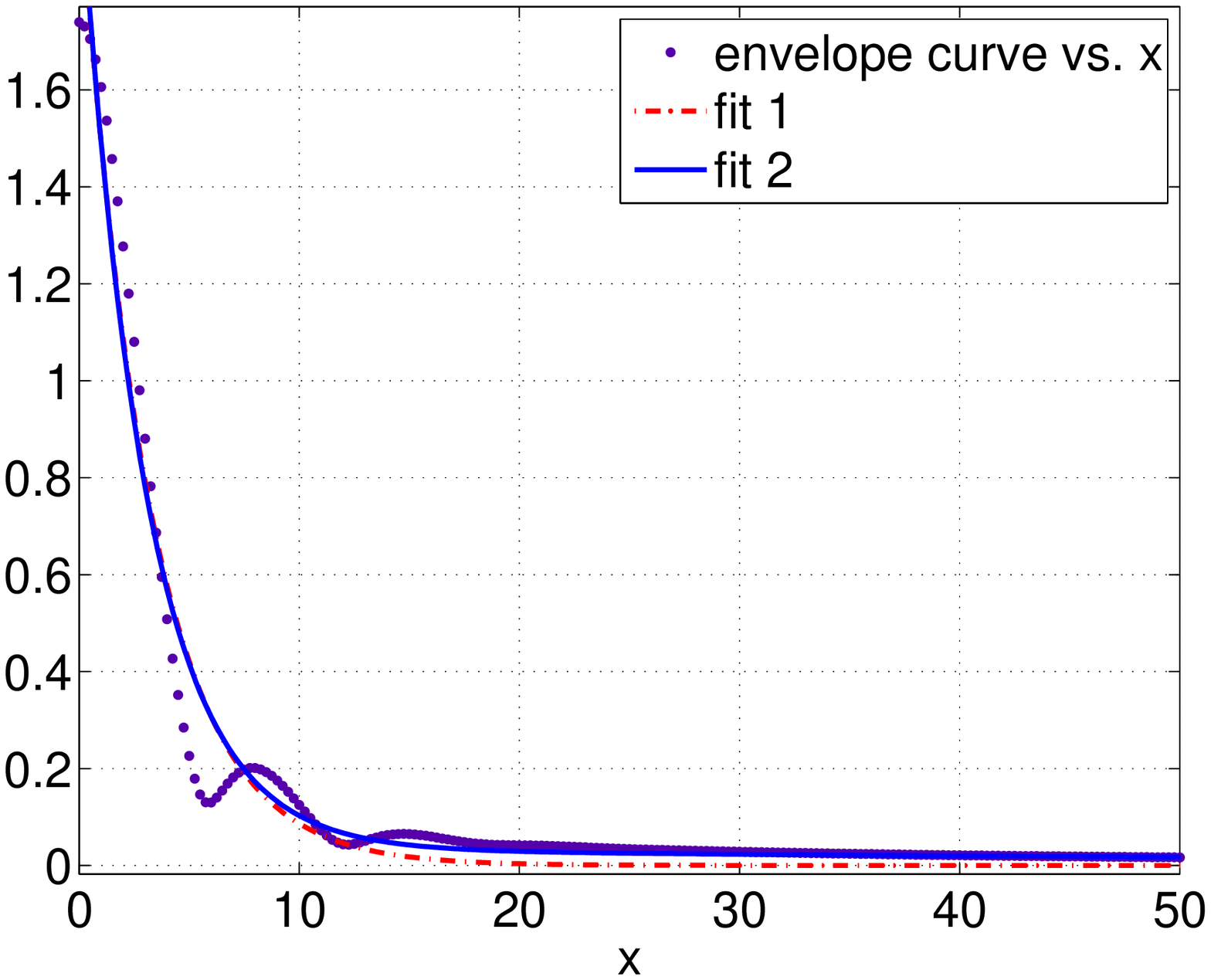}}
\caption{Envelope of the absolute value of the computed profile for $r=1,m=2, q=2, c_{s}=1.01$ and $\gamma=1.5$, see Table \ref{gbentab12bpp}.}
\label{gbenfig12bfit}
\end{figure}

It may be worth illustrating the general case with a couple of examples in order to confirm the results on decay of \cite{BonaC1998}, where Chen and Bona prove exponential decay when $r$ is an integer and algebraic decay of order $2r+1$ otherwise. The first situation is illustrated by computing the profile with $r=1,m=2, q=2, c_{s}=1.01$ (thus $\gamma_{max}(c_{s})\approx 2.009975$) and $\gamma=1.5$. {A magnification of the phase plot. near the origin, of the solitary-wave profile obtained from the normalized equation (\ref{gben8b}) is shown in Figure \ref{gbenfig3ppzz}.}
When fitting the part of the envelope curve  of the absolute value of the profile corresponding to the positive real semiaxis to exponentials, we obtain the results shown in the upper part of Table \ref{gbentab12pp} and in Figure \ref{gbenfig12bfit}. We have tried two fits, one with a function of the form $ae^{bx}$ (fit1) and one with $ae^{bx}+ce^{dx}$ (fit2). As inferred from Table \ref{gbentab12pp}, the second one looks better. This is also observed in Figure \ref{gbenfig12bfit}, for $x$ large, so that the asymptotic decay is evident. The goodness of the fits improves if we consider the data starting from a larger value of $x$. For example, the lower part of Table \ref{gbentab12pp} shows the results corresponding to fitting the part of the curve for $x\geq 120$. In this case, the second fit is much better.

\begin{table}
\smallskip
\begin{center}
\begin{tabular}{|c|c|c|}
\hline
 &fit1&fit2\\
\hline
Parameters&$a=2.049 (2.019,2.078)$&$a=2.034 (2.006,2.062)$\\
(95\% confidence bounds)&$b=-0.3164 (-0.3231,-0.3098)$&$b=-0.3339 (-0.3419,-0.3259)$\\
&&$c=0.03559 (0.02674,0.04444)$\\
&&$d=-0.01438 (-0.01871,-0.01006)$\\
g.o.f.&$SSE=0.842$&$SSE=0.7153$\\
&$R^{2}=0.9703$&$R^{2}=0.9748$\\
\hline
Parameters&$a=0.04114$&$a=-3.737E-05$\\
(95\% confidence bounds)&$(0.03942,0.04286)$&$ (-3.95E-05,-3.54E-05)$\\
&$b=-0.01575 $&$b=0.0141 $\\
&$(-0.01602,-0.01547)$&$(0.01391,0.01428)$\\
&&$c=0.02231 (0.02223,0.02239)$\\
&&$d=-0.01091 (-0-01094,-0.01088)$\\
g.o.f.&$SSE=4.126E-05$&$SSE=1.935E-08$\\
&$R^{2}=0.9725$&$R^{2}=1$\\
\hline
\end{tabular}
\end{center}
\caption{Envelope of the absolute value of the computed profile for $r=1,m=2, q=2, c_{s}=1.01$ and $\gamma=1.5$. Fitting curves $f(x)=ae^{bx}$ (fit1) and $f(x)=ae^{bx}+ce^{dx}$ (fit2) and goodness of fit parameters (with 95\% confidence bounds). Upper data: $x\geq 0$, lower data: $x\geq 120$.}
\label{gbentab12pp}
\end{table}

\begin{figure}[htbp]
\centering
\subfigure
{\includegraphics[width=6cm]{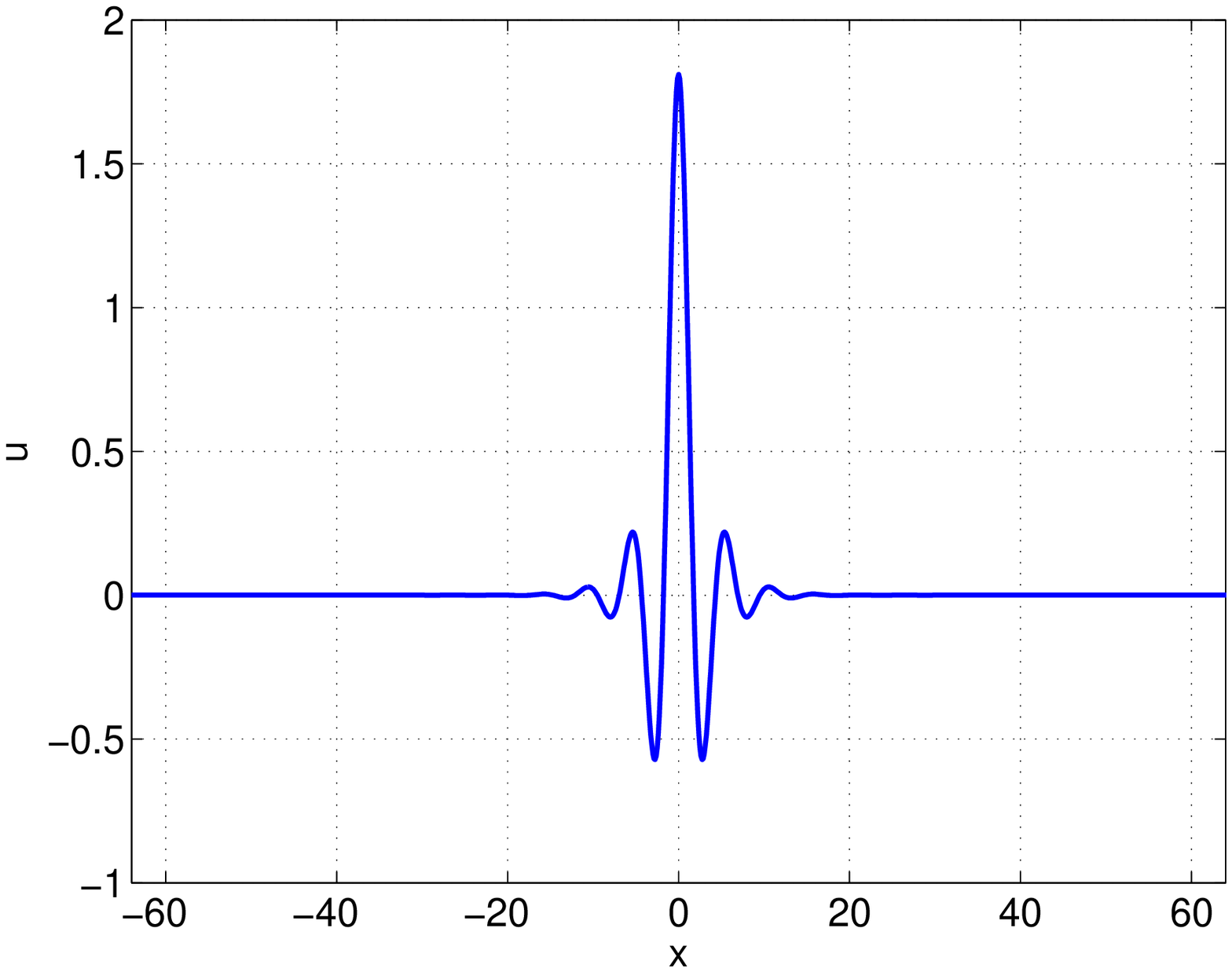}}
\subfigure
{\includegraphics[width=6cm]{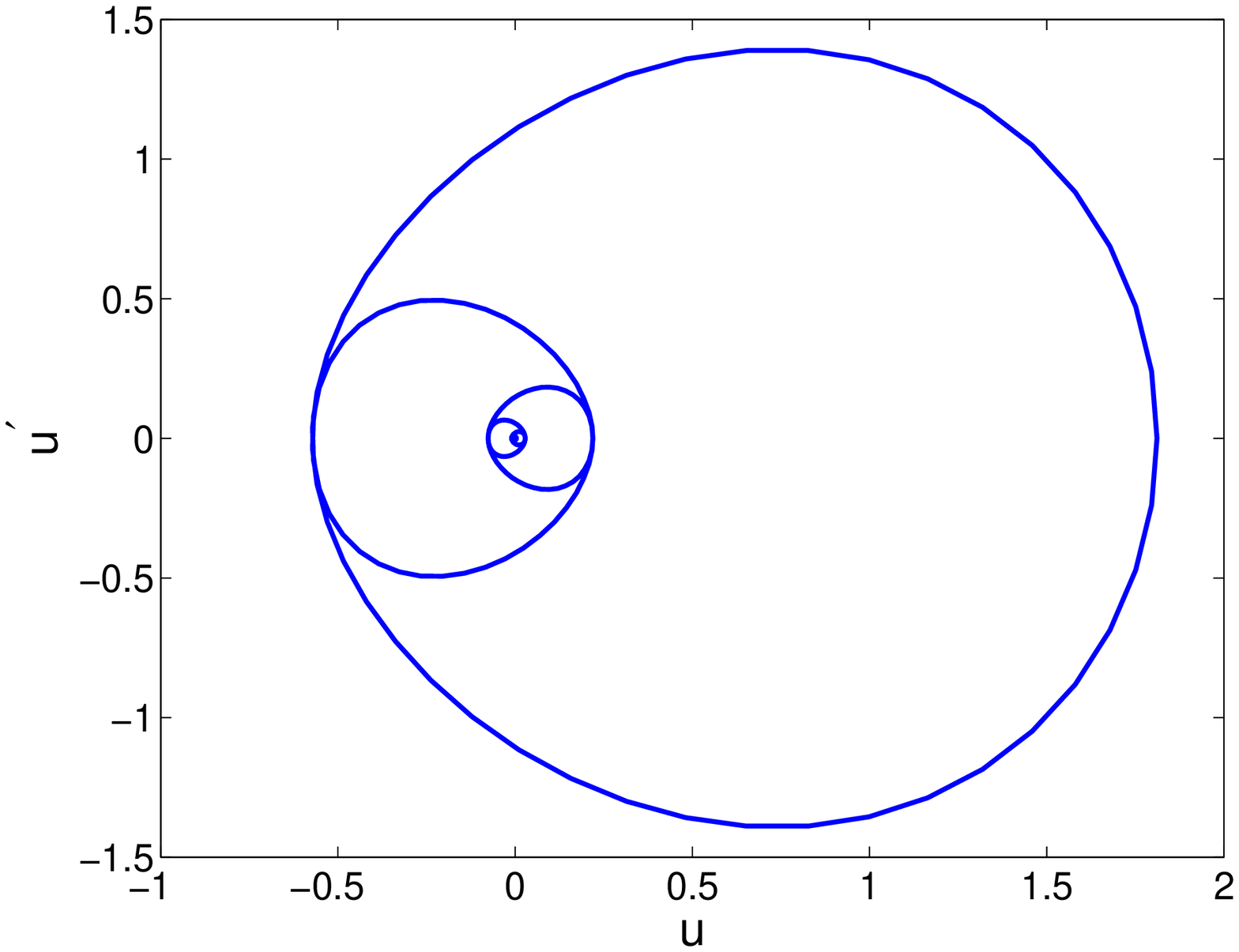}}
\caption{(a) Computed profile for $r=3/2, m=2, q=2, c_{s}=1.01$ and $\gamma=1.5$; (b) Phase plot.}
\label{gbenfig12b}
\end{figure}
\begin{figure}[htbp]
\centering
{\includegraphics[width=0.8\textwidth]{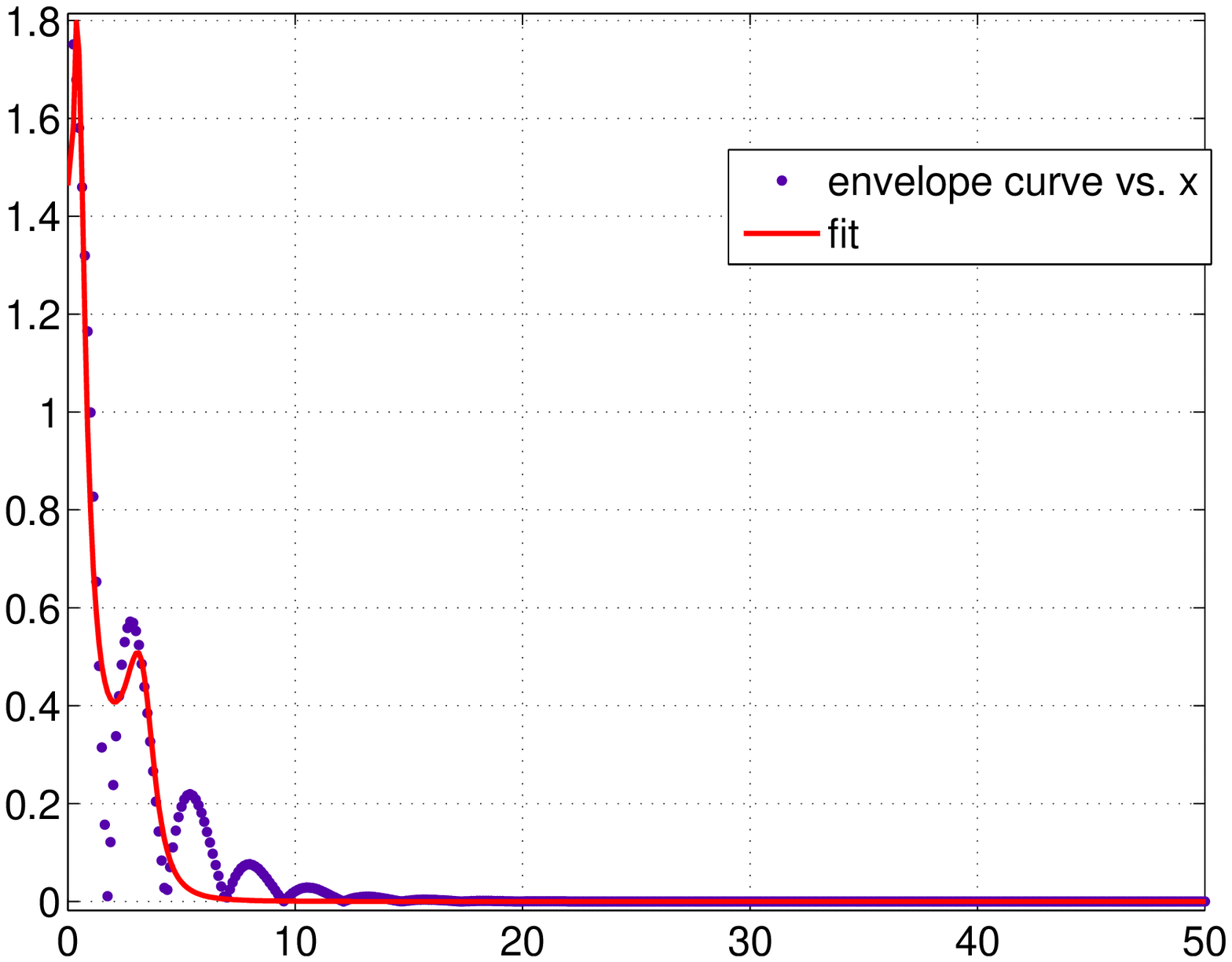}}
\caption{Envelope of the absolute value of the computed profile for $r=3/2,m=2, q=2, c_{s}=1.01$ and $\gamma=1.5$, see Table \ref{gbentab12bpp}.}
\label{gbenfig12bfit2}
\end{figure}

Algebraic decay is illustrated for an approximate profile corresponding to $r=3/2,m=2,q=2,c_{s}=1.01$ (thus $\gamma_{max}(c_{s})\approx 1.759136$) and $\gamma=1.5$. This profile, along with the phase plot,  is shown in Figure \ref{gbenfig12b}. The fit with a rational function of the form $f(x)=p_{1}/(x^{4}+q_{1}x^{3}+q_{2}x^{2}+q_{3}x+q_{4})$ gives the results shown in Table \ref{gbentab12bpp} and Figure \ref{gbenfig12bfit2}. (Recall that the theoretical results in \cite{BonaC1998} imply an algebraic decay of order $1/|x|^{4}$.)


\begin{table}
\smallskip
\begin{center}
\begin{tabular}{|c|c|}
\hline
 Parameters&g.o.f.\\
\hline
$p_{1}= 21.53$&$SSE=1.34$\\
$q_{1}= -10.1$&$R^{2}=0.949$\\
$q_{2}= 31.68$&\\
$q_{3}= -10.93 $&\\
$q_{4}=12.55$&\\
\hline
\end{tabular}
\end{center}
\caption{Envelope of the absolute value of the computed profile for $r=3/2,m=2, q=2, c_{s}=1.01$ and $\gamma=1.5$. Fitting curve $f(x)=p_{1}/(x^{4}+q_{1}x^{3}+q_{2}x^{2}+q_{3}x+q_{4})$ and goodness of fit parameters (with 95\% confidence bounds).}
\label{gbentab12bpp}
\end{table}

%

\subsubsection{Generation of multi-pulses}
Multi-pulse solitary waves can also be generated as in \cite{DougalisDM2012}. In this case superpositions of solitary wave solutions of the gKdV equation are taken as initial iterates and the Petviashvili method with acceleration (using MPE) is used for the numerical generation.

\begin{figure}[htbp]
\centering
\subfigure[]
{\includegraphics[width=6cm]{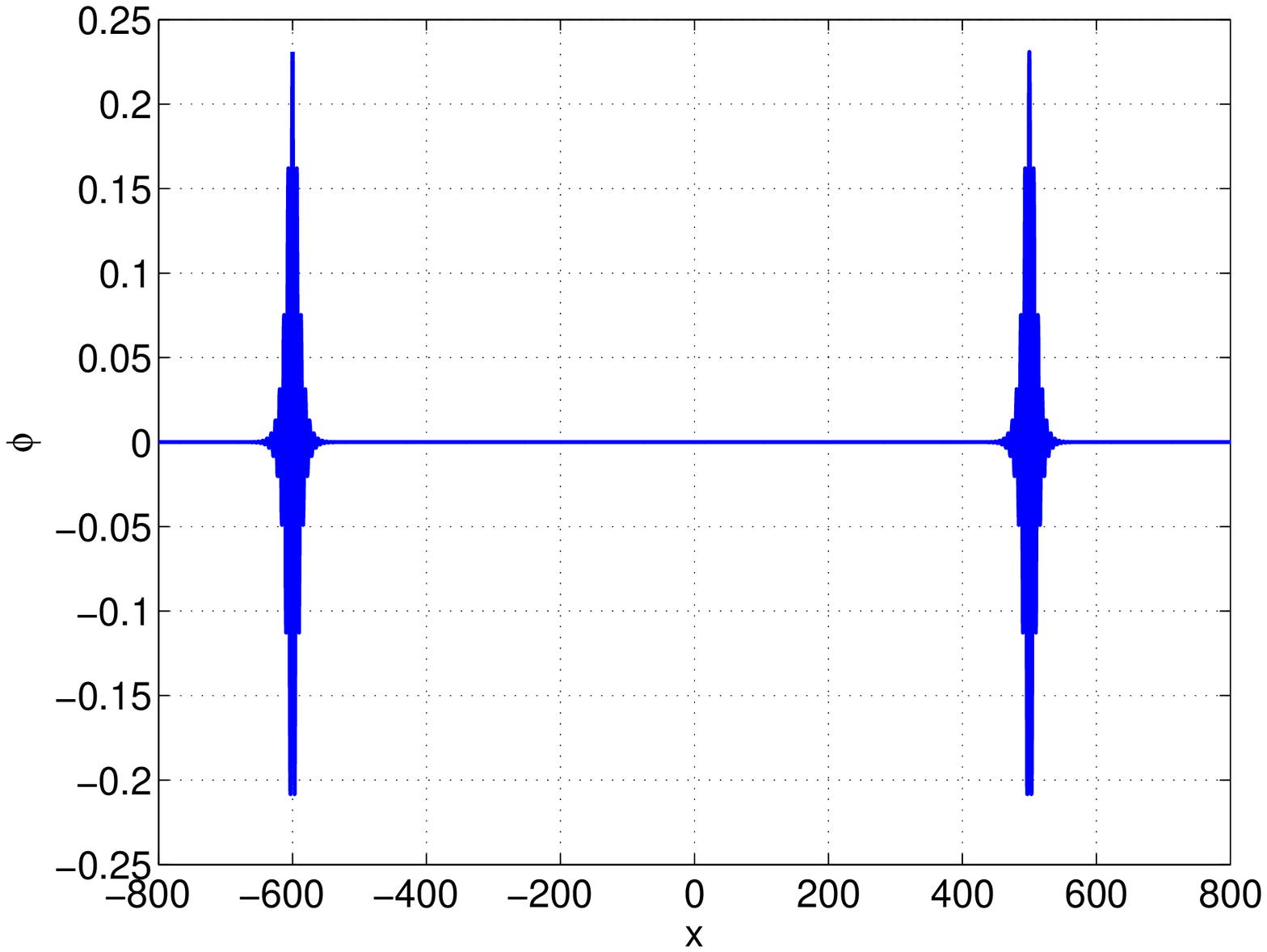}}
\subfigure[]
{\includegraphics[width=6cm]{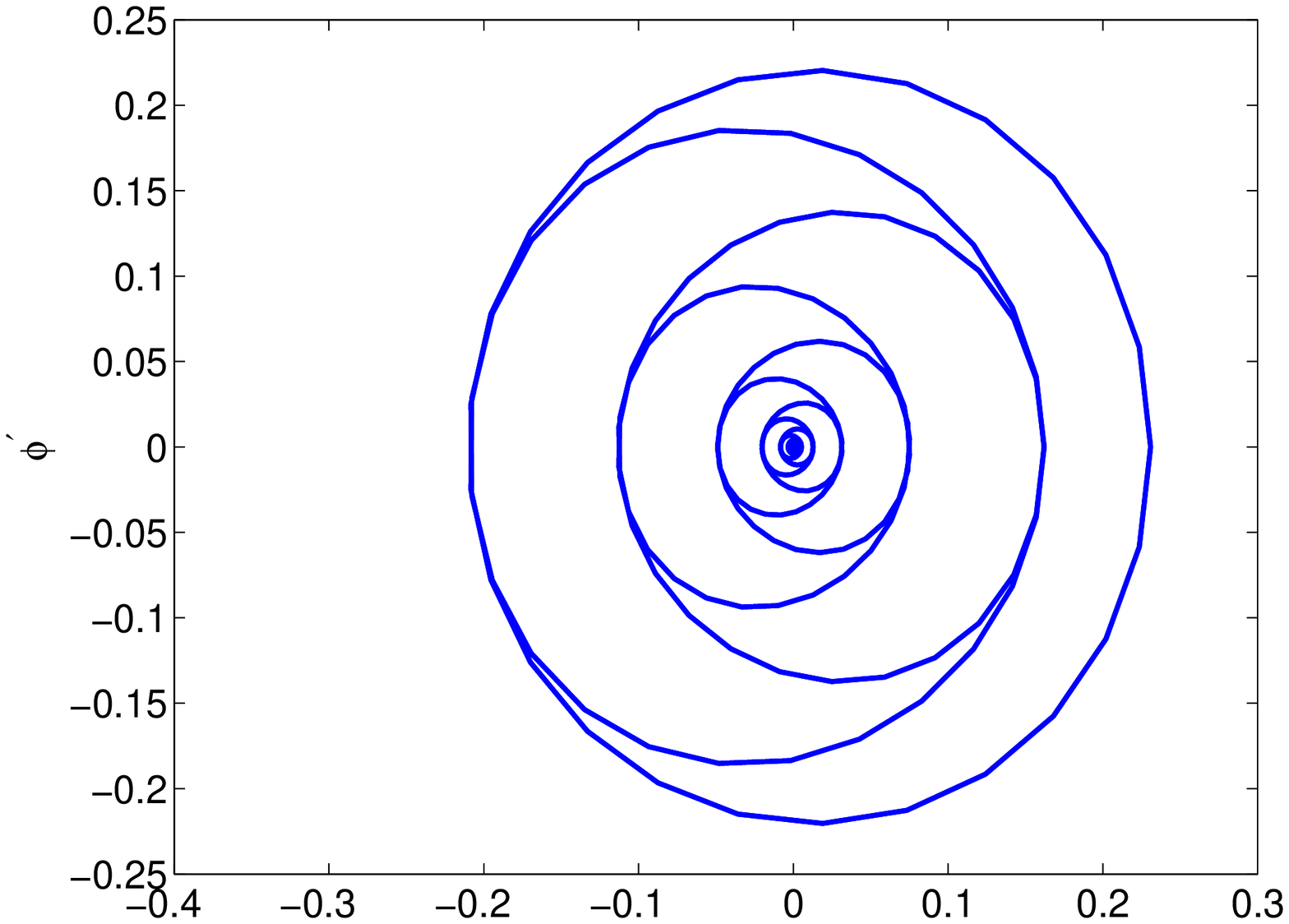}}
\subfigure[]
{\includegraphics[width=6cm]{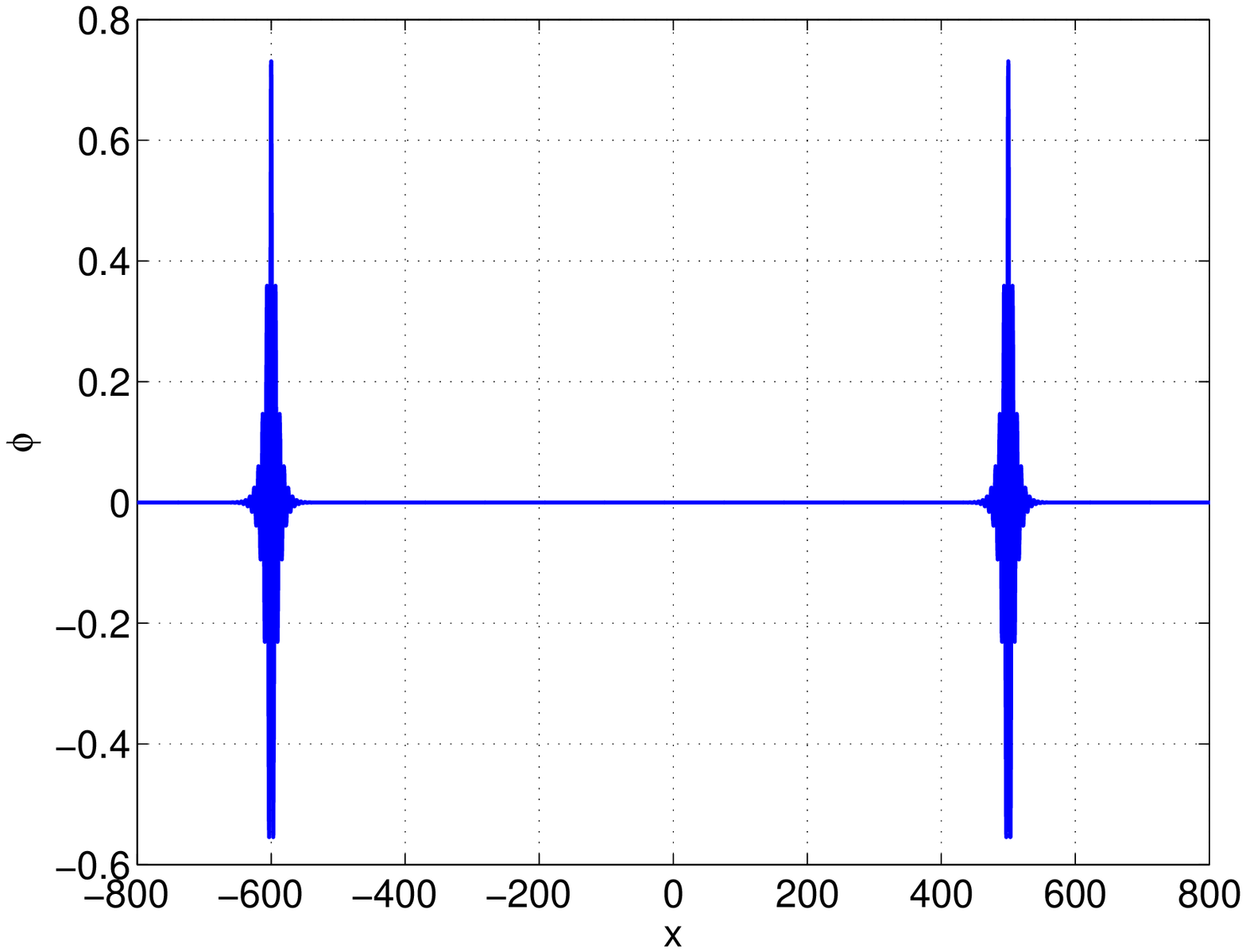}}
\subfigure[]
{\includegraphics[width=6cm]{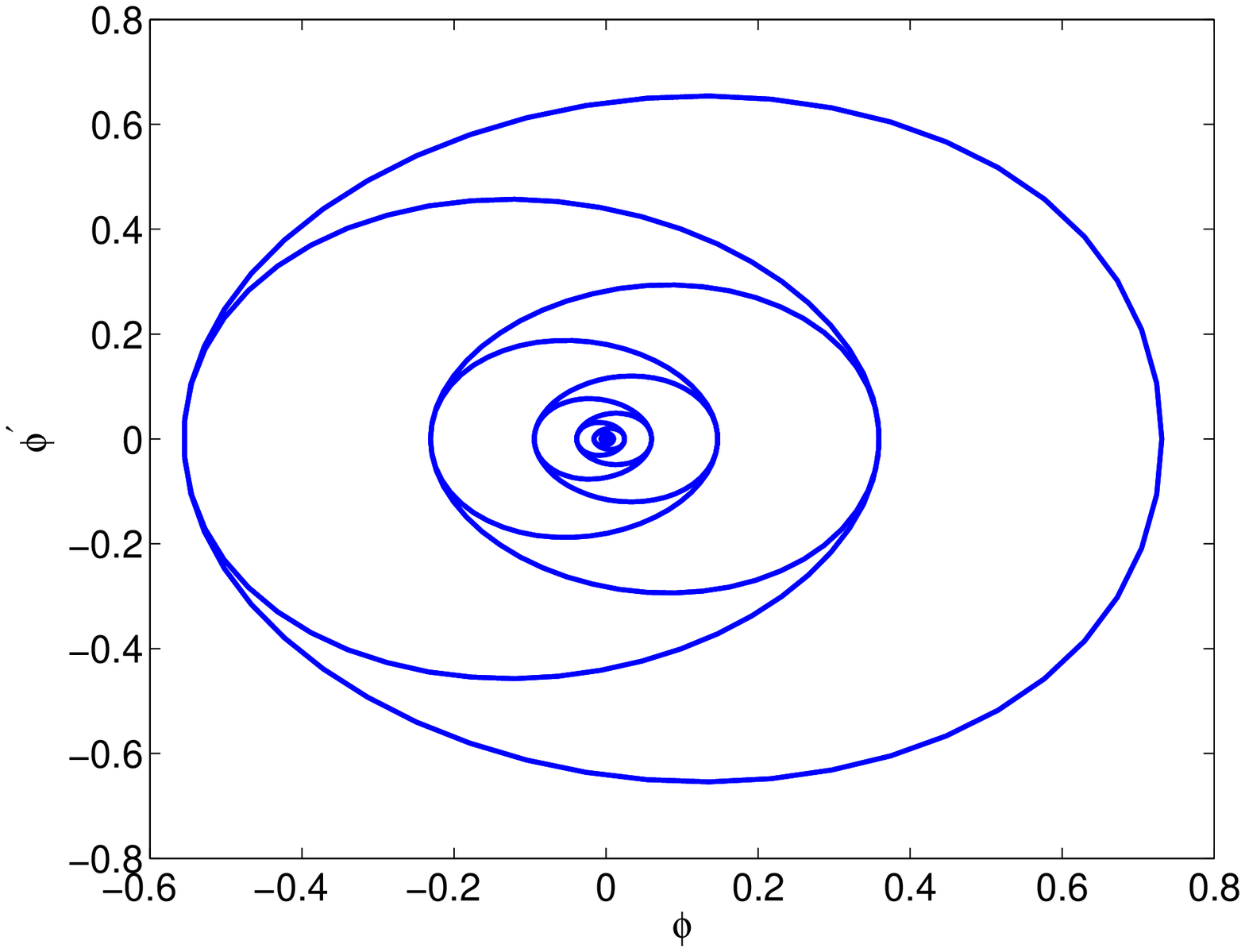}}
\caption{Computed two-pulse profiles (left) and corresponding phase plots (right) for (\ref{gben8b}) with $r=1/2,m=1, c_{s}=1.1$ and $\widetilde{\gamma}=0.99$; (a),(b) $q=2$; (c),(d) $q=6$.}
\label{gbenfig12c}
\end{figure}
\begin{figure}[htbp]
\centering
\subfigure[]
{\includegraphics[width=6cm]{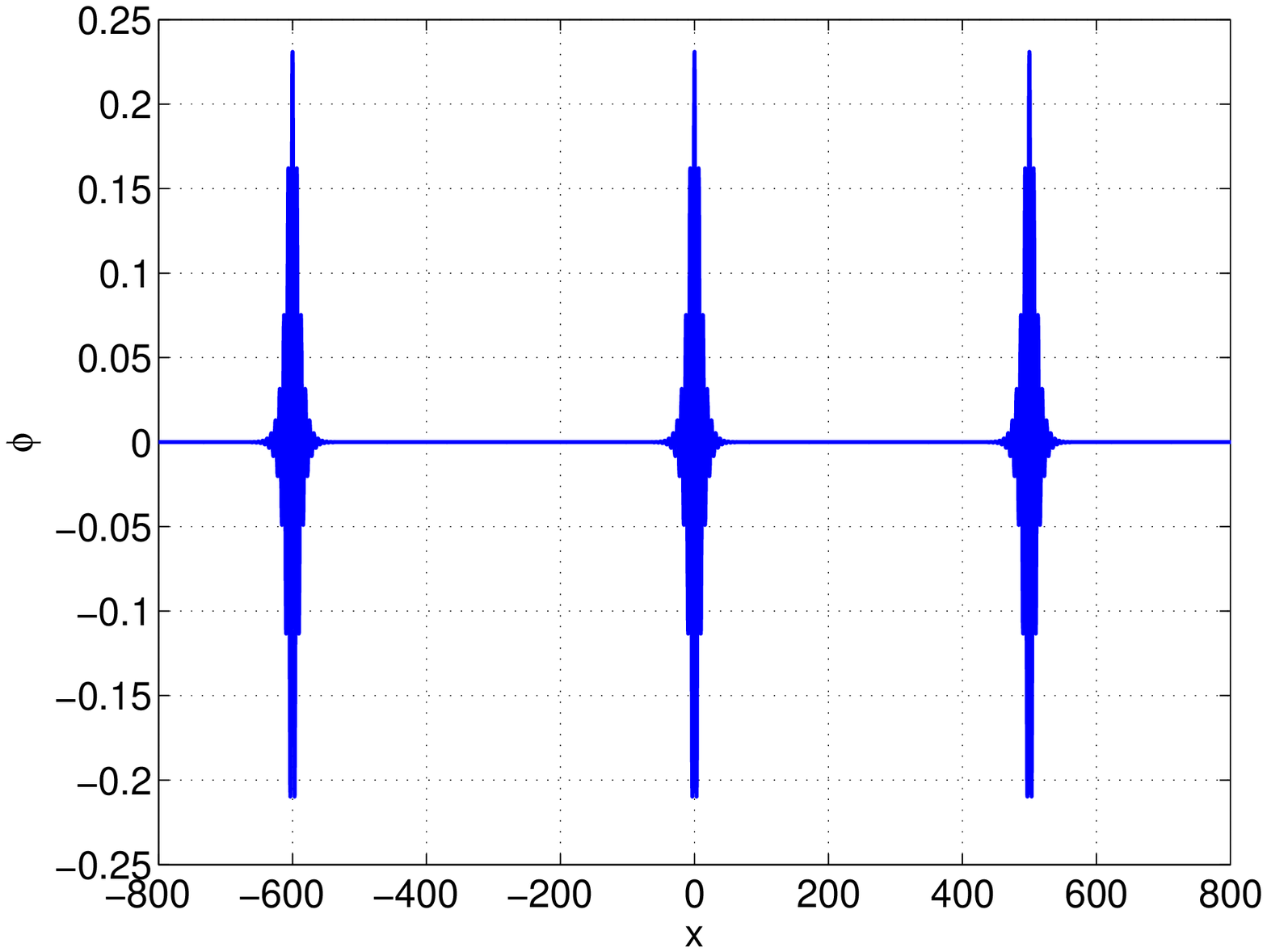}}
\subfigure[]
{\includegraphics[width=6cm]{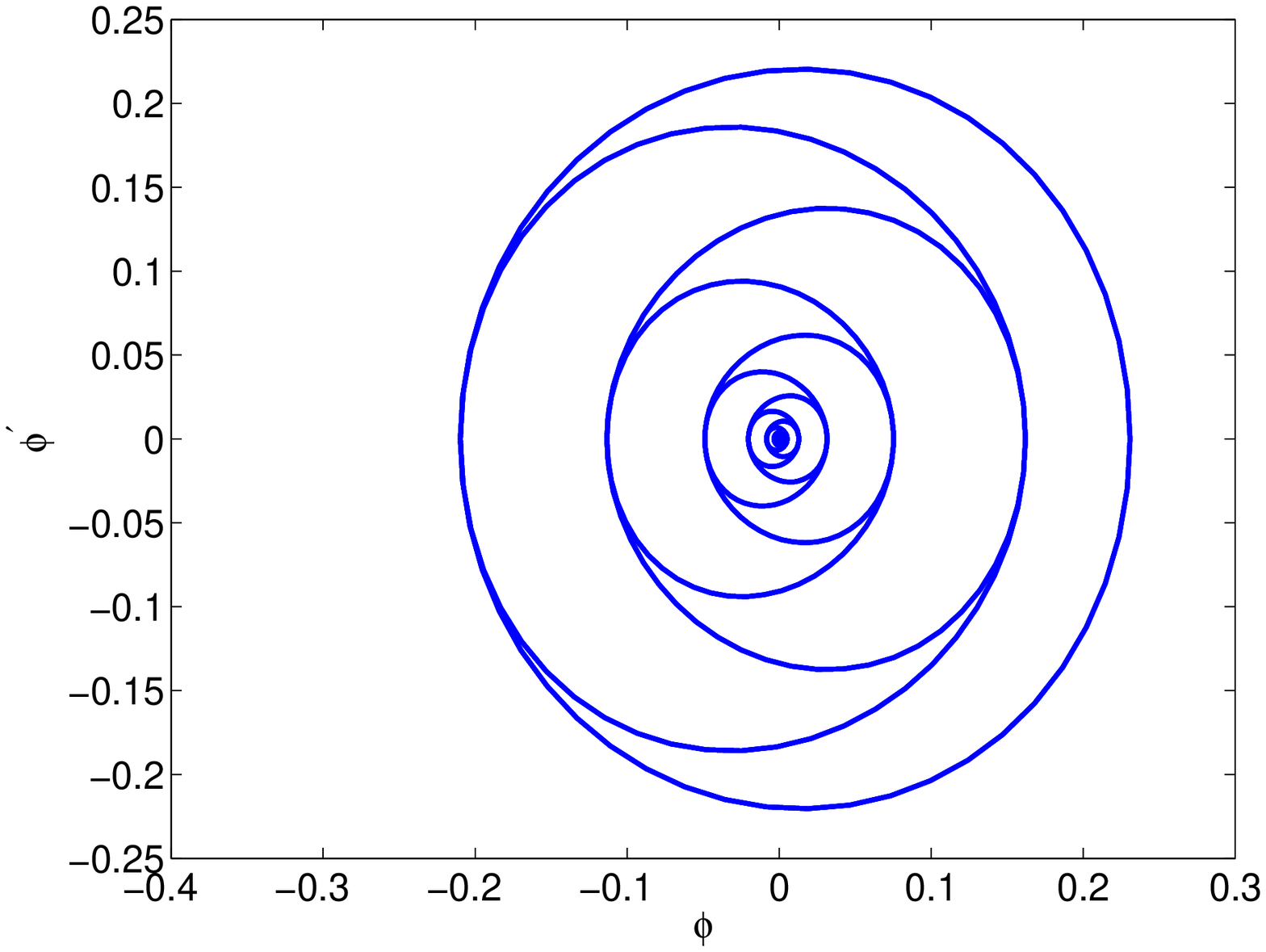}}
\subfigure[]
{\includegraphics[width=6cm]{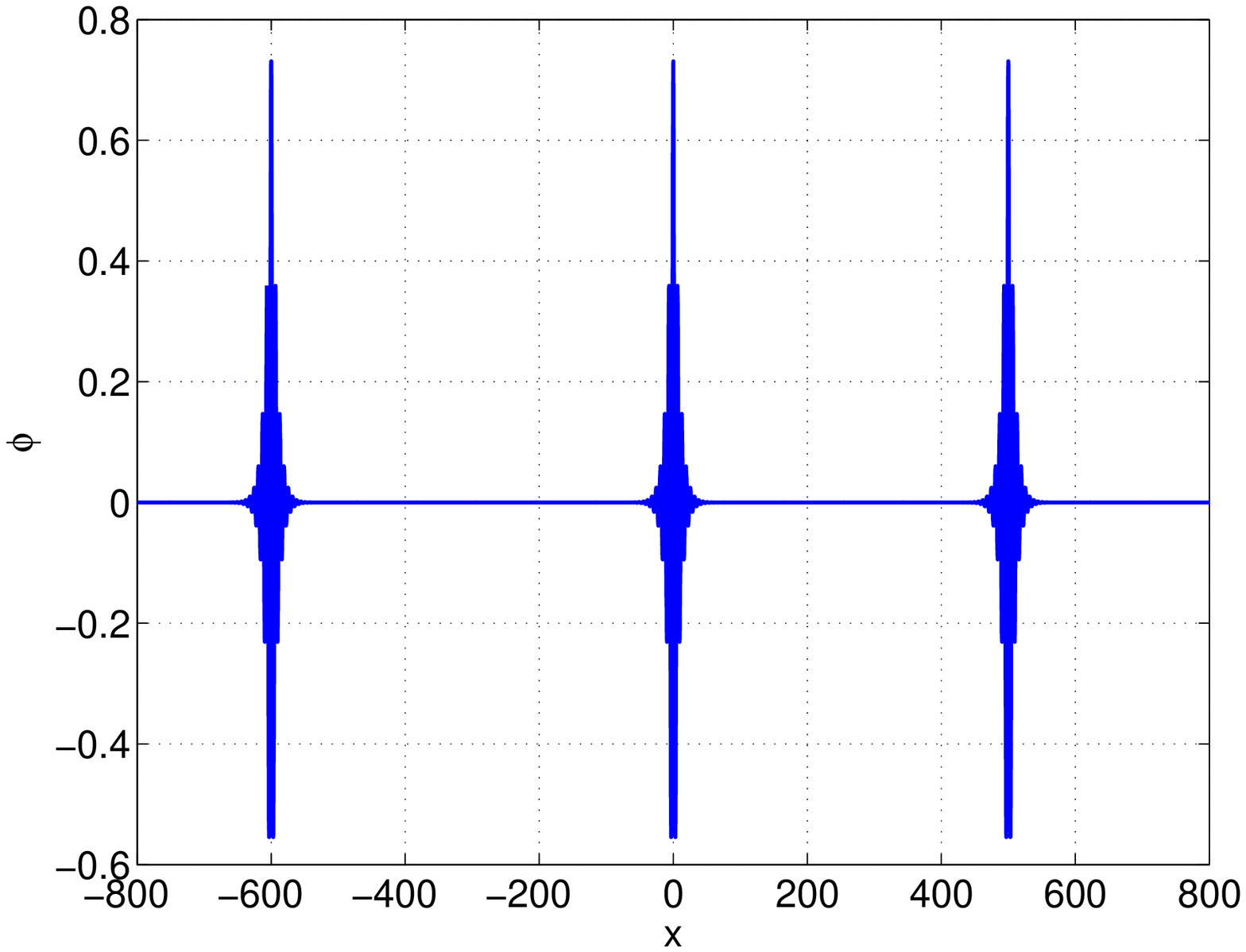}}
\subfigure[]
{\includegraphics[width=6cm]{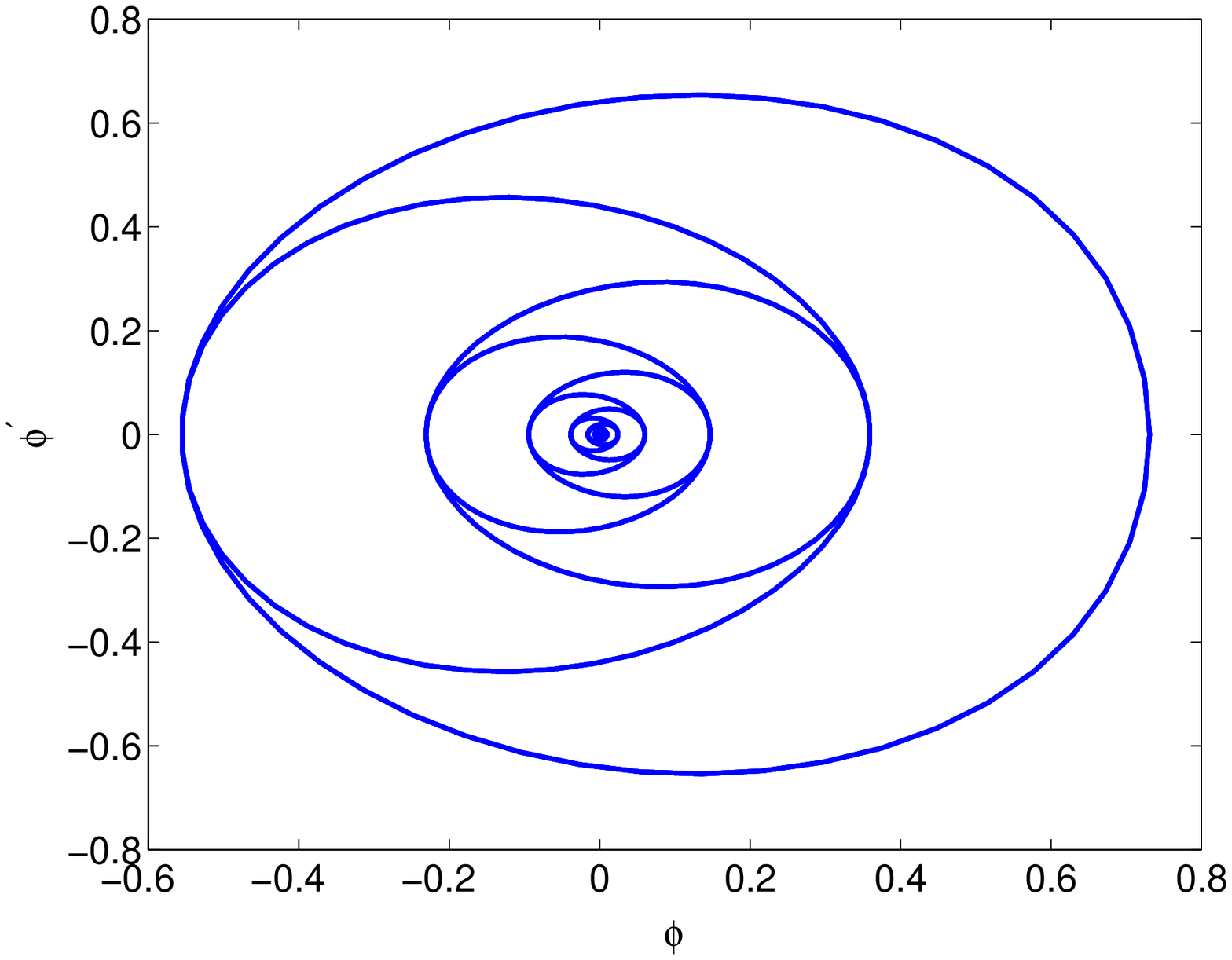}}
\caption{Computed three-pulse profiles (left) and corresponding phase plots (right) for (\ref{gben8b}) with $r=1/2,m=1, c_{s}=1.1$ and $\widetilde{\gamma}=0.99$; (a),(b) $q=2$; (c),(d) $q=6$.}
\label{gbenfig12d}
\end{figure}
\begin{figure}[htbp]
\centering
\subfigure[]
{\includegraphics[width=6cm]{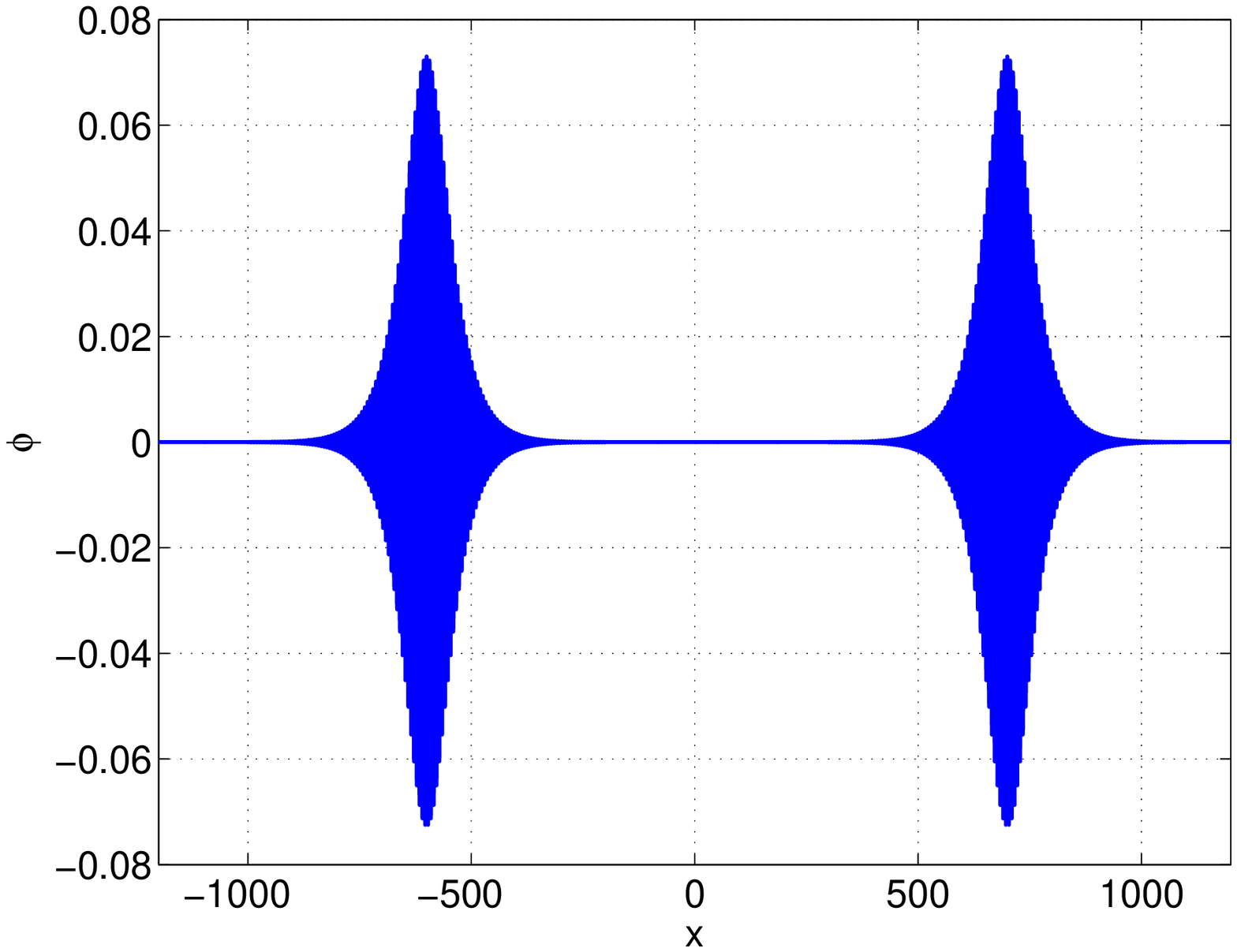}}
\subfigure[]
{\includegraphics[width=6cm]{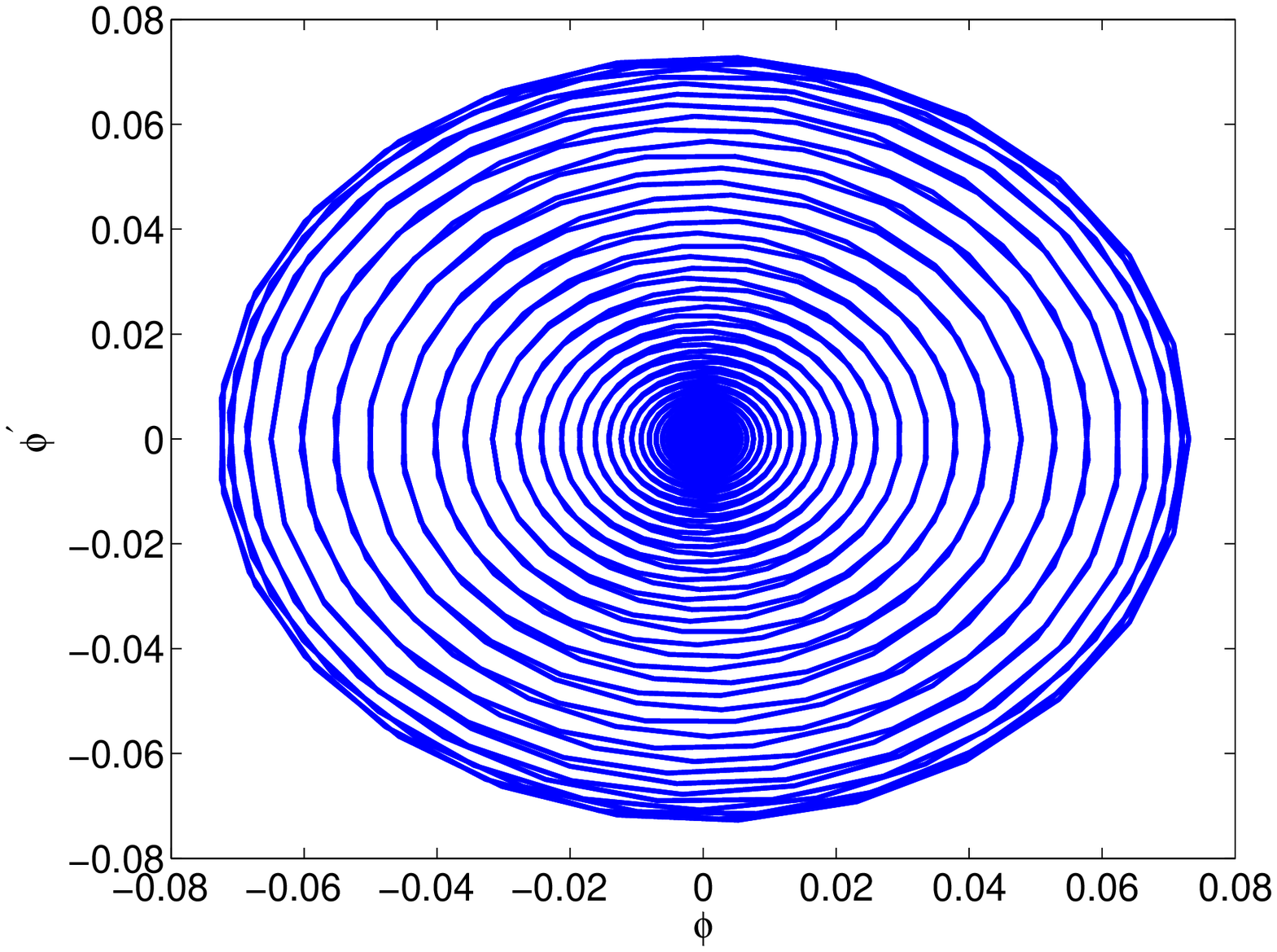}}
\subfigure[]
{\includegraphics[width=6cm]{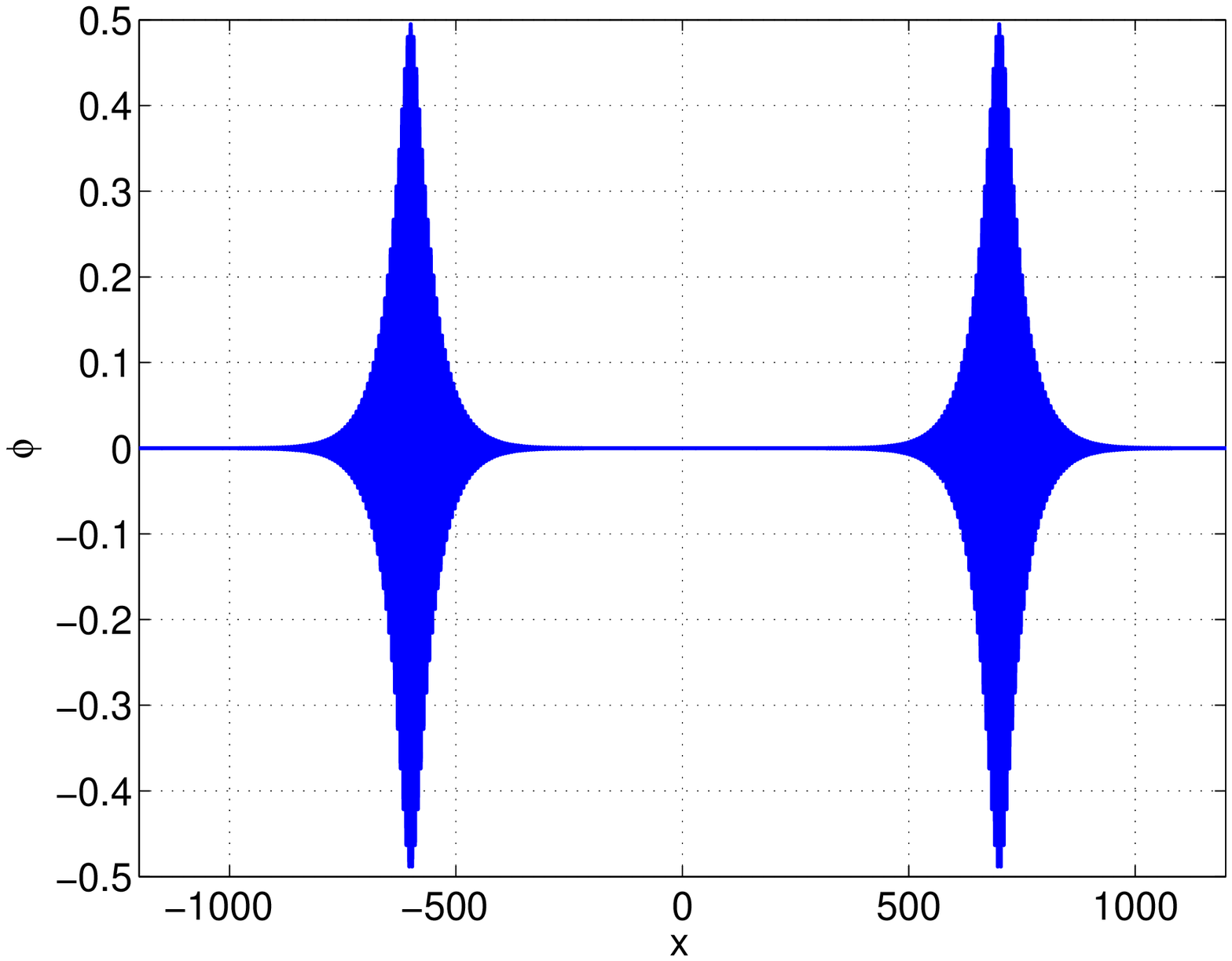}}
\subfigure[]
{\includegraphics[width=6cm]{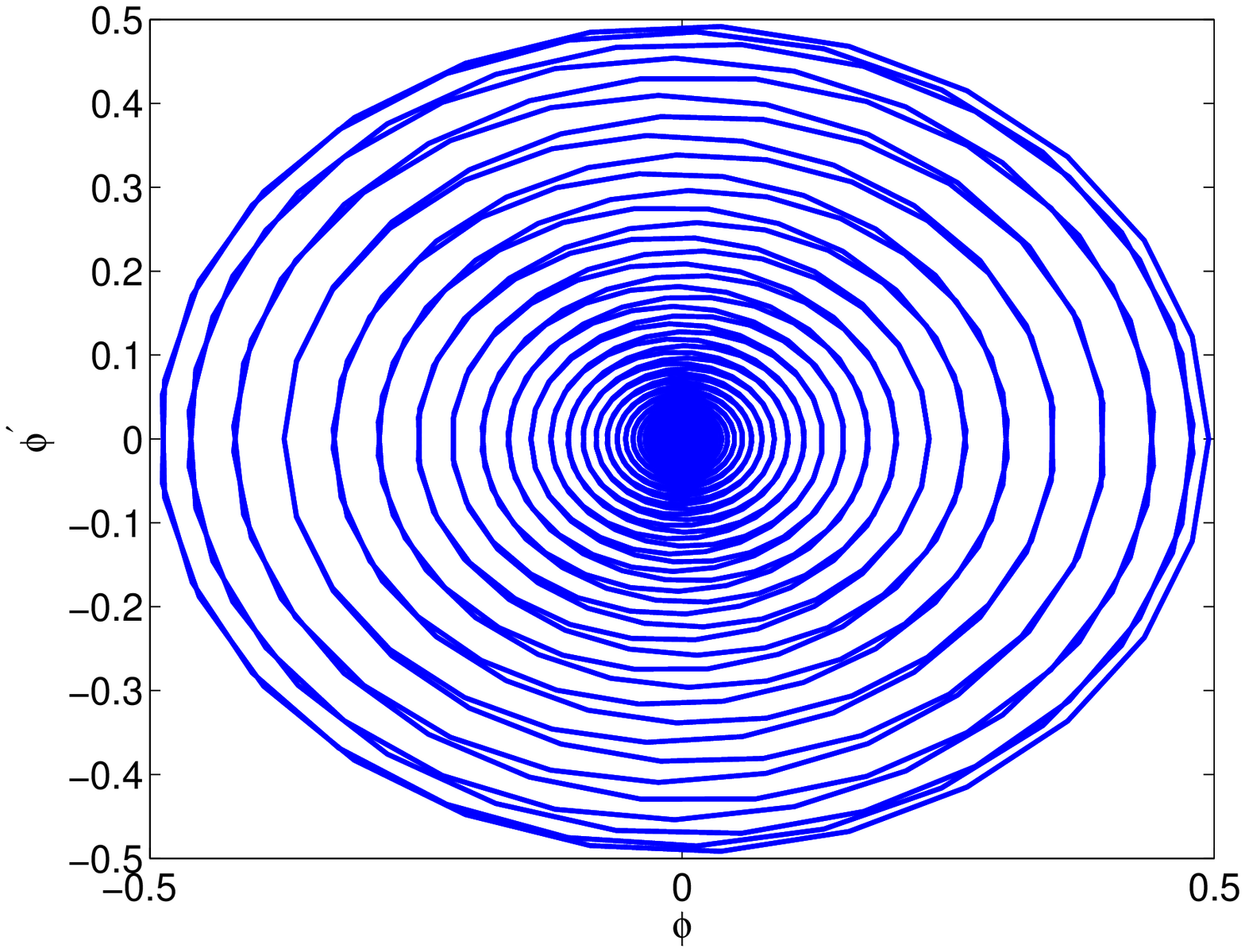}}
\caption{Computed two-pulse profiles (left) and corresponding phase plots (right) for (\ref{gben8b}) with $r=1,m=2, c_{s}=1.1$ and $\widetilde{\gamma}=0.999$; (a),(b) $q=2$; (c),(d)  $q=6$.}
\label{gbenfig12e}
\end{figure}

For simplicity we display some multi-pulses for the normalized equation (\ref{gben8b}). The initial profiles are superpositions of solitary-wave solutions of the gKdV equation, centered at the positions where the final computed waves are placed in the figures.
In the case of the gBenjamin equation ($r=1/2, m=1$) approximate two-pulse profiles are shown in Figure \ref{gbenfig12c} for $\widetilde{\gamma}=0.99$ and $q=2, 6$. The same parameters are used to generate numerically three-pulse profiles in Figure \ref{gbenfig12d}. The general case is illustrated by the two-pulse profiles in Figure \ref{gbenfig12e} for the values $r=1, m=2$,  $\widetilde{\gamma}=0.999$ and $q=2, 6$.

\subsection{A numerical method for the periodic initial-value problem of  Benjamin type equations}
\label{sec23}
In this section we consider a numerical method to study the stability of the solitary waves. The scheme approximates the periodic initial value problem for (\ref{gben1}) on an interval $(-l,l)$. The problem is discretized in space with the standard Fourier-Galerkin spectral method, while time discretization is carried out with a fourth-order, diagonally implicit composition method obtained from the implicit midpoint rule, \cite{FrutosS1992}, as described below. The method will be initially used to verify the accuracy of the numerical profiles, generated in Section \ref{sec3}.
\subsubsection{Semi-discrete scheme}
\label{sec41}
The spatial discretization is standard. For $N\geq 1$ integer we consider
\begin{eqnarray*}
S_{N}=span\{e^{ikx}, -N\leq k\leq N\},
\end{eqnarray*}
and denote by $(\cdot,\cdot)$ the usual $L^{2}$ inner product in $(-l,l)$ with corresponding norm $||\cdot ||$. The semidiscrete Fourier-Galerkin approximation to (\ref{gben1}) is a map $u^{N}:[0,\infty)\rightarrow S_{N}$ such that, for all $\chi\in S_{N}$,
\begin{eqnarray}
&&(u_{t}^{N},\chi)+((-\mathcal{L} u_{x}^{N}+F(u^{N})_{x}),\chi)=0,t>0,\label{gben41}\\
&&u^{N}(x,0)=P_{N}u_{0}(x),\nonumber
\end{eqnarray}
where $P_{N}$ is the $L^{2}$ orthogonal projection onto $S_{N}$. When $\chi=e^{ikx}, k=1,\ldots,N, k\neq 0$, (\ref{gben41}) becomes the initial value problem for the Fourier coefficients of $u^{N}$ given by
\begin{eqnarray}
&&\widehat{u^{N}}_{t}(k,t)+(ik)(-(\delta |k|^{2m}-\gamma |k|^{2r})\widehat{u^{N}}(k,t)+\widehat{F(u^{N})}(k,t))=0, t>0,\label{gben43}\\
&&\widehat{u^{N}}(k,0)=\widehat{u_{0}}(k).\nonumber
\end{eqnarray}
A convergence proof for the semidiscrete problem (\ref{gben41}) will appear elsewhere.
\subsubsection{Fully discrete scheme}
\label{sec42}
{
The ode system (\ref{gben43}) is stiff in general, so an implicit time stepping integrator is typically required. Considered here is a simply diagonally implicit Runge-Kutta (SDIRK) method of Butcher tableau
\begin{equation}\label{RKM2}
\begin{tabular}{c | ccccc}
& $b_{1}/2$ & & &&\\
&  $b_{1}$ & $b_{2}/2$ &&&\\
&$b_{1}$&$b_{2}$&$b_{3}/2$&&\\ \hline
 & $b_{1}$&$b_{2}$&$b_{3}$
\end{tabular},
\end{equation}
with $b_{1}=\displaystyle\frac{1}{2-2^{1/3}}, b_{2}=1-2b_{1}, b_{3}=b_{1}$. The fully discretization is formulated as follows. Given $0<t^{*}<\infty$, a step size $\Delta t$ and $M$ such that $t^{*}=M\Delta t$, we consider a discretization of the interval $[0,t^{*}]$ with the points $t_{m}=m\Delta t, m=0,\ldots,M$. The fully discrete solution corresponding to (\ref{RKM2}) is defined as the sequence $\{U^{m}\}_{m=0}^{M}$ of elements of $S_{N}$, with $U^{0}=P_{N}u_{0}$, satisfying, for every $\chi\in S_{N}$ and $m=0,\ldots,M$
\begin{eqnarray}
&&\left(\frac{U^{m,1}-U^{m}}{\Delta t},\chi\right)-\left(\frac{b_{1}}{2}G(U^{m,1}),\chi\right)=0\nonumber\\
&&\left(\frac{U^{m,i}-U^{m}}{\Delta t},\chi\right)-\left(\sum_{j=1}^{i-1}b_{j}G(U^{m,j})+\frac{b_{i}}{2}G(U^{m,i}),\chi\right)=0\nonumber\\
&&i=2,3,\label{fully}\\
&&\left(\frac{U^{m+1}-U^{m}}{\Delta t},\chi\right)-\left(\sum_{j=1}^{3}b_{j}G(U^{m,j}),\chi\right)=0,\nonumber
\end{eqnarray}
where ${G}:S_{N}\rightarrow S_{N}$ is given by 
\begin{eqnarray*}
({G}(u),\chi)=((\mathcal{L} u_{x}-P_{N}F(u)_{x}),\chi), \quad \forall \chi\in S_{N}.
\end{eqnarray*}
Note that a step of length $\Delta t$ with (\ref{RKM2}) is a composition of the implicit midpoint rule of lengths $b_{i}\Delta t, i=1,2,3$. This leads to the alternative formulation of (\ref{fully}) given by, \cite{FrutosS1992,SanzSC1994}
\begin{eqnarray}
&&\left(\frac{Y^{1}-U^{m}}{b_{1}\Delta t},\chi\right)=\left({G}\left(\frac{Y^{1}+U^{m}}{2}\right),\chi\right)\nonumber\\
&&\left(\frac{Y^{i}-Y^{i-1}}{b_{i}\Delta t},\chi\right)=\left({G}\left(\frac{Y^{i-1}+Y^{i}}{2}\right),\chi\right)\nonumber\\
&&i=2,3,\label{fully2}\\
&&U^{m+1}=Y^{3}.\nonumber
\end{eqnarray} 
The method (\ref{fully2}) is implemented in the Fourier space and the fixed-point iteration is used to solve numerically the systems for the intermediate stages.

The SDIRK method (\ref{RKM2}) belongs to the family of high-order, symplectic, symmetric somposition methods based on the implicit midpoiny rule and constructed by Yoshida in \cite{Yoshida1990}. It was applied to wave problems in \cite{FrutosS1992}, where some properties are emphasized, namely:
\begin{itemize}
\item The nondissipative, symplectic and symmetric character of the method ensures its
suitability for long time integration, the preservation of quadratic invariants and an accurate simulation of the evolution of the energy, \cite{SanzSC1994,HairerLW2004}.
\item Since $b_{2}<0$, the integration of dissipative equations would require some stability step size restriction.
\item
With the aim of preserving the convergence of order four, some implementation details concerning the fixed-point iterations for the intermediate stages are necessary.
\end{itemize}
}
\subsection{Numerical experiments of validation}
\label{sec24}
We first check the accuracy of the solitary-wave profiles computed in Section \ref{sec3} as traveling-wave solutions of the initial-value problem.
We will consider one example for the generalized Benjamin equation ($r=1/2, m=1$). (The behaviour does not change when other values of the coefficients in (\ref{gben1}) are considered.)
IThe computed solitary-wave profile will be taken as initial condition for the evolution code with spatial step size $h=2l/N=0.125$ and will be integrated up to a final time $t^{*}=100$ with three time steps $\Delta t=6.25E-03, 3.125E-03, 1.5625E-03$. The evolution of the amplitude, speed and phase errors will be tracked, as well as the discrete versions of the invariants, i.e. the quantities
\begin{eqnarray}
&&I_{h}(U)=h\sum_{j=0}^{M-1}U_{j}^{2},\quad U=(U_{0},\ldots,U_{M-1})^{T},\label{gben46}\\
&&E_{h}(U)=h\sum_{j=0}^{M-1}\left(\delta |D^{m}U)_{j}|^{2}-\gamma U_{j}(\mathbb{H}^{2r}U)_{j}-\frac{2}{(q+1)(q+2)}U_{j}^{q+2}\right),\label{gben47}
\end{eqnarray}
where $\mathbb{H}$ is the discrete version of $H$, computed in terms of the discrete Fourier components as
\begin{eqnarray*}
\widehat{\mathbb{H}U}_{k}=|k|\widehat{U}_{k},\quad -N\leq k\leq N,
\end{eqnarray*}
while $D$ is the pseudospectral differentiation matrix.

\begin{figure}[htbp]
\centering
{\includegraphics[width=0.8\textwidth]{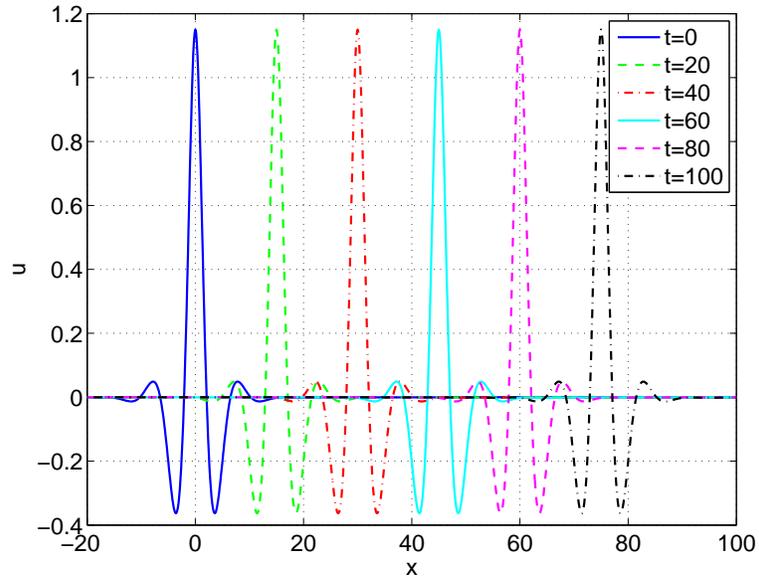}}
\caption{Numerical approximation of gBenjamin equation from the profile computed with $\gamma=1.5, \delta=1, c_{s}=0.75, q=2$ and $h=0.125, \Delta t=1.5625E-03$, at several times.}
\label{gbenfig4321}
\end{figure}
\begin{figure}[htbp]
\centering
\subfigure[]
{\includegraphics[width=6cm]{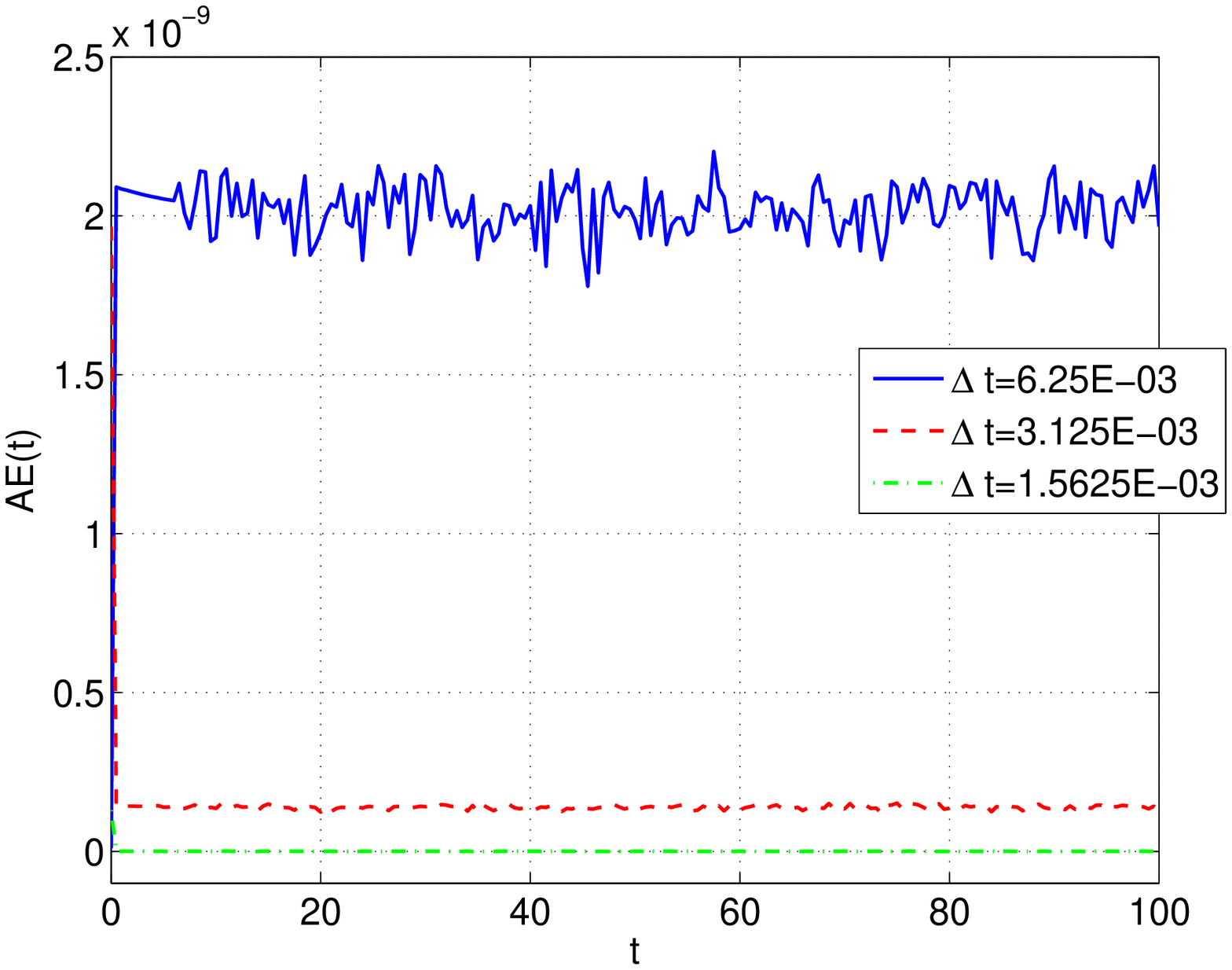}}
\subfigure[]
{\includegraphics[width=6cm]{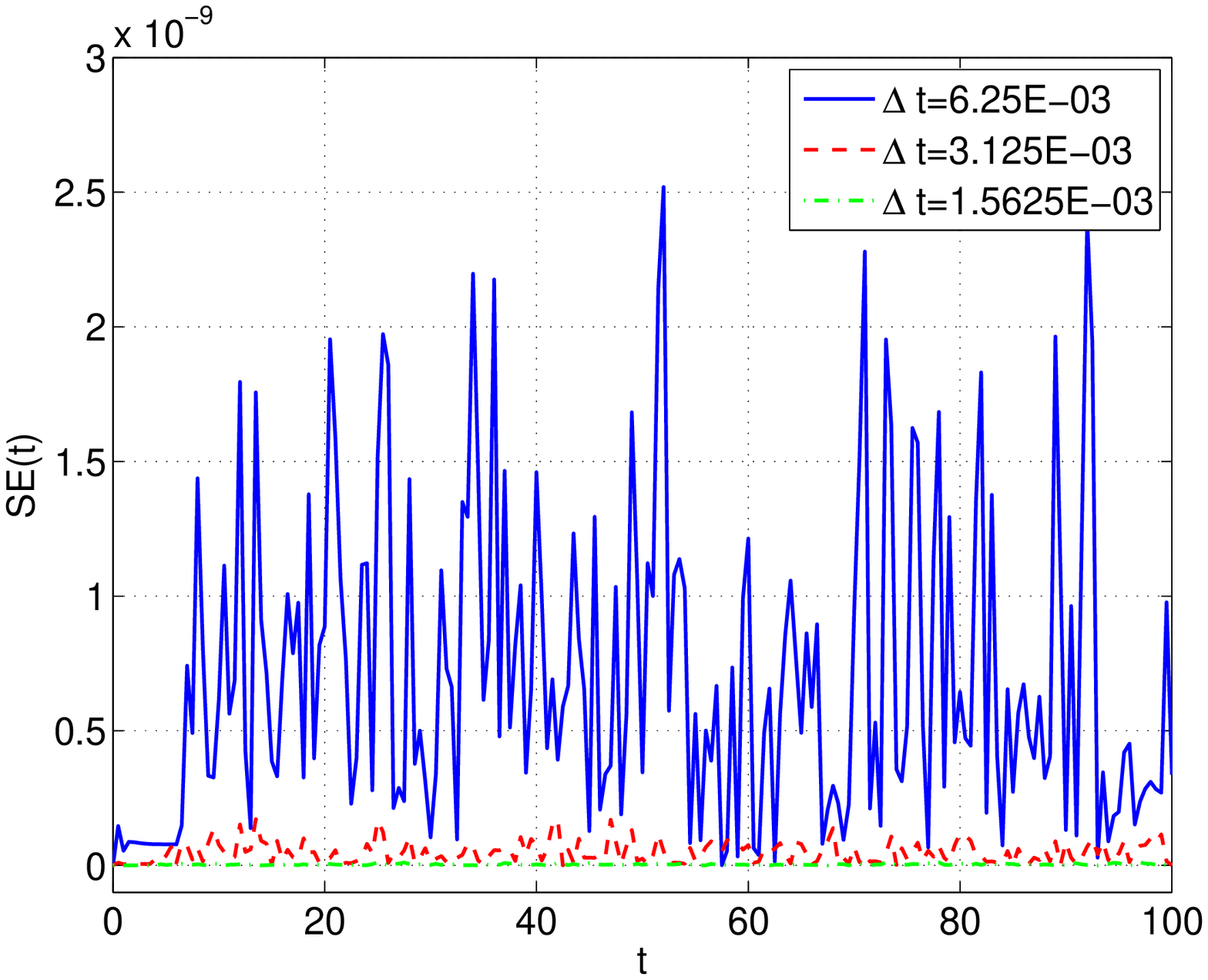}}
\subfigure[]
{\includegraphics[width=6cm]{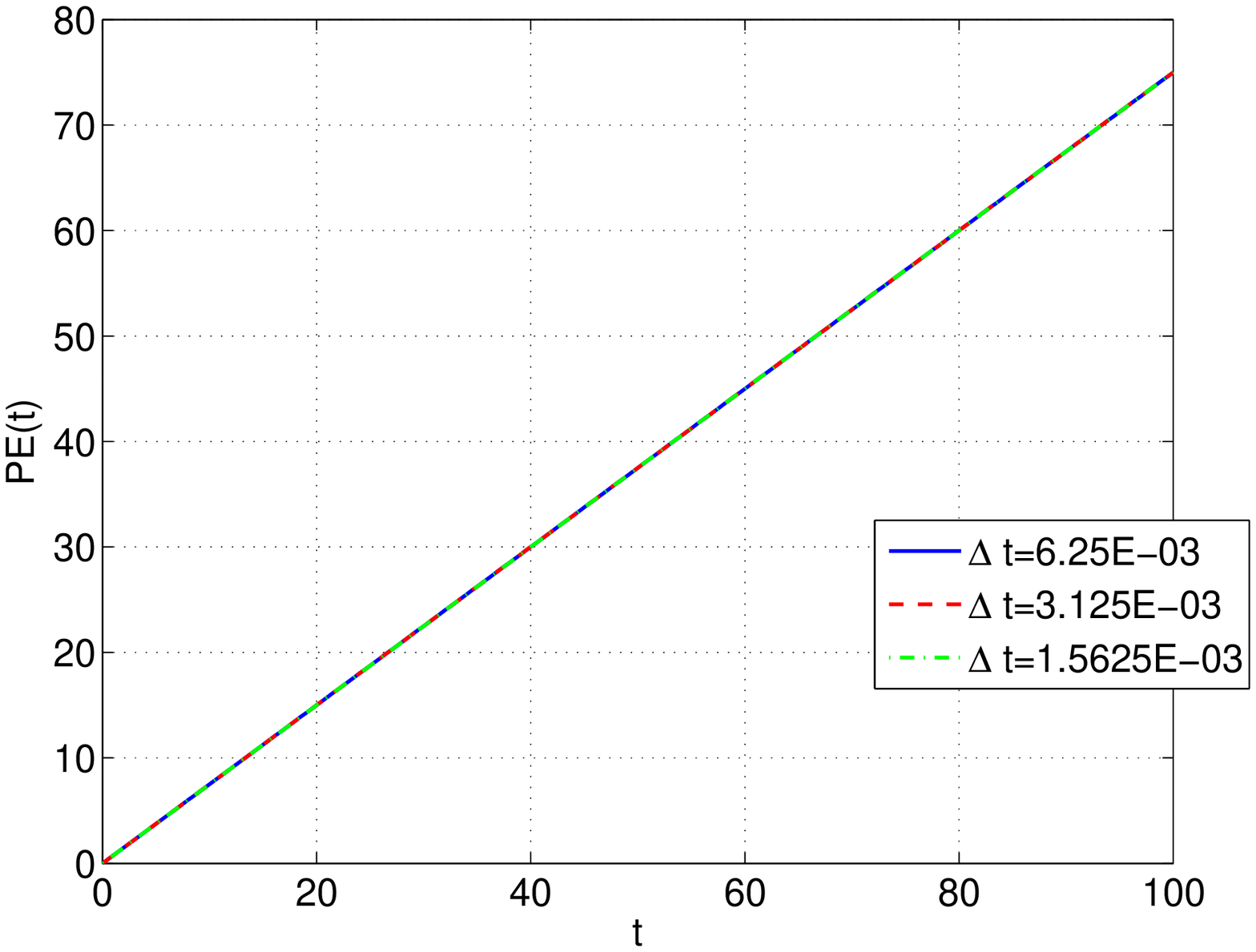}}
\caption{Amplitude (a), speed (b) and phase (c) errors vs. time, for the numerical approximation of gBenjamin equation from the initial solitary-wave profile computed with $\gamma=1.5, \delta=1, c_{s}=0.75, q=2$ and $h=0.125$.}
\label{gbenfig432234}
\end{figure}
\begin{figure}[htbp]
\centering
\subfigure[]
{\includegraphics[width=6cm]{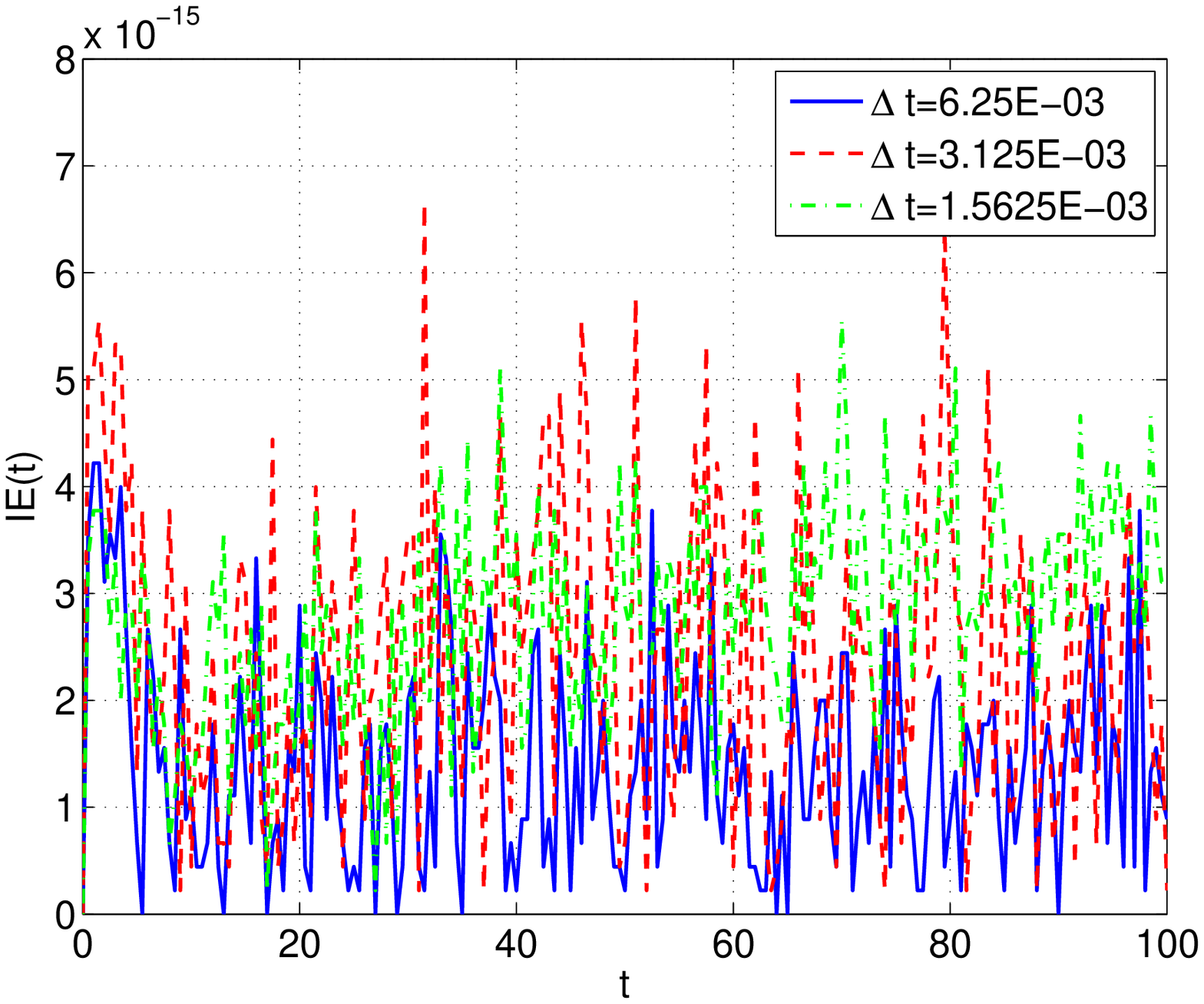}}
\subfigure[]
{\includegraphics[width=6cm]{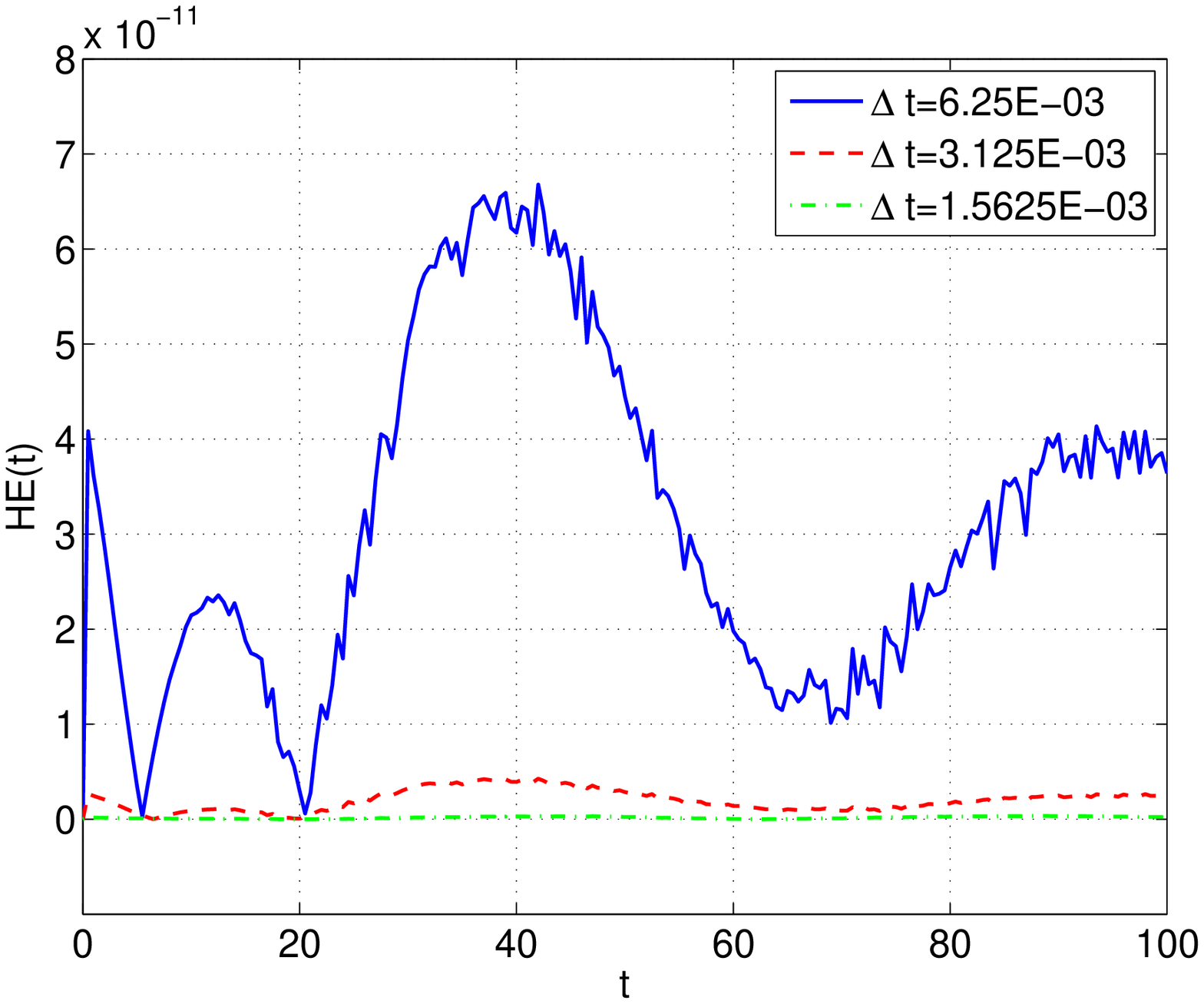}}
\caption{(a) Momentum ($|I_{h}(U^{n})-I_{h}(U^{0})|$) and (b) Energy ($|E_{h}(U^{n})-E_{h}(U^{0})|$) errors vs. time, where $I_{h}, E_{h}$ are given by (\ref{gben46}), (\ref{gben47}) resp. from the generalized Benjamin initial solitary-wave profile computed with $\gamma=1.5, \delta=1, c_{s}=0.75, q=2$ and $h=0.125$.}
\label{gbenfig43256}
\end{figure}

Figures \ref{gbenfig4321}-\ref{gbenfig43256} summarize the results corresponding to the generalized Benjamin equation for the values $q=2, \delta=1, \gamma=1.5$ and $c_{s}=0.75$.  Figure \ref{gbenfig4321} illustrates the evolution of the initial profile at several time stances with $\Delta t=1.5625E-03$. The numerical solution does not appear to develop relevant enough spurious oscillations to perturb the traveling wave evolution. This level of accuracy is checked and confirmed by the rest of the figures. Thus Figures \ref{gbenfig432234}(a)-(c) show, respectively, the evolution of the errors in amplitude, speed and phase (with respect to those of the computed initial profile) for different values of $\Delta t$. The errors in amplitude and speed remain small and bounded in time while the numerical profile is affected by a phase shift that grows linearly with time (cf. as in the gKdV case, \cite{BonaDKM1995}). Finally, Figure \ref{gbenfig43256}(a) confirms the performance of the projection technique for the preservation of (\ref{gben46}) while the good behaviour in the evolution of the error of the energy (\ref{gben47}) is shown in Figure \ref{gbenfig43256}(b). We observed that 
%
this behaviour does not change when other values of the coefficients in (\ref{gben1}) are considered. 

\section{Numerical study of stability and interactions of solitary waves}
\label{sec3}
In this section we discuss several experiments concerning the stability of solitary waves. They are essentially of three types: small and larger perturbations of an approximate solitary wave and interactions of waves. The study focuses on results with $q=2,3$. We have also included some experiments corresponding to $q=4$. These last computations suggest instability and may point to singularity formation of some sort if the initial perturbation is large enough.
\subsection{Structure of dispersive tails}
\label{sec31}
The behaviour of small-amplitude solutions of (\ref{gben1}) will be described first. This will be useful to identify small-amplitude dispersive tails generated by the nonlinear interactions in the experiments to follow. We consider a solitary wave of speed $c_{s}$. In a frame moving with the wave $y=x-c_{s}t$, small-amplitude solutions $u(y,t)$ evolve approximately according to the linear dispersive partial differential equation
\begin{eqnarray*}
u_{t}-c_{s}u_{y}-\mathcal{L}u_{y}=0.
\end{eqnarray*}
For plane wave solutions $u(y,t)=e^{i(ky-\omega(k)t)}, k\in \mathbb{R}$ we have
\begin{eqnarray*}
\omega(k)=-c_{s}k-k(\delta |k|^{2m}-\gamma |k|^{2r}).
\end{eqnarray*}
The local phase speed is
\begin{eqnarray*}
v(k)=\frac{\omega(k)}{k}=-c_{s}-(\delta |k|^{2m}-\gamma |k|^{2r})=-c_{s}+\phi(|k|^{2}),
\end{eqnarray*}
where $\phi(x)=\gamma x^{r}-\delta x^{m}, x>0$. The function $\phi$ attains its  maximum at $x^{*}=(r\gamma/m\delta)^{1/(m-r)}$. Thus if $x>0$
$$\phi(x)\leq \phi(x^{*})=\left(\frac{r\gamma}{m\delta}\right)^{\frac{r}{m-r}}\gamma\left(1-\frac{r}{m}\right).
$$ Then if $k\neq 0$ and since $\gamma<\gamma_{max}(c_{s})$,
\begin{eqnarray*}
v(k)\leq -c_{s}+\phi(x^{*})<-c_{s}+\left(\frac{r\gamma_{max}}{m\delta}\right)^{\frac{r}{m-r}}\gamma_{max}\left(1-\frac{r}{m}\right)=-c_{s}+c_{s}=0.
\end{eqnarray*}
Therefore $v(k)<0$ for all  wavenumbers $k\neq 0$ and the solution component $u(y,t)=e^{i(ky-\omega(k)t)}$ travels behind the solitary pulse to the left in the frame of reference of the pulse. Also, $\phi$ is decreasing for $x>x^{*}$, tending to $-\infty$ as $x\rightarrow\infty$. Thus among the components with $|k|^{2}>x^{*}$ those corresponding to longer wavelength (smaller $k$) are faster than those of shorter wavelength.

In this frame of reference, the associated group velocity is
\begin{eqnarray*}
\omega^{\prime}(k)=-c_{s}+\psi(|k|^{2}),
\end{eqnarray*}
where $\psi(x)=(2r+1)\gamma x^{r}-(2m+1)\delta x^{m}, x>0$. The function $\psi$ attains its maximum at $x^{*}=(r(2r+1)\gamma/m(2m+1)\delta)^{1/(m-r)}$. Thus if $x>0$
$$\psi(x)\leq \psi(x^{*})=\left(\frac{r(2r+1)\gamma}{m(2m+1)\delta}\right)^{\frac{r}{m-r}}(2r+1)\gamma\left(1-\frac{r}{m}\right).
$$ 
The discussion in the general case is a bit more complicated. We consider first the gBenjamin case ($r=1/2, m=1$) for which the function
$$ \mathbb{F}(x)=\psi(x)-c_{s}=2\gamma\sqrt{x}-3\delta x-c_{s}, \quad x>0,$$ determines the group velocity ($\omega^{\prime}(k)=\mathbb{F}(|k|^{2})$). We study the sign of $\mathbb{F}$. Note that $\mathbb{F}(x)=0$ when
$$\mathbb{G}(x)=9\delta^{2}x^{2}+(6\delta c_{s}-\gamma^{2})x+c_{s}^{2}=0.$$ This leads to the following cases (recall that $\gamma<\gamma_{max}(c_{s})=2\sqrt{\delta c_{s}}$):
\begin{enumerate}
\item If $\sqrt{3\delta c_{s}}<\gamma<2\sqrt{\delta c_{s}}$ then $\mathbb{G}$ has two positive roots $0<x_{-}<x_{+}$,
\begin{eqnarray*}
x_{\pm}&=&\frac{1}{2}\left(-\left(\frac{6\delta c_{s}-4\gamma^{2}}{9\delta^{2}}\right)\pm\frac{4\gamma}{9\delta^{2}}\sqrt{\Delta}\right),\\
\Delta&=&\gamma^{2}-3\delta c_{s},
\end{eqnarray*}
with $\mathbb{G}(x)<0$ if $x_{-}<x<x_{+}$ and $\mathbb{G}(x)>0$ if $0<x<x_{-}$ or $x<x_{+}$. This implies:
\begin{itemize}
\item For those $k$ with $|k|^{2}\in (x_{-},x_{+}), \omega^{\prime}(k)>0$. Furthermore, $\psi(x)>0\Leftrightarrow x<x_{c}=(4\gamma^{2})/(9\delta^{2})$. Since $x_{c}>x_{+}$ if $|k|^{2}\in (x_{-},x_{+})$ then $\psi(|k|^{2})>c_{s}, \psi(|k|^{2})>0$ so those components travel to the right and in front of the pulse.
\item When $0<|k|^{2}<x_{-}$ or $x_{+}<|k|^{2}<x_{c}$  then $\psi(|k|^{2})<c_{s}$, i.~e. $\omega^{\prime}(k)<0$. Furthermore, when $0<|k|^{2}<x_{-}$ then $\psi(|k|^{2})>0$ so the components travel to the right and behind the wave.
\item When $|k|^{2}>x_{c}$ then $\psi(|k|^{2})<0, \psi(|k|^{2})<c_{s}$ and the components travel behind the pulse and to the left with respect to the moving frame of reference.
\end{itemize}
\item If $\gamma\leq \sqrt{3\delta c_{s}}$ then $\mathbb{G}(x)>0$ for all $x$ and consequently $\mathbb{F}(x)<0, x>0$, that is $\psi(x)<c_{s}$. This means $\omega^{\prime}(k)<0$ and those components travel behind the solitary pulse to the right if $|k|^{2}<x_{c} (\psi(|k|^{2})>0)$ or to the left if $|k|^{2}>x_{c} (\psi(|k|^{2})<0)$, always with respect to the reference moving with the solitary wave.
\end{enumerate}
Finally, we can compute $\omega^{\prime}(k)-v(k)=\varphi(|k|^{2})$, where
$
\varphi(x)=\psi(x)-\phi(x)=2r\gamma x^{r}-2m\delta x^{m}.
$ In the case of the gBenjamin equation ($r=1/2, m=1$), $\varphi(x)=\gamma \sqrt{x}-2\delta x$. Now $\varphi(x)=0\Leftrightarrow x=x_{p}=(\frac{\gamma}{2\delta})^{2}$. Therefore if $0<|k|^{2}<x_{p}$ the group velocity exceeds the phase speed and if $|k|^{2}>x_{p}$ the phase speed exceeds the group velocity.

\begin{figure}[htbp]
\centering
{\includegraphics[width=0.8\textwidth]{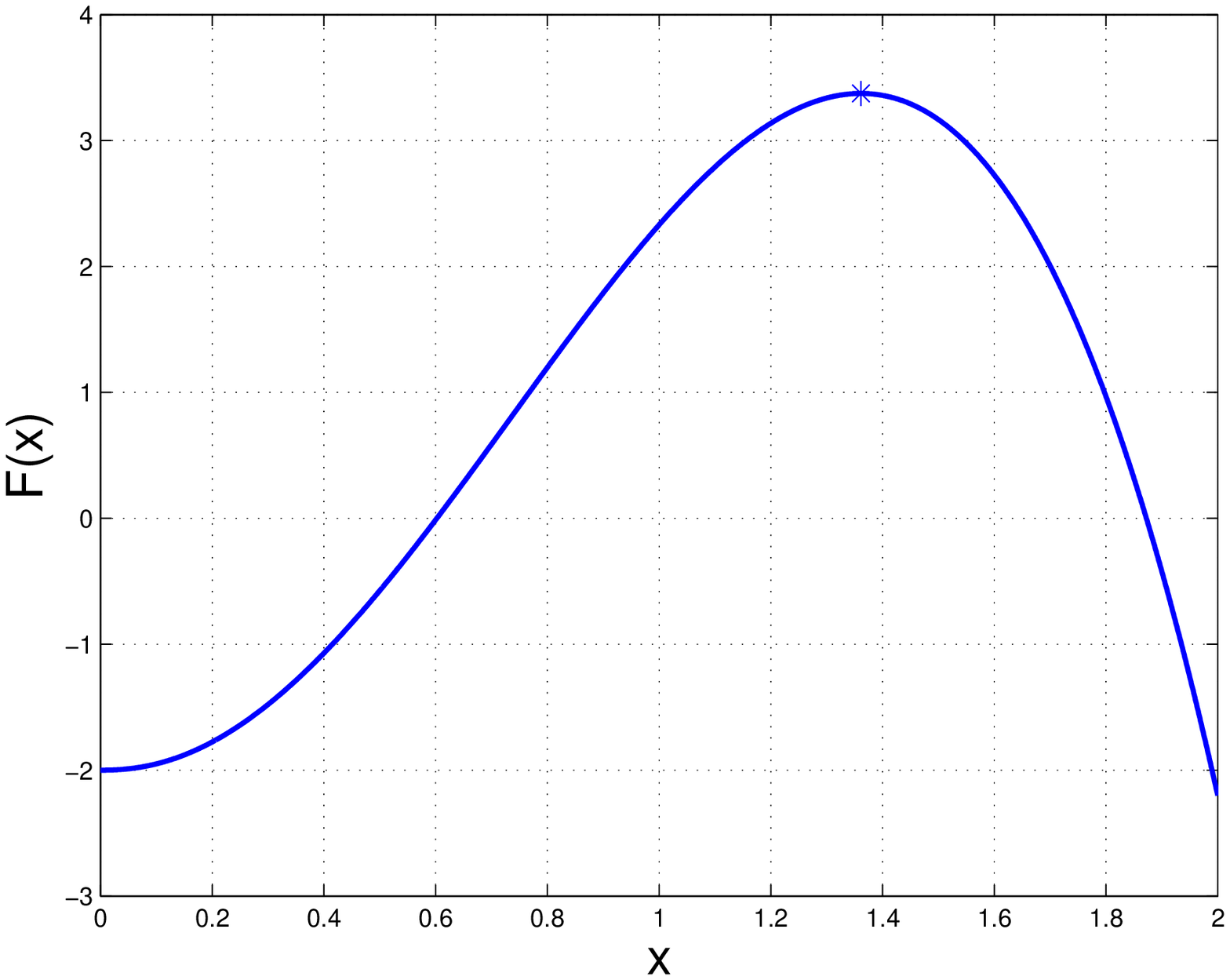}}
\caption{$\mathbb{F}(x)$ for  $r=2.3, m=3, c_{s}=2, \delta=1$ with $\gamma=\gamma_{max}-10^{-4}$. Case  $\gamma\in (\gamma_{*},\gamma_{max})$.}
\label{gbenfig513}
\end{figure}

\begin{figure}[htbp]
\centering
{\includegraphics[width=0.8\textwidth]{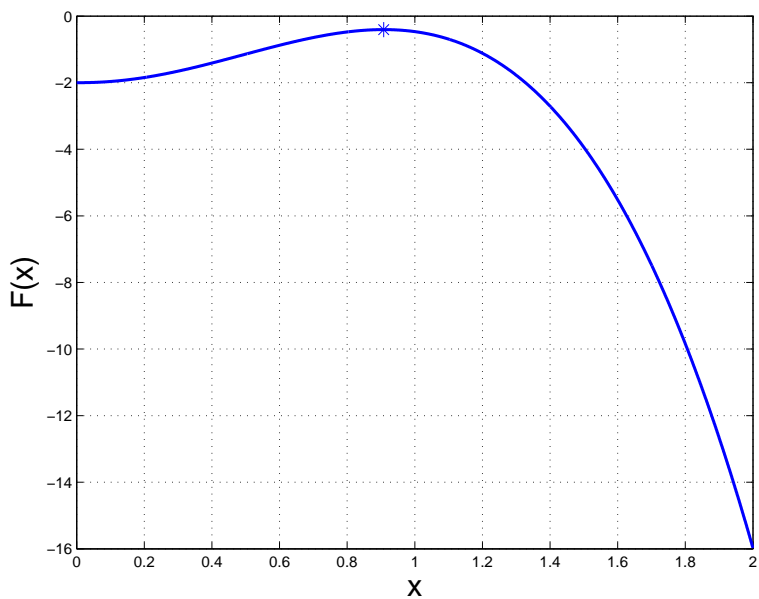}}
\caption{$\mathbb{F}(x)$ for  $r=2.3, m=3, c_{s}=2, \delta=1$ with $\gamma=\gamma_{max}-10^{-1}$. Case  $\gamma<\gamma_{*}$.}
\label{gbenfig514}
\end{figure}

For the general case, when $r=s/t\in \mathbb{Q}, m\in\mathbb{N}$ then $$\mathbb{F}(x)=\psi(x)-c_{s}=0\Leftrightarrow
(2r+1)\gamma x^{s}=(c_{s}+(2m+1)\delta x^{m})^{t},$$ and the discussion is more complicated, but some numerical experiments suggest a similar behaviour:Two ranges of $\gamma$, determined by some $\gamma_{*}=\gamma_{*}(c_{s},r,m,\delta)$ and $\gamma_{max}(c_{s},r,m,\delta)$, for which if $\gamma\in (\gamma_{*},\gamma_{max})$ we are in the case (1), and if $0<\gamma<\gamma_{*}$ in case (2). This is illustrated in Figures \ref{gbenfig513} and \ref{gbenfig514}. They display the form of $\mathbb{F}(x)$ for the case $r=2.3, m=3, c_{s}=2, \delta=1$ with $\gamma=\gamma_{max}-10^{-4}$ and $\gamma=\gamma_{max}-10^{-1}$ respectively. In the first one $\mathbb{F}$ has two zeros (case (1)) while in the second one $\mathbb{F}(x)<0$ for $x>0$ (case (2)). Since $\psi$ attains its maximum at $x^{*}$, the maximum value of $\mathbb{F}(x)$ is $\mathbb{F}(x^{*})=\psi(x^{*})-c_{s}$ which vanishes when $\psi(x^{*})=c_{s}$, that is
$$
\gamma^{\frac{m}{m-r}}\left(\frac{r(2r+1)}{m(2m+1)\delta}\right)^{\frac{r}{m-r}}(2r+1)\left(1-\frac{r}{m}\right)=c_{s}
$$ This leads to
\begin{eqnarray}
\gamma_{*}=\left(\frac{c_{s}}{(2r+1)\left(1-\frac{r}{m}\right)\left(\frac{r(2r+1)}{m(2m+1)\delta}\right)^{\frac{r}{m-r}}}\right)^{\frac{m-r}{m}}.\label{gben52}
\end{eqnarray}
In the case of  Figures \ref{gbenfig513} and \ref{gbenfig514}, $\gamma_{*}\approx 1.142677$ and $\gamma_{*}\approx 1.184416$ respectively; when $r=1/2, m=1$ $\gamma_{*}=\sqrt{3\delta c_{s}}$ as shown before.

%

\subsection{Small perturbations of solitary waves}
\label{sec32}

Here the numerical solitary wave profiles $\phi_{s}$ with speed $c_{s}$ are slightly perturbed in amplitude, the perturbed profiles $A\phi_{s}$, $A$ small, are taken as initial conditions $u(x,0)=A\phi_{s}(x)$ of the code, and the corresponding evolution of the numerical approximation is monitored. According to the number of parameters of the problem (i.~e. $c_{s}, \gamma$ and $q$) we have considered experiments mainly for the gBenjamin equation and illustrated the general case with some examples as well. In all cases $\delta=1, \Delta t=1.5625E-03, h=0.125$. 


We first consider the gBenjamin equation ($r=1/2, m=1$). The perturbation factor $A=1.1$ is fixed. The first group of experiments corresponds to the values $q=2, c_{s}=0.75$ (for which $\gamma_{max}(c_{s})\approx 1.732051$) and $\gamma=1.5, 1.7$. This experiment aims at exploring the influence of taking $\gamma$ closer to $\gamma_{max}$.

\begin{figure}[htbp]
\centering
\subfigure[]
{\includegraphics[width=6cm]{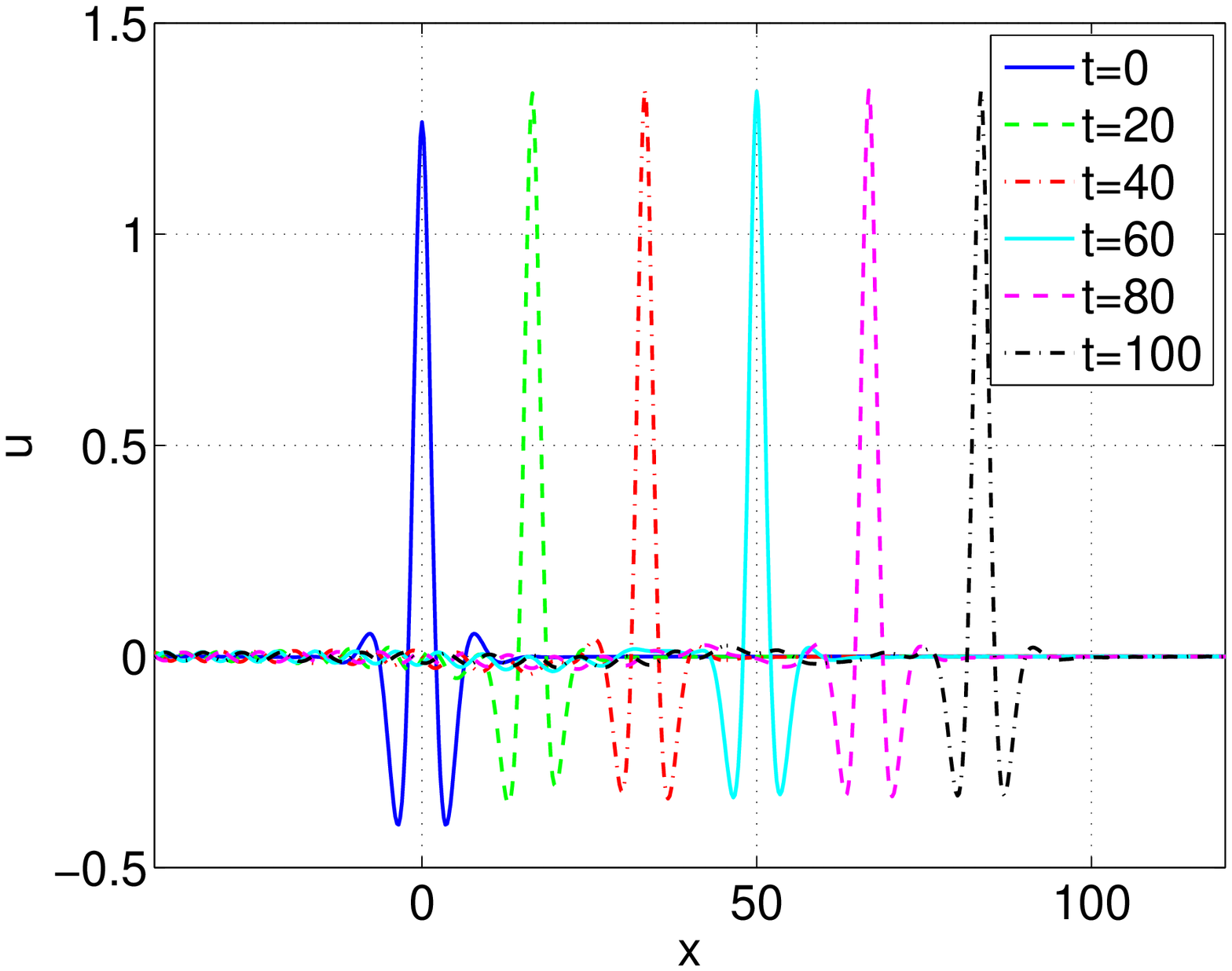}}
\subfigure[]
{\includegraphics[width=6cm]{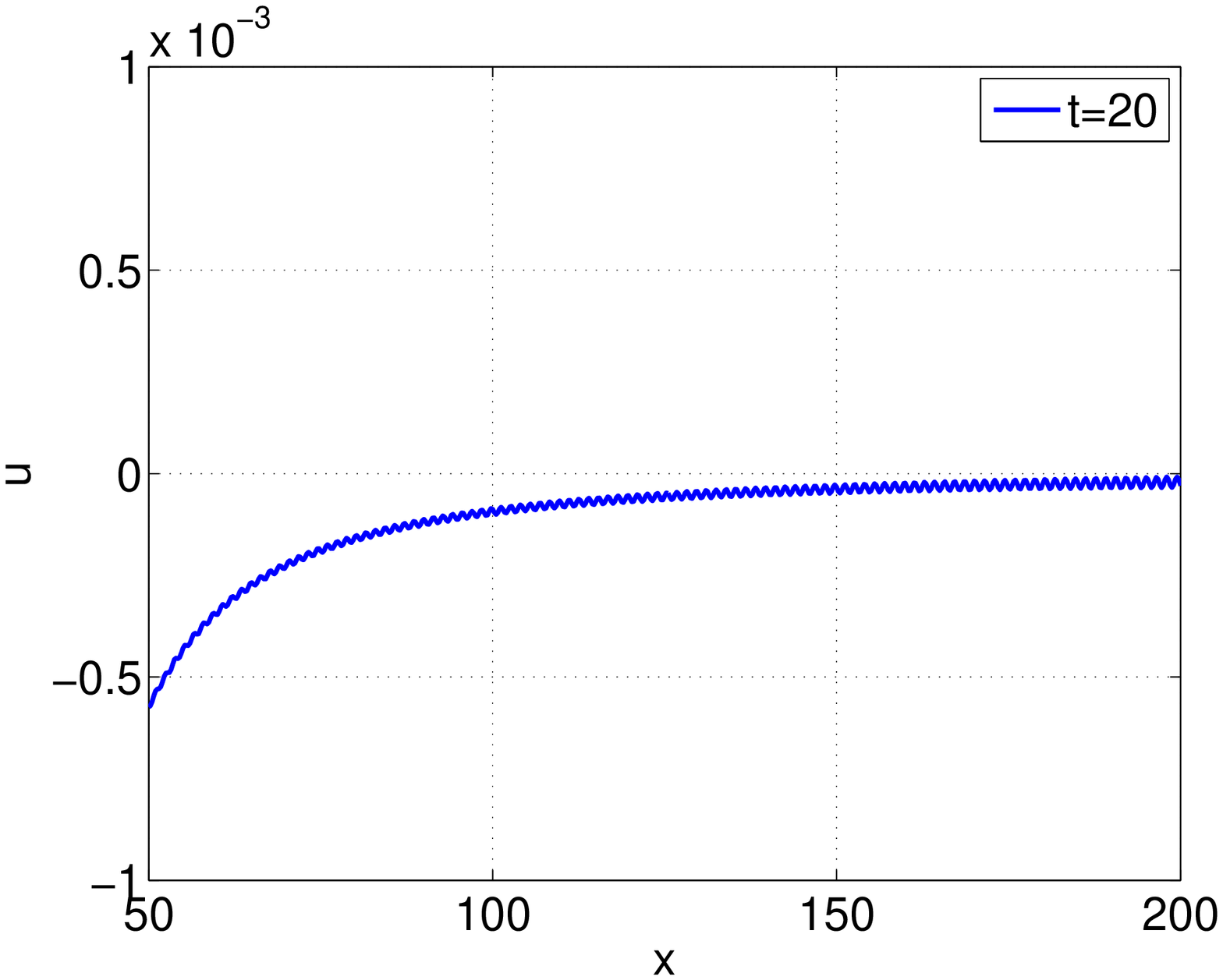}}
\caption{$q=2, c_{s}=0.75$,  $\gamma=1.5$, $A=1.1$. (a) Evolution of the numerical approximation; (b) Magnification when $t=20$.}
\label{gbenfig511_1}
\end{figure}

\begin{figure}[htbp]
\centering
\subfigure[]
{\includegraphics[width=6cm]{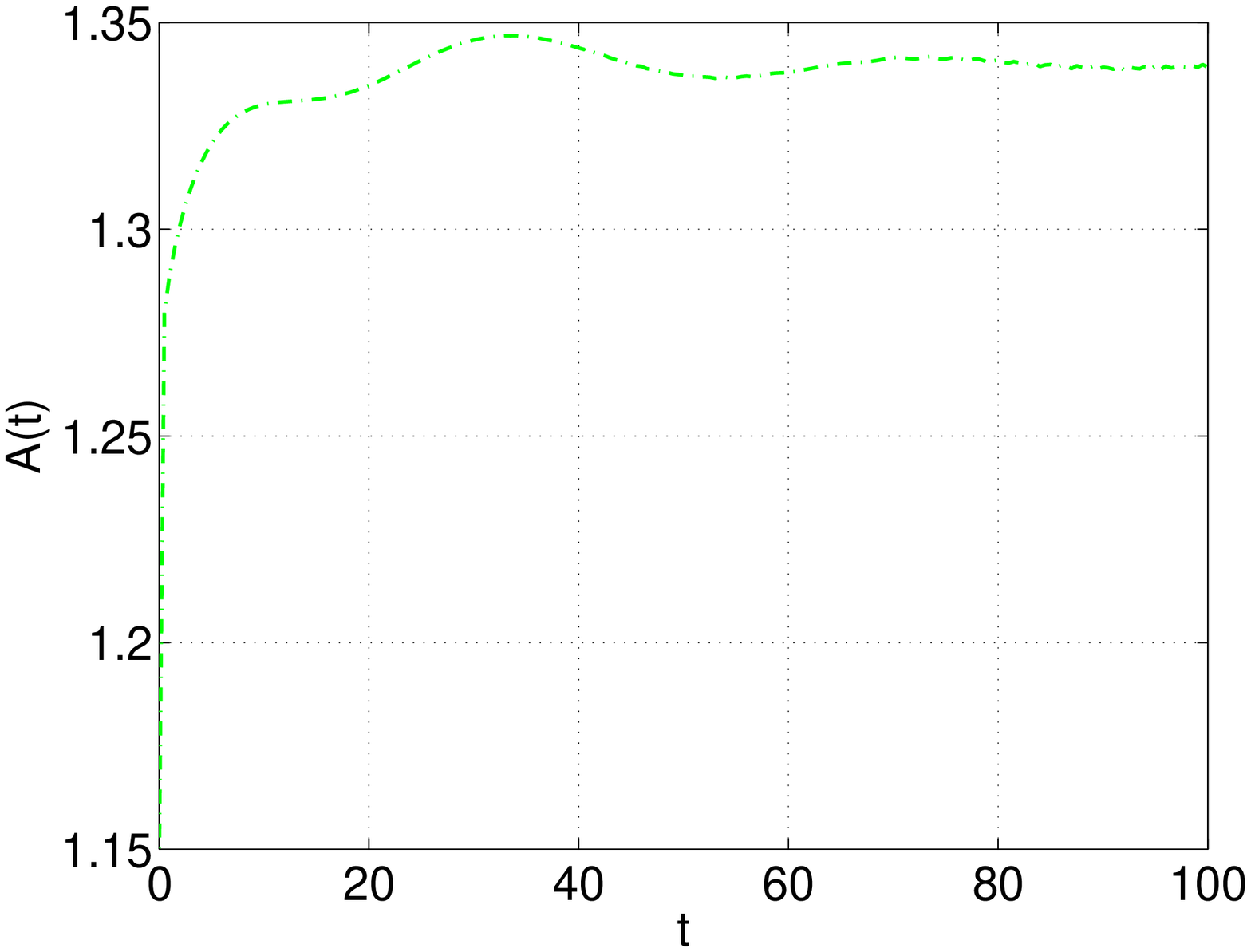}}
\subfigure[]
{\includegraphics[width=6cm]{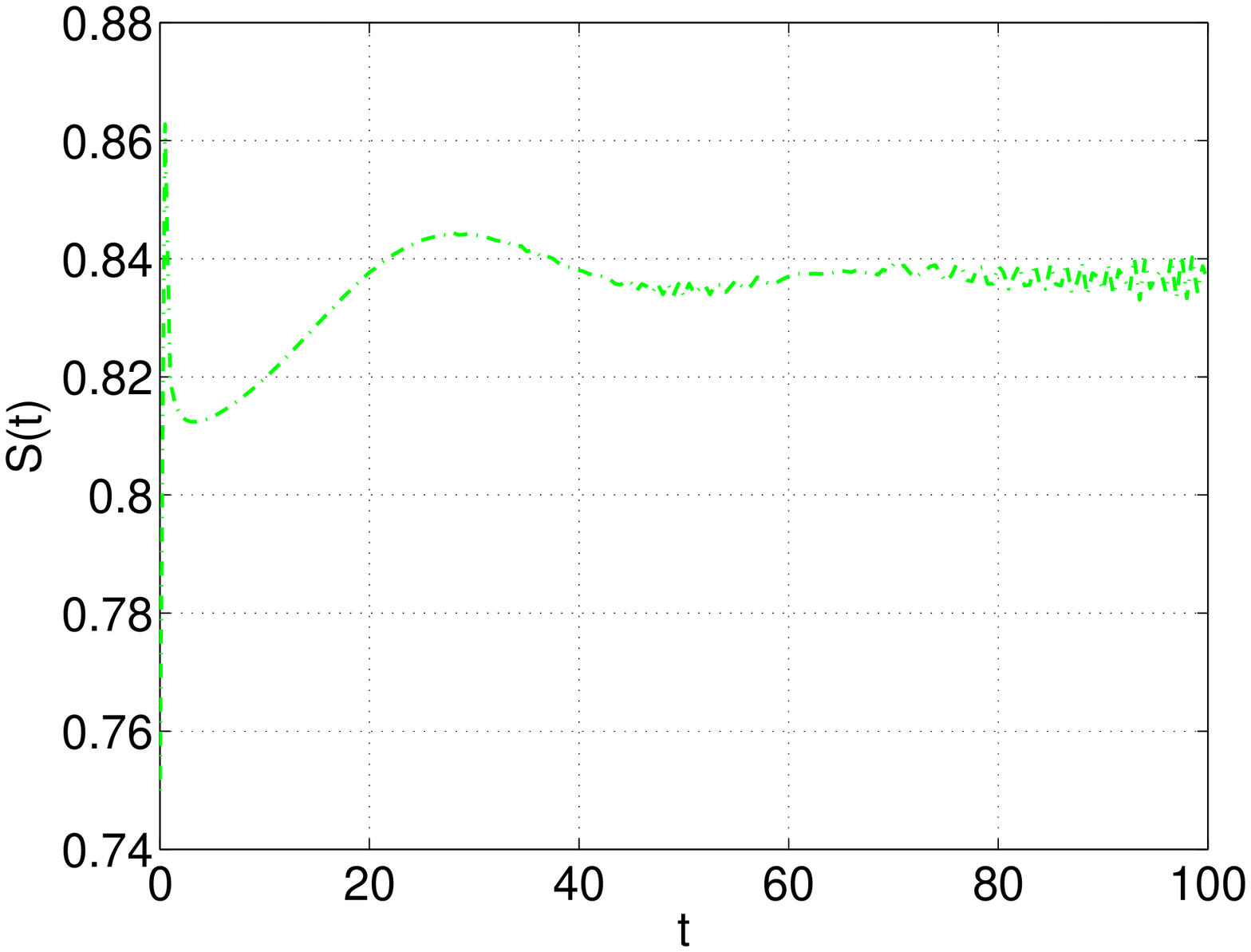}}
\caption{$q=2, c_{s}=0.75$,  $\gamma=1.5$, $A=1.1$. (a) Evolution of the amplitude $A(t)$ of the main numerical pulse; (b) Evolution of speed $S(t)$.}
\label{gbenfig511_2}
\end{figure}

The case $\gamma=1.5$ is illustrated by Figures \ref{gbenfig511_1}-\ref{gbenfig511_2}. Figure \ref{gbenfig511_1}(a) shows the numerical solution at several times; it consists of a solitary pulse with slightly larger amplitude, a long dispersive tail behind this pulse and a very small tail traveling in front of it. (See Figure \ref{gbenfig511_1}(b); this was predicted in Section \ref{sec51}. Here $\gamma=\sqrt{3\delta c_{s}}$, so the oscillations are very small.) The initial profile (without perturbation) has an approximate amplitude of $u_{m}\approx 1.150295$. Thus the perturbation increases the amplitude so about $Au_{m}\approx 1.265325$. The evolution of the amplitude of the main pulse is observed in Figure \ref{gbenfig511_2}(a); compared to the initial amplitude (with and without perturbation) this is larger. (Fitting the values from $t=20$ to a constant gives an approximate amplitude of $1.340417$.) Furthermore, the results shown in Figure \ref{gbenfig511_2}(b) suggest an approximate speed of the main pulse of $c_{m}\approx 0.836720$, if the data are fitted to a constant for $t\geq 40$.)

\begin{figure}[htbp]
\centering
\subfigure[]
{\includegraphics[width=6cm]{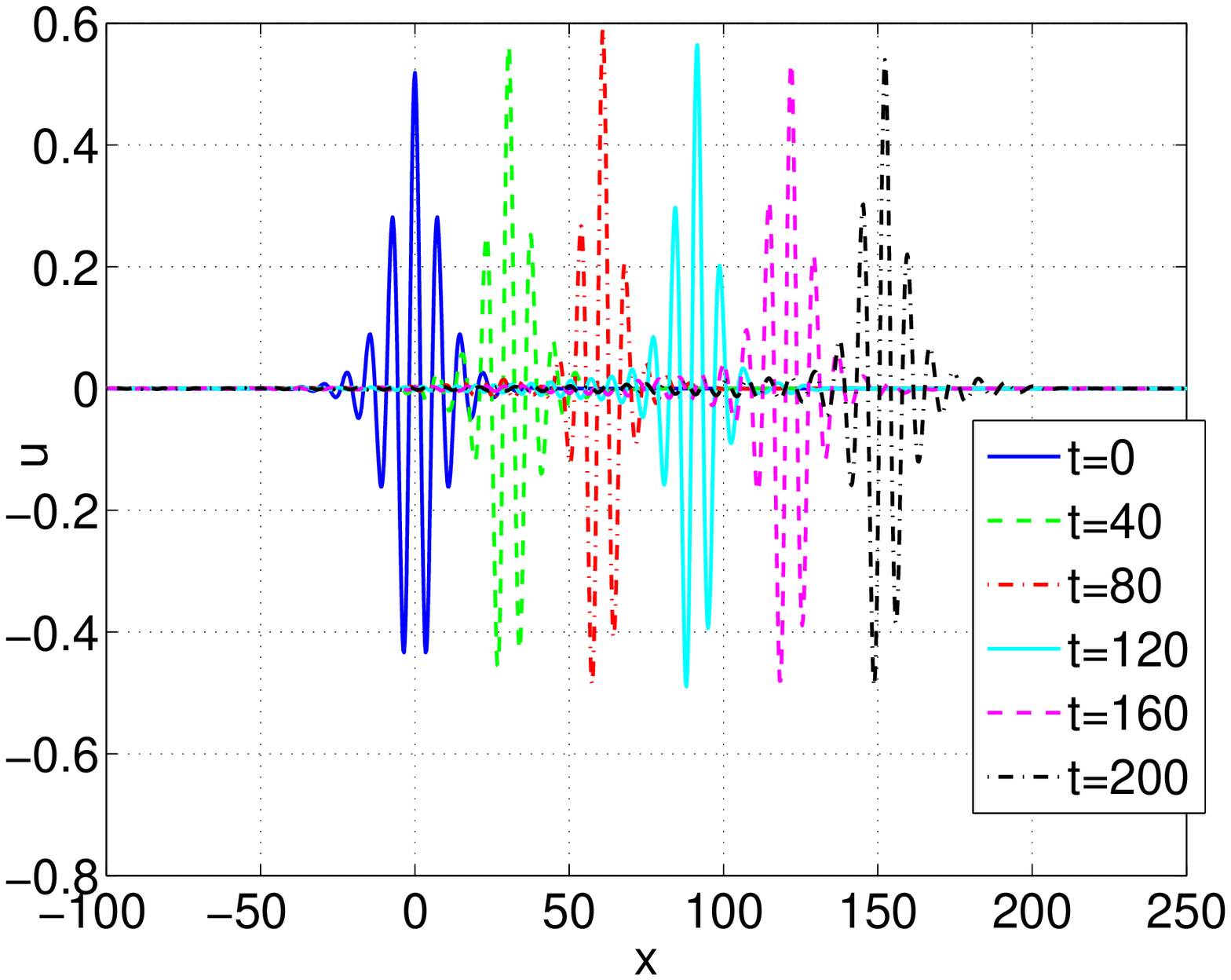}}
\subfigure[]
{\includegraphics[width=6cm]{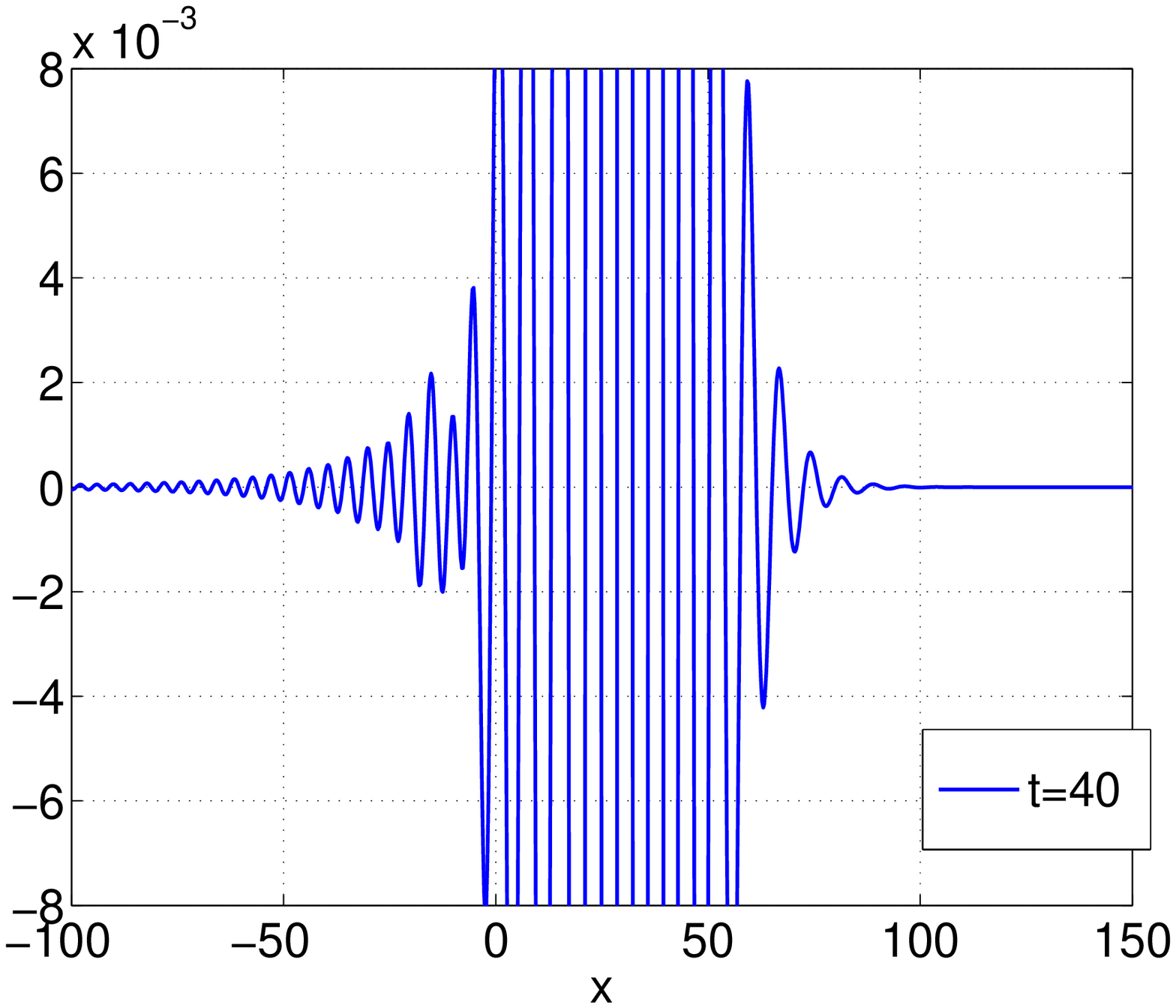}}
\subfigure[]
{\includegraphics[width=6cm]{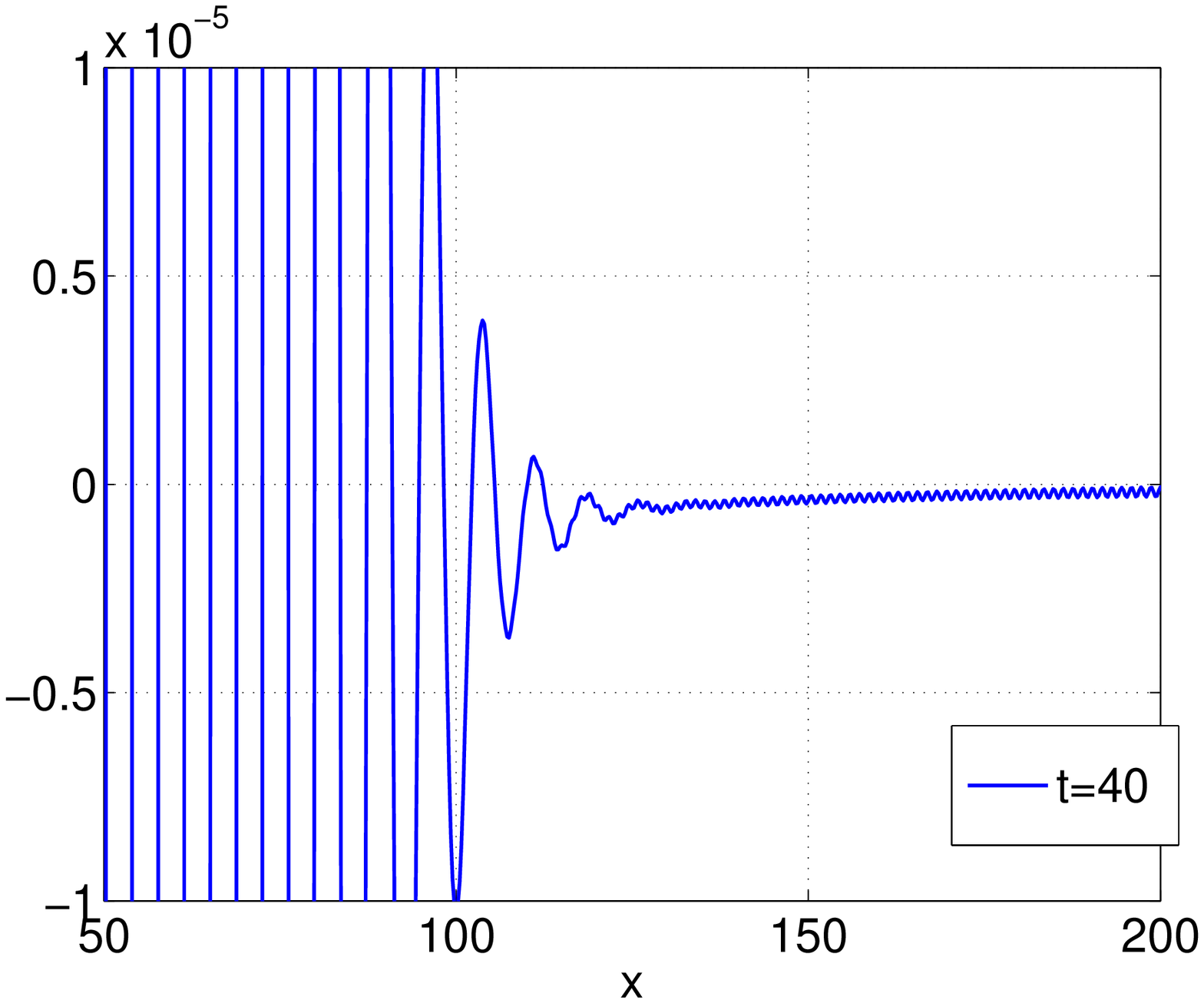}}
\caption{$q=2, c_{s}=0.75$,  $\gamma=1.7$, $A=1.1$. (a) Evolution of the numerical approximation; (b), (c)  Magnifications when $t=40$.}
\label{gbenfig511_4}
\end{figure}

\begin{figure}[htbp]
\centering
\subfigure[]
{\includegraphics[width=6cm]{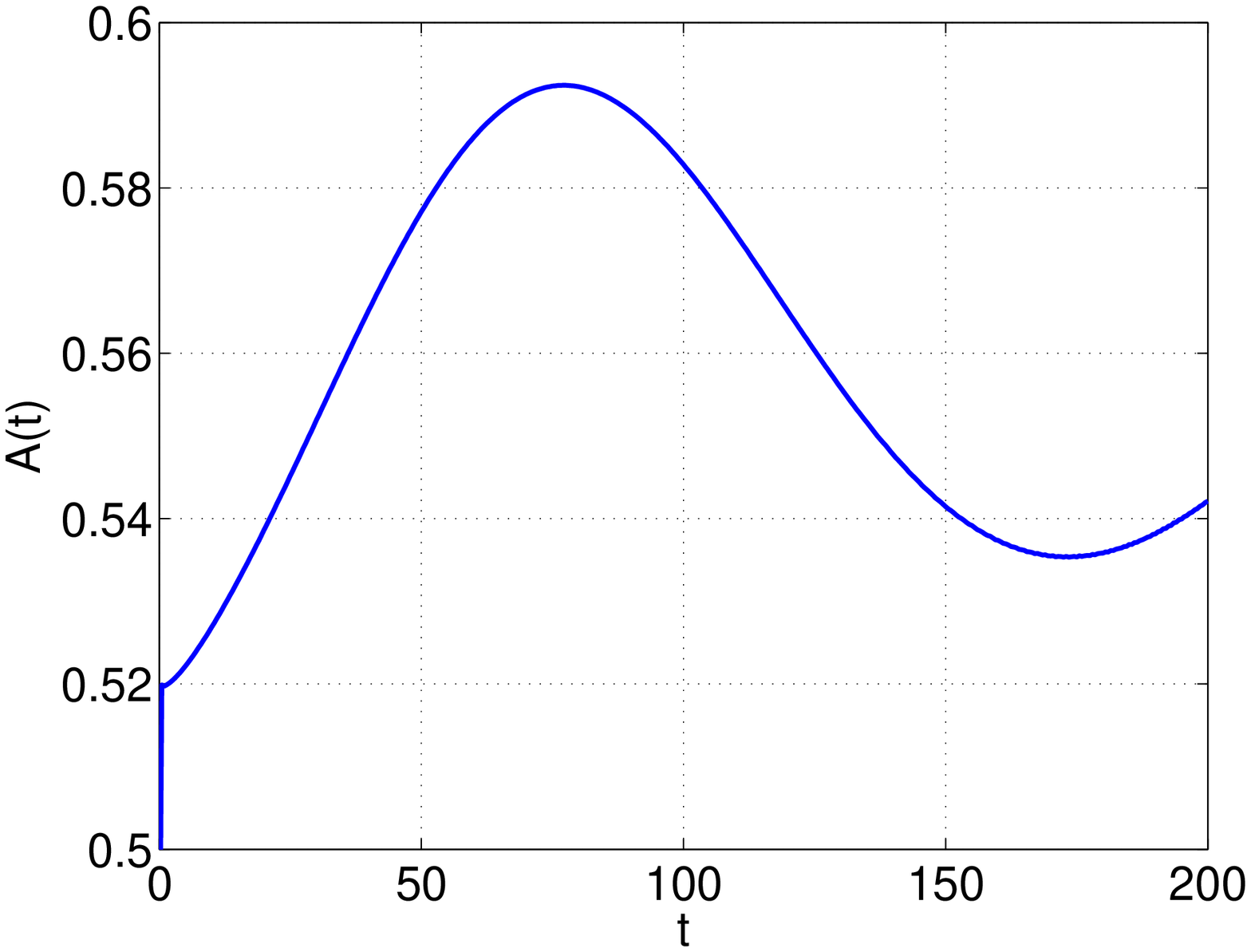}}
\subfigure[]
{\includegraphics[width=6cm]{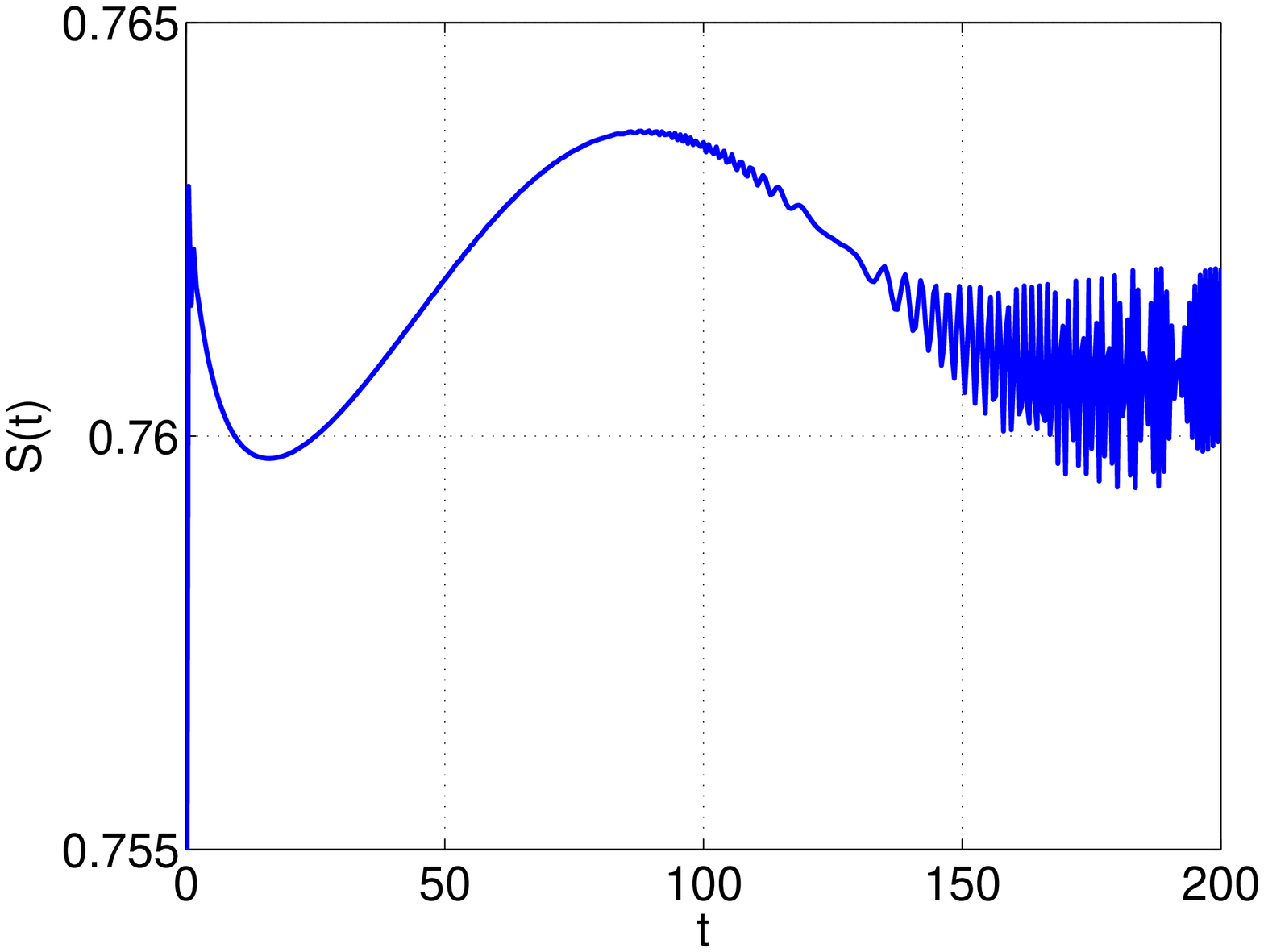}}
\caption{$q=2, c_{s}=0.75$,  $\gamma=1.7$, $A=1.1$. (a) Evolution of the amplitude of the main numerical pulse; (b) Evolution of speed.}
\label{gbenfig511_5}
\end{figure}

When $\gamma$ is closer to $\gamma_{max}$ some slight differences are observed. The main one is that the stabilization of the pulse requires a longer computation. This is observed in Figure \ref{gbenfig511_5} depicting the amplitude and speed of the emerging pulse; the pulse is higher and faster than its initial value as is also observed in Figure \ref{gbenfig511_4}(a). Thus, if the profile were exact, it would correspond to a smaller value of $\gamma$. 
The formation of tails behind and in front of this main pulse is observed in Figures \ref{gbenfig511_4}(b), (c). In this case the initial pulse has an amplitude of $u_{m}\approx 4.716103E-01$ while the emerging wave increases up to $u_{m}\approx 5.766060E-01$.

\begin{figure}[htbp]
\centering
\subfigure[]
{\includegraphics[width=6cm]{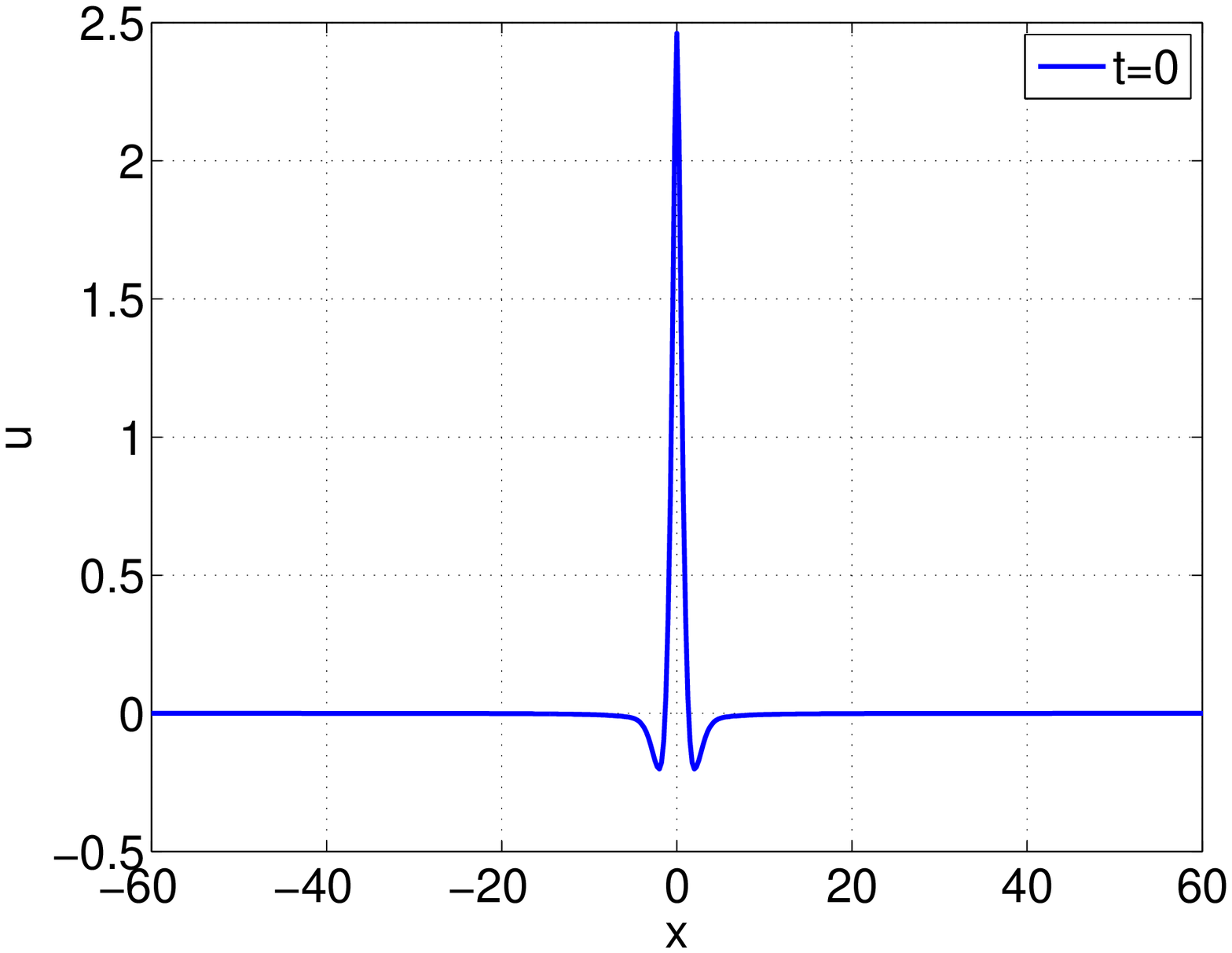}}
\subfigure[]
{\includegraphics[width=6cm]{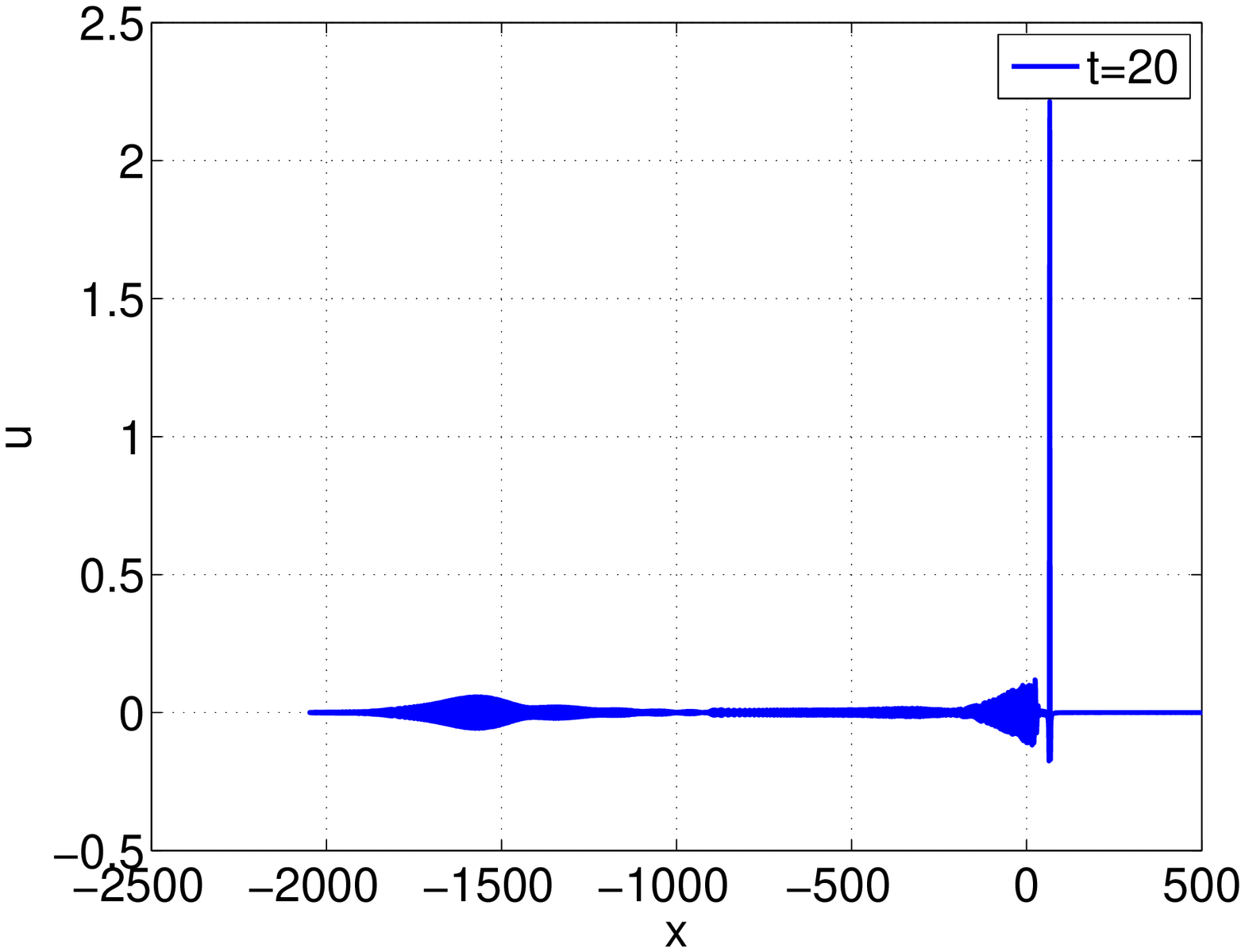}}
\subfigure[]
{\includegraphics[width=6cm]{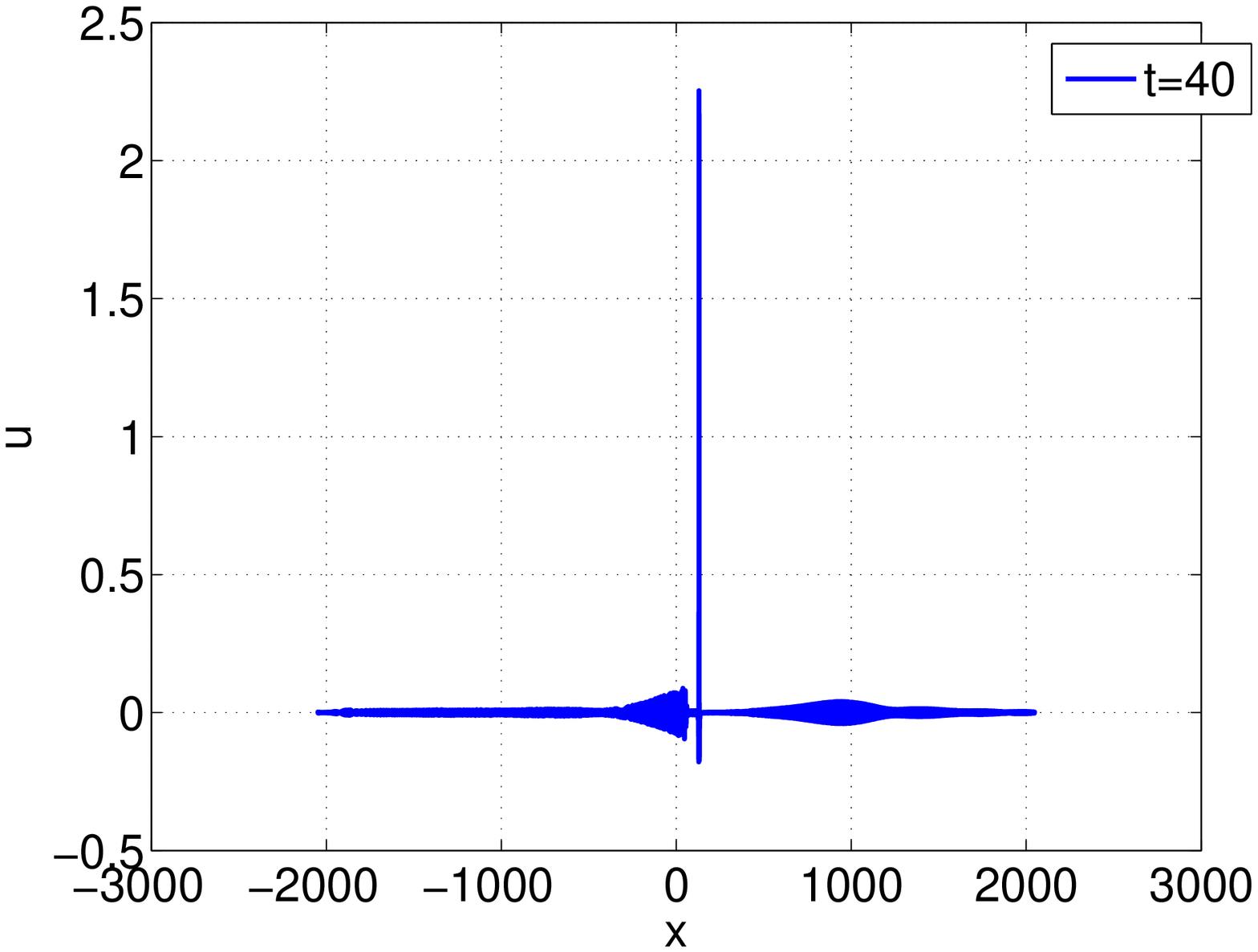}}
\subfigure[]
{\includegraphics[width=6cm]{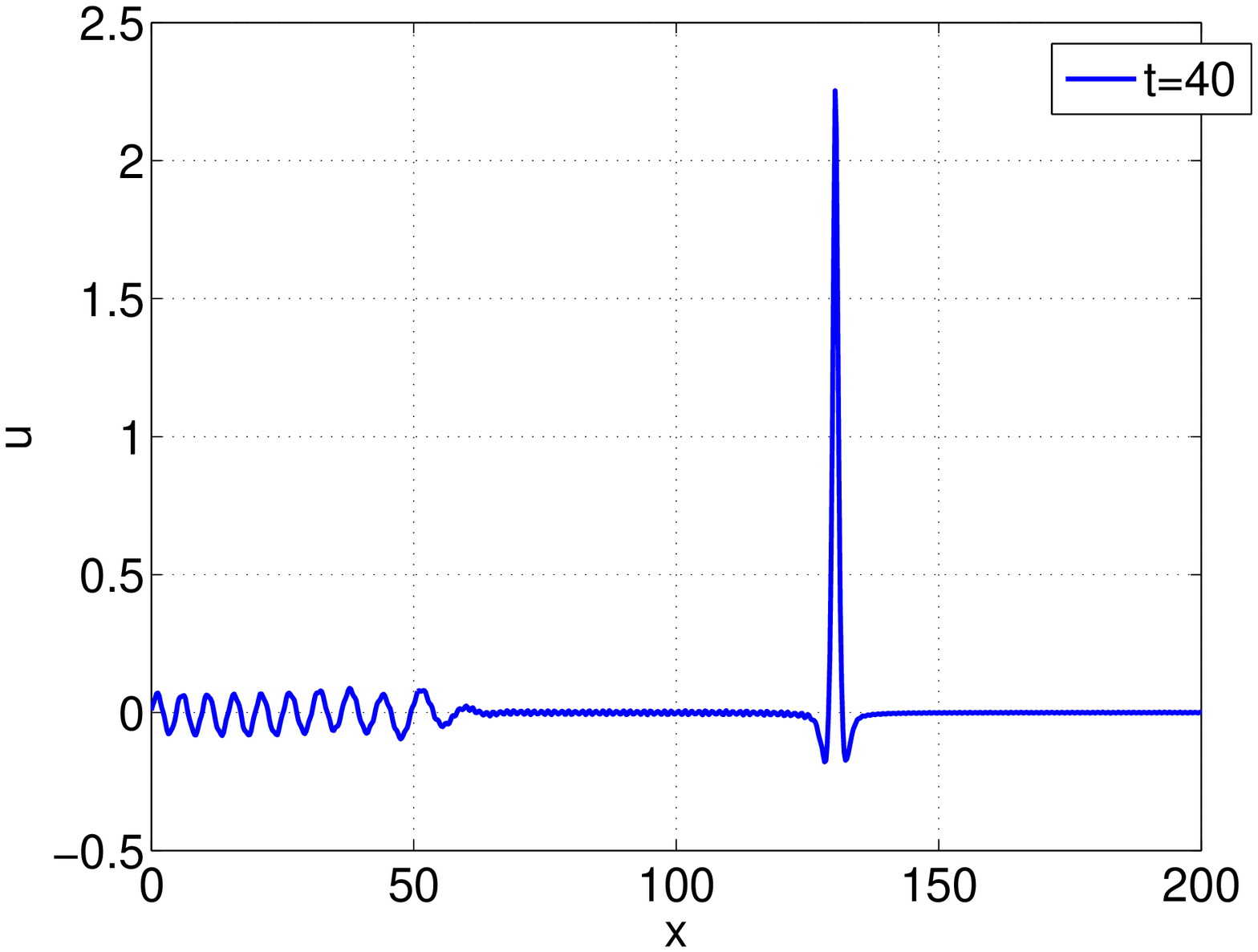}}
\caption{$q=4, c_{s}=3$,  $\gamma=2$, $A=1.1$. (a)-(c) Numerical approximation at times $t=0, 20, 40$. (d) Magnification of (c).}
\label{gbenfig516_1}
\end{figure}

For such a small perturbation $A=1.1$, changing the nonlinearity to $q=4$ does not seem to affect the formation of a perturbed solitary pulse, at least for moderately large times.  If we suspect a singularity formation for $q\geq 4$, then this size of the perturbation does seem however to be large enough to generate it.
Increasing the speed results to generating dispersive bullets. This is observed in Figure \ref{gbenfig516_1}, which corresponds to $q=4, c_{s}=3$ (then $\gamma_{max}(c_{s})\approx 3.464102$) and $\gamma=2$. In this case, the initial amplitude is $u_{m}\approx 2.237447$ and the emerging wave has $u_{m}\approx 2.295051$ with $c_{s}\approx 3.254965$. This figure points to an instability of some sort for $q=4$ since the size of the ripples does not diminish as time evolves and then they might not be of dispersive nature. (Recall the orbital instability result established in \cite{Angulo2003}; see also \cite{Angulo2009}.)

\begin{figure}[htbp]
\centering
\subfigure[]
{\includegraphics[width=6cm]{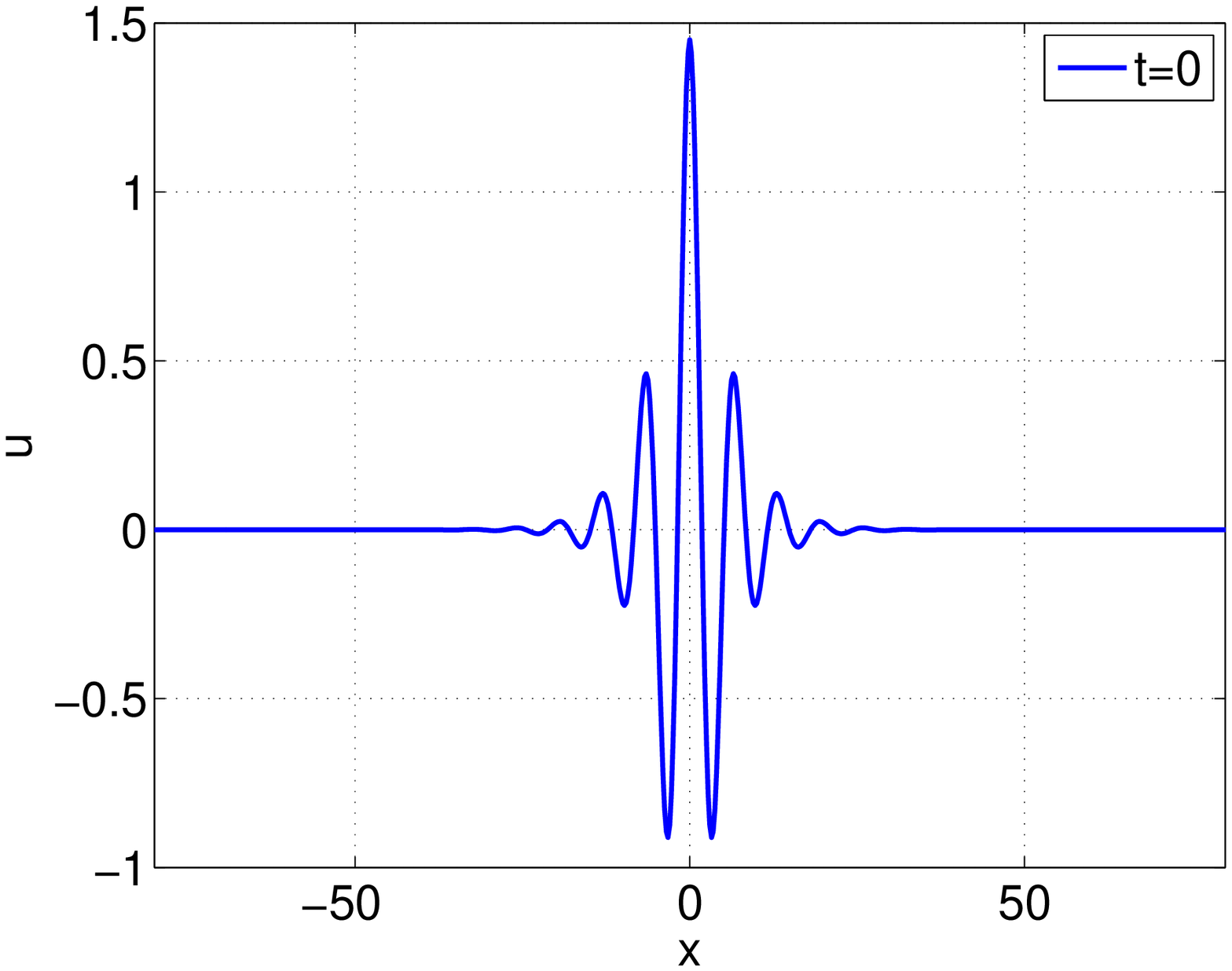}}
\subfigure[]
{\includegraphics[width=6cm]{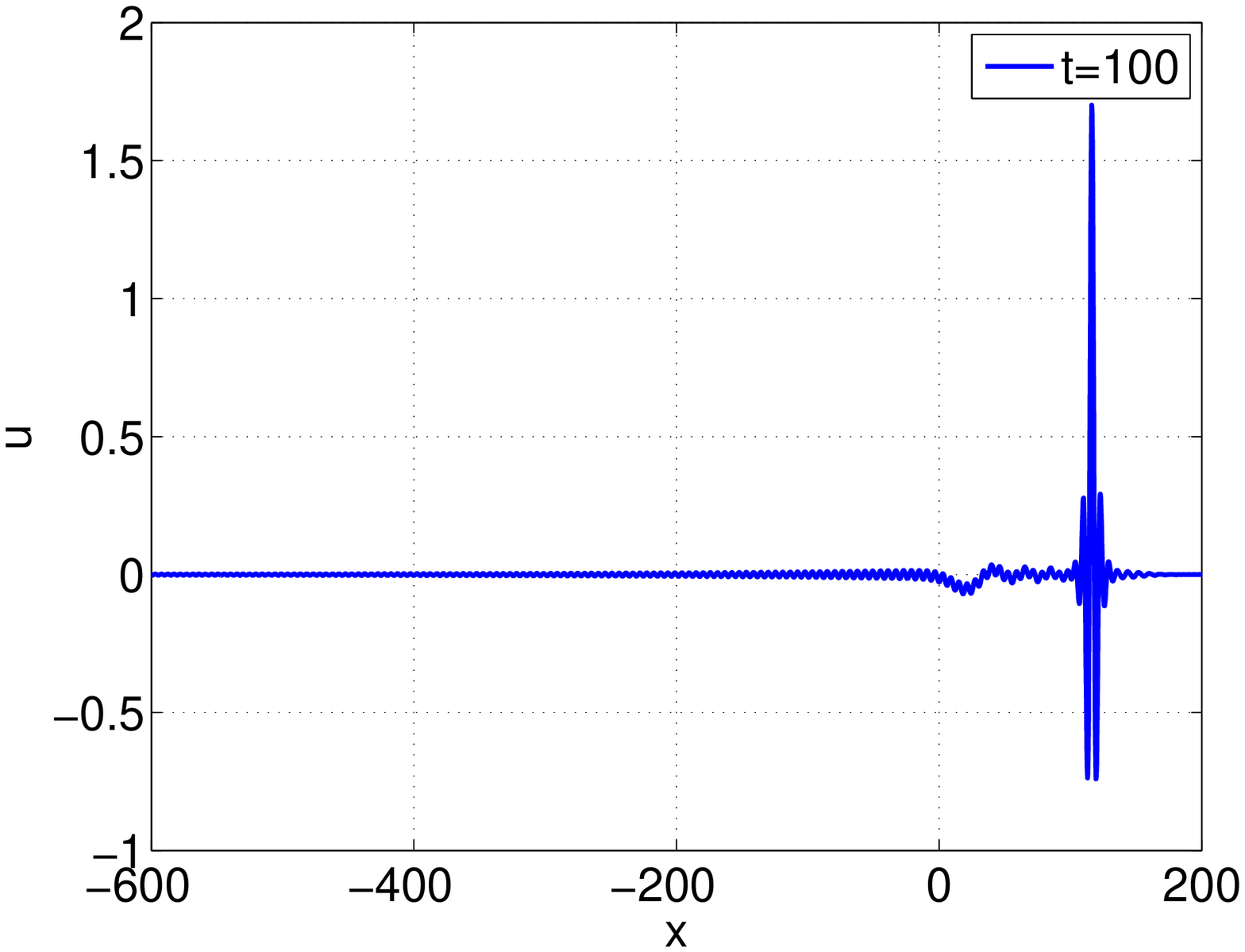}}
\subfigure[]
{\includegraphics[width=6cm]{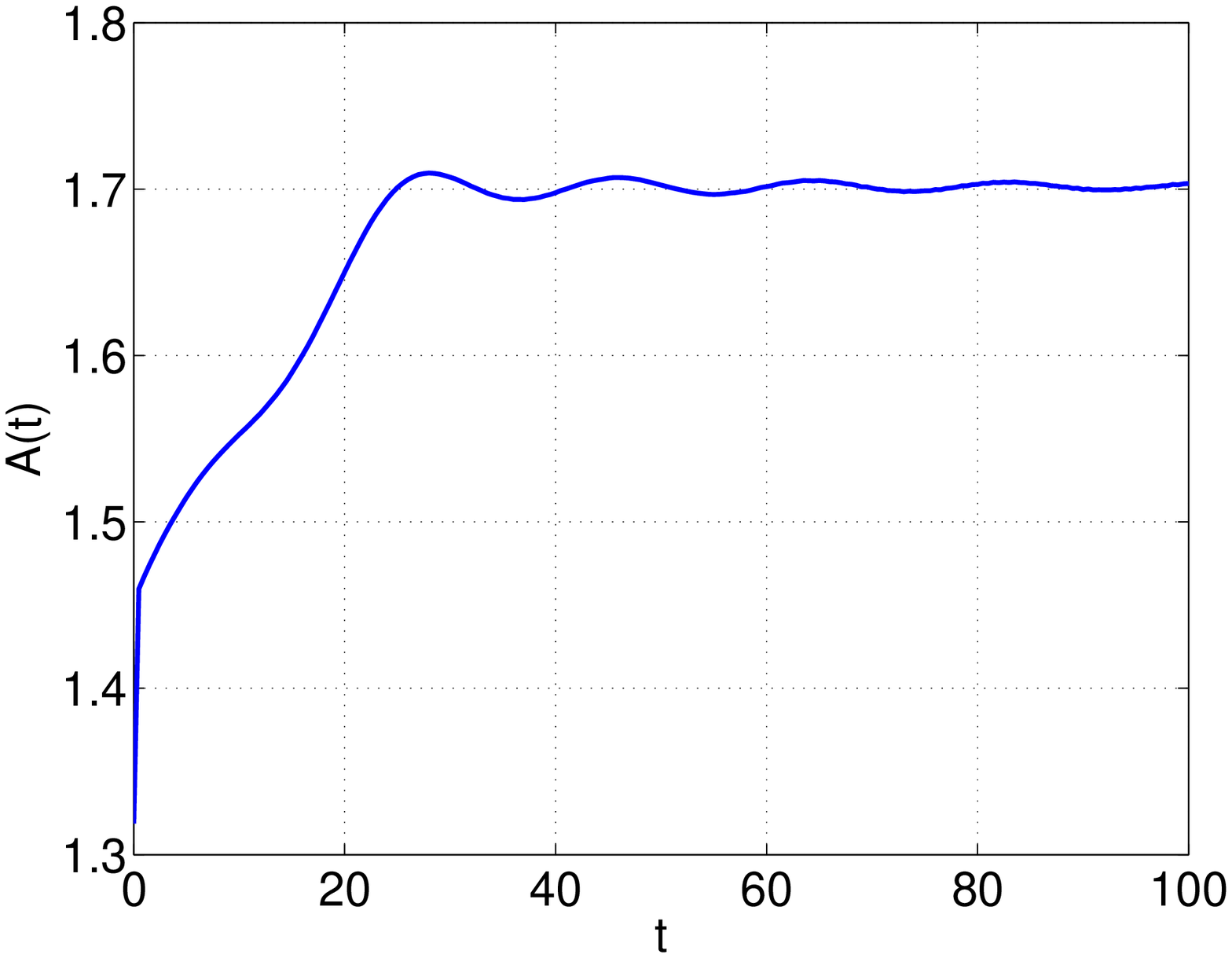}}
\subfigure[]
{\includegraphics[width=6cm]{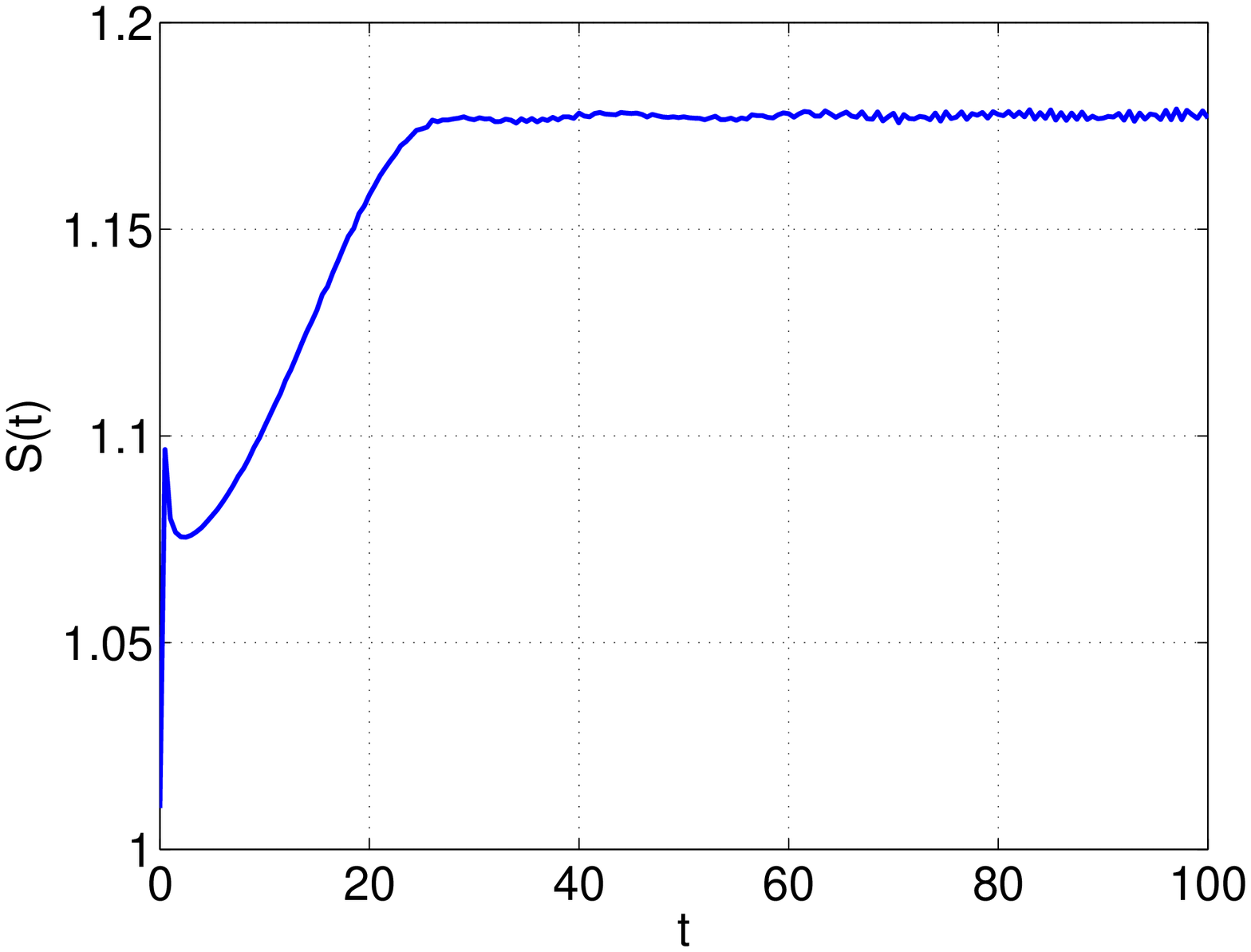}}
\caption{$r=1, m=2, q=2, c_{s}=1.01$,  $\gamma=1.8$, $A=1.1$. (a)  Numerical approximation at $t=0$; (b) numerical approximation at $t=100$; (c) Evolution of the amplitude of the main numerical pulse; (d) Evolution of speed.}
\label{gbenfig513_1}
\end{figure}
\begin{figure}[htbp]
\centering
\subfigure[]
{\includegraphics[width=6cm]{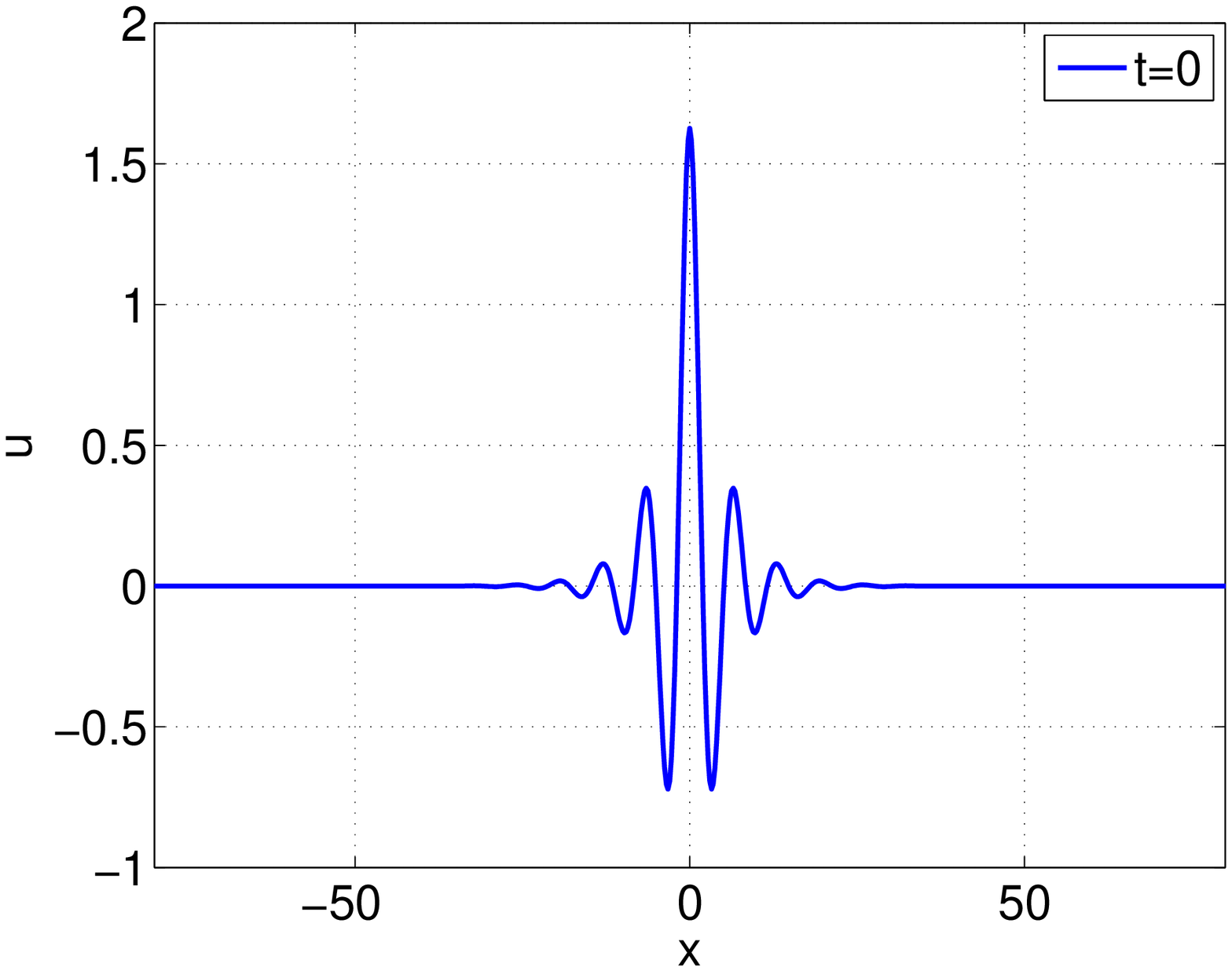}}
\subfigure[]
{\includegraphics[width=6cm]{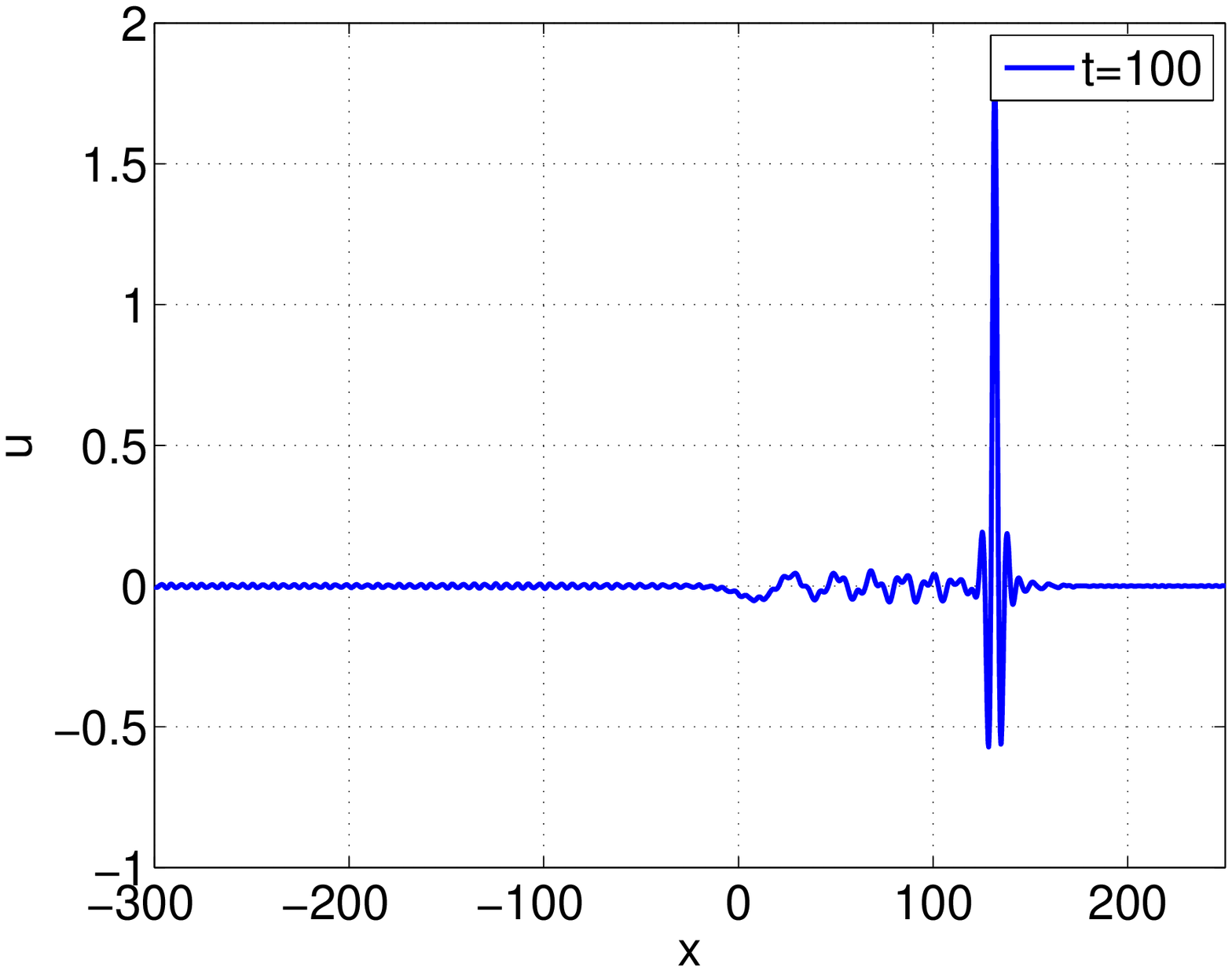}}
\subfigure[]
{\includegraphics[width=6cm]{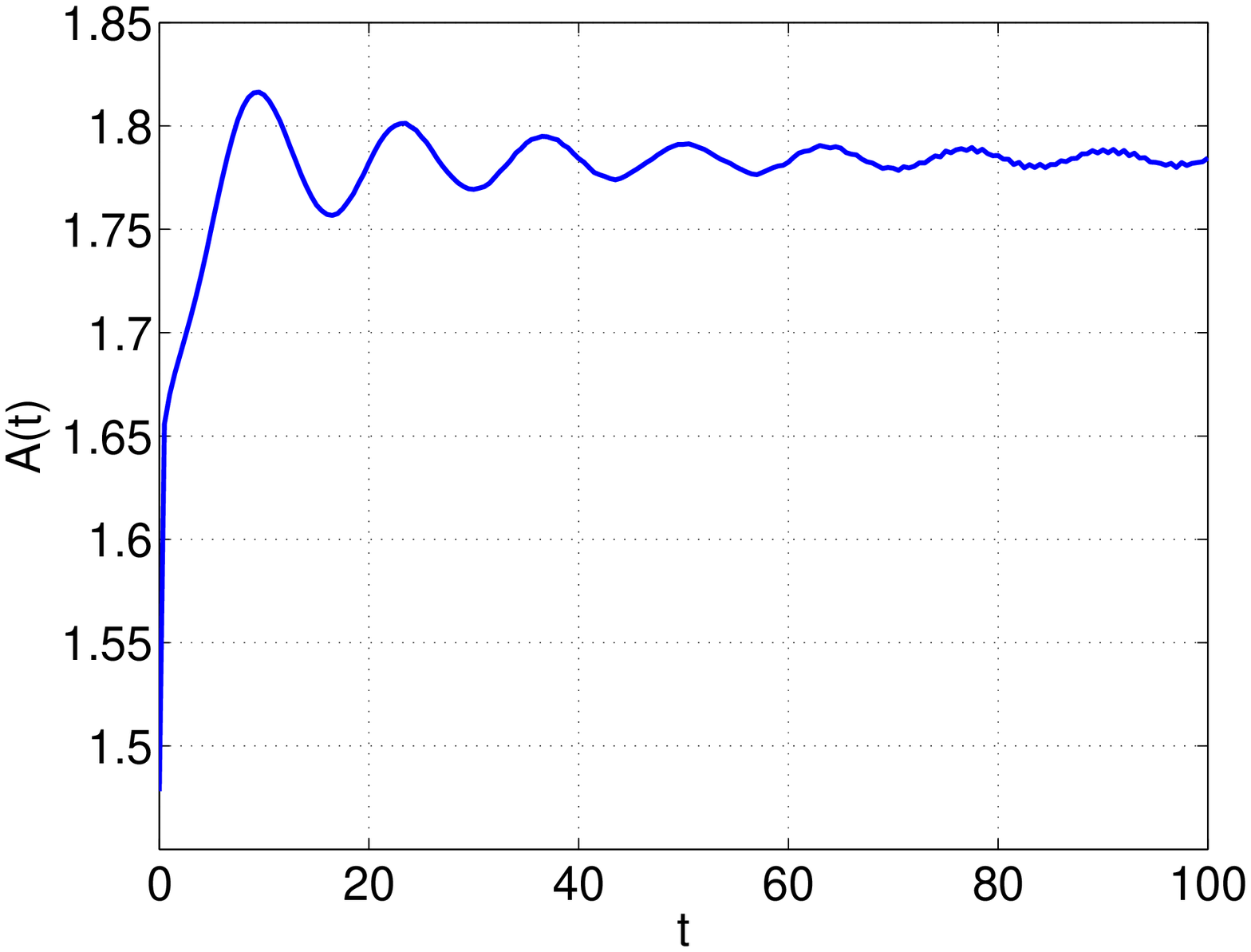}}
\subfigure[]
{\includegraphics[width=6cm]{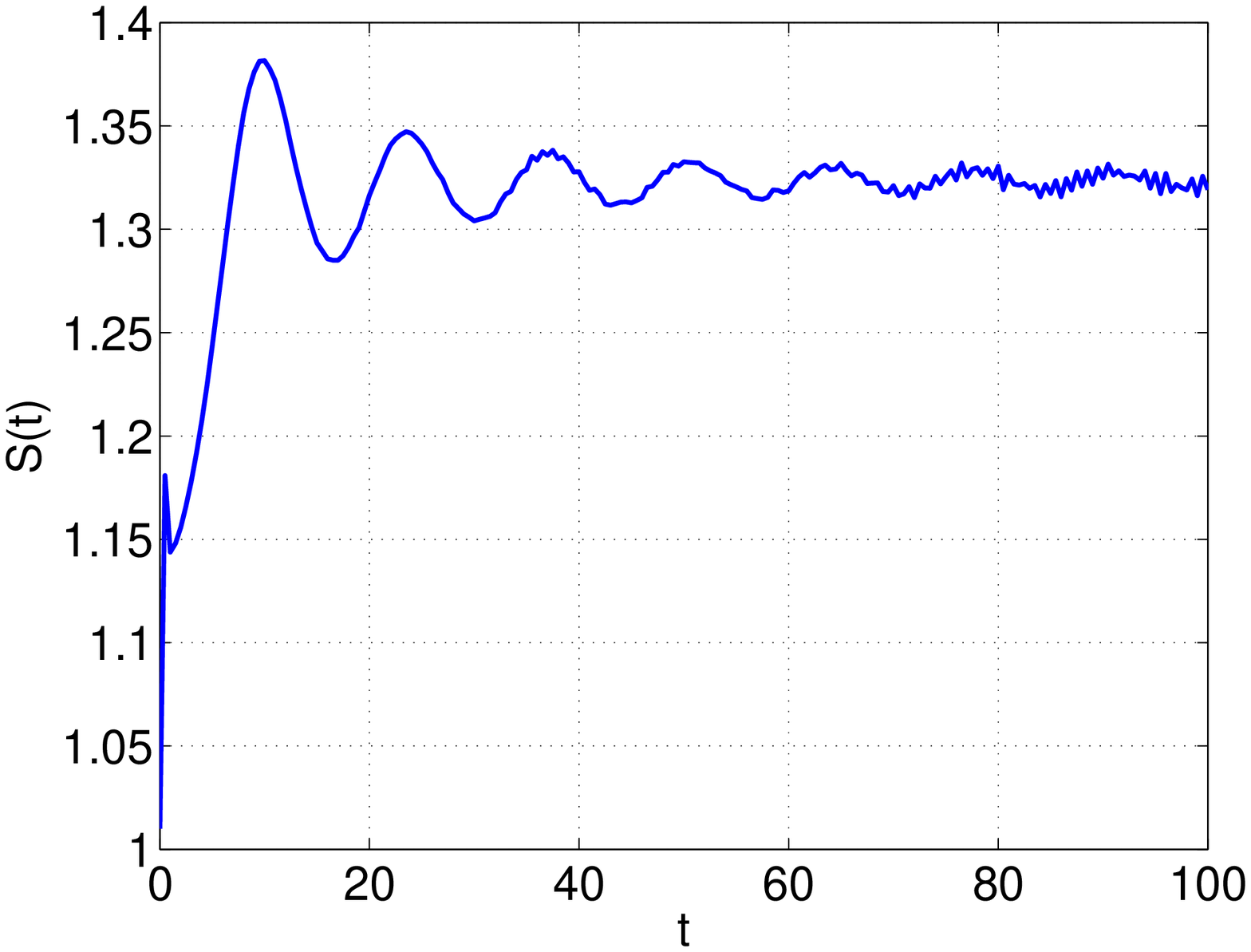}}
\caption{$r=1, m=2, q=3, c_{s}=1.01$,  $\gamma=1.8$, $A=1.1$. (a)  Numerical approximation at $t=0$; (b) numerical approximation at $t=100$; (c) Evolution of the amplitude of the main numerical pulse; (d) Evolution of speed.}
\label{gbenfig513_2}
\end{figure}

The general case does not seem to give a different behaviour, including the case $q=4$. By way of illustration we display the results corresponding to $m=2, r=1$ by fixing $\delta=1, c_{s}=1.01$ (thus $\gamma_{max}(c_{s})\approx 2.009975$), $\gamma=1.8$ and $q=2,3$. For $A=1.1$ these are given in Figures \ref{gbenfig513_1}, \ref{gbenfig513_2} respectively. The amplitudes of the initial profiles (without perturbation) are $u_{m}\approx 1.318697$ and $u_{m}\approx 1.478188$, and for the emerging waves $u_{m}\approx 1.699632$ and $u_{m}\approx 1.784304$ respectively, with corresponding approximate speeds of $c_{s}\approx 1.176478$ and $c_{s}\approx 1.323700$. 

%
%
\subsection{Large perturbations of solitary waves. Resolution property}
\label{sec33}
%

\begin{figure}[htbp]
\centering
\subfigure[]
{\includegraphics[width=6cm]{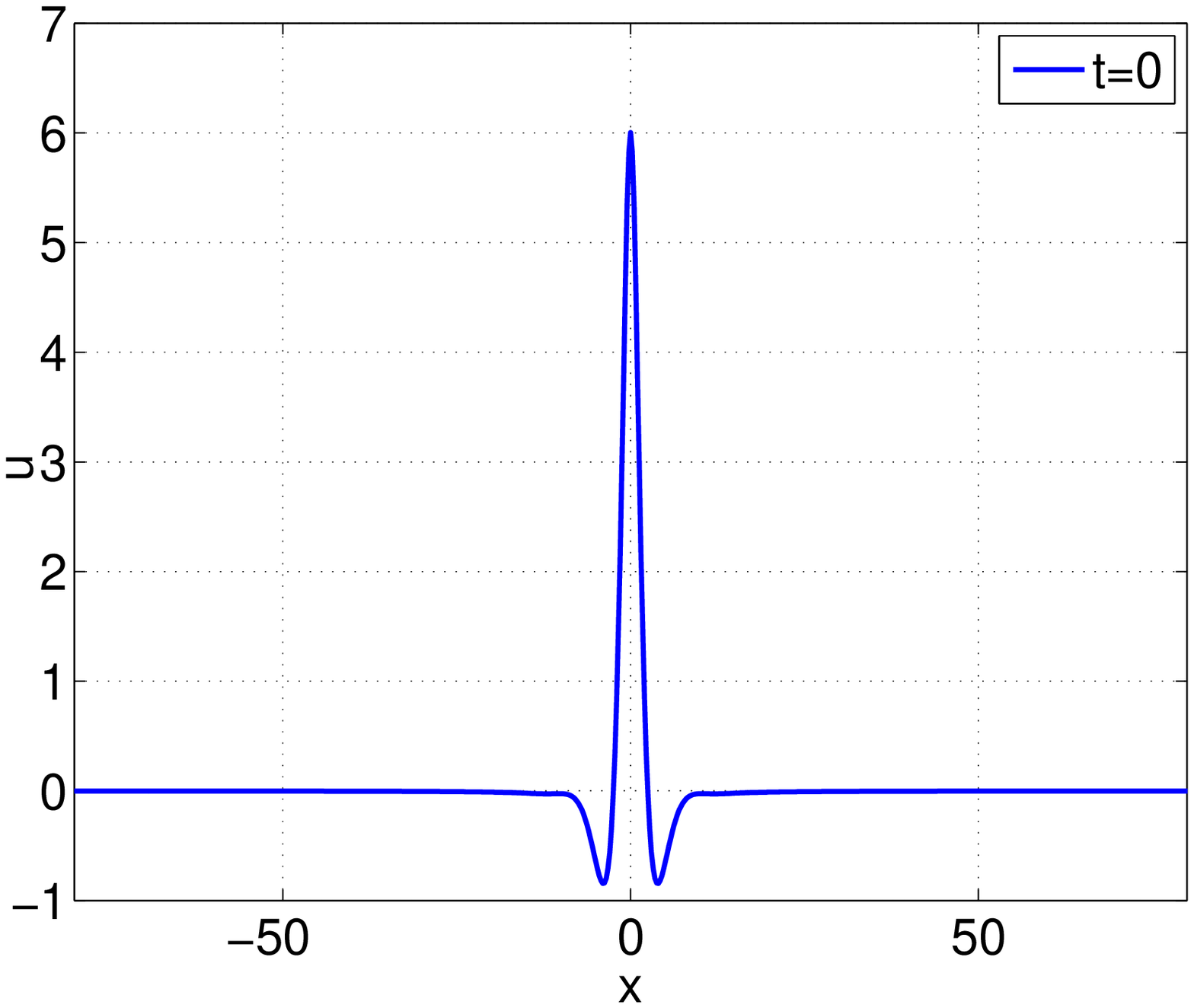}}
\subfigure[]
{\includegraphics[width=6cm]{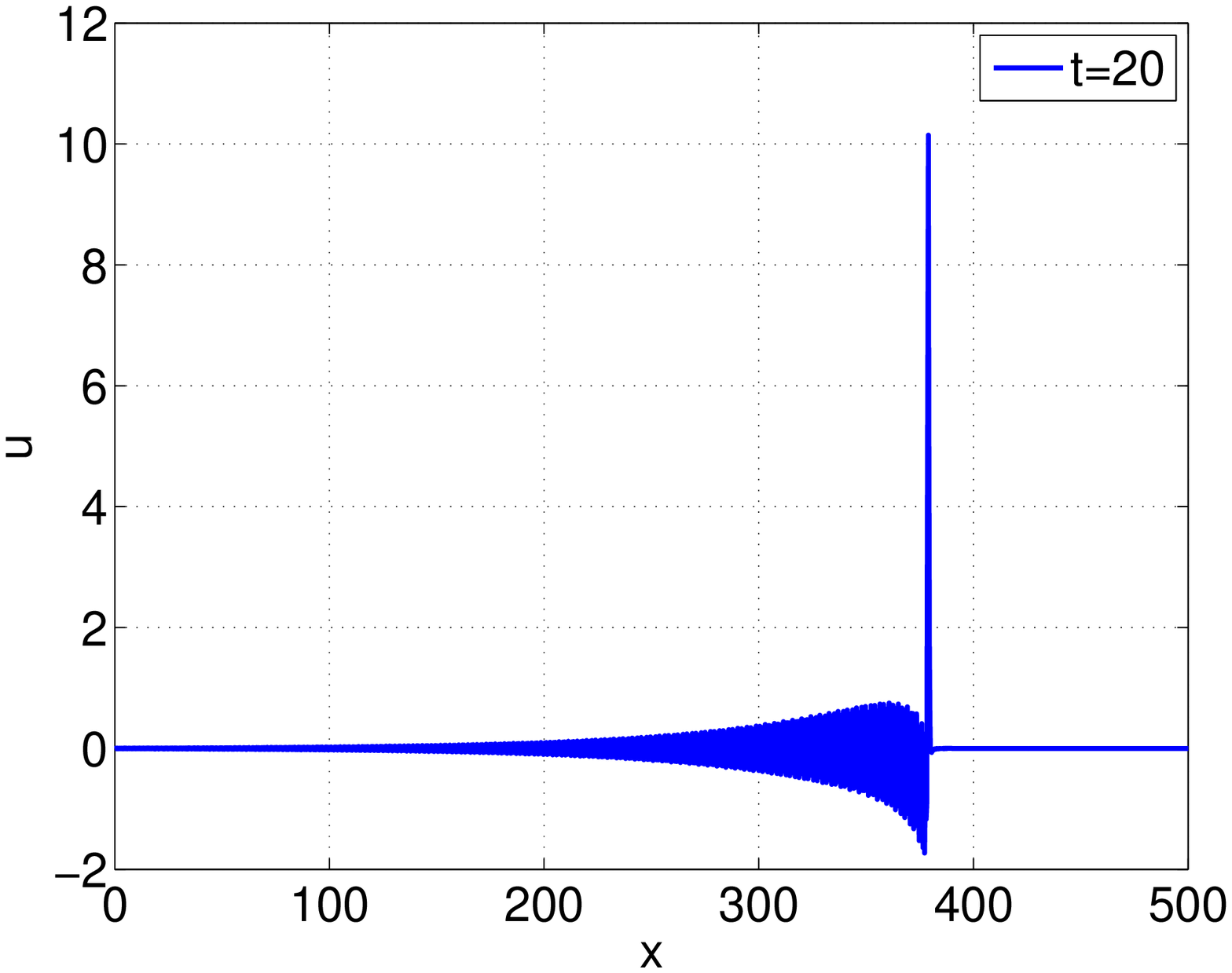}}
\subfigure[]
{\includegraphics[width=6cm]{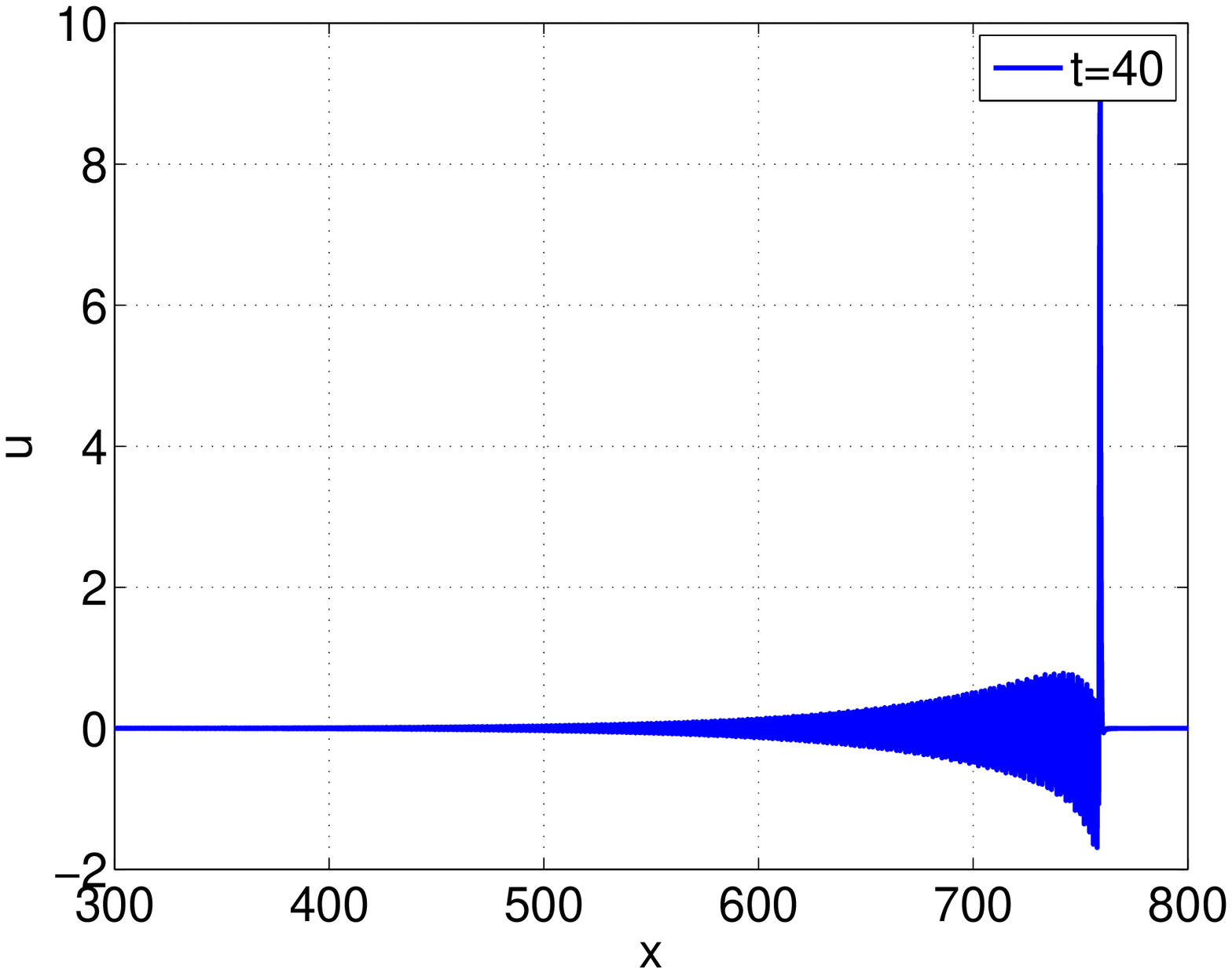}}
\subfigure[]
{\includegraphics[width=6cm]{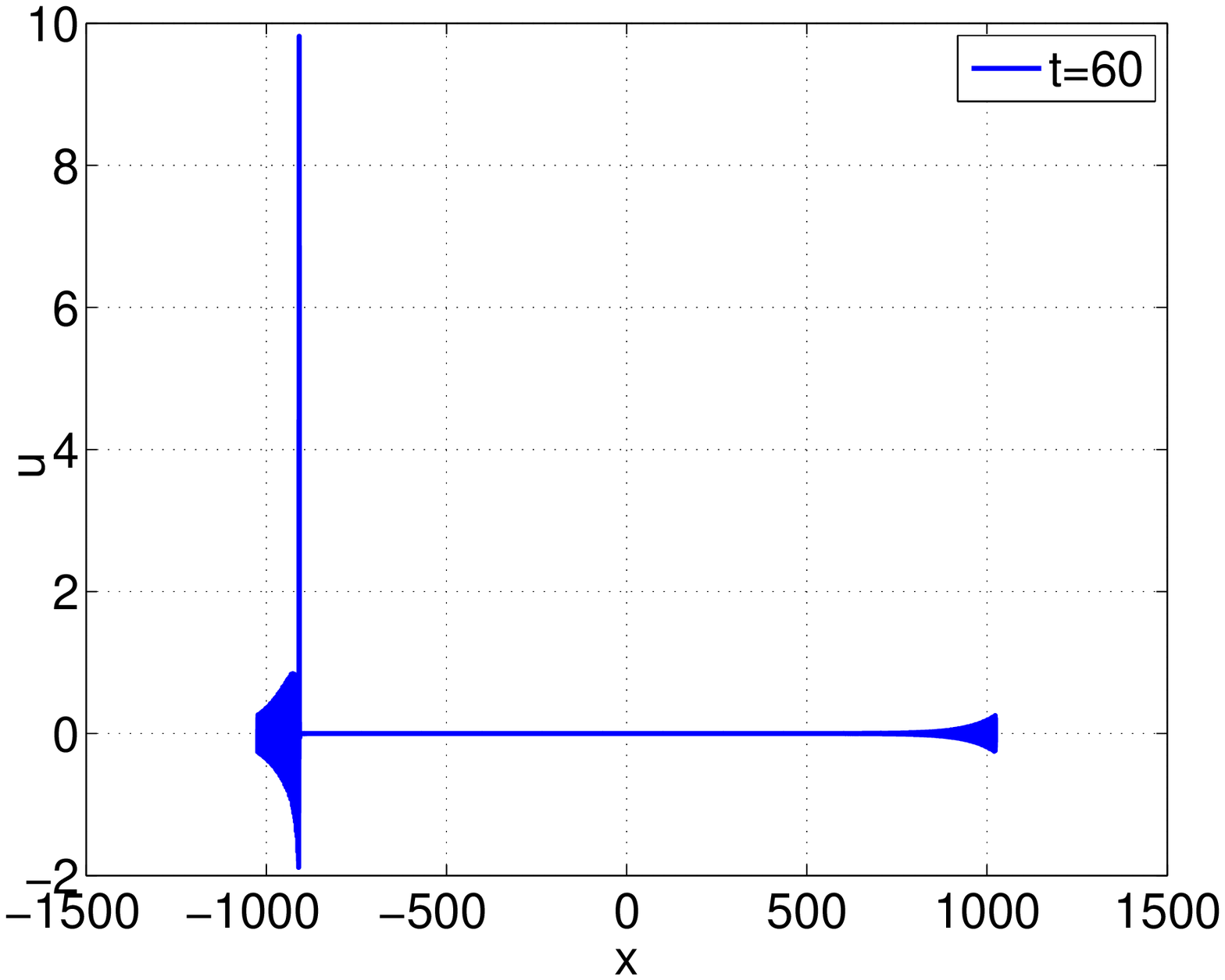}}
\subfigure[]
{\includegraphics[width=6cm]{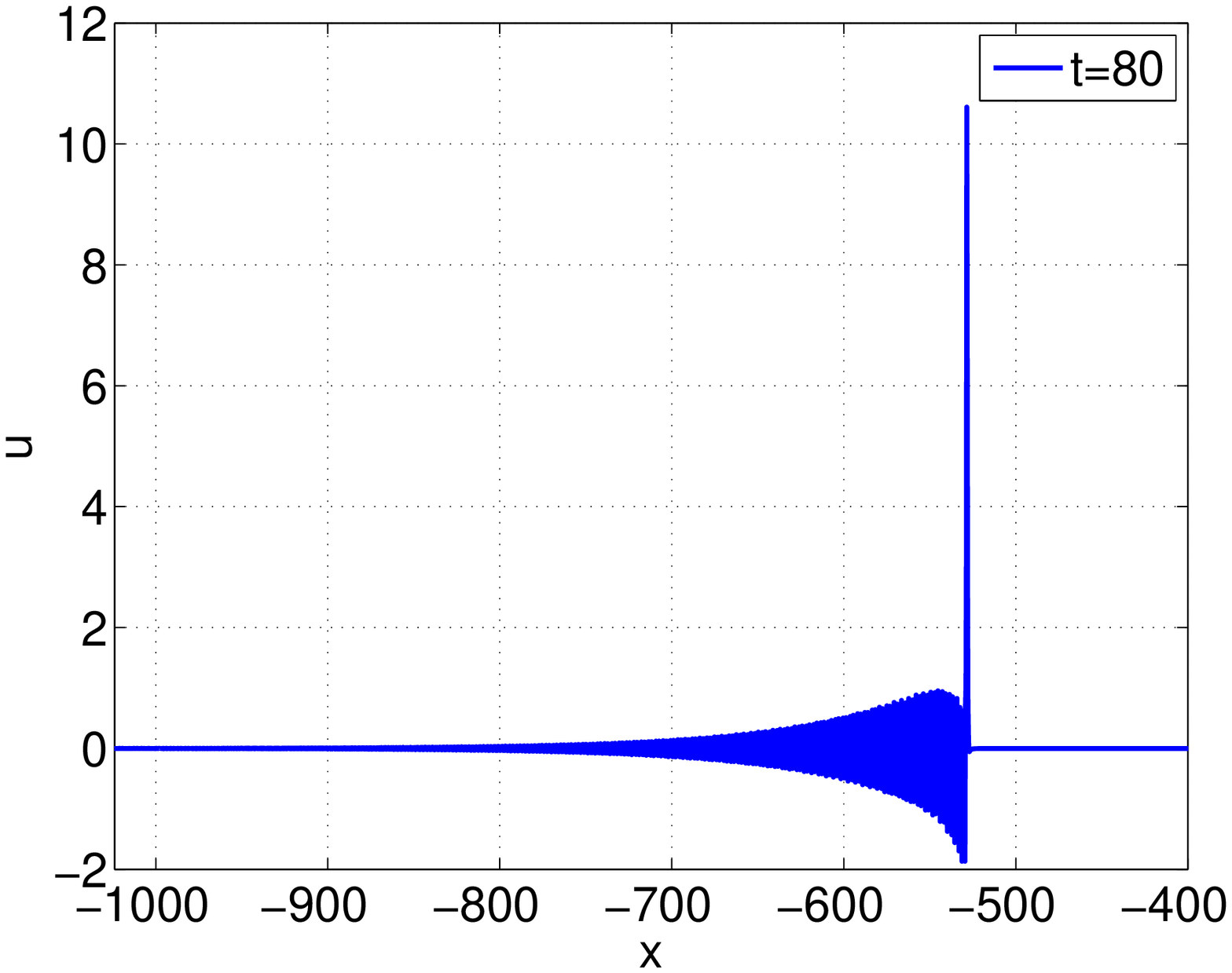}}
\subfigure[]
{\includegraphics[width=6cm]{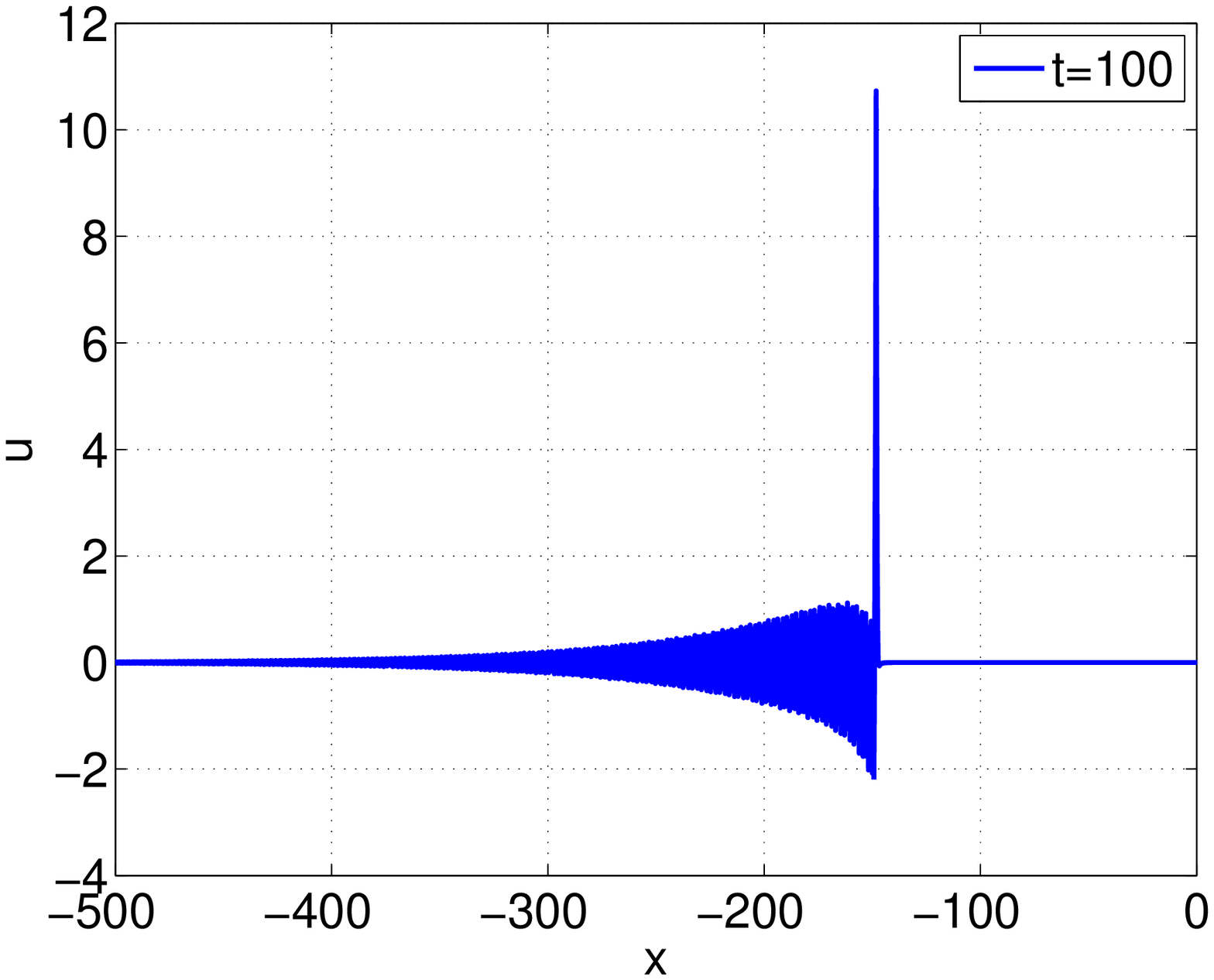}}
\caption{$r=1/2, m=1, q=2, c_{s}=0.75$,  $\gamma=1.2$, $A=4$. (a)-(f)  Numerical approximation at $t=0,20,40,60,80,100$.}
\label{gbenfig512_2}
\end{figure}

\begin{figure}[htbp]
\centering
\subfigure[]
{\includegraphics[width=6cm]{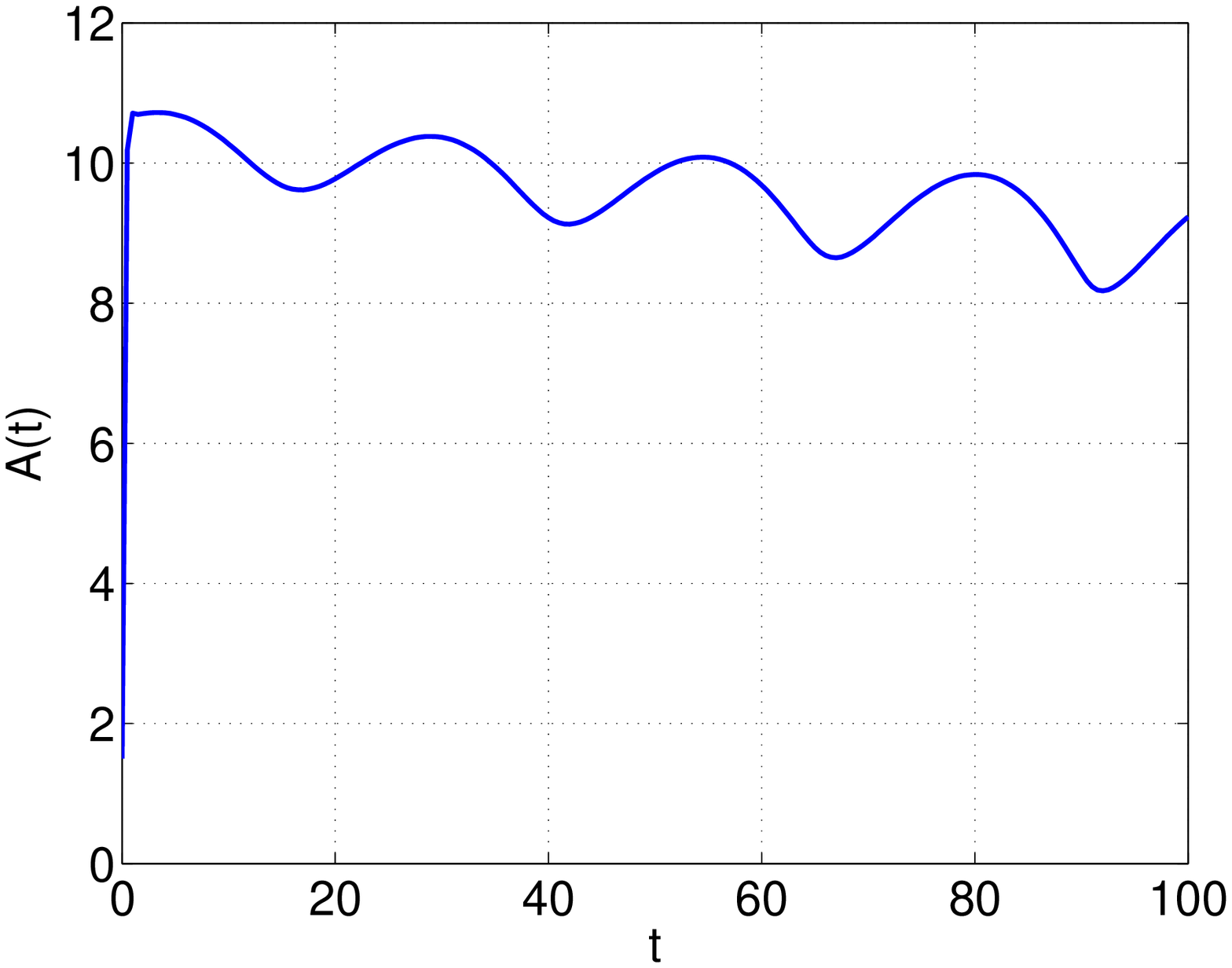}}
\subfigure[]
{\includegraphics[width=6cm]{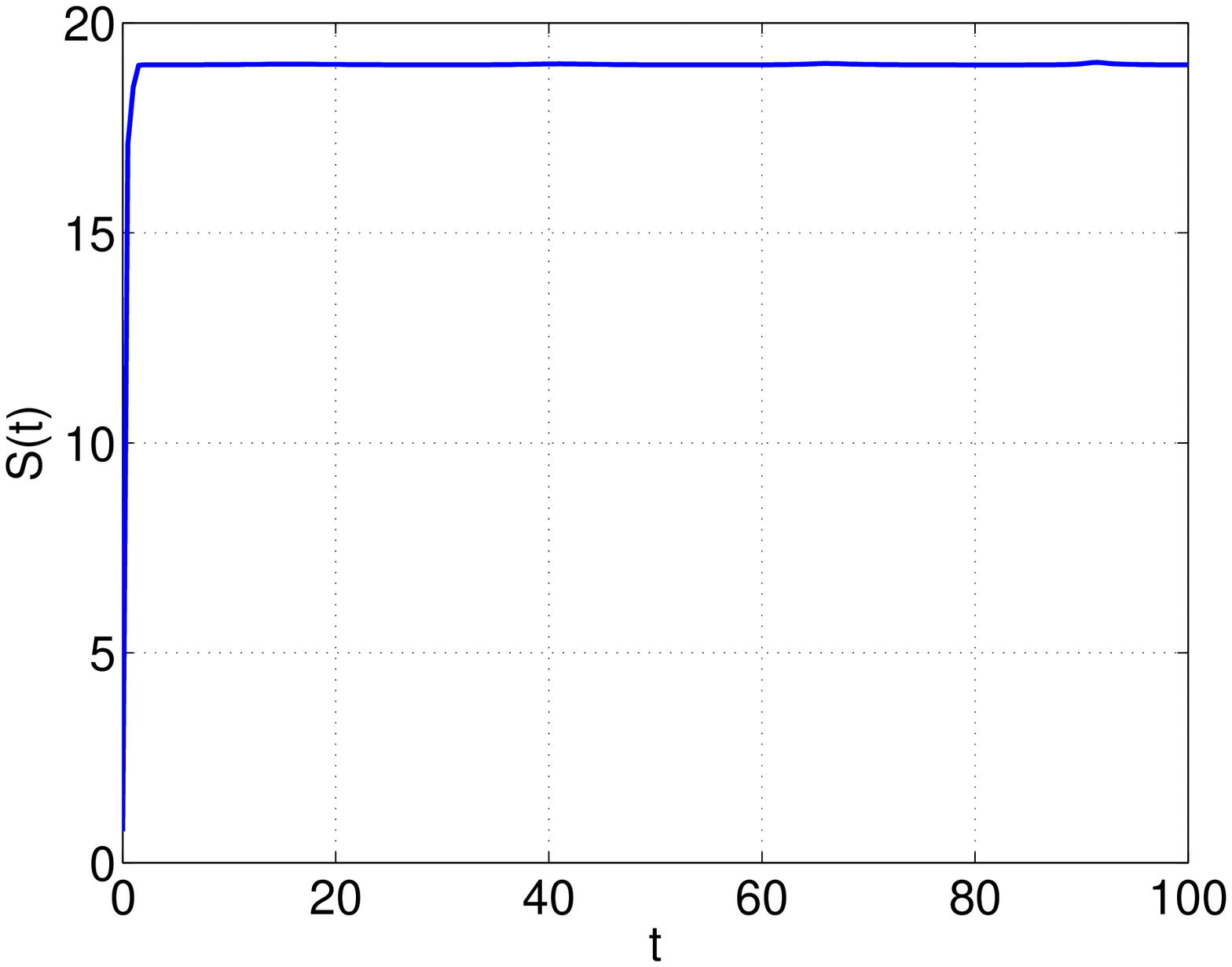}}
\caption{$r=1/2, m=1, q=2, c_{s}=0.75$,  $\gamma=1.2$, $A=4$. (a)  Evolution of the amplitude of the main numerical pulse; (b) Evolution of speed.}
\label{gbenfig512_2b}
\end{figure}

\begin{figure}[htbp]
\centering
\subfigure[]
{\includegraphics[width=6cm]{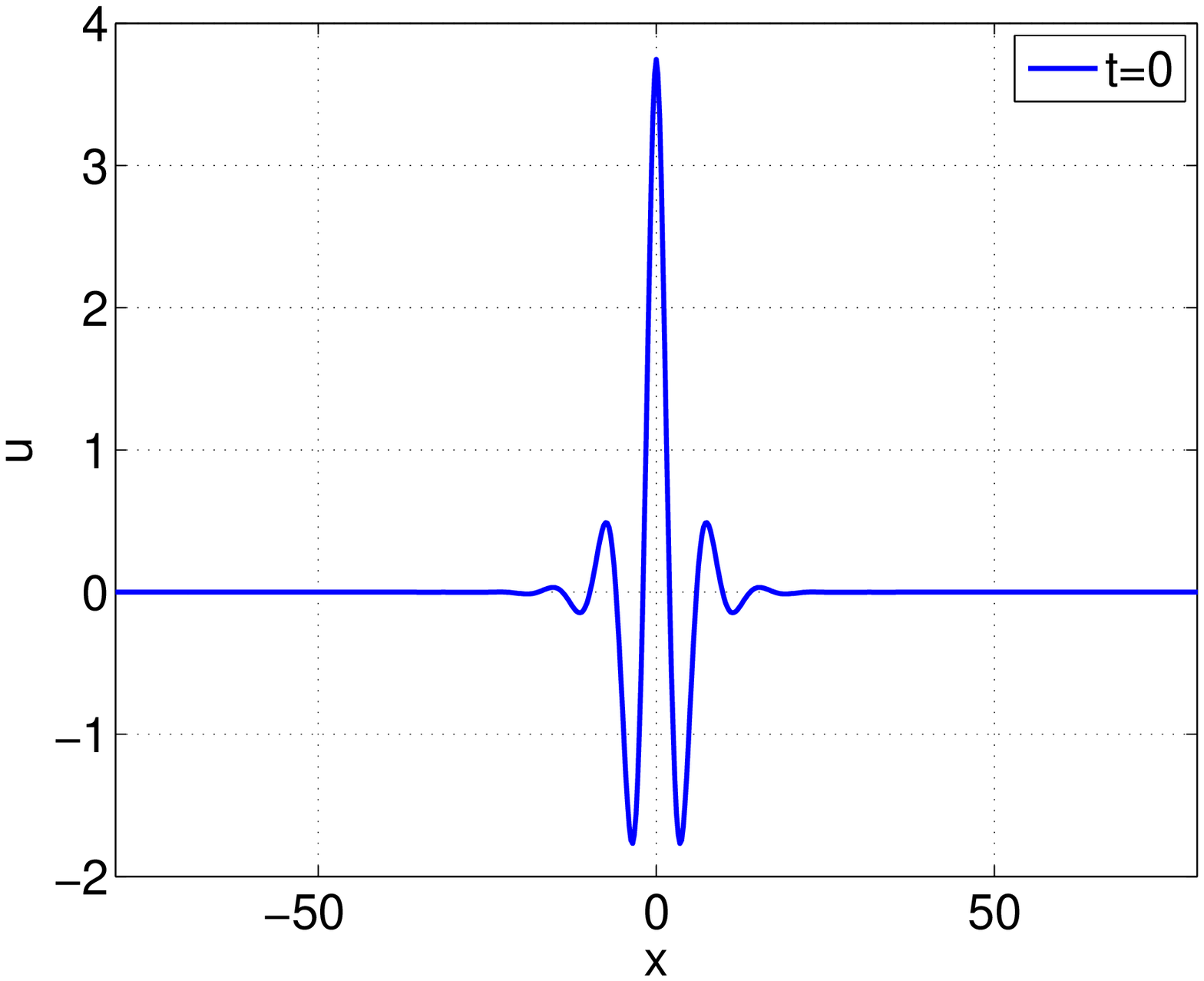}}
\subfigure[]
{\includegraphics[width=6cm]{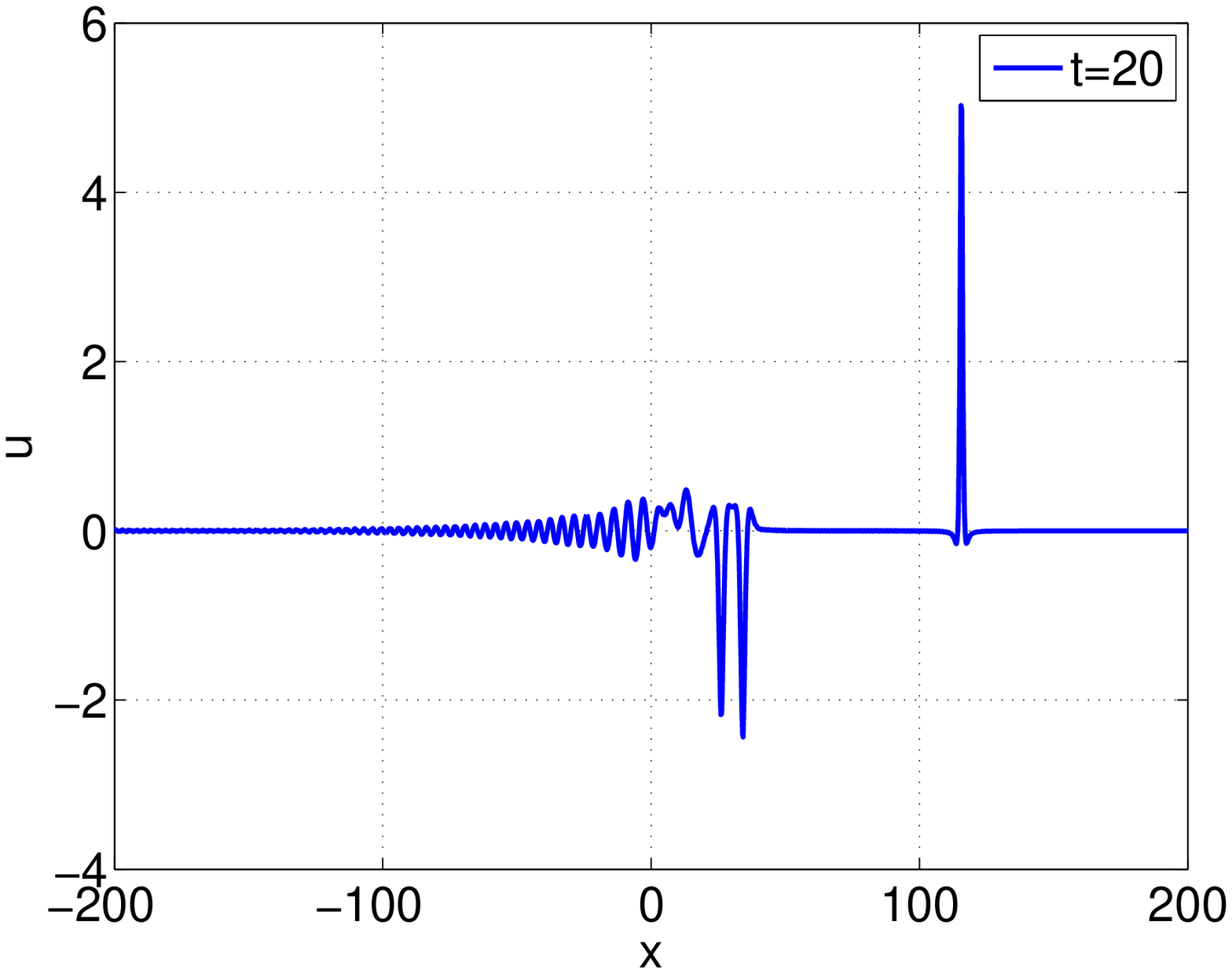}}
\subfigure[]
{\includegraphics[width=6cm]{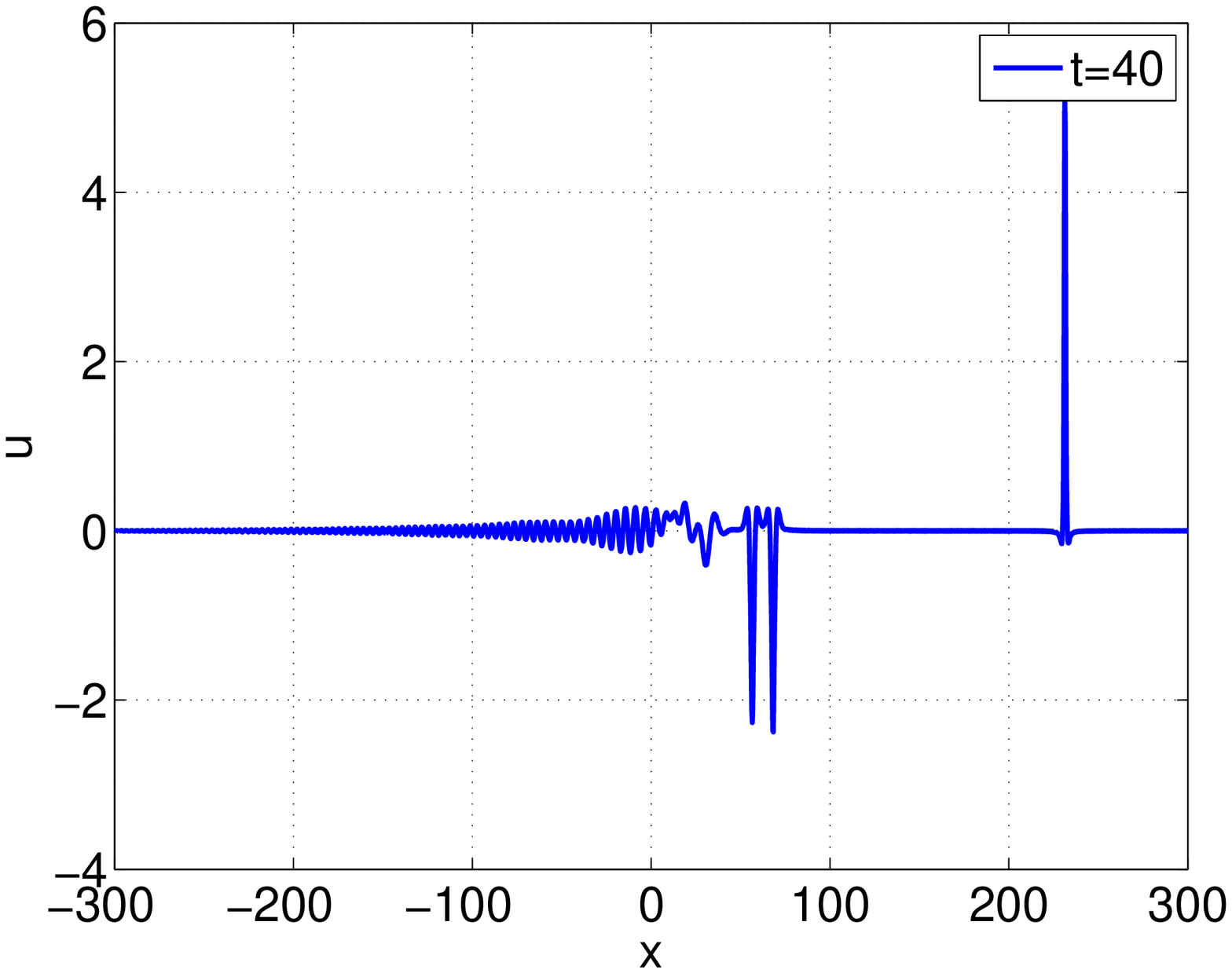}}
\subfigure[]
{\includegraphics[width=6cm]{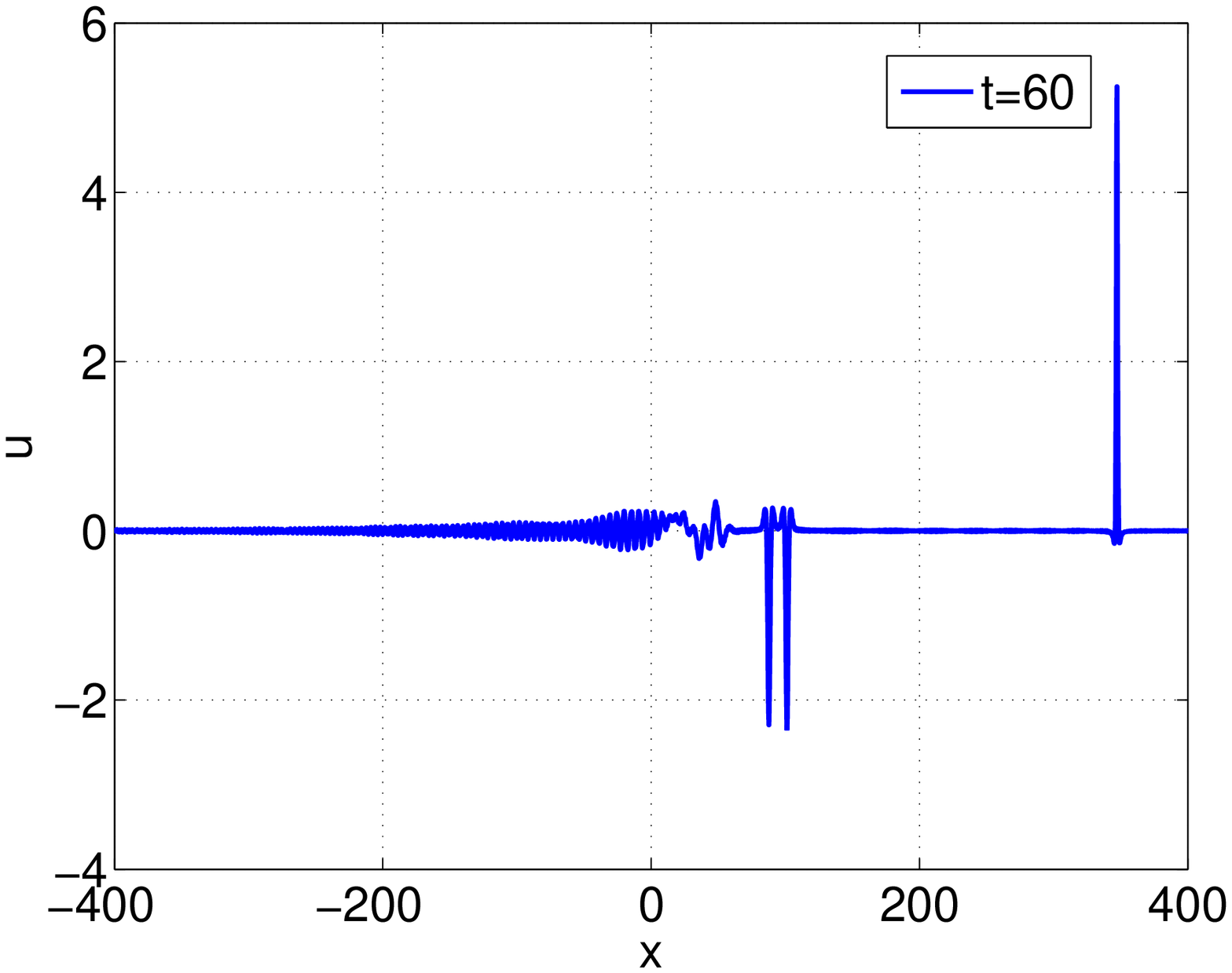}}
\subfigure[]
{\includegraphics[width=6cm]{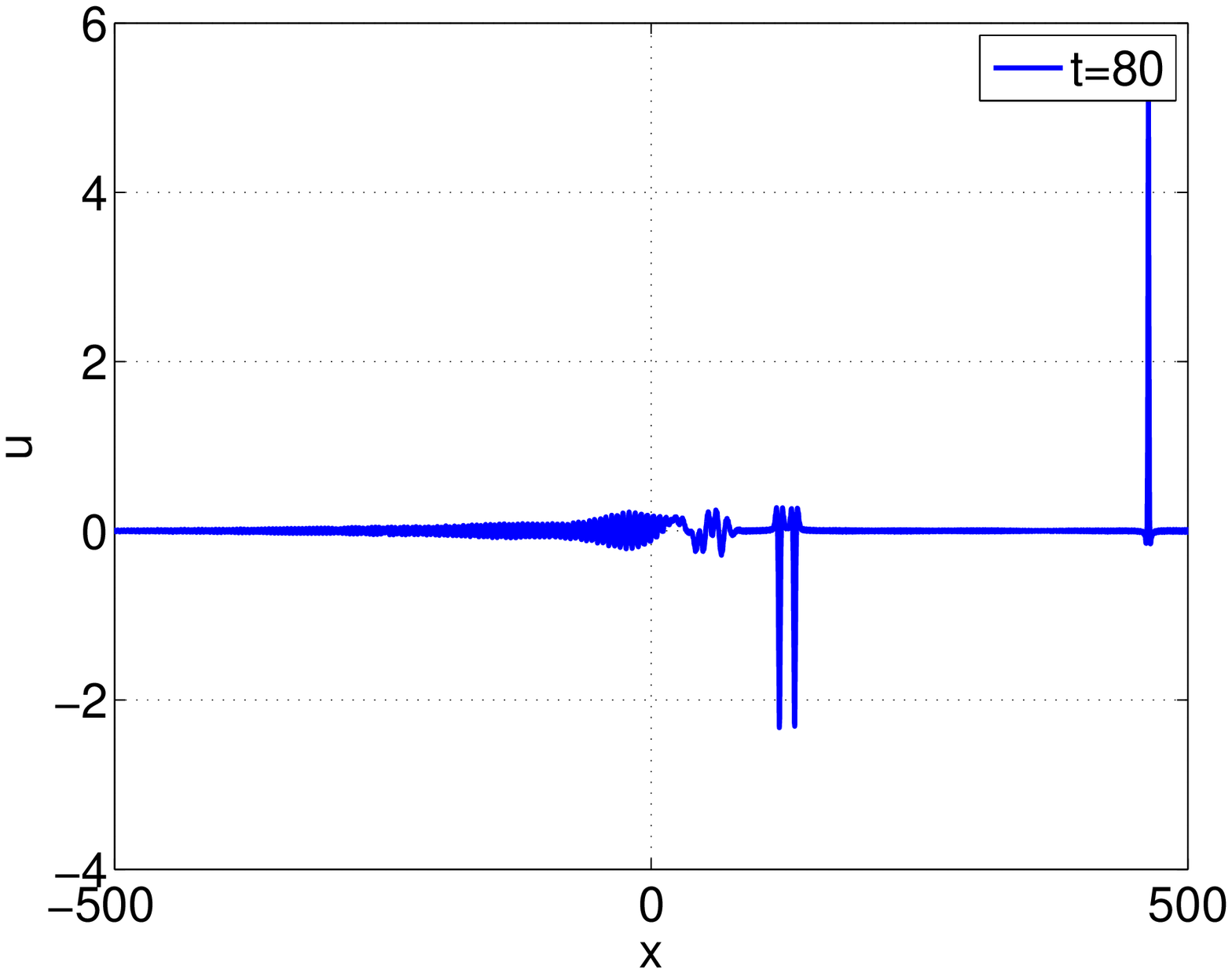}}
\subfigure[]
{\includegraphics[width=6cm]{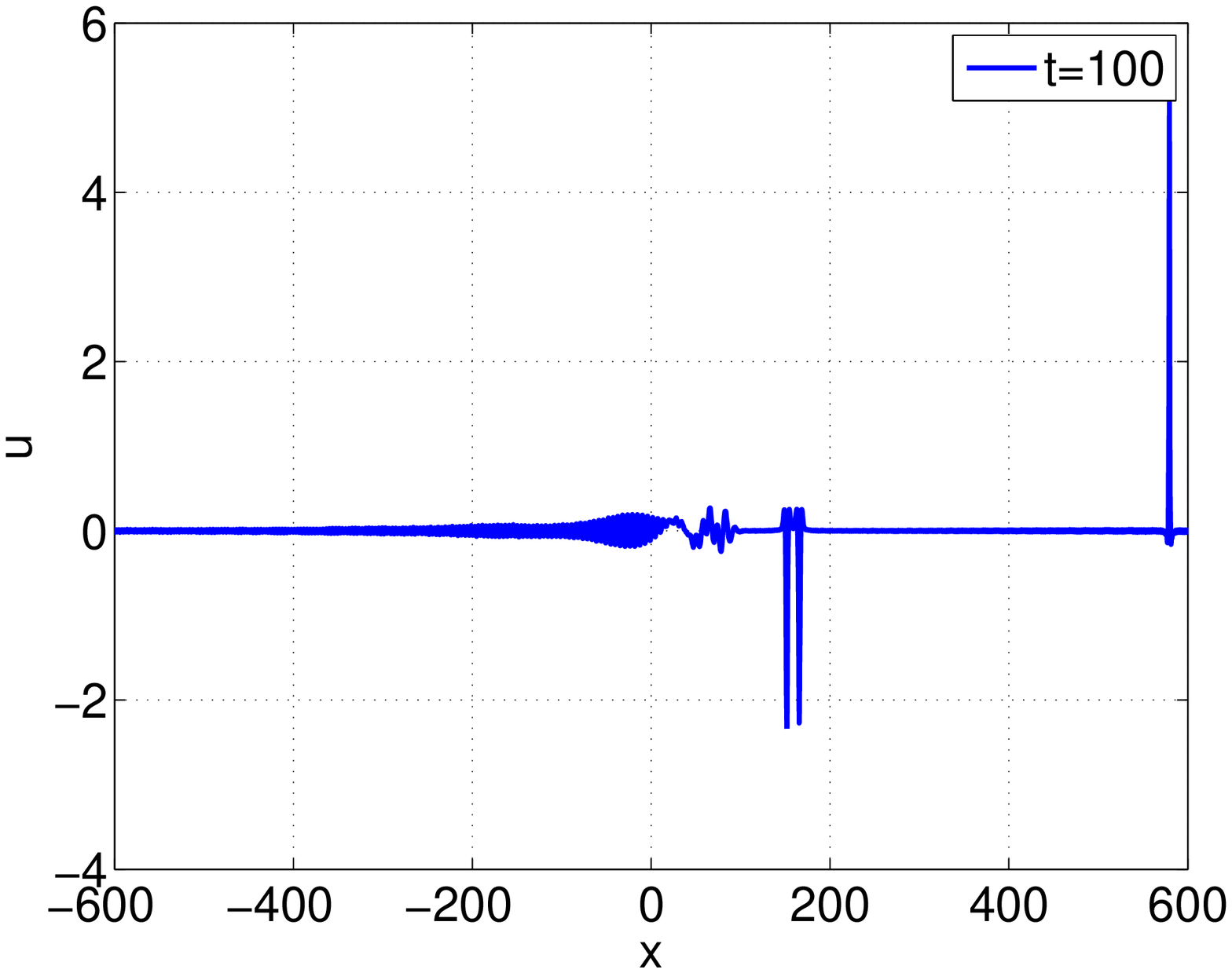}}
\caption{$r=1/2, m=1, q=2, c_{s}=0.75$,  $\gamma=1.6$, $A=4$. (a)-(f)  Numerical approximation at $t=0,20,40,60,80,100$.}
\label{gbenfig512_1}
\end{figure}

\begin{figure}[!htbp]
\centering
\subfigure[]
{\includegraphics[width=6cm]{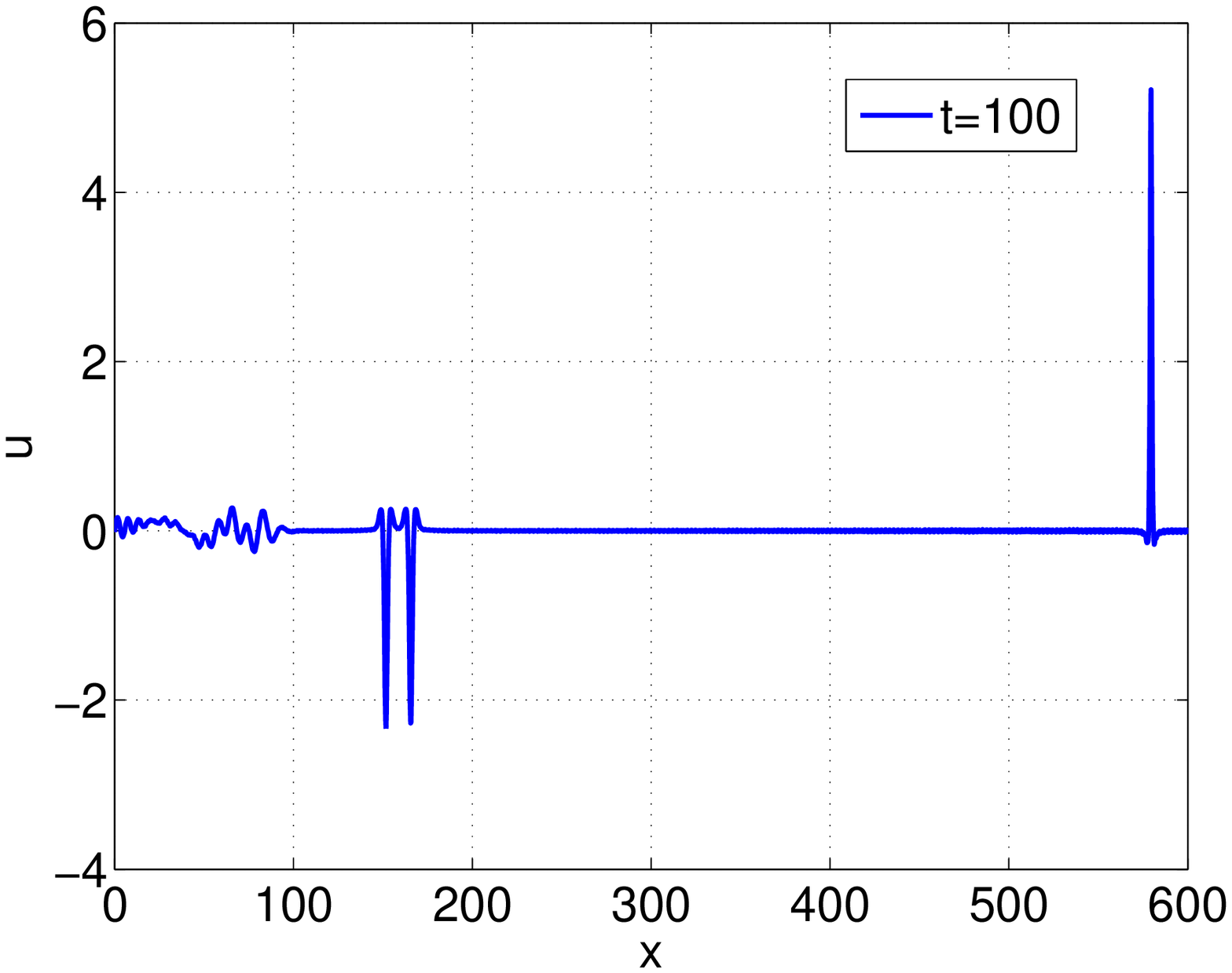}}
\subfigure[]
{\includegraphics[width=6cm]{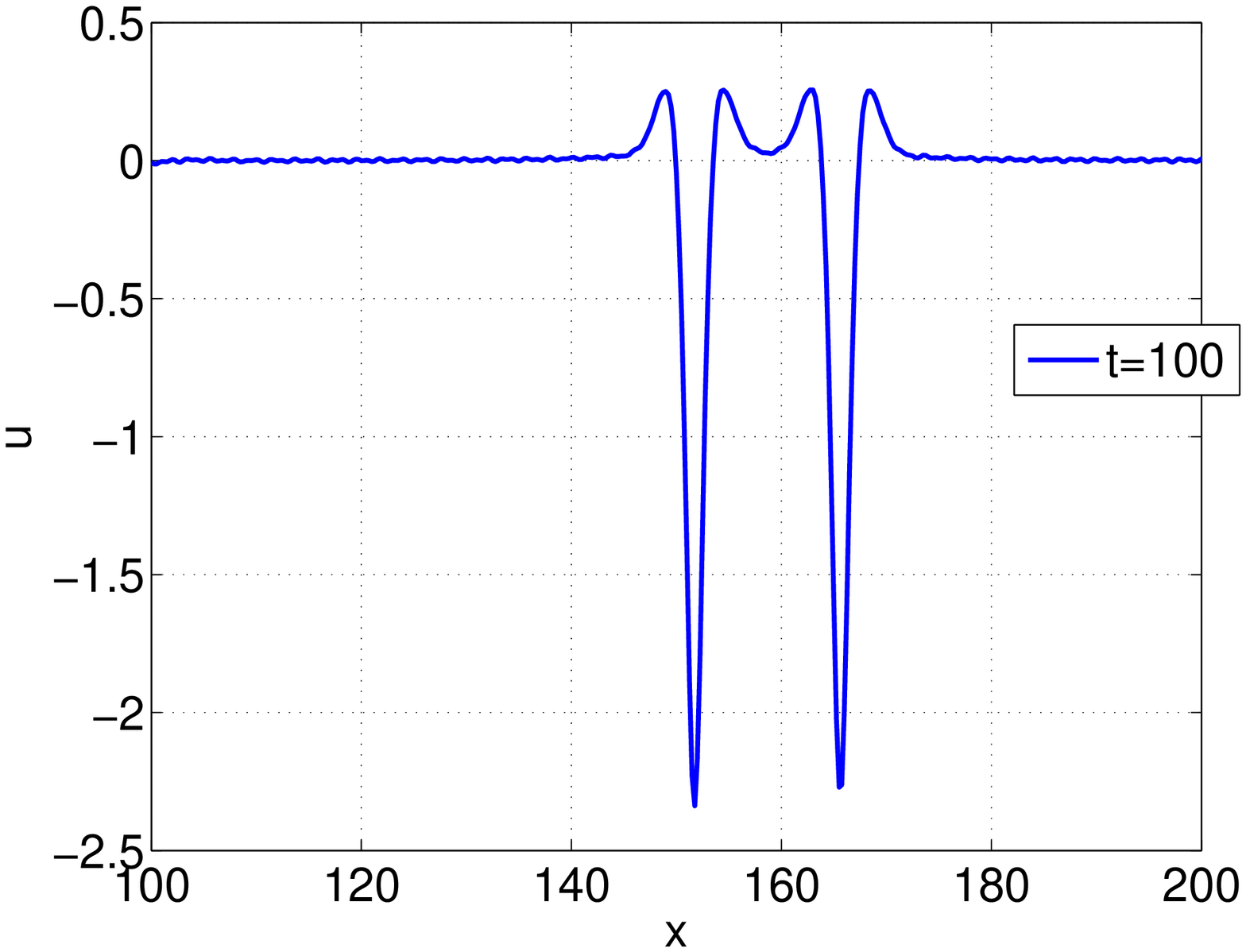}}
\caption{$r=1/2, m=1, q=2, c_{s}=0.75$,  $\gamma=1.6$, $A=4$. Magnifications of Figure \ref{gbenfig512_1}(f).}
\label{gbenfig512_1m}
\end{figure}

\begin{figure}[!htbp]
\centering
\subfigure[]
{\includegraphics[width=6cm]{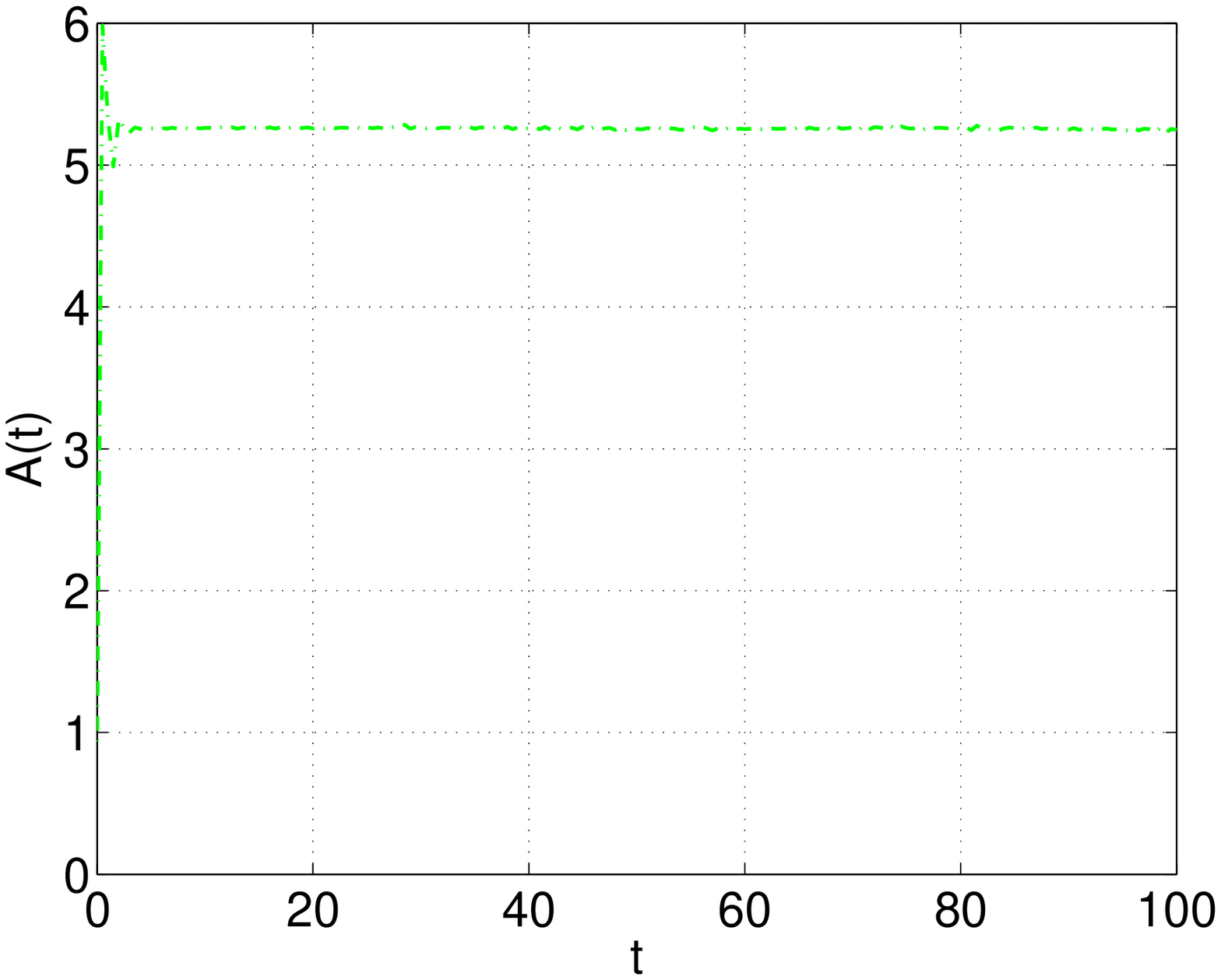}}
\subfigure[]
{\includegraphics[width=6cm]{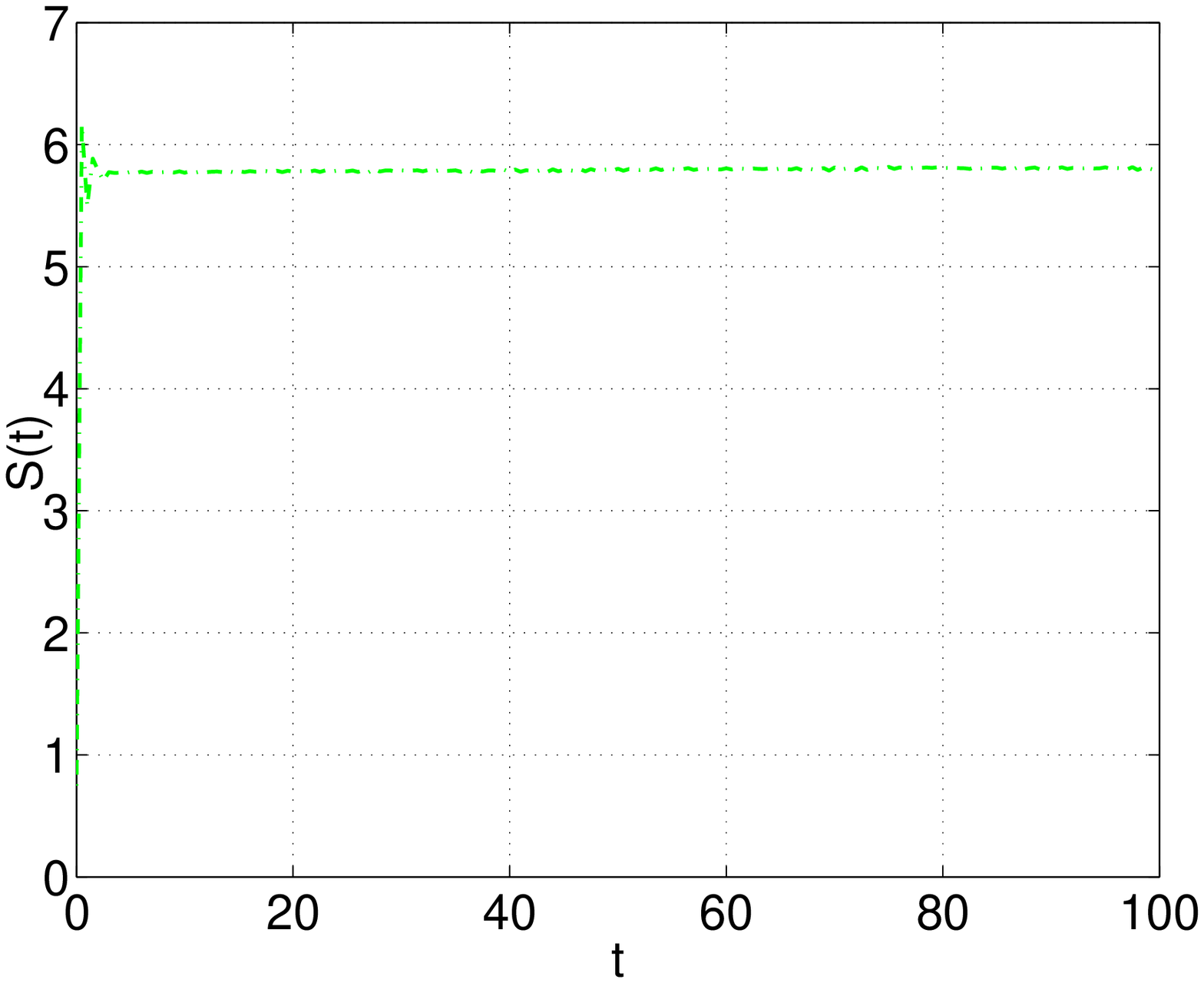}}
\caption{$r=1/2, m=1, q=2, c_{s}=0.75$,  $\gamma=1.6$, $A=4$. (a)  Evolution of the amplitude of the main numerical pulse; (b) Evolution of speed.}
\label{gbenfig512_1b}
\end{figure}

\begin{figure}[!htbp]
\centering
\subfigure[]
{\includegraphics[width=6cm]{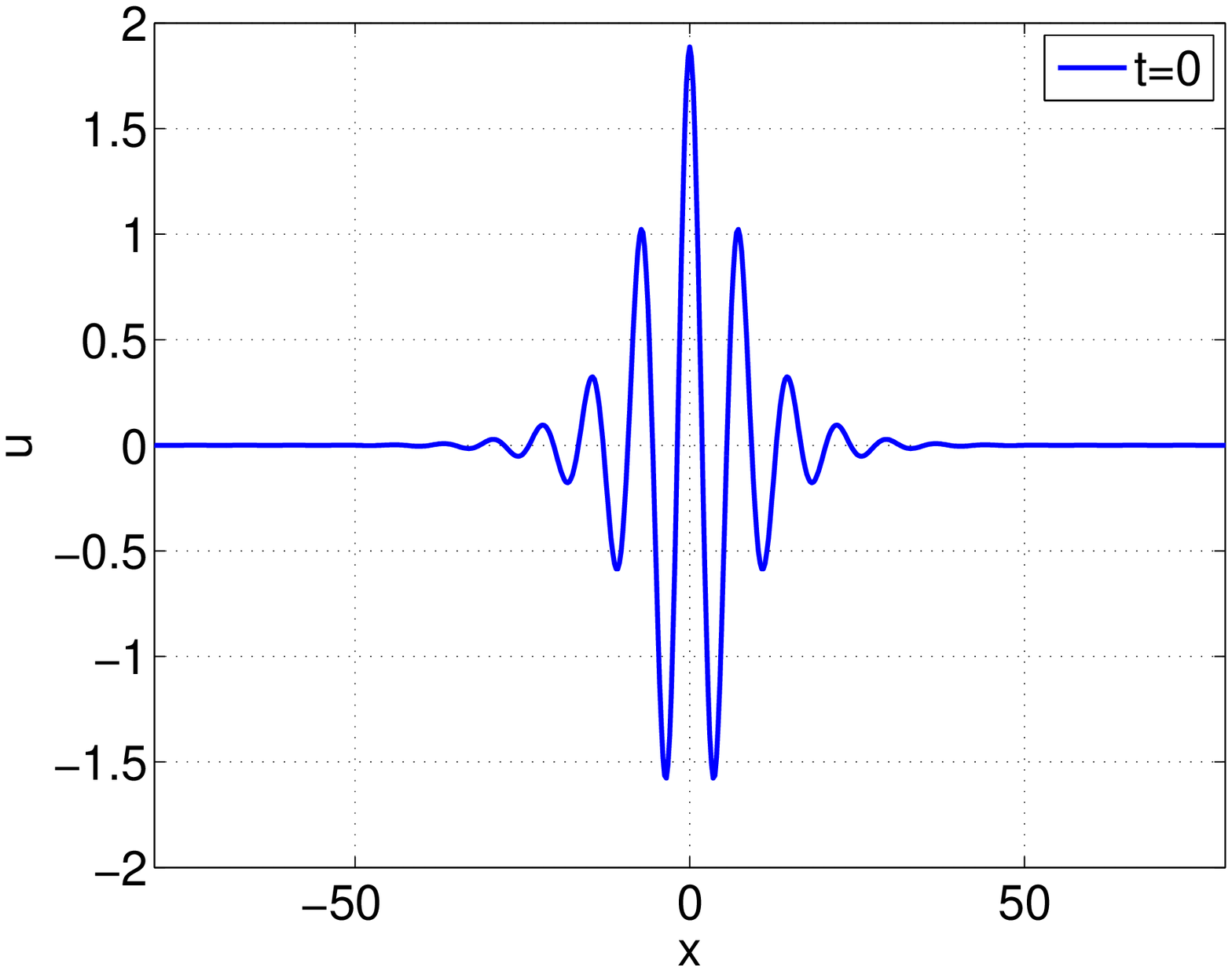}}
\subfigure[]
{\includegraphics[width=6cm]{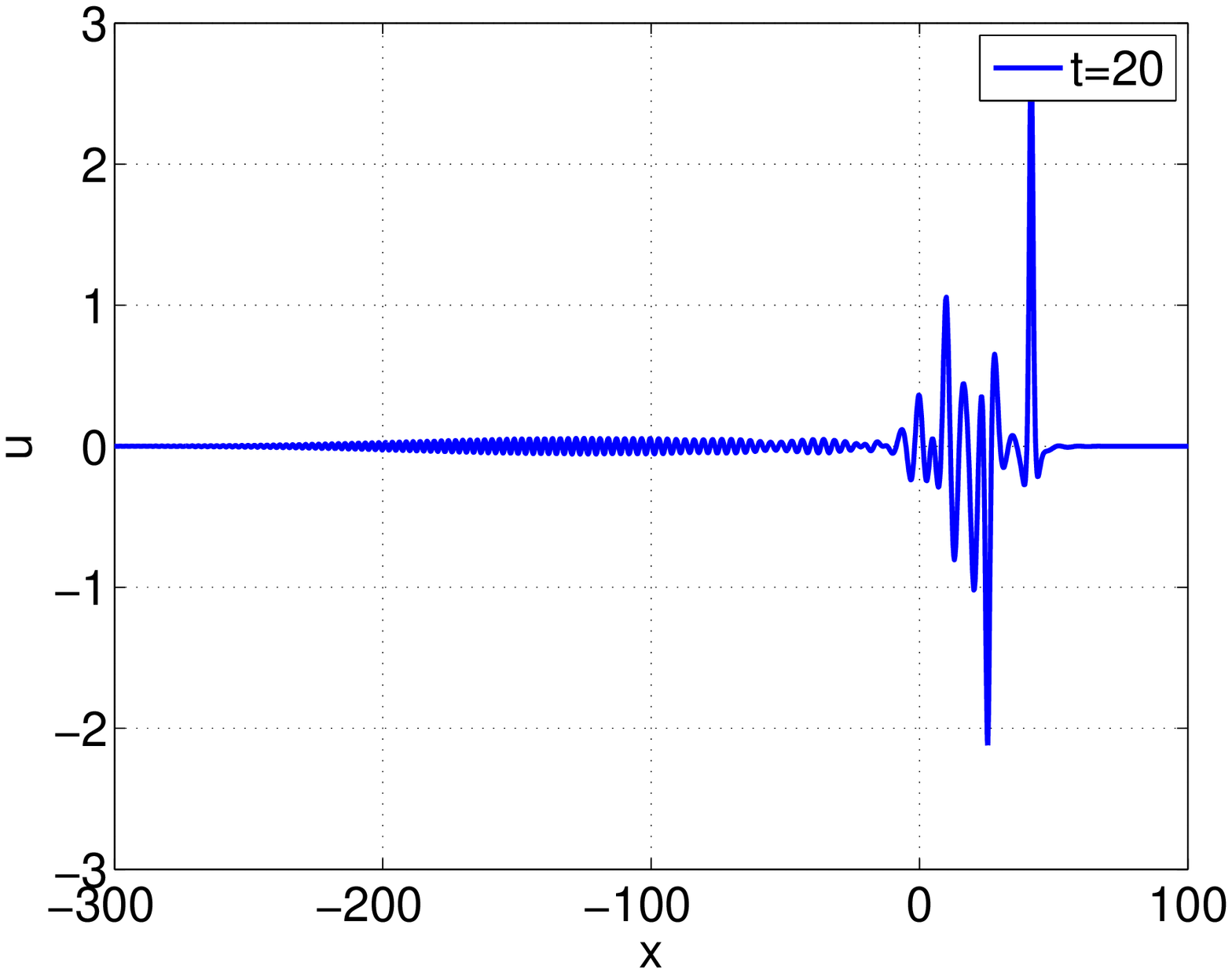}}
\subfigure[]
{\includegraphics[width=6cm]{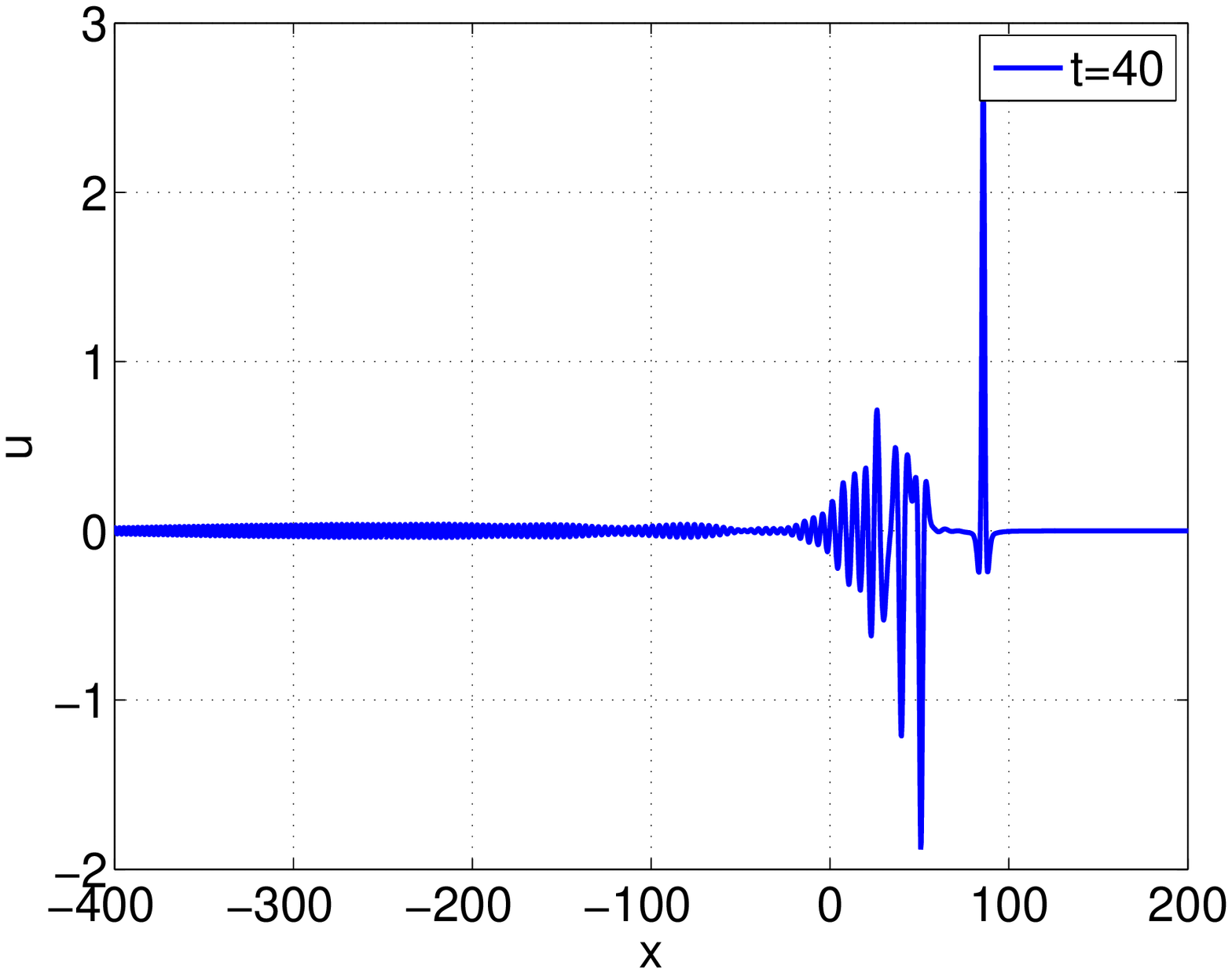}}
\subfigure[]
{\includegraphics[width=6cm]{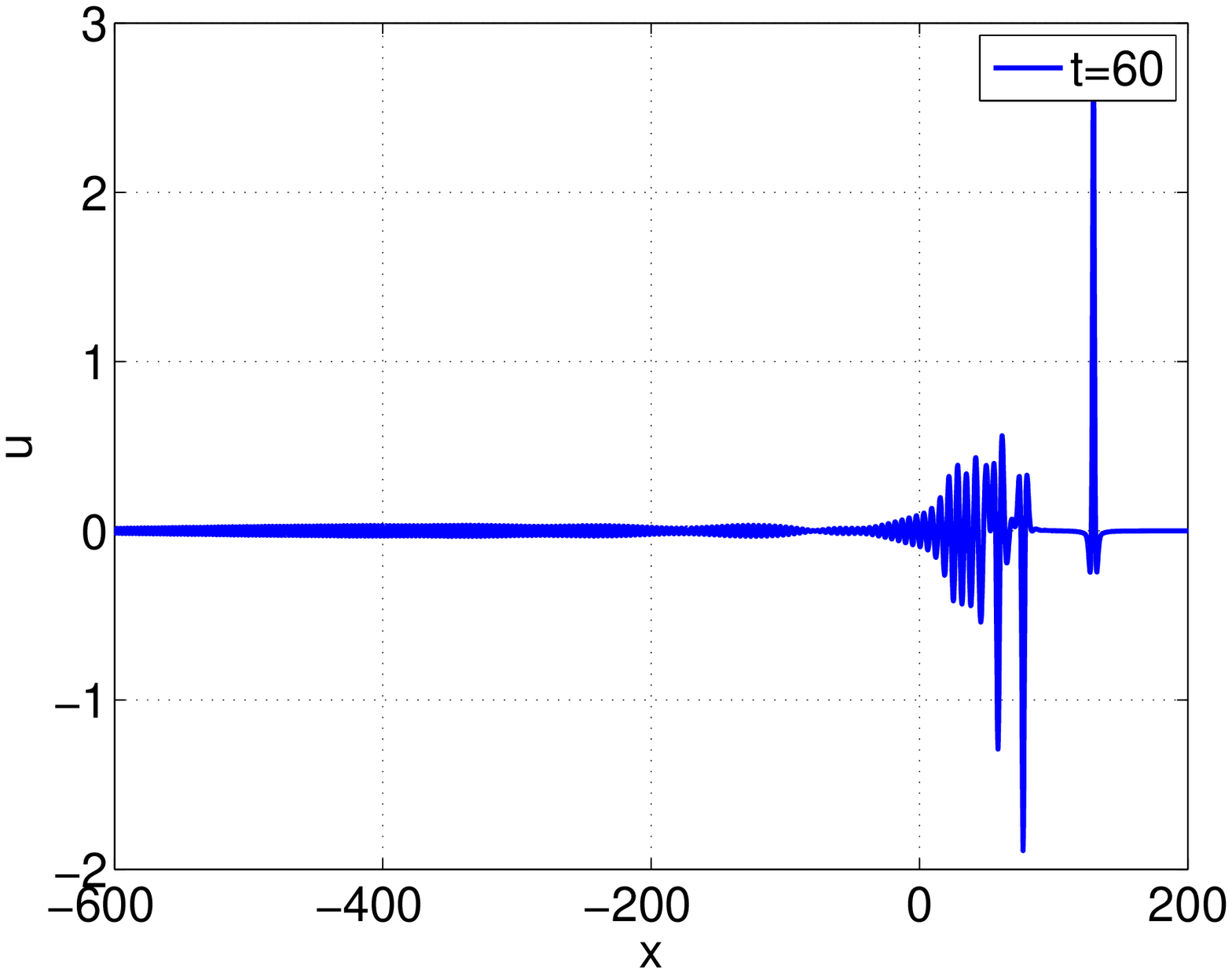}}
\subfigure[]
{\includegraphics[width=6cm]{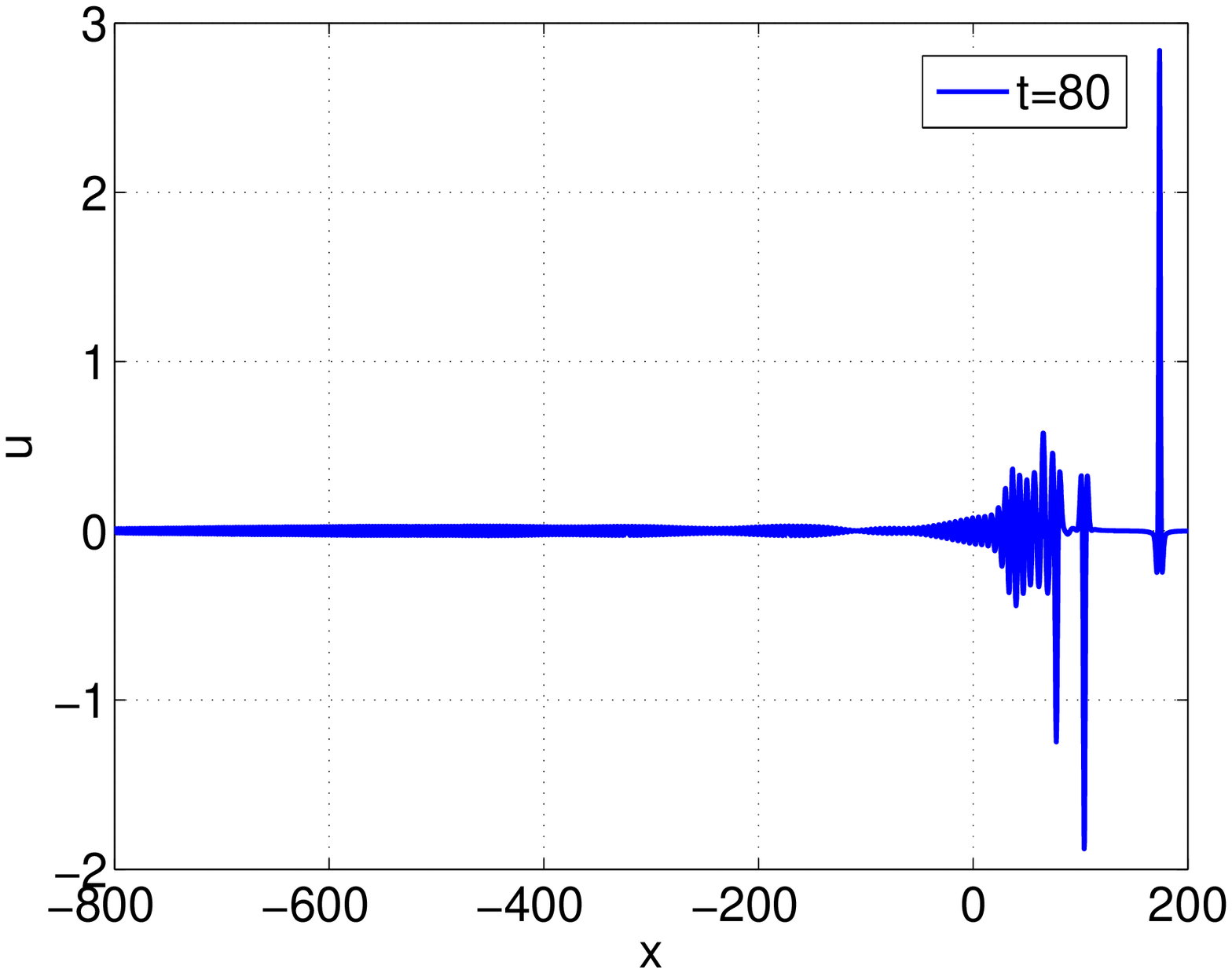}}
\subfigure[]
{\includegraphics[width=6cm]{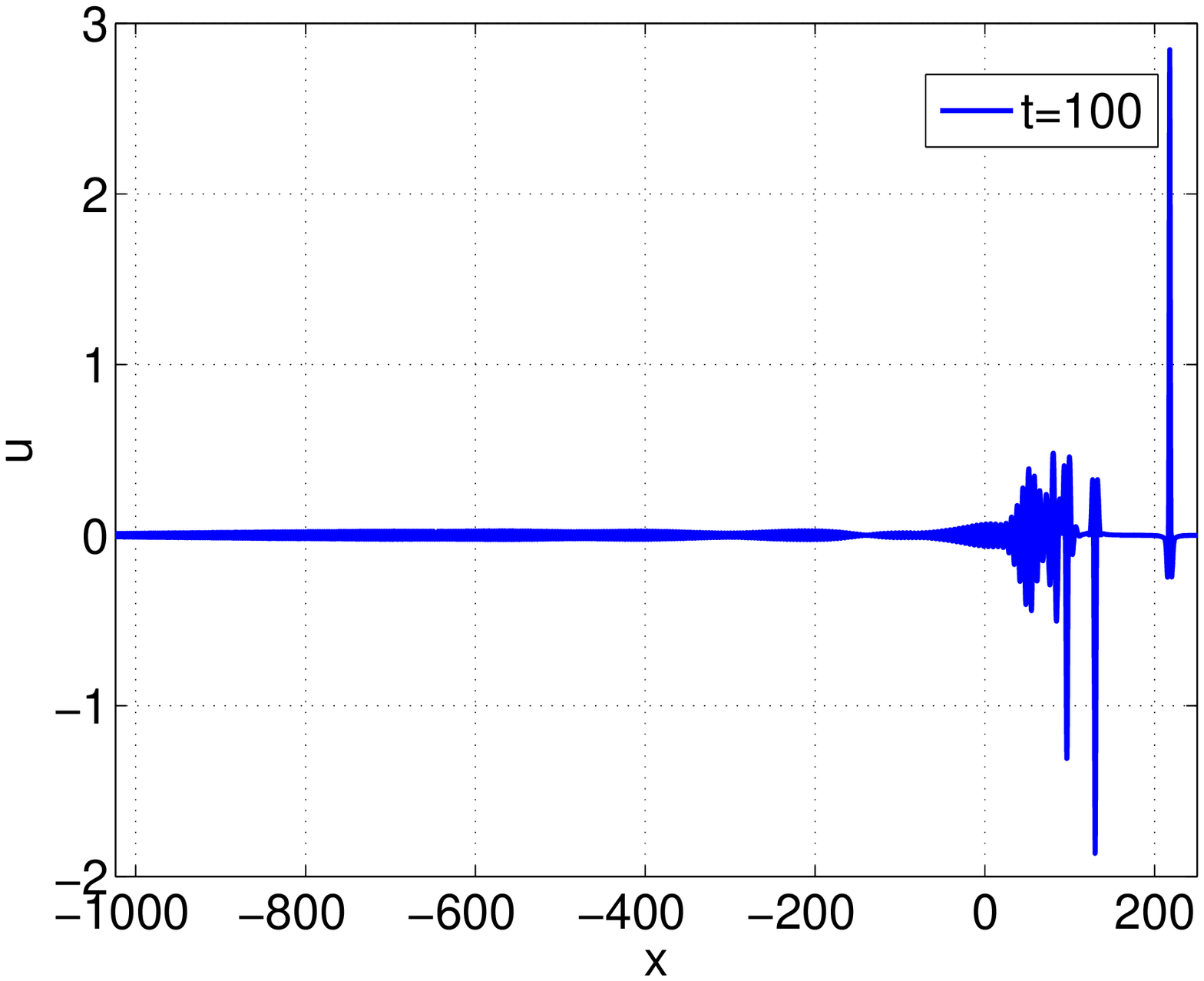}}
\caption{$r=1/2, m=1, q=2, c_{s}=0.75$,  $\gamma=1.7$, $A=4$. (a)-(f)  Numerical approximation at $t=0,20,40,60,80,100$.}
\label{gbenfig512_3}
\end{figure}

In some cases, increasing the amplitude perturbation parameter $A$ leads to the generation of more solitary pulses, i.~e. a resolution property. Our examples concern nonlinearities up to $q=3$. For the first group of experiments we have considered the gBenjamin equation with $A=4, q=2, c_{s}=0.75$ and three values of $\gamma=1.2, 1.6, 1.7$. The results are presented in Figures \ref{gbenfig512_2}-\ref{gbenfig512_3}. 
We look again for some connection between taking $\gamma$ closer to $\gamma_{max}(c_{s})\approx 1.732051$ and the resulting form of the approximations. 


When $\gamma=1.2$, the perturbation does not appear to generate more than one solitary pulse up to $t^{*}=100$ (Figures \ref{gbenfig512_2}, \ref{gbenfig512_2b}).
A different behaviour is observed in the other two cases. When $\gamma=1.6$  (Figures \ref{gbenfig512_1}, \ref{gbenfig512_1b}) the evolution of the perturbed initial profile leads to the formation of a main pulse followed by two additional presumed solitary waves (perhaps a two-pulse) of depression. A magnification of these structures can be seen in Figure \ref{gbenfig512_1m}. This pair is followed by a dispersive tail along with probably a nonlinear structure that may hide the formation of additional profiles. (The initial profile has amplitude $u_{m}\approx 9.365112E-01$ which gives $Au_{m}\approx 3.746045$ after perturbation. The main pulse has an amplitude $u_{m}\approx5.259017$ with a speed of $c_{s}\approx 5.795548$, see Figure \ref{gbenfig512_1b}. The following (two?-) pulse has a maximum negative excursion of about $-2.338071$.)

\begin{figure}[!htbp]
\centering
\subfigure[]
{\includegraphics[width=6cm]{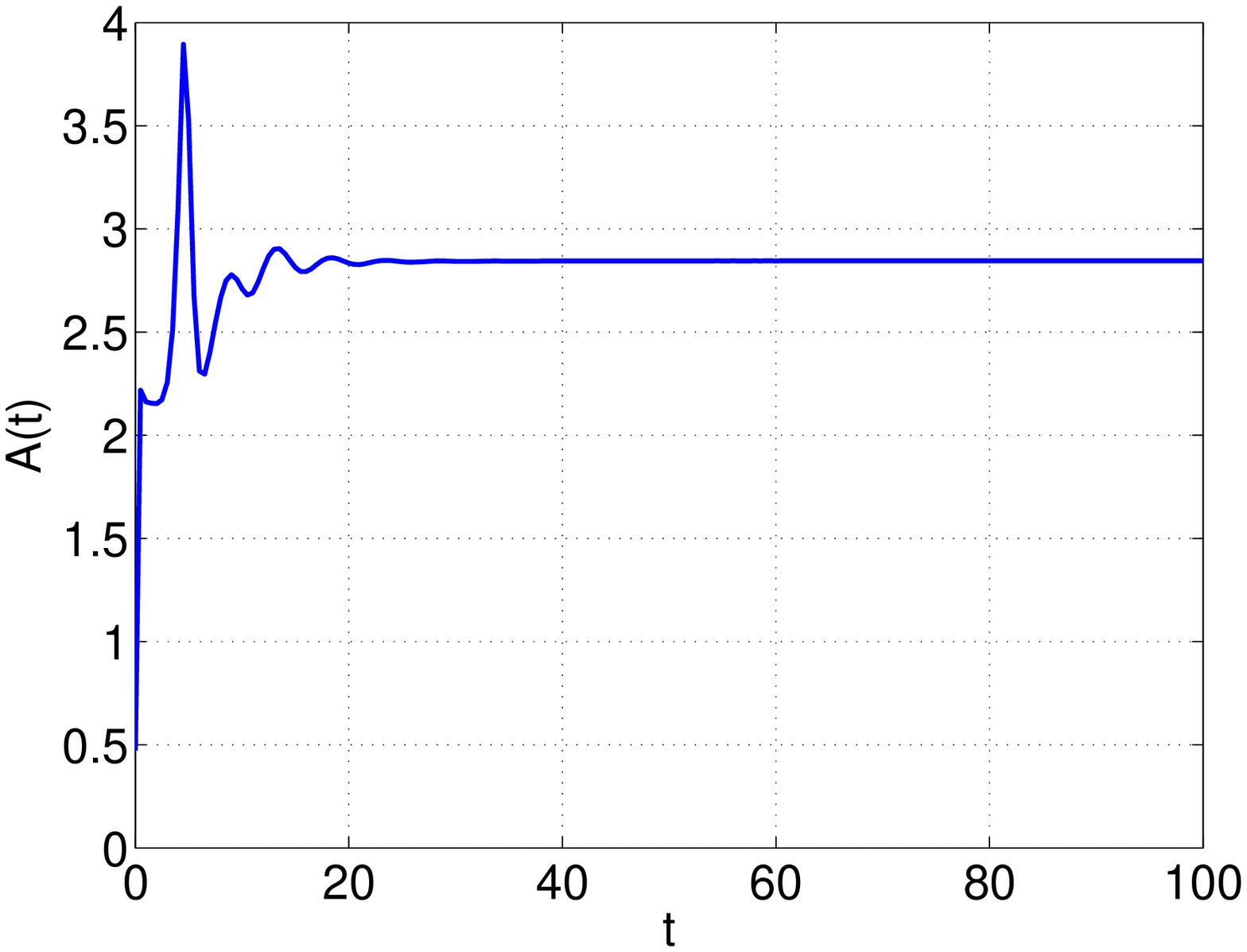}}
\subfigure[]
{\includegraphics[width=6cm]{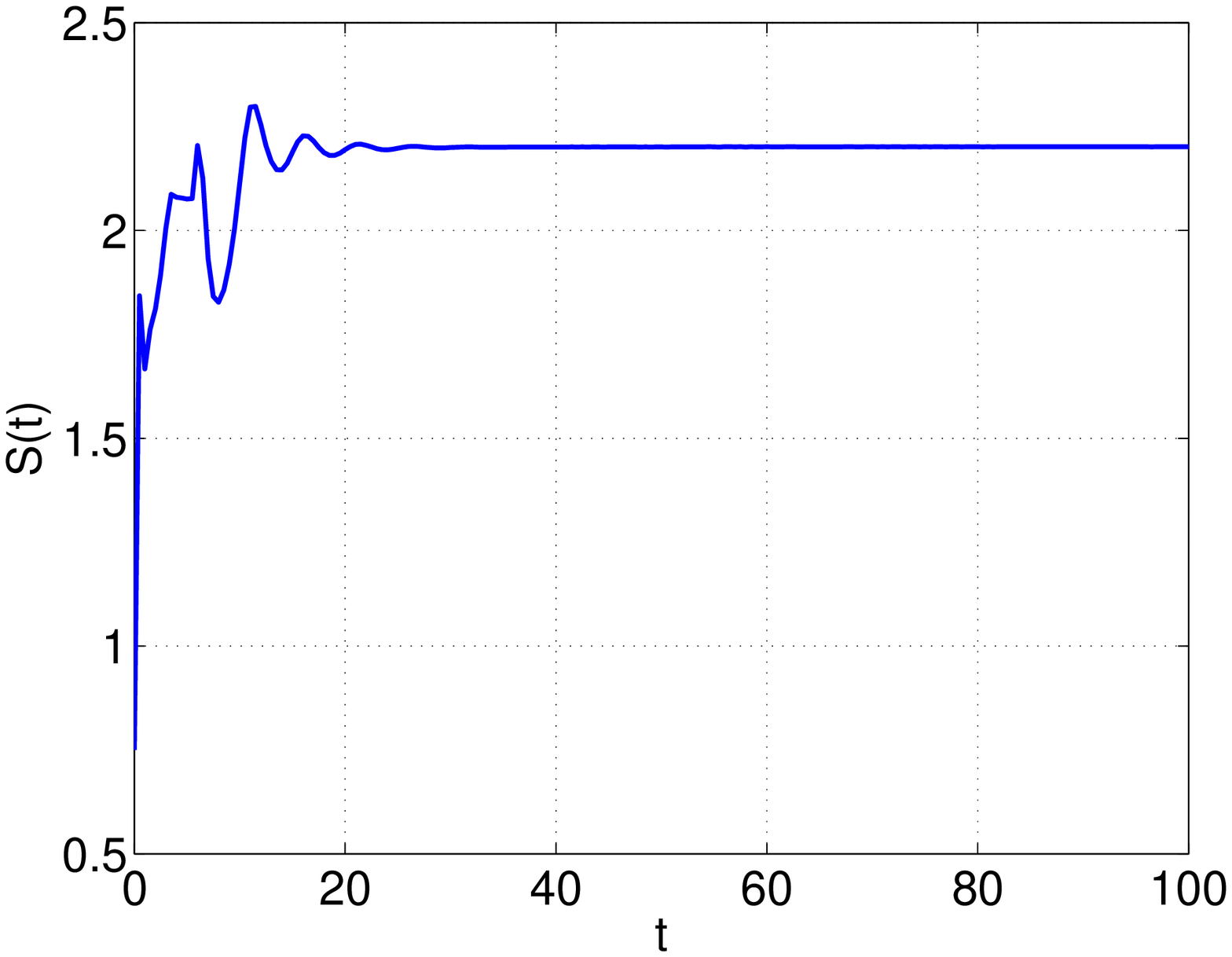}}
\caption{$r=1/2, m=1, q=2, c_{s}=0.75$,  $\gamma=1.7$, $A=4$. (a)  Evolution of the amplitude of the main numerical pulse; (b) Evolution of speed.}
\label{gbenfig512_3b}
\end{figure}

In the case $\gamma=1.7$ (Figures \ref{gbenfig512_3}, \ref{gbenfig512_3b}) the initial profile has amplitude $u_{m}\approx 4.716103E-01$ ($Au_{m}\approx 1.886441$). By $t=100$ two clearly separated solitary pulses emerge with amplitudes $u_{m}\approx 2.846051$ and (negative) $u_{m}\approx -1.884761$. The form of the tail behind the second wave suggests perhaps the emergence of another solitary wave at a larger time. The speed of the main pulse is $c_{s}\approx 2.200950$. We conclude that the resolution property for large perturbations in the amplitude seems to be related to the amplitude of the initial profile, that is how small $\gamma$ is, or, equivalently, how close to the gKdV the gBenjamin equation is.

\begin{figure}[!htbp]
\centering
\subfigure[]
{\includegraphics[width=6cm]{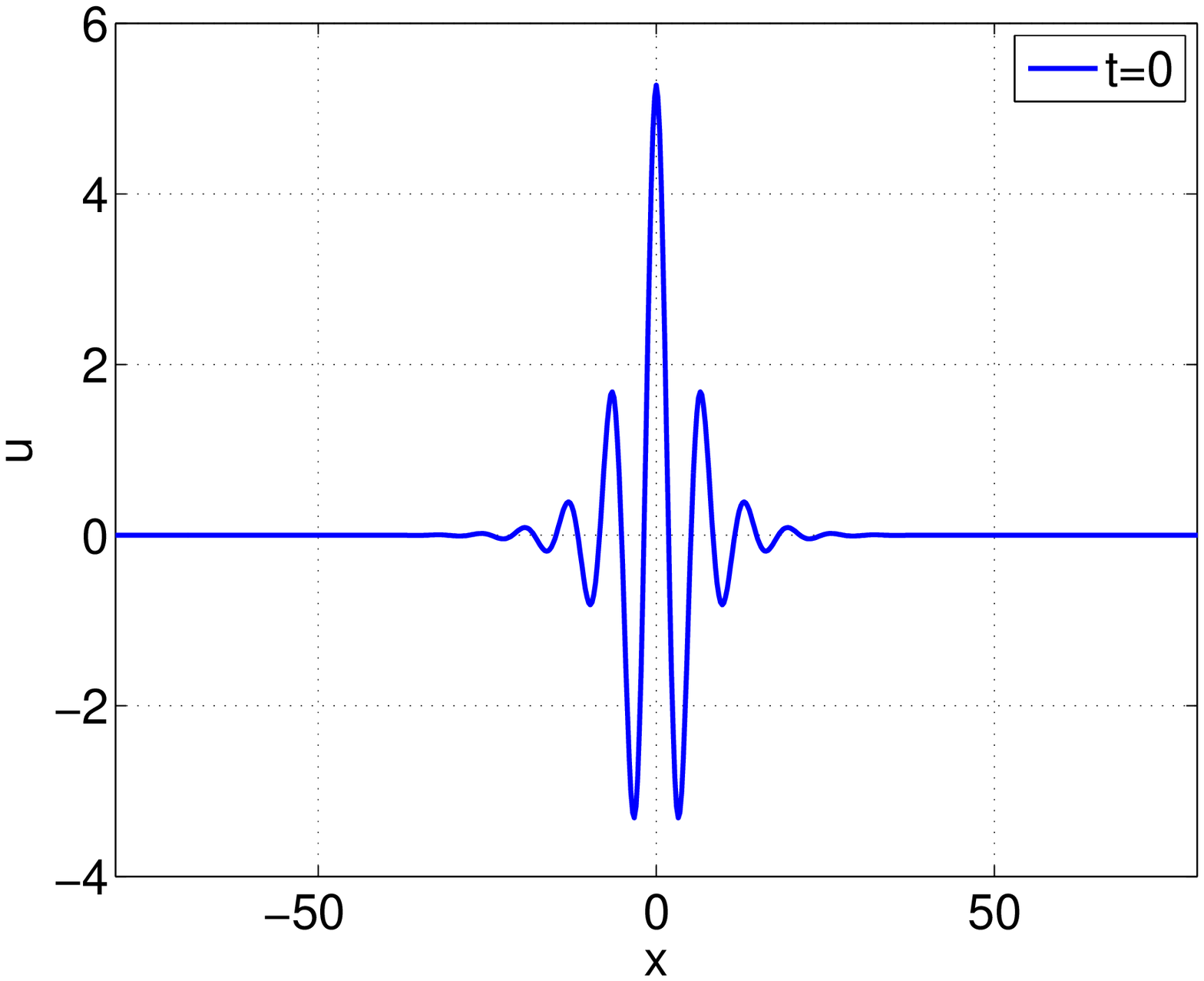}}
\subfigure[]
{\includegraphics[width=6cm]{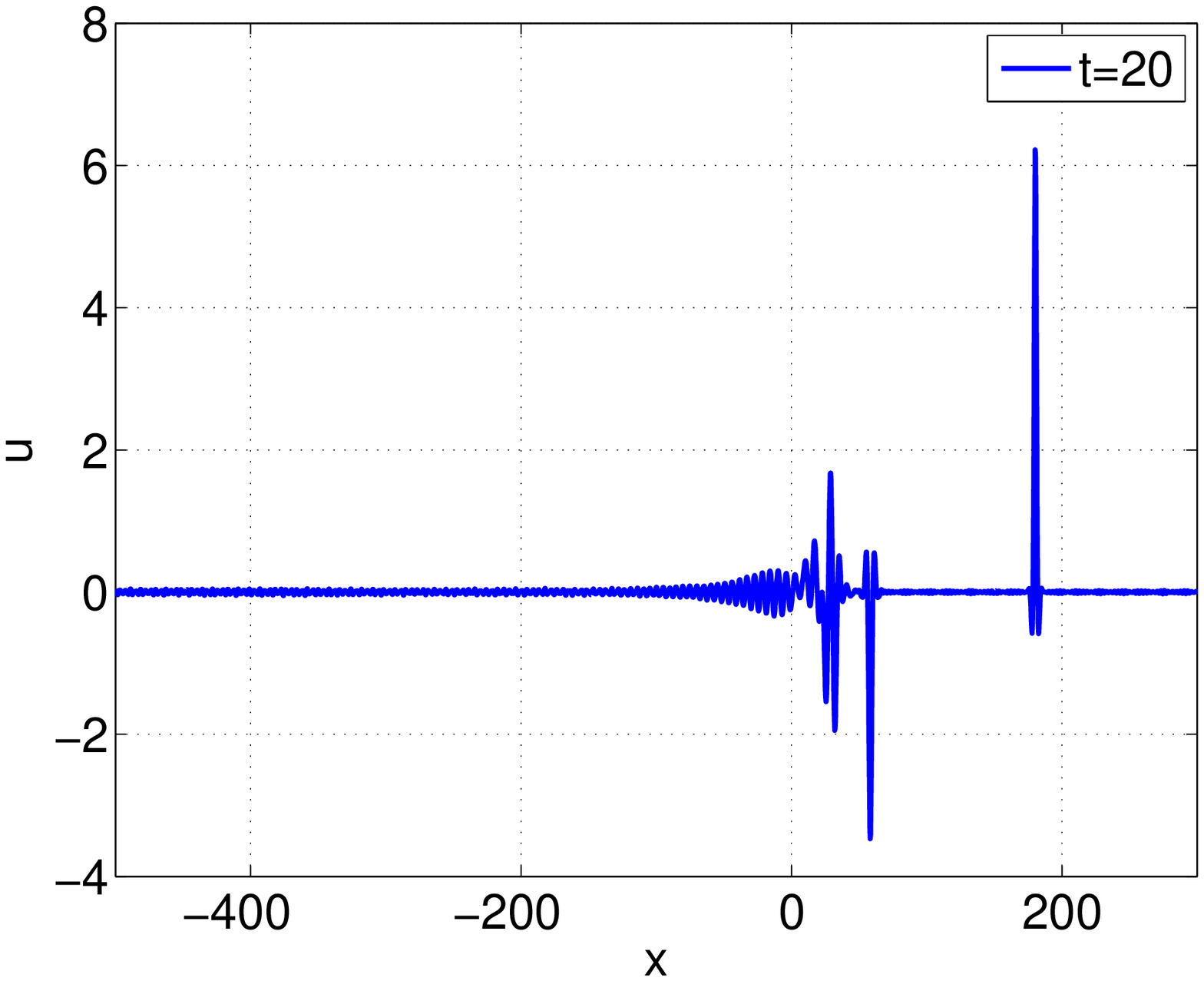}}
\subfigure[]
{\includegraphics[width=6cm]{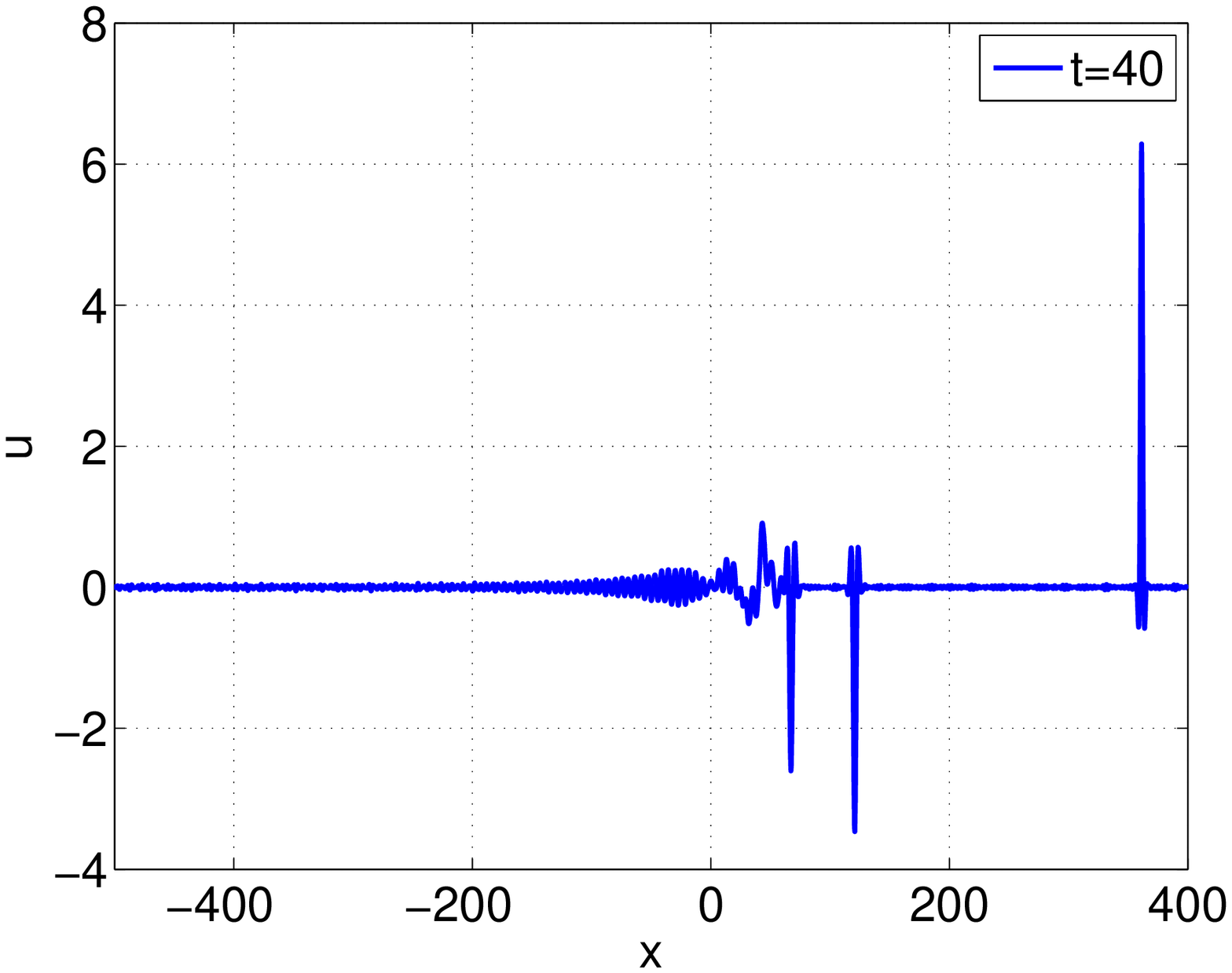}}
\subfigure[]
{\includegraphics[width=6cm]{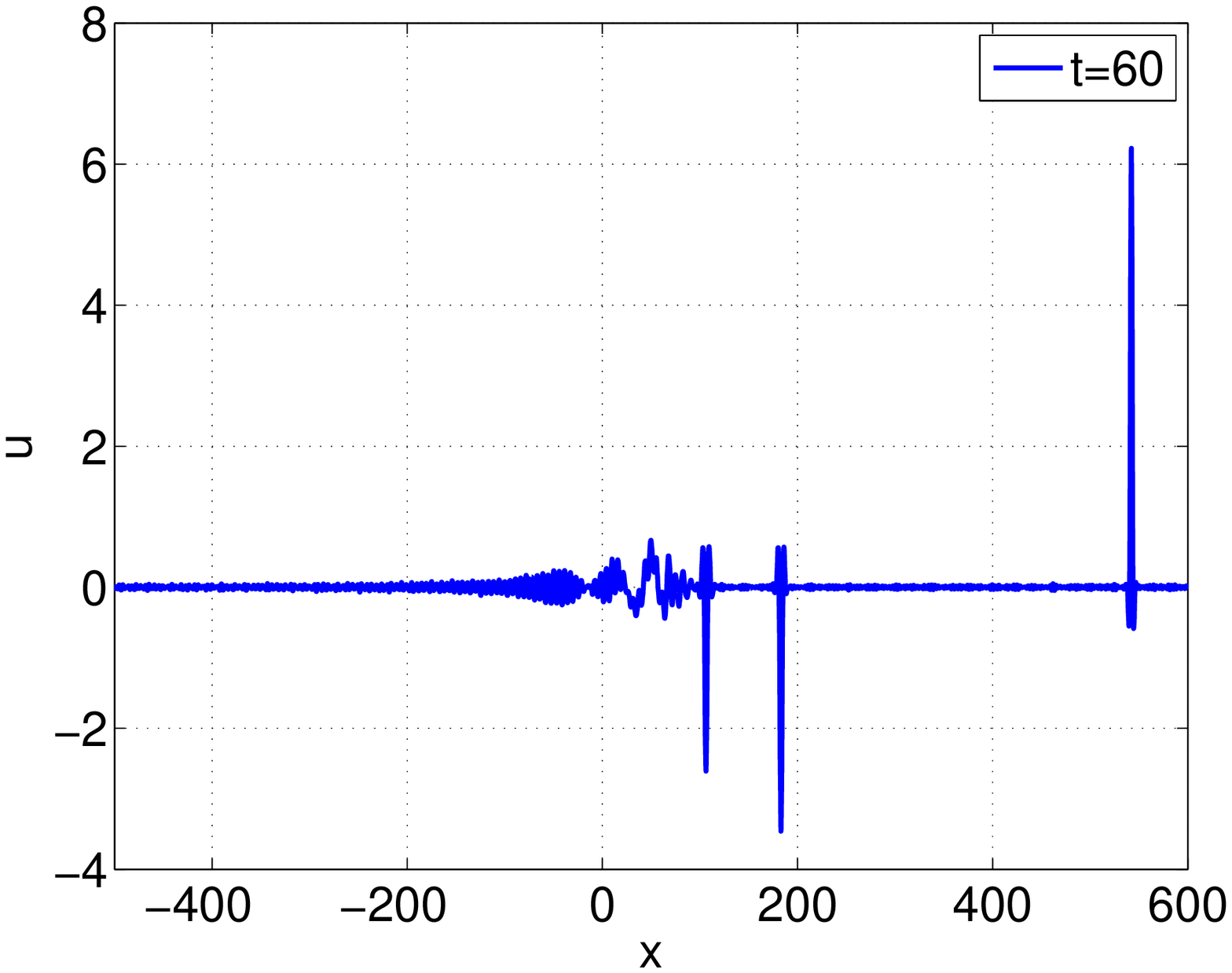}}
\subfigure[]
{\includegraphics[width=6cm]{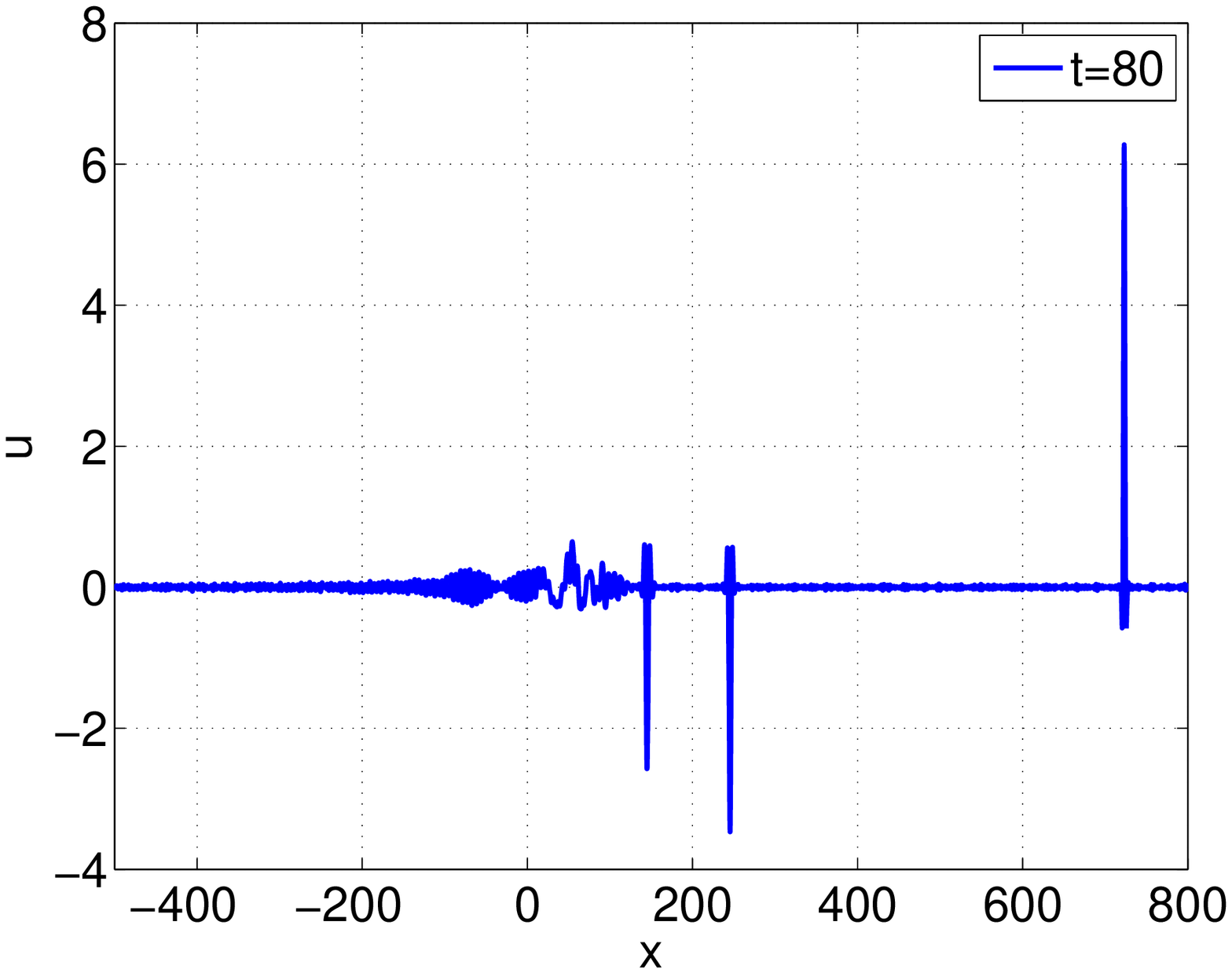}}
\subfigure[]
{\includegraphics[width=6cm]{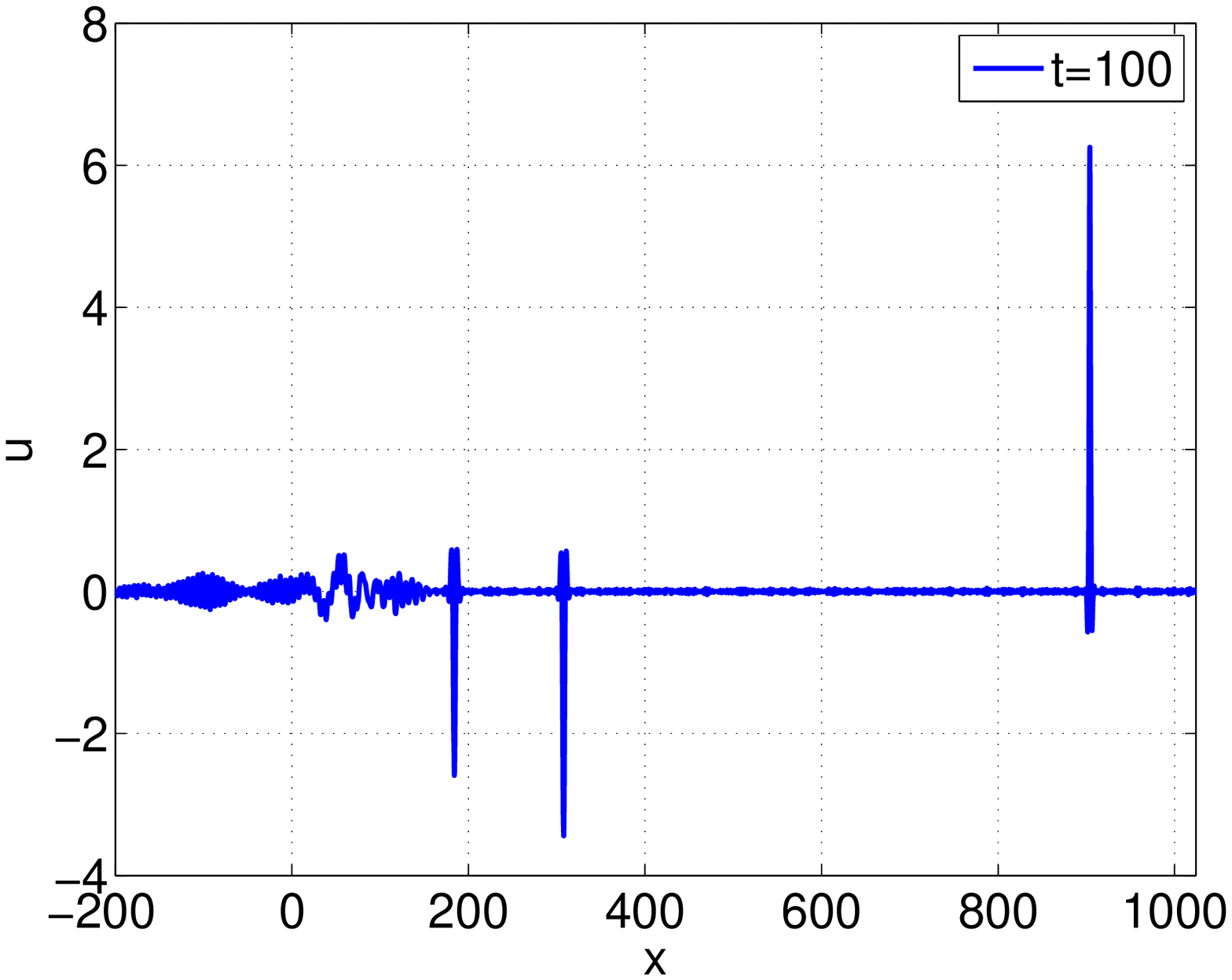}}
\subfigure[]
{\includegraphics[width=6cm]{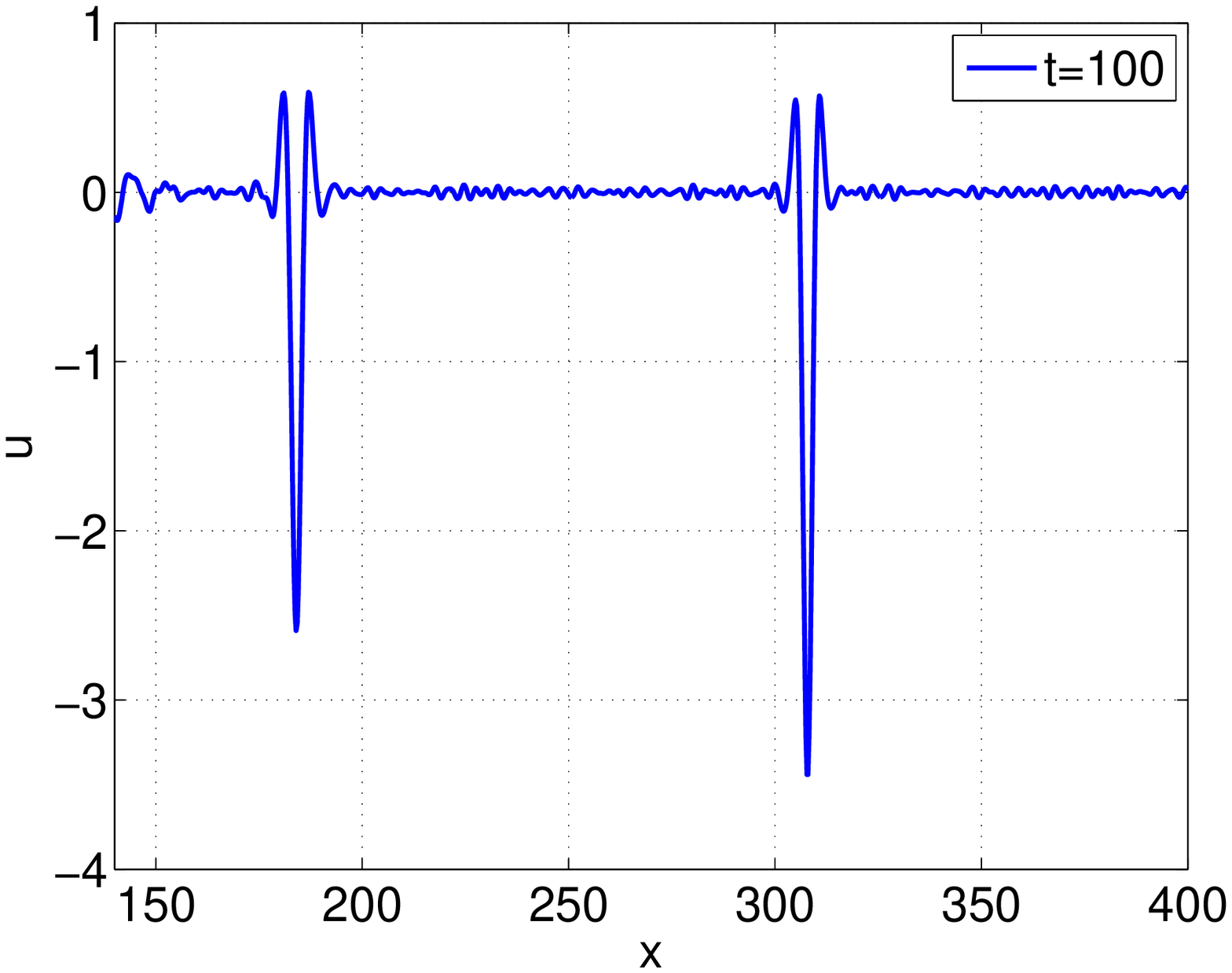}}
\caption{$r=1, m=2, q=2, c_{s}=1.01$,  $\gamma=1.8$, $A=4$. (a)-(f)  Numerical approximation at $t=0,20,40,60,80,100$. (g) Magnification of (f).}
\label{gbenfig513_2bb}
\end{figure}
\begin{figure}[!htbp]
\subfigure[]
{\includegraphics[width=6cm]{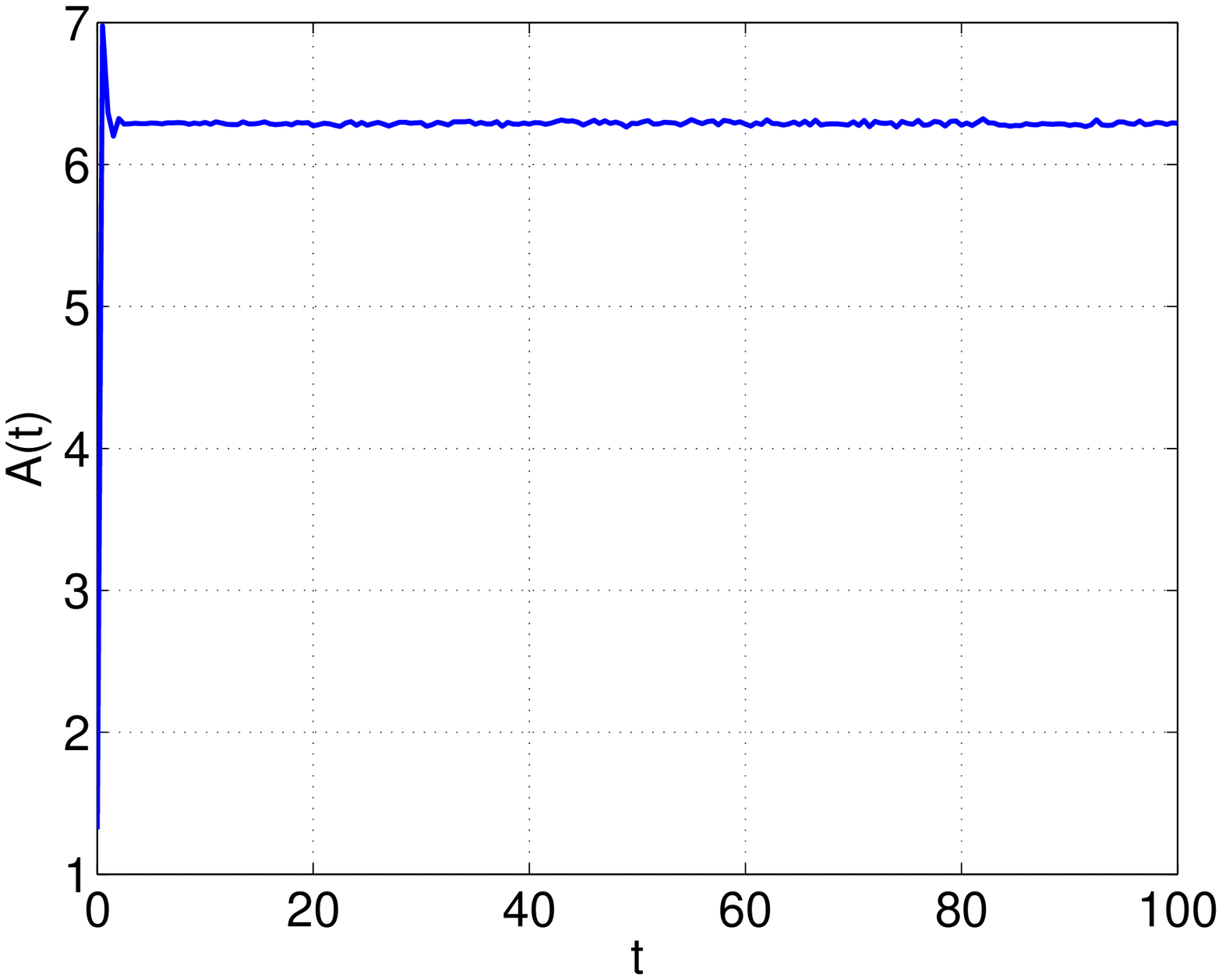}}
\subfigure[]
{\includegraphics[width=6cm]{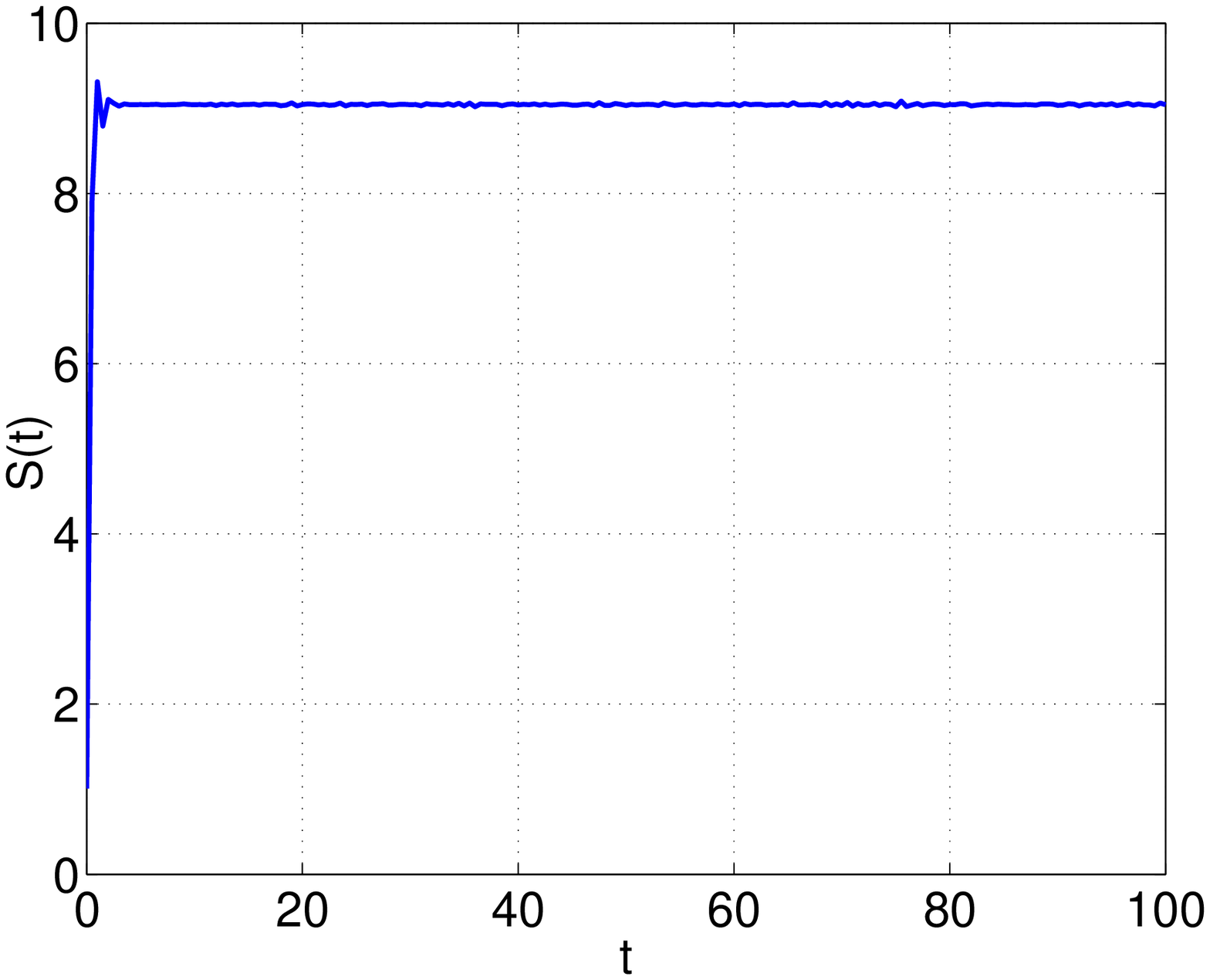}}
\caption{$r=1, m=2, q=2, c_{s}=1.01$,  $\gamma=1.8$, $A=4$. (a) Evolution of the amplitude of the main numerical pulse; (b) Evolution of speed.}
\label{gbenfig513_2b}
\end{figure}
\begin{figure}[htbp]
\centering
\subfigure[]
{\includegraphics[width=6cm]{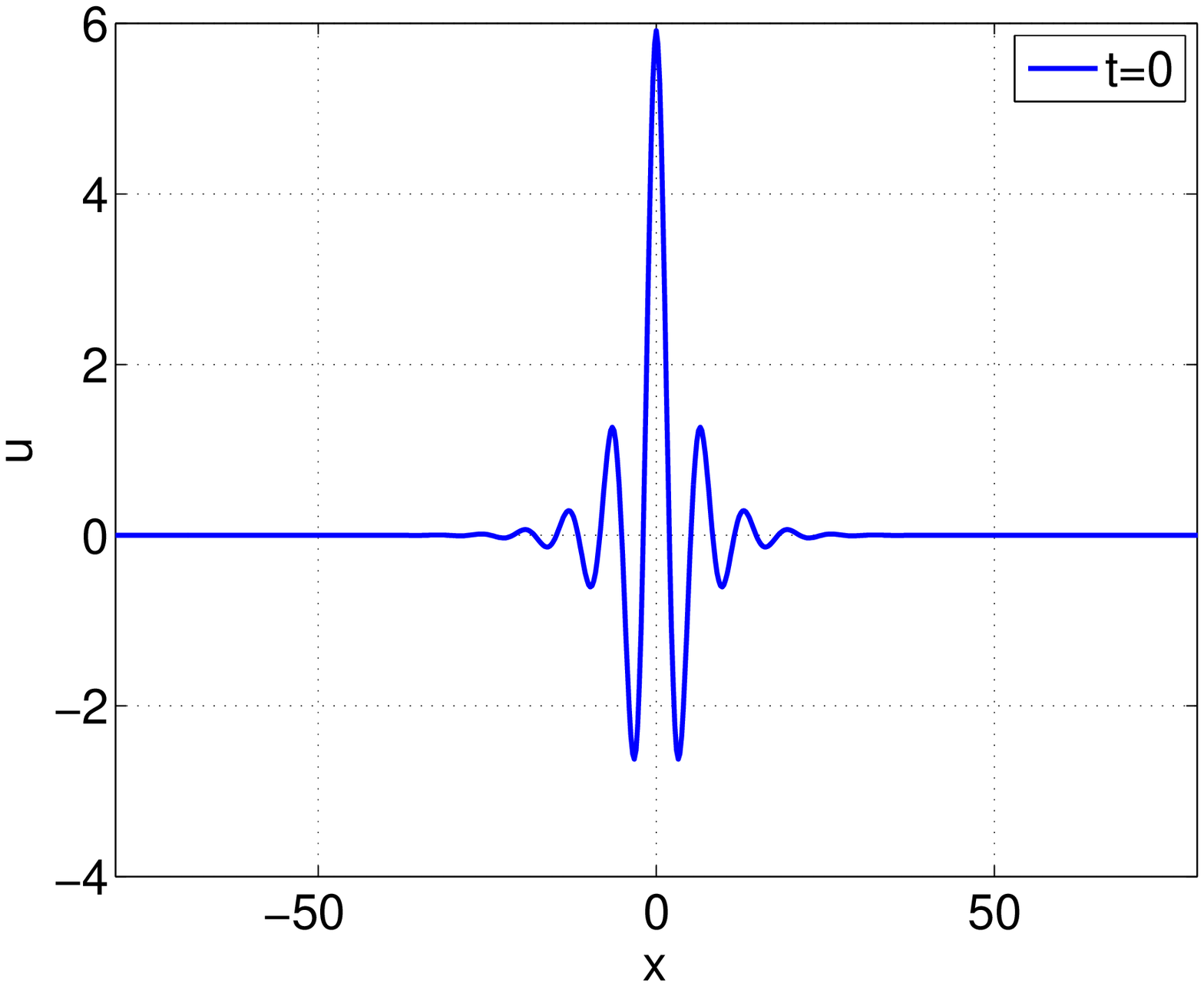}}
\subfigure[]
{\includegraphics[width=6cm]{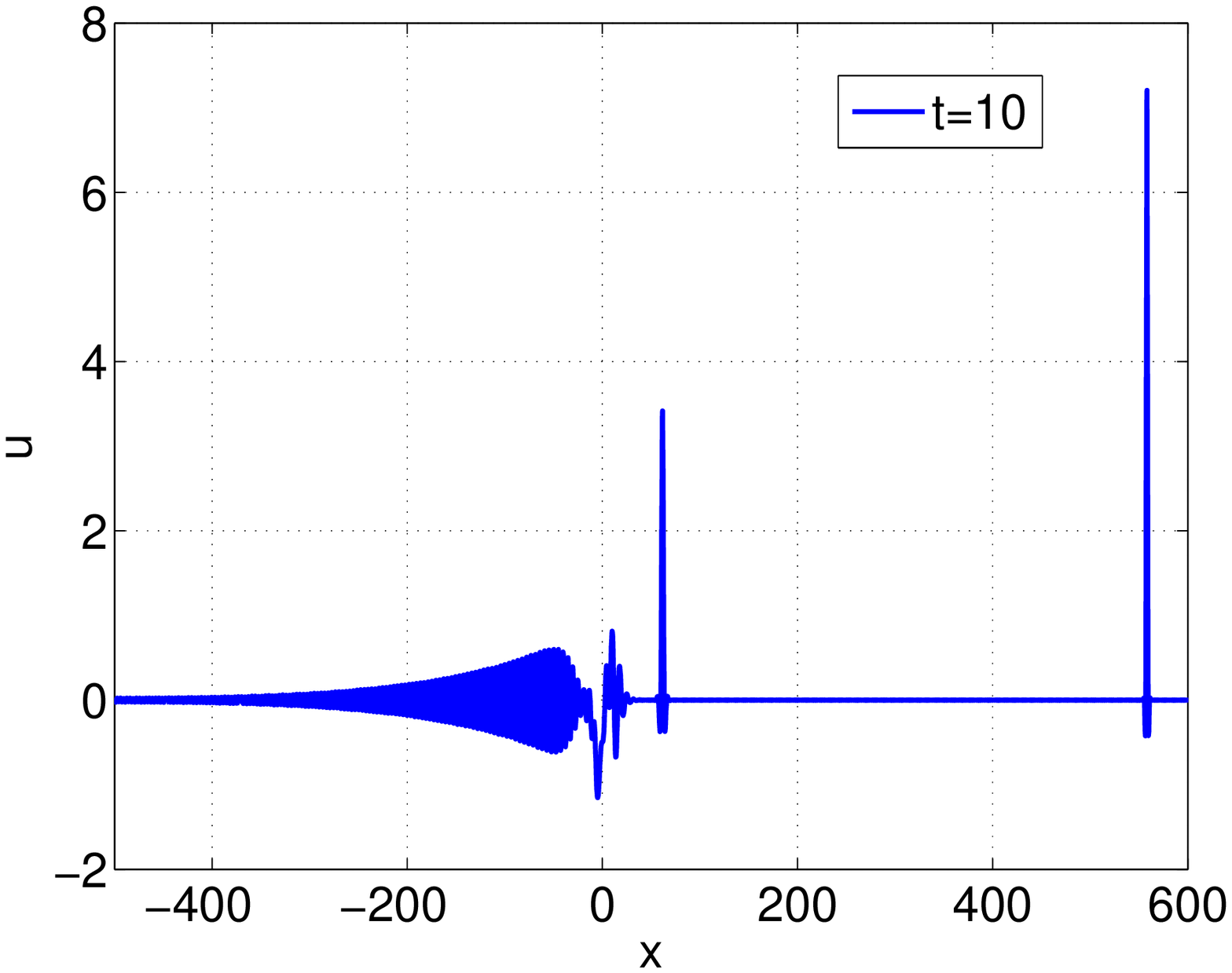}}
\subfigure[]
{\includegraphics[width=6cm]{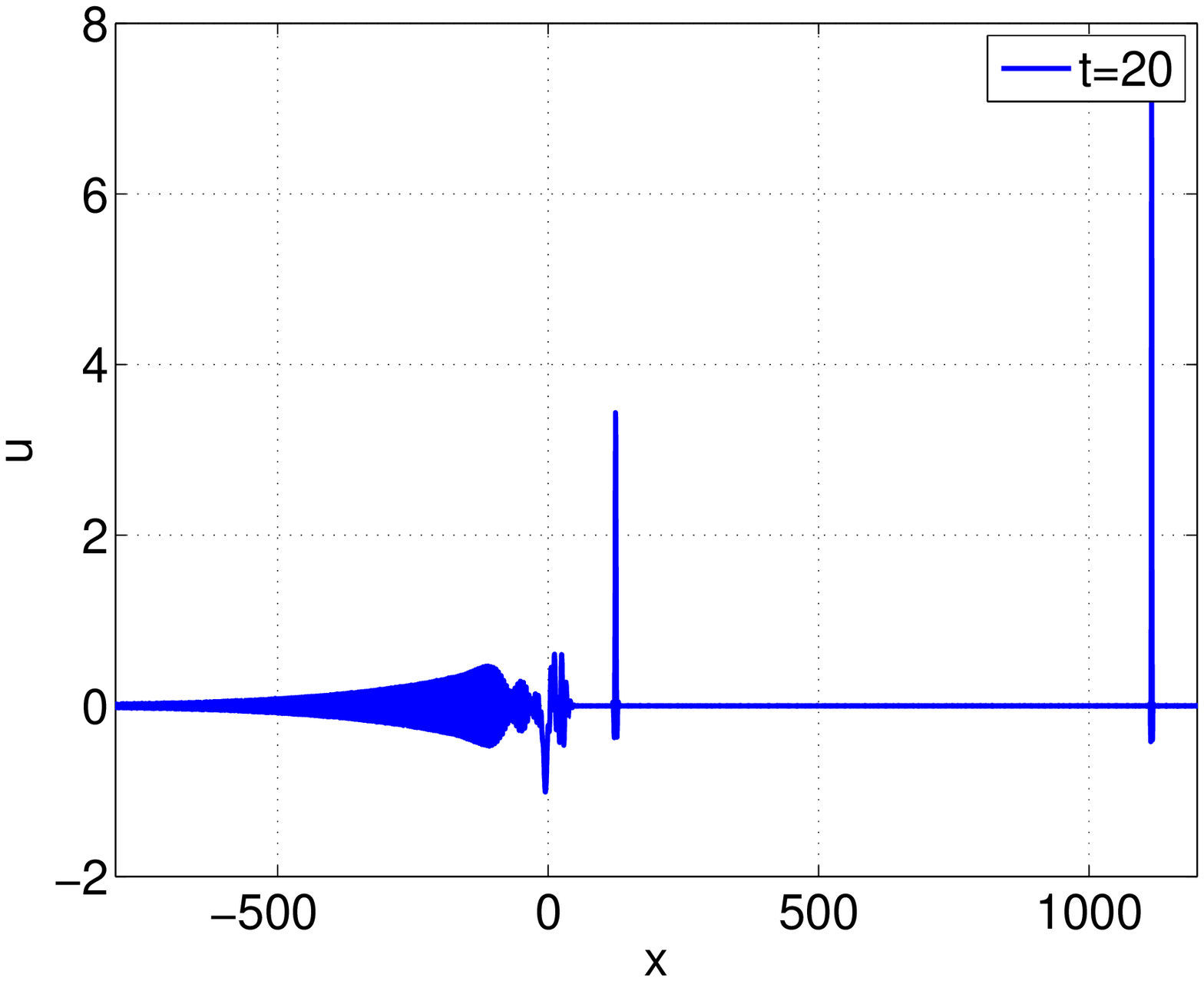}}
\subfigure[]
{\includegraphics[width=6cm]{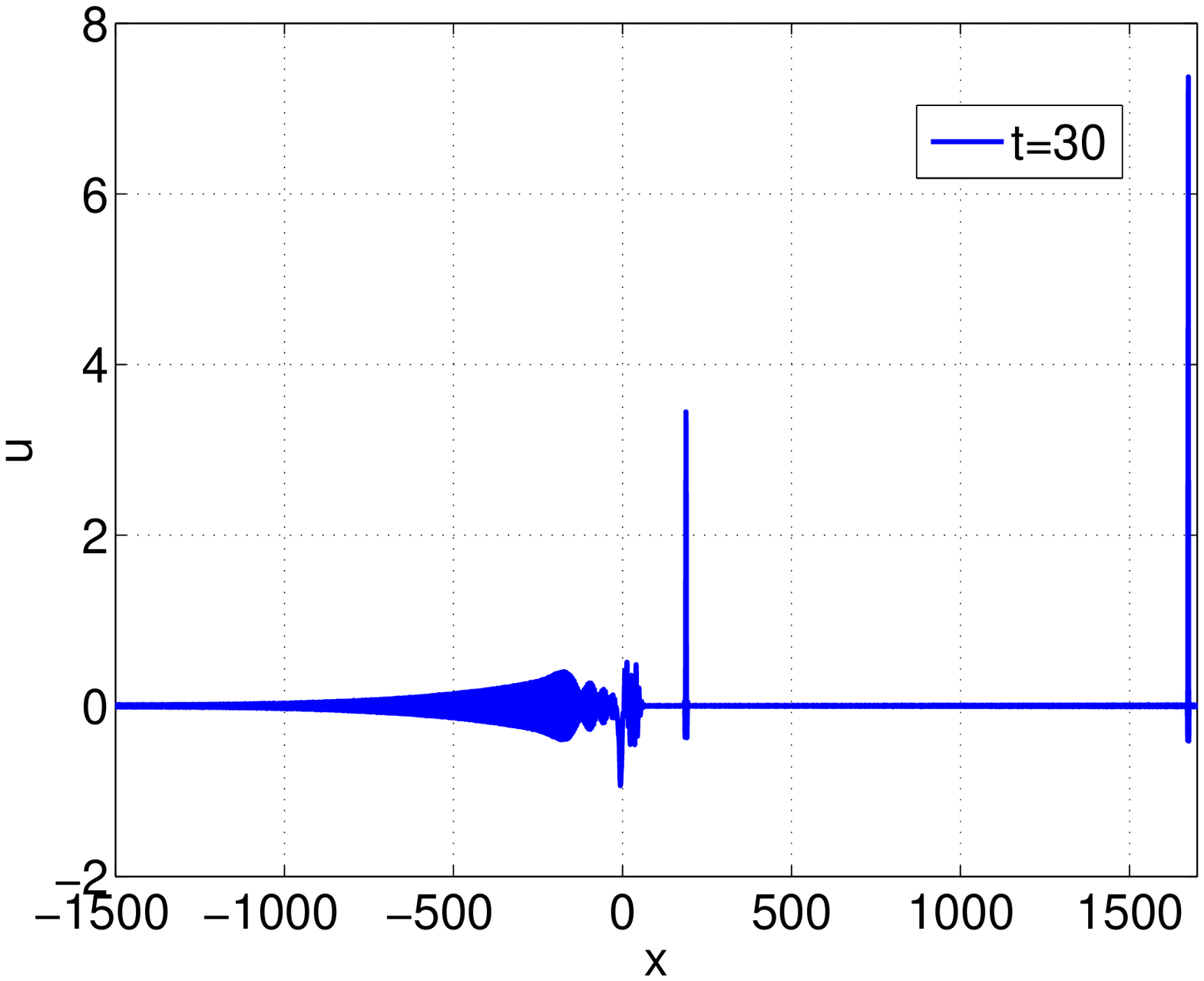}}
\subfigure[]
{\includegraphics[width=6cm]{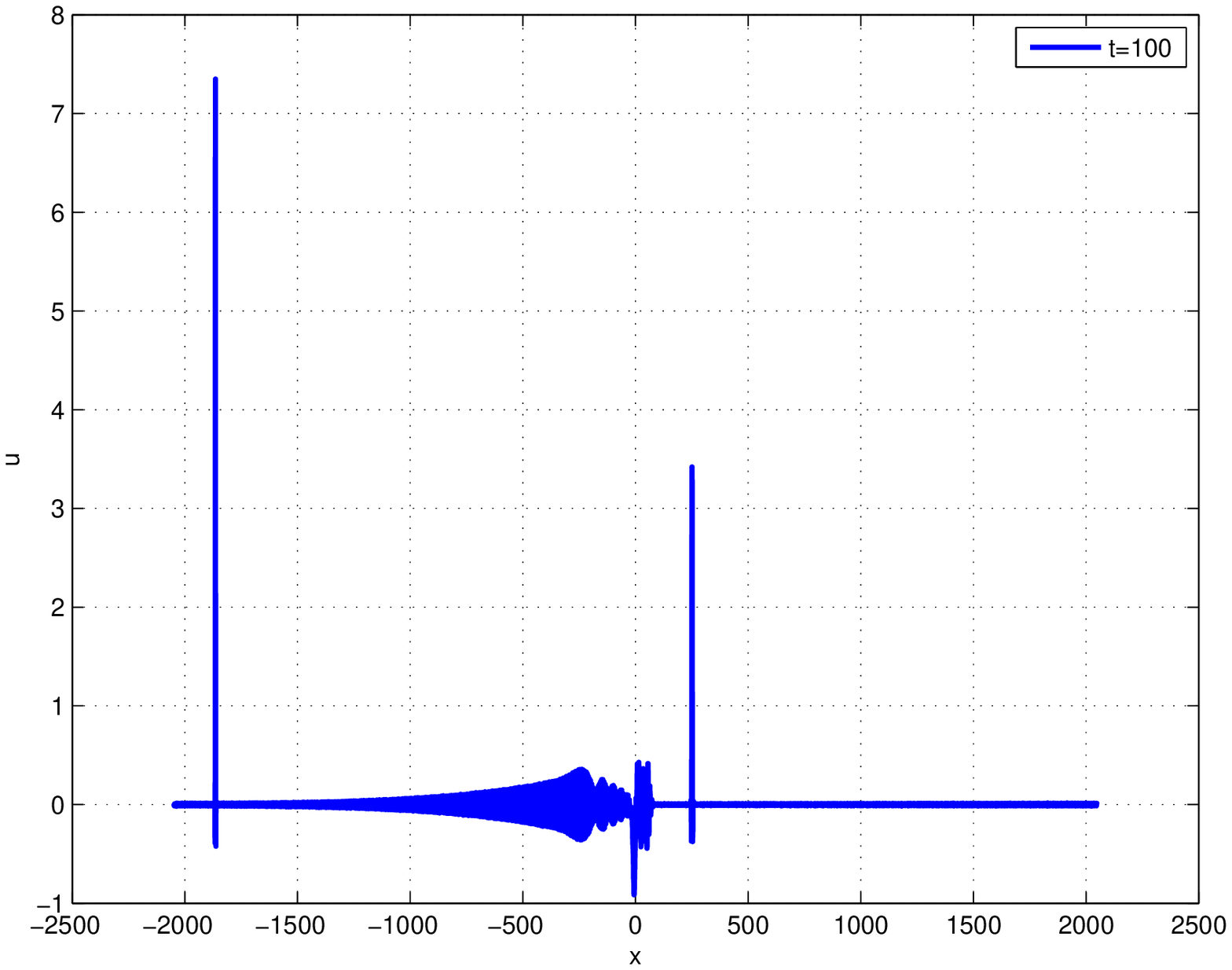}}
\subfigure[]
{\includegraphics[width=6cm]{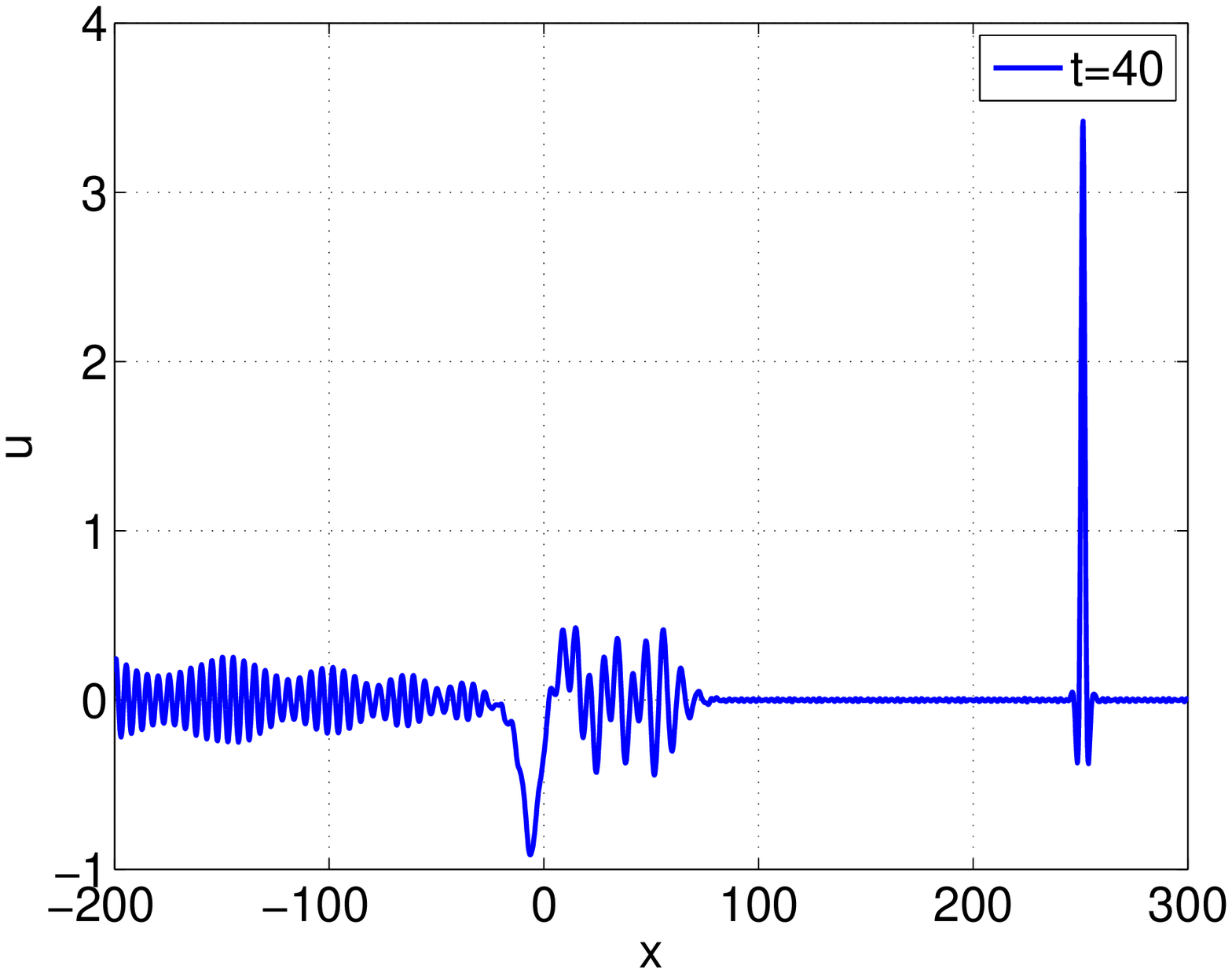}}
\subfigure[]
{\includegraphics[height=4cm,width=6.3cm]{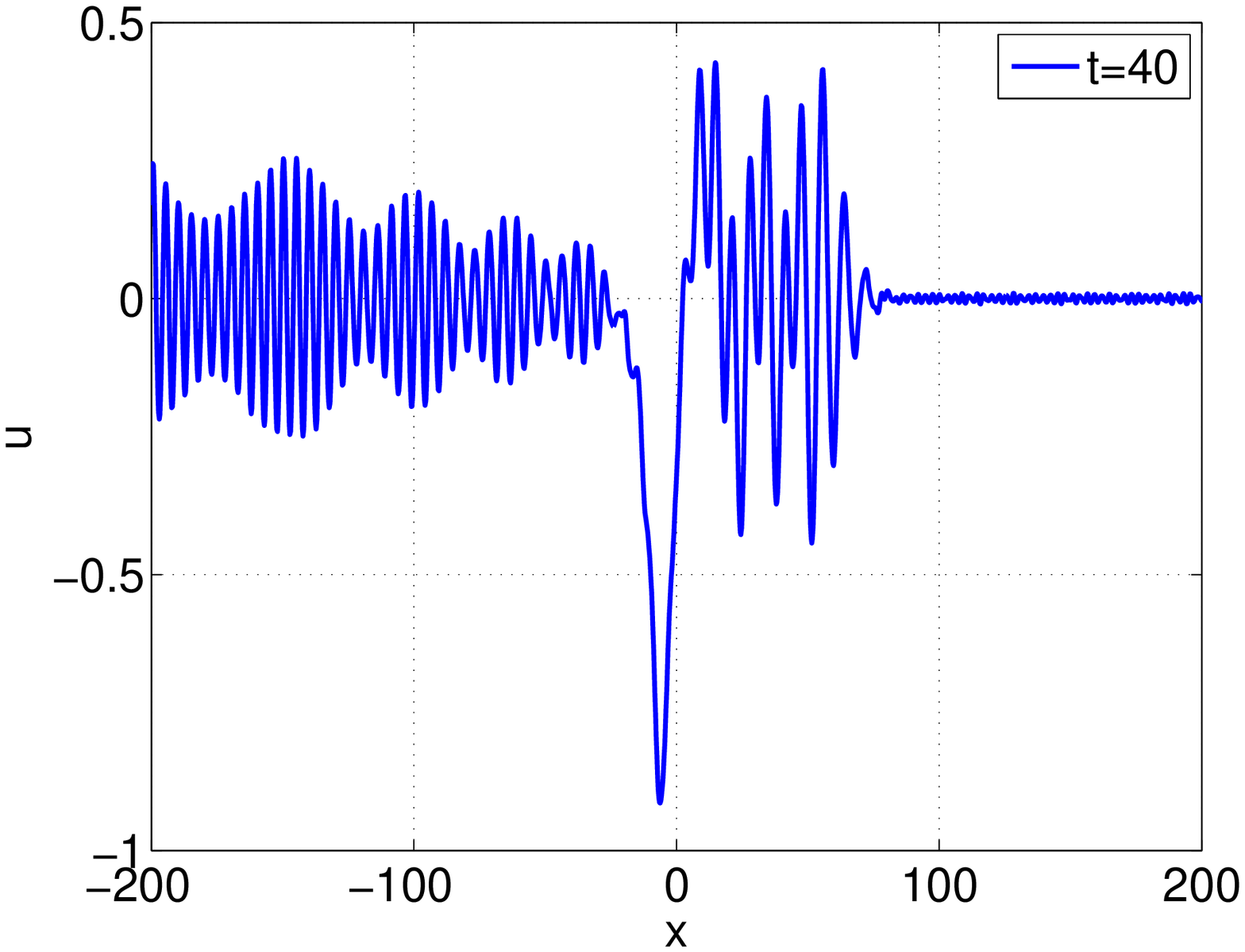}}
\subfigure[]
{\includegraphics[height=4cm,width=6cm]{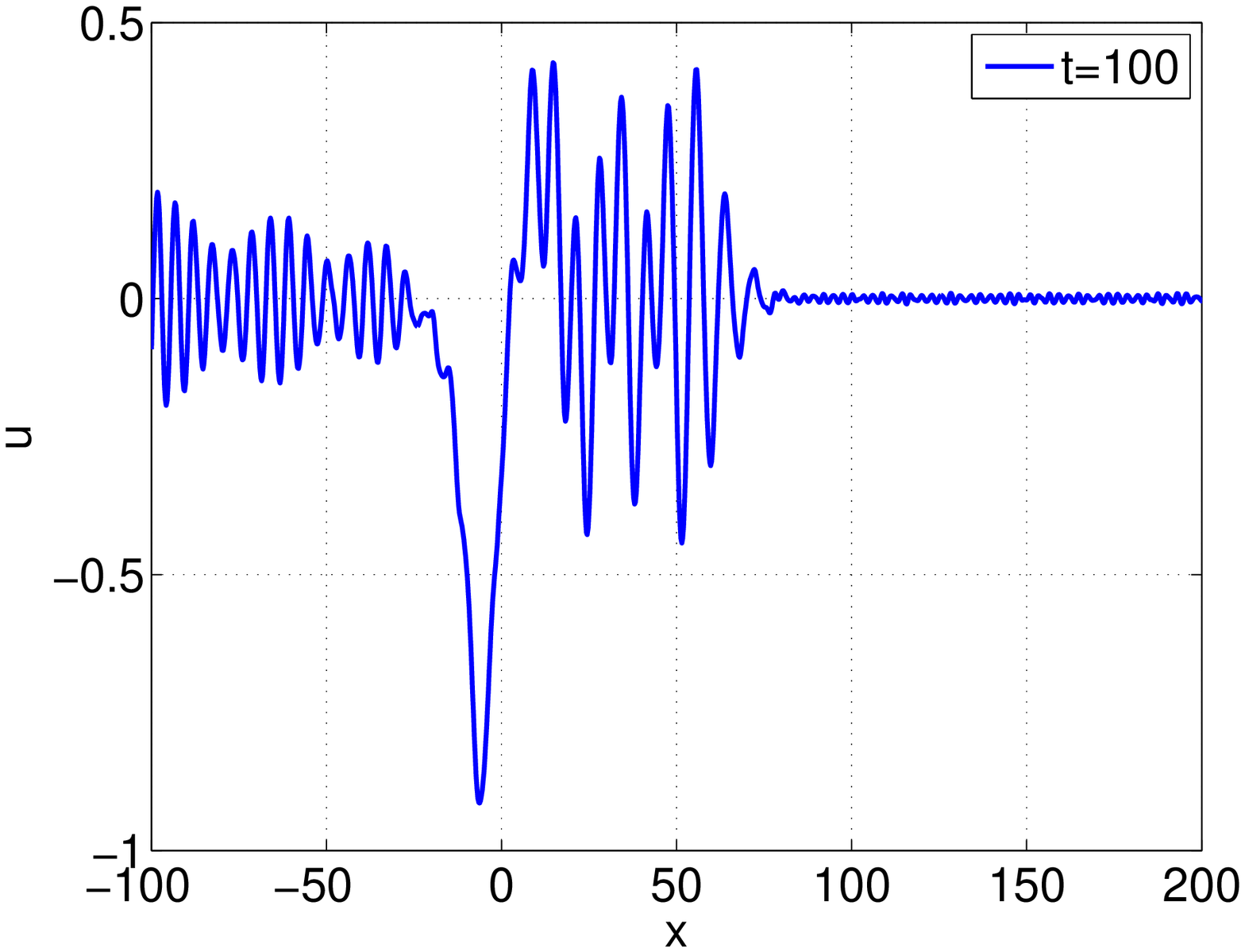}}
\caption{$r=1, m=2, q=3, c_{s}=1.01$,  $\gamma=1.8$, $A=4$. (a)-(e)  Numerical approximation at $t=0,10,20,30,100$; (f) is a magnification of the second pulse at $t=40$; (g) and (h) are magnifications of the tail behind the second pulse at $t=40$ and $t=100$ respectively.}
\label{gbenfig513_4}
\end{figure}
\begin{figure}[htbp]
\centering
\subfigure[]
{\includegraphics[width=6cm]{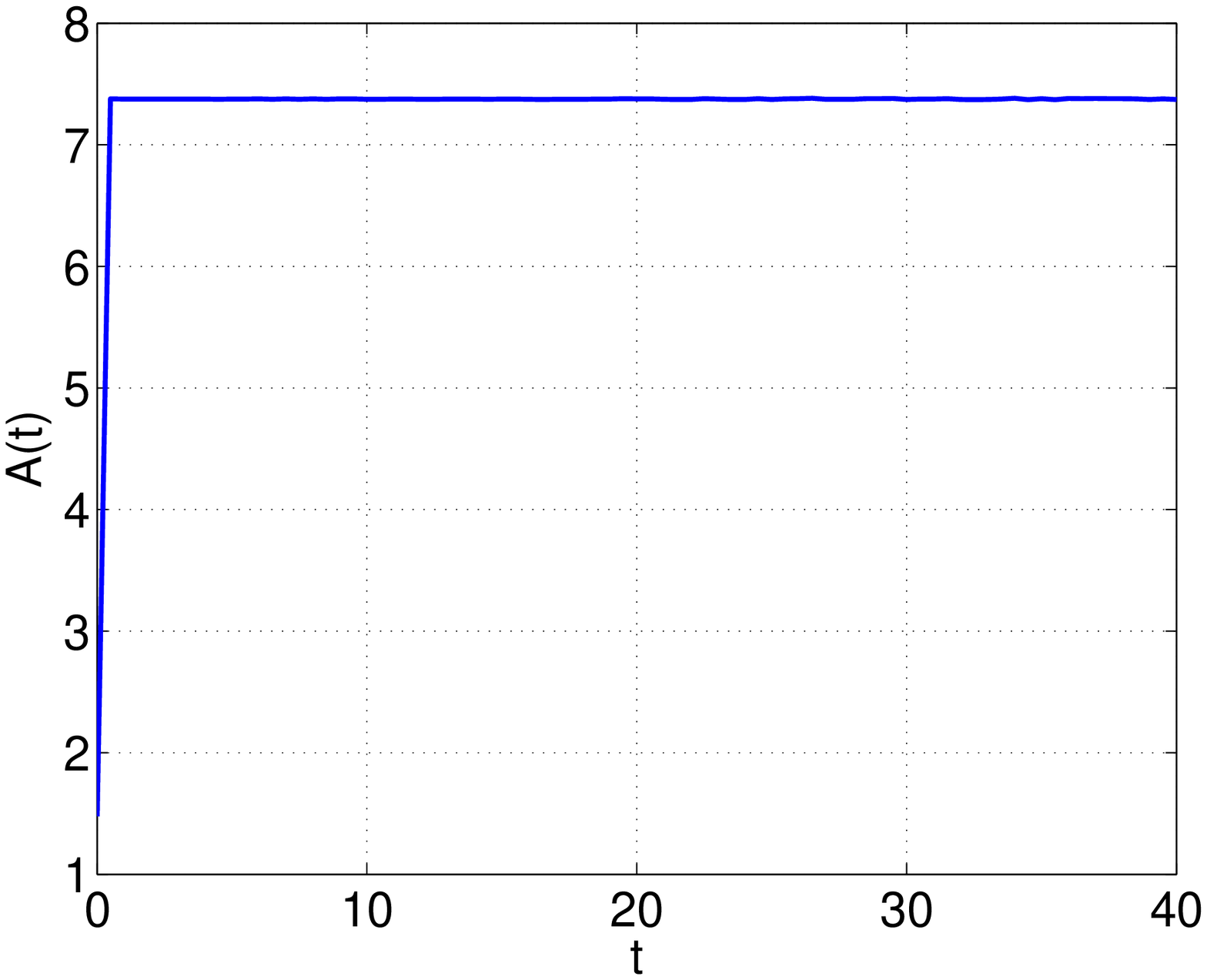}}
\subfigure[]
{\includegraphics[width=6cm]{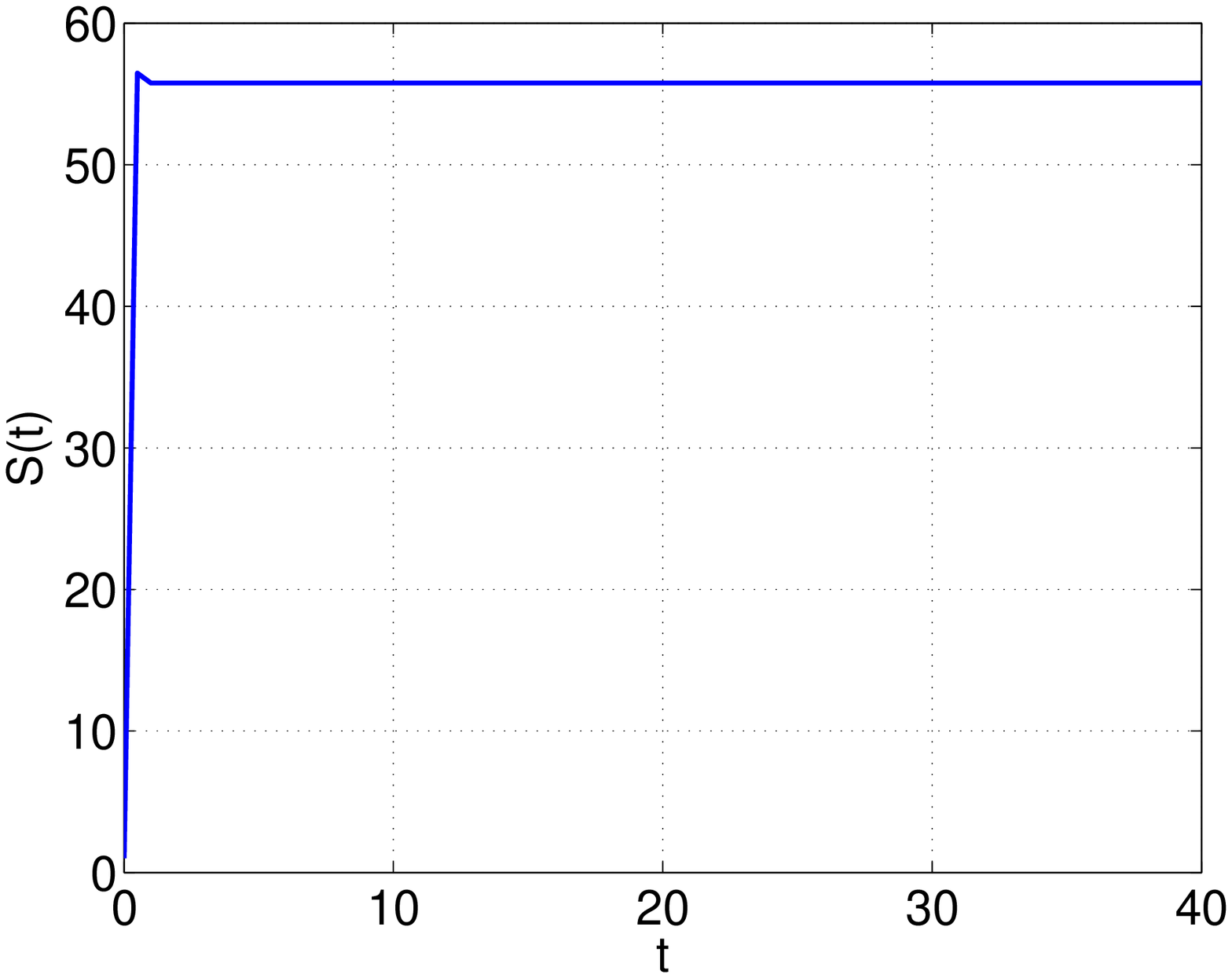}}
\caption{$r=1, m=2, q=3, c_{s}=1.01$,  $\gamma=1.8$, $A=4$.  (a) Evolution of the amplitude of the main numerical pulse; (b) Evolution of speed.}
\label{gbenfig513_4b}
\end{figure}

Now a few observations on the general case. This is illustrated with the values $r=1, m=2$ with $c_{s}=1.01, \delta =1, \gamma=1.8$ and $q=2, 3$ (Figures \ref{gbenfig513_2bb}, \ref{gbenfig513_2b} and \ref{gbenfig513_4}, \ref{gbenfig513_4b} respectively). The behaviour looks similar to that of the gBenjamin equation.
%
In the case $q=2$ two clear pulses of depression are generated behind the main pulse and in front of a tail with probably linear and nonlinear structures. The perturbed initial profile has amplitude $Au_{m}\approx 5.274788$ ($u_{m}\approx 1.318697$) and a maximum negative excursion of about $-3.312797$ (without the perturbation factor $A$ this is $-8.281993E-01$ approximately). The main wave has amplitude $u_{m}\approx 6.289751$ and speed $c_{s}\approx 9.044469$. The maximum negative excursion of the second pulse is about $-3.440279$ and of the third pulse about $-2.589379$.


The evolution in the case $q=3$ is generated from an initial solitary-wave profile with amplitude $u_{m}\approx 1.478188$. Two solitary pulses of elevation are generated with approximate amplitudes $u_{m}\approx 7.376428$ and $u_{m}\approx 3.421021$. The speed of the first one is $c_{s}\approx 57.997822$ (cf. Figure \ref{gbenfig513_4b}) . The form of the tail behind the second pulse is shown magnified in Figures \ref{gbenfig513_4}(g), (h). The tail appears to be very slow, compared to the high-speed main pulse. The figures suggest a nonlinear structure in front with a chain of ripples behind. Some numerical or real instability cannot be discarded.


\subsection{Interactions of solitary waves}
\label{sec34}
As in the classical Benjamin equation, solitary-wave collisions are expected to be inelastic. We will illustrate this inelastic behaviour with some examples.

The results concerning the gBenjamin equation are collected into two groups:
\begin{itemize}
\item[(G1)] $\delta=1,\gamma=1.7$ and two initial profiles of speeds $c_{s}^{(1)}=2, c_{s}^{(2)}=1$ centered at $x_{0}^{(1)}=-50, x_{0}^{(2)}=0$. The final time is $t=200$ and three values $q=2,3,4$ have been considered. In this case, $\gamma_{max}(c_{s}^{(1)})\approx 2.828427$ and $\gamma_{max}(c_{s}^{(2)})=2$. The corresponding amplitudes before ($u_{b}$) and after ($u_{a}$) the interaction are displayed in Table \ref{gbentab541}.
\item[(G2)] $\delta=1,\gamma=1.999$ and two initial profiles of speeds $c_{s}^{(1)}=3, c_{s}^{(2)}=1$ centered at $x_{0}^{(1)}=-200, x_{0}^{(2)}=0$. The final time is $t=200$ and three values $q=2,3,4$ have also been considered. Now $\gamma_{max}(c_{s}^{(1)})\approx 3.464102$ and $\gamma$ is closer to  $\min\{\gamma_{max}(c_{s}^{(1)}), \gamma_{max}(c_{s}^{(2)})\}=2$. The corresponding amplitudes are also displayed in Table \ref{gbentab541}.
\end{itemize}

\begin{table}
\smallskip
\begin{center}
\begin{tabular}{|c|c|c|}
\hline
 &(G1)&(G2)\\
\hline\hline
$q=2$&$u_{b}^{(1)}\approx 2.656731$,$u_{a}^{(1)}\approx 2.666032$&$u_{b}^{(1)}\approx 3.312262$
,$u_{a}^{(1)}\approx 3.312365$\\
&$c_{sb}^{(1)}=2, c_{sa}^{(1)}\approx 2.009556$
&$c_{sb}^{(1)}=3, c_{s}^{(1)}\approx 3.000601$\\
&$u_{b}^{(2)}\approx 1.379312$,$u_{a}^{(2)}\approx 1.326514$&$u_{b}^{(2)}\approx 8.944152E-02$
,$u_{a}^{(2)}\approx 8.940440E-02$\\
\hline
$q=3$&$u_{b}^{(1)}\approx 2.233158$,$u_{a}^{(1)}\approx 2.757020$&$u_{b}^{(1)}\approx 2.589268$
,$u_{a}^{(1)}\approx 2.596668$\\
&$c_{sb}^{(1)}=2, c_{sa}^{(1)}\approx 2.043637$
&$c_{sb}^{(1)}=3, c_{s}^{(1)}\approx 3.026974$\\
&$u_{b}^{(2)}\approx 1.437308$,$u_{a}^{(2)}\approx 1.361308$&$u_{b}^{(2)}\approx 5.784419E-01$
,$u_{a}^{(2)}\approx 5.714085E-01$\\
\hline
$q=4$&$u_{b}^{(1)}\approx 2.001032$,$u_{a}^{(1)}\approx 1.525542 $&\\
&$c_{sb}^{(1)}=2, c_{sa}^{(1)}\approx1.104477$
&instability\\
&$u_{b}^{(2)}\approx 1.427771$,$u_{a}^{(2)}$: instability&\\
\hline
\end{tabular}
\end{center}
\caption{Interactions of solitary waves, gBenjamin equation. Amplitudes before ($u_{b}$) and after ($u_{a}$) the collision; speeds ($c_{sb}, c_{sa}$) of the tallest wave.}\label{gbentab541}
\end{table}



\begin{figure}[htbp]
\centering
\subfigure[]
{\includegraphics[width=6cm]{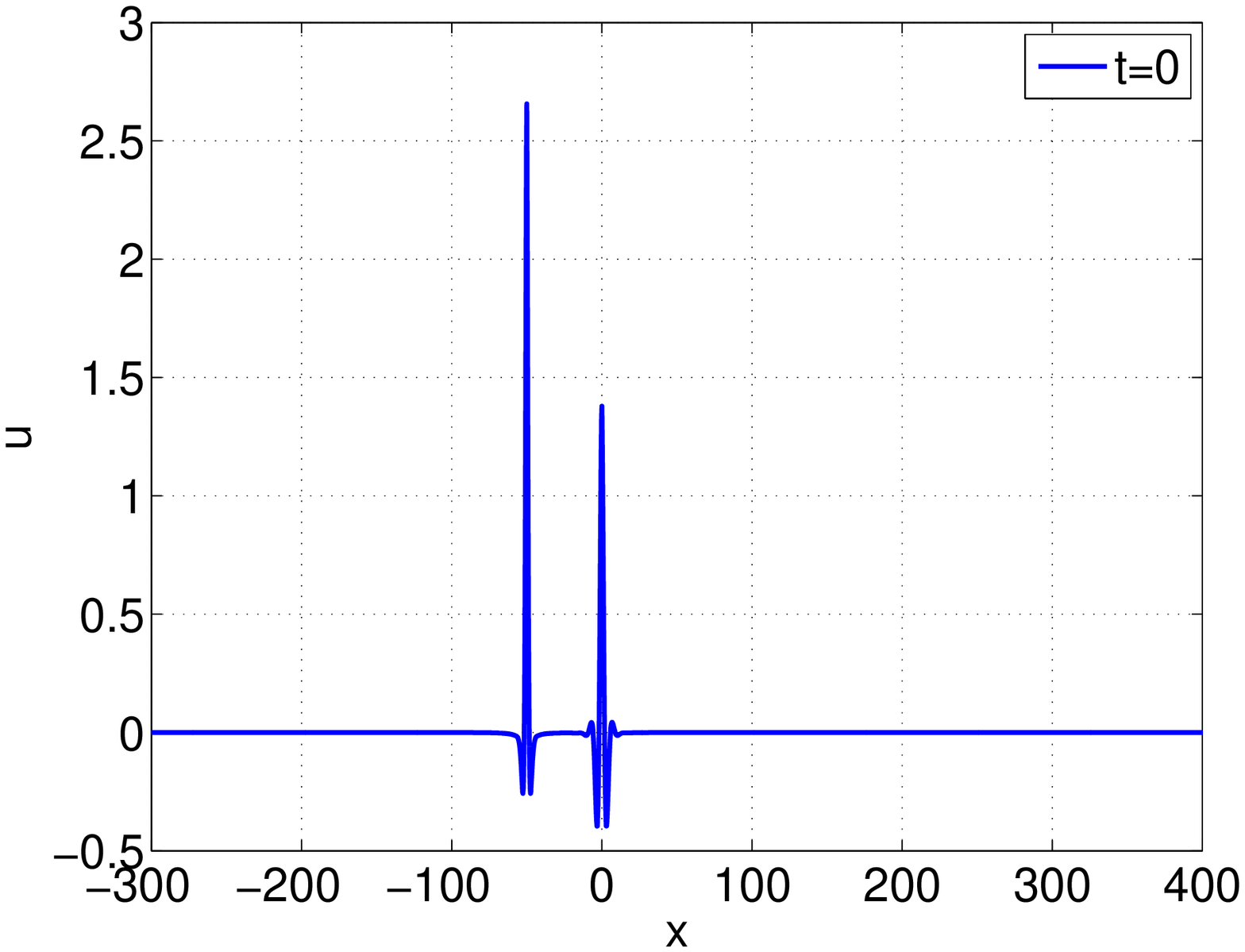}}
\subfigure[]
{\includegraphics[width=6cm]{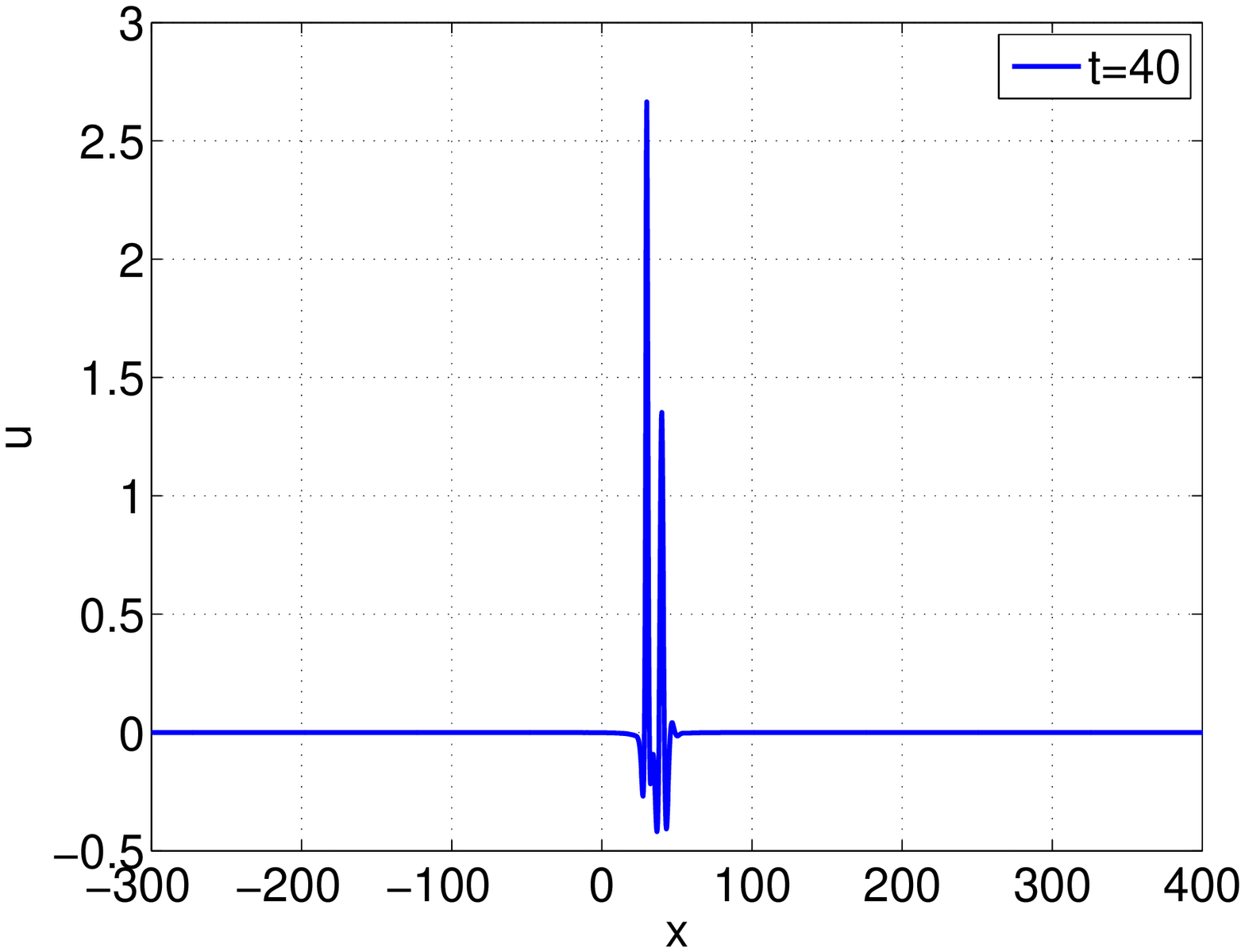}}
\subfigure[]
{\includegraphics[width=6cm]{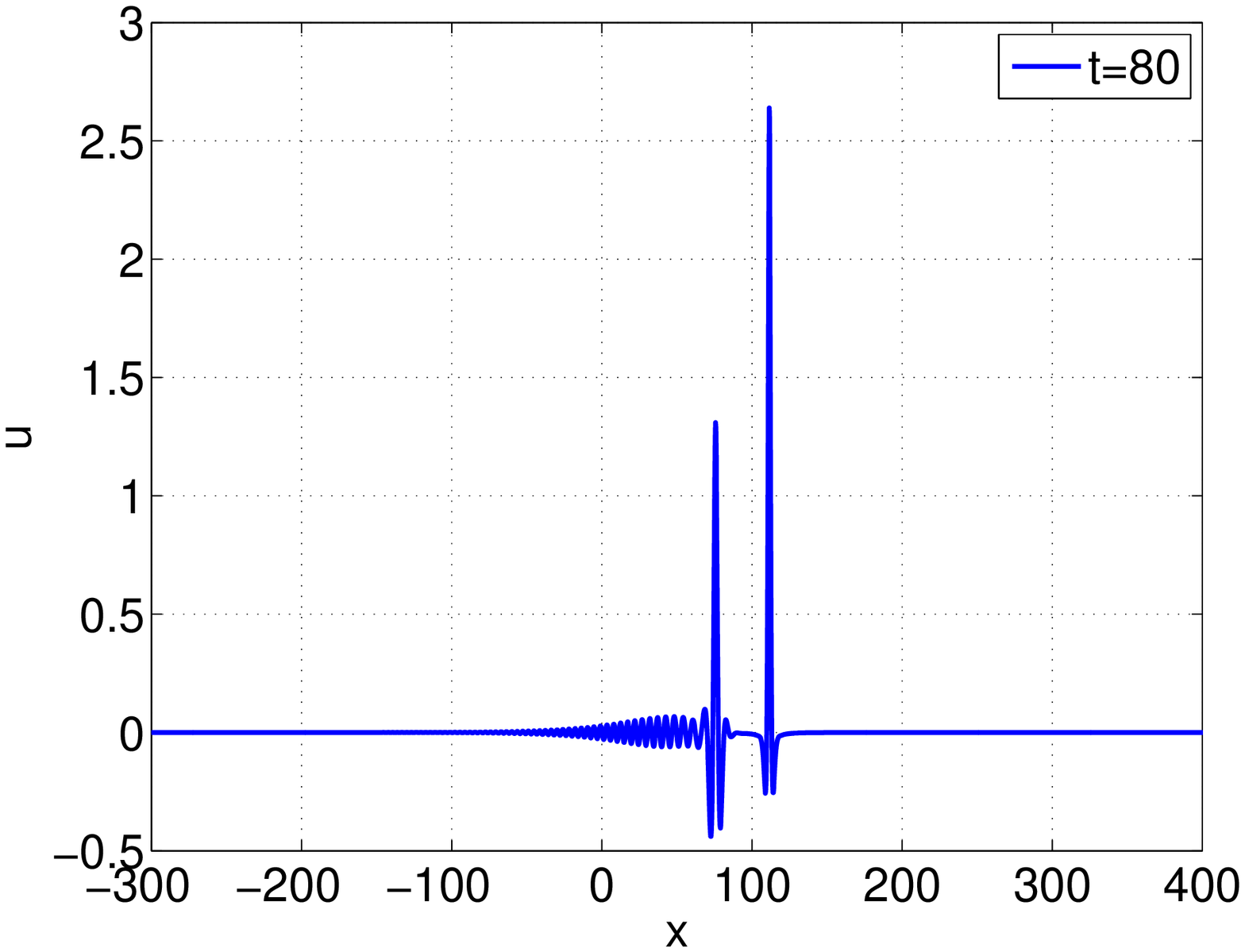}}
\subfigure[]
{\includegraphics[width=6cm]{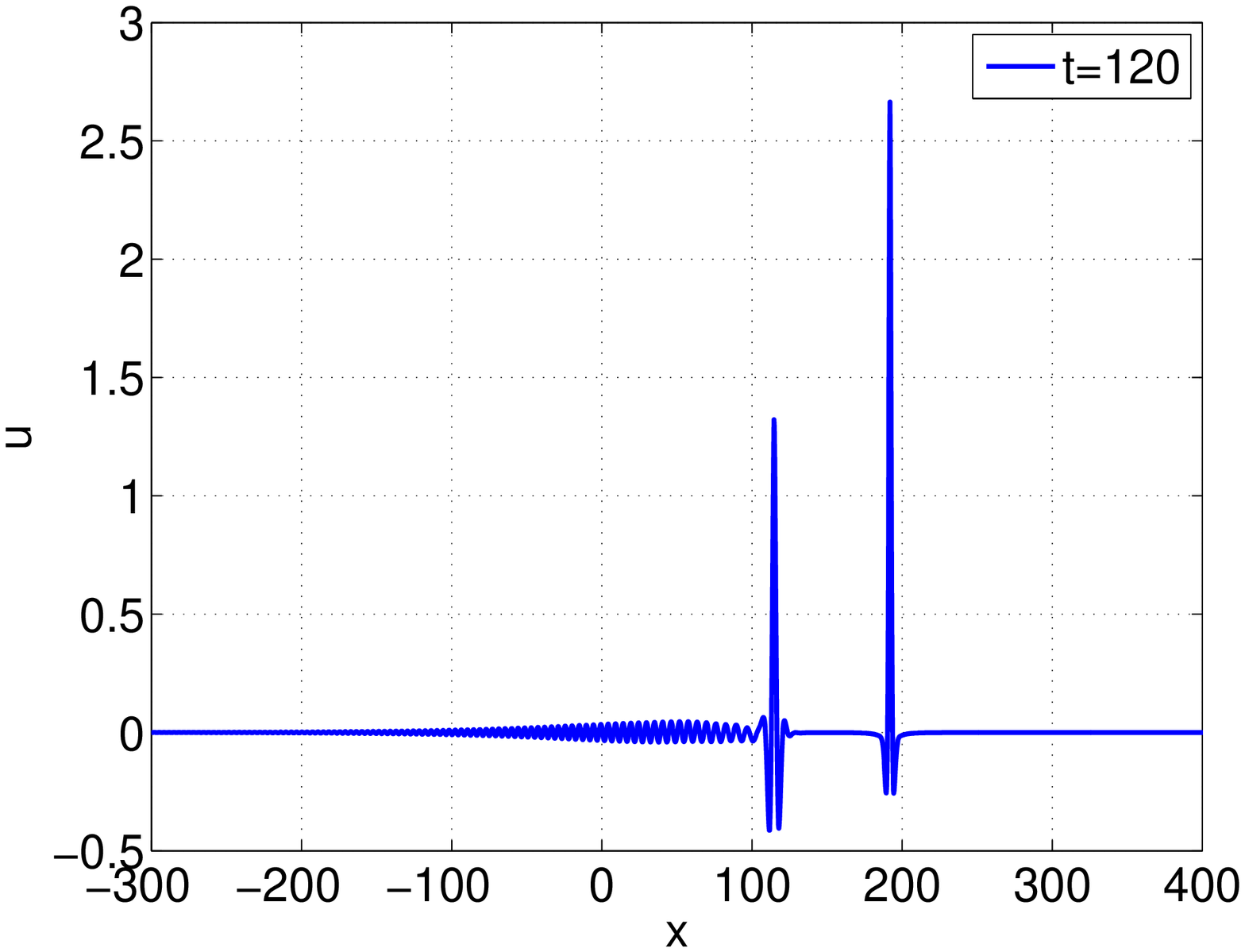}}
\subfigure[]
{\includegraphics[width=6cm]{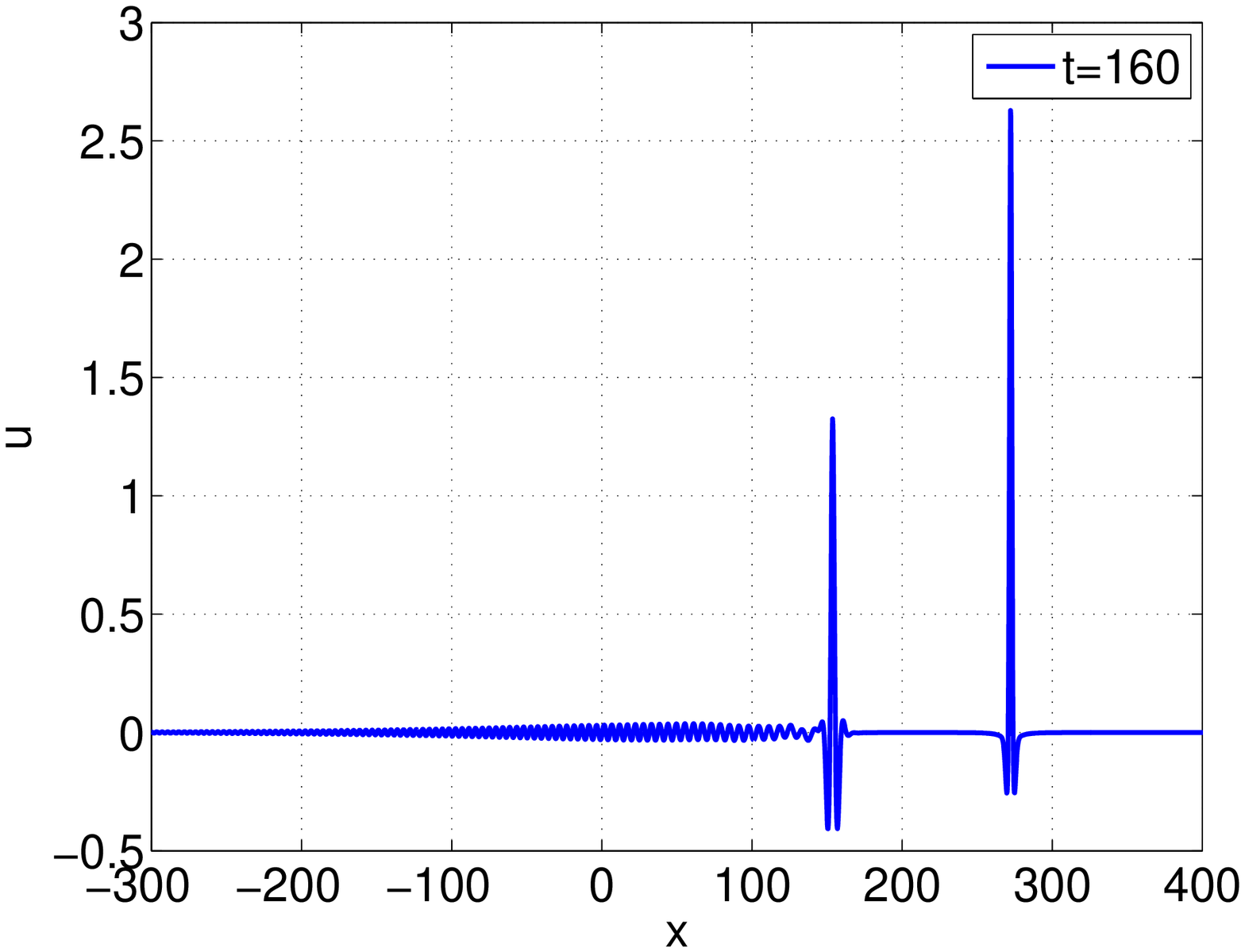}}
\subfigure[]
{\includegraphics[width=6cm]{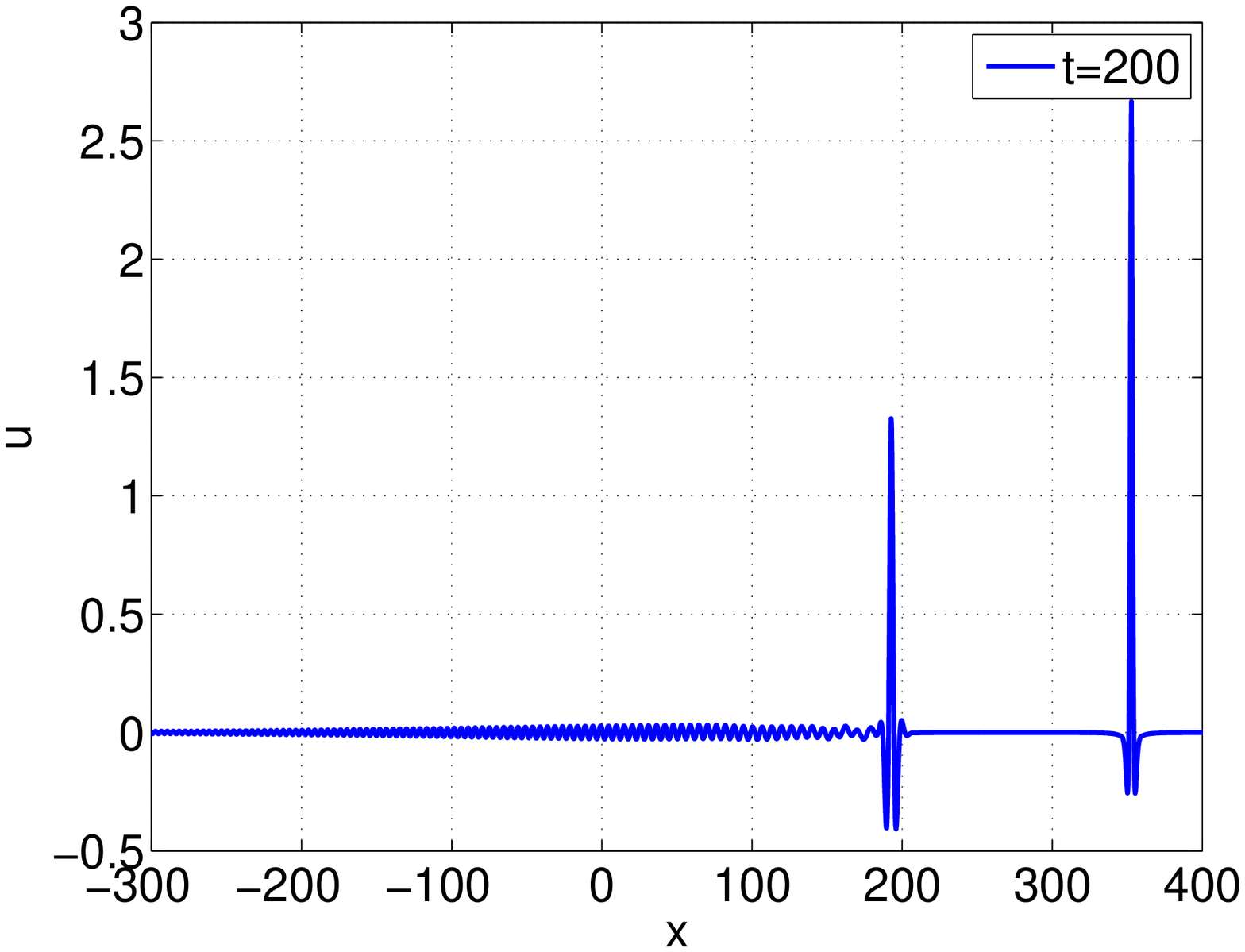}}
\caption{Interactions of solitary waves. Case (G1) $q=2$ of gBenjamin equation. (a)-(f) Numerical approximation at $t=0,40,80,120,180,200$.}
\label{gbenfig54G11}
\end{figure}

\begin{figure}[htbp]
\centering
\subfigure[]
{\includegraphics[width=6cm]{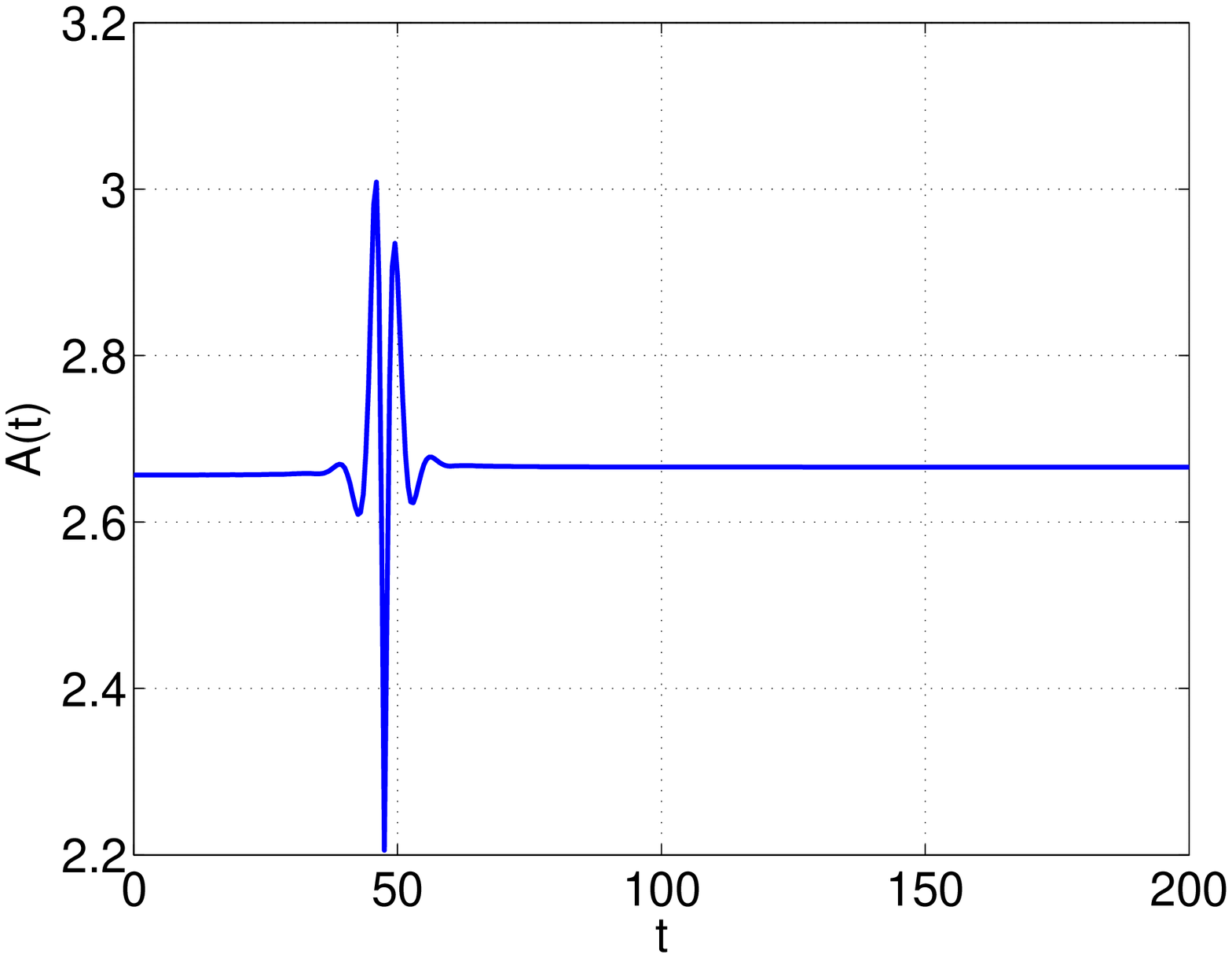}}
\subfigure[]
{\includegraphics[width=6cm]{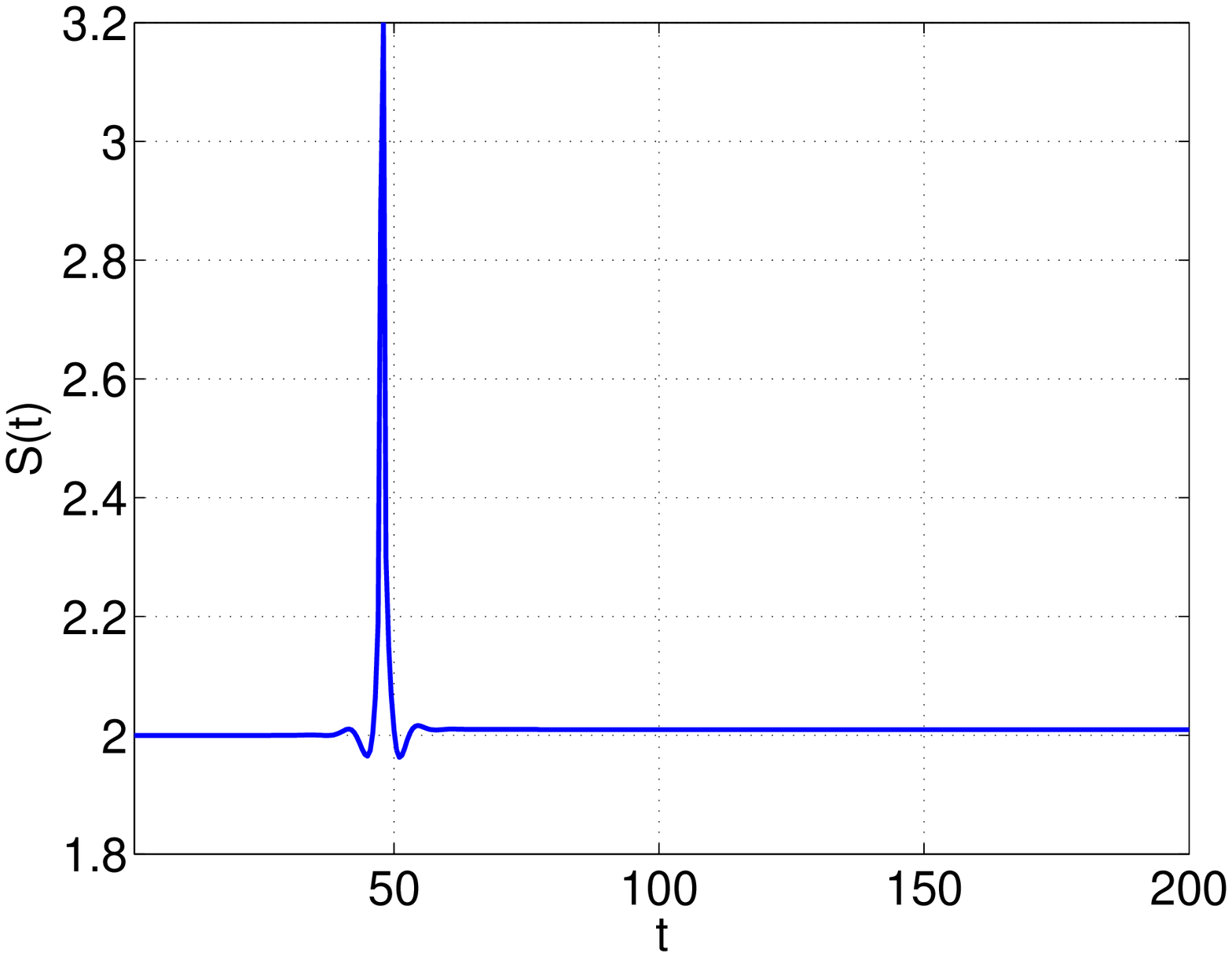}}
\caption{Interactions of solitary waves. Case (G1) $q=2$ of gBenjamin equation. (a) Evolution of the amplitude of the main numerical pulse; (b) Evolution of speed.}
\label{gbenfig54G11b}
\end{figure}

We start by describing (G1). The cases $q=2,3$ give similar results. (We show only the first one in Figures \ref{gbenfig54G11}, \ref{gbenfig54G11b}.) After the collision two new solitary wave profiles emerge followed by a dispersive tail. The larger wave increases slightly its amplitude and speed. The increment is slightly larger in the case $q=3$.
The smaller wave after the interaction has a slightly reduced amplitude. The reduction is approximately $3.83\%$ when $q=2$ and $5.29\%$ when $q=3$. The increment of speed of the tallest wave is also larger when $q=3$ ($2.18\%$) than when $q=2$ ($0.48\%$).

\begin{figure}[!htbp]
\centering
\subfigure[]
{\includegraphics[width=6cm]{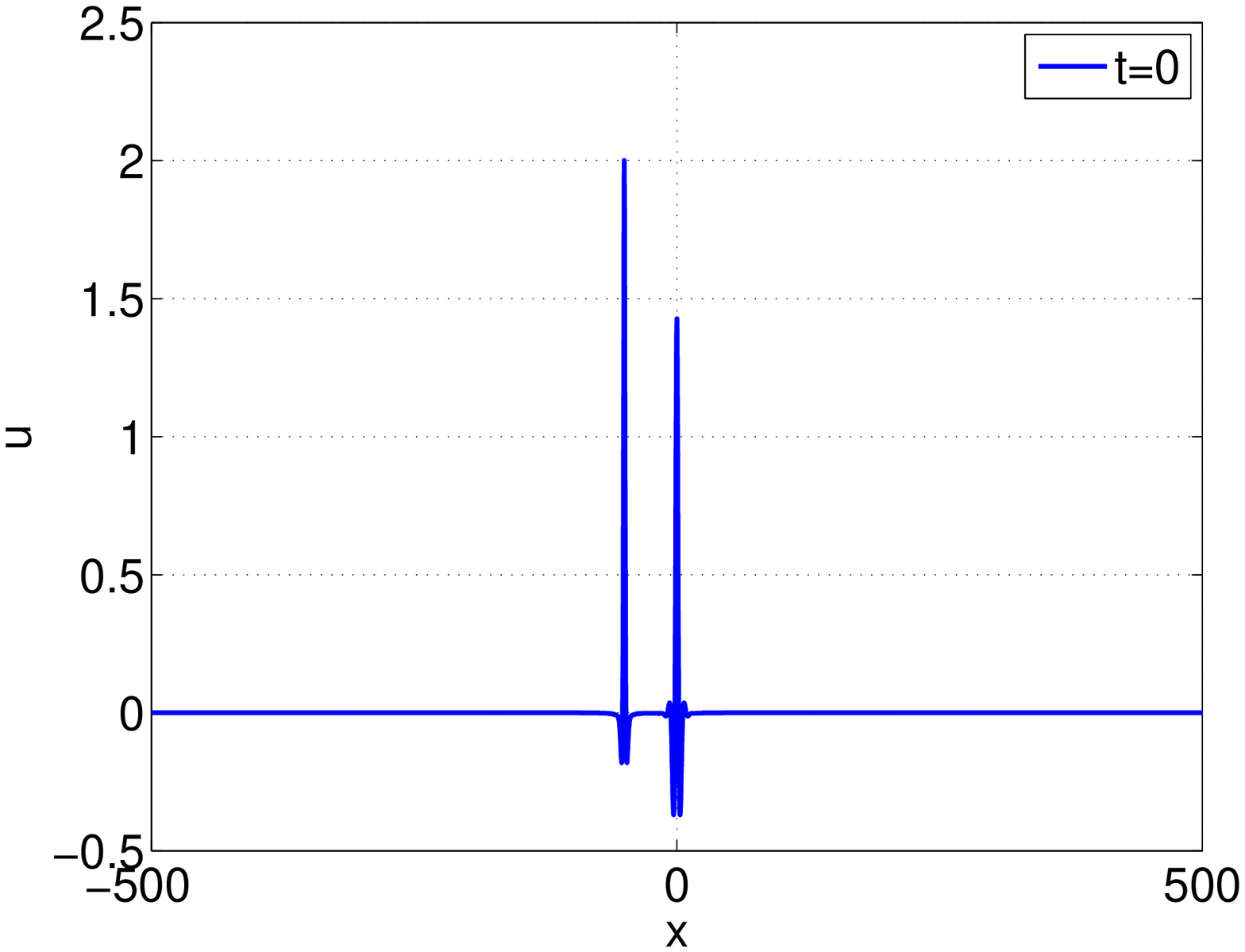}}
\subfigure[]
{\includegraphics[width=6cm]{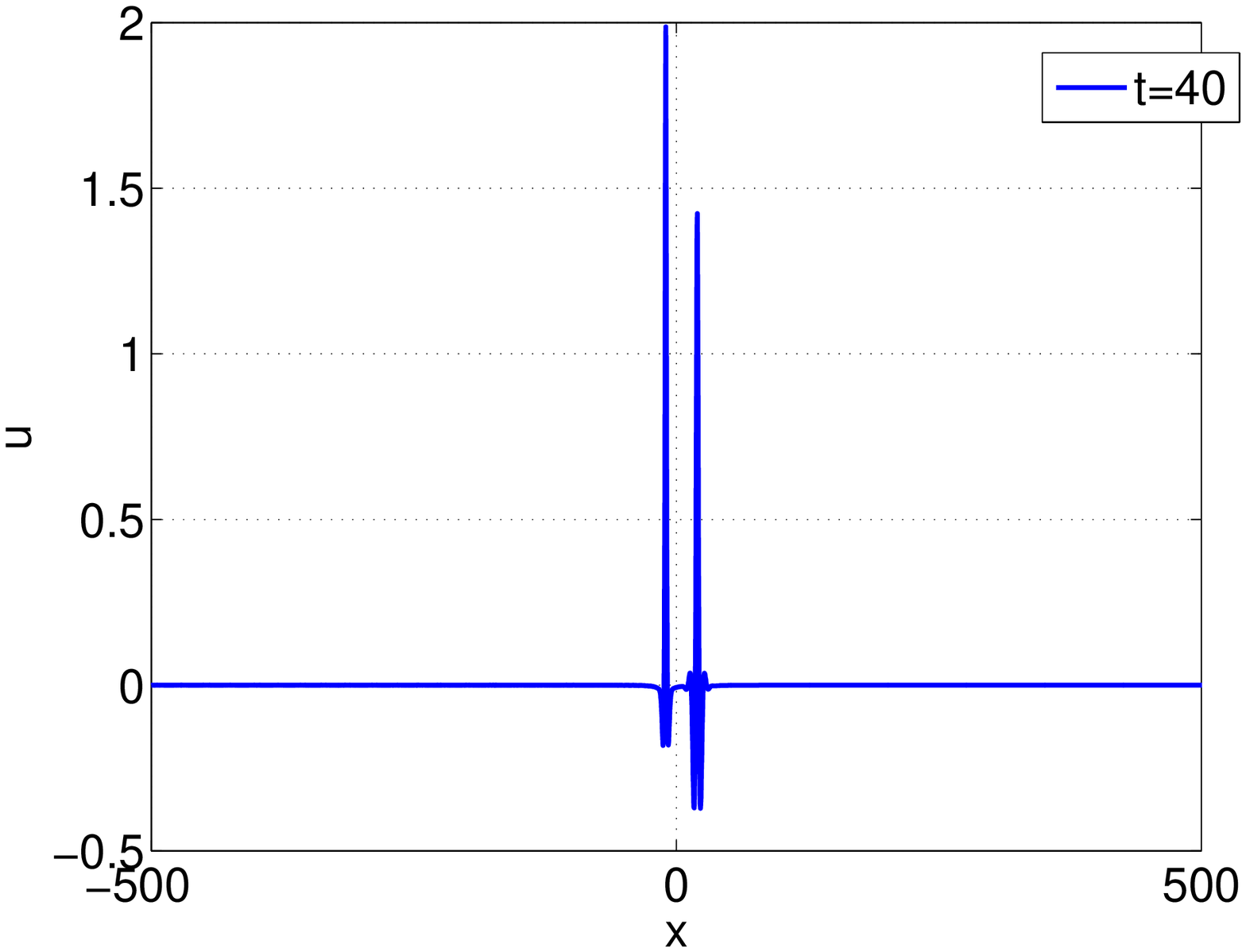}}
\subfigure[]
{\includegraphics[width=6cm]{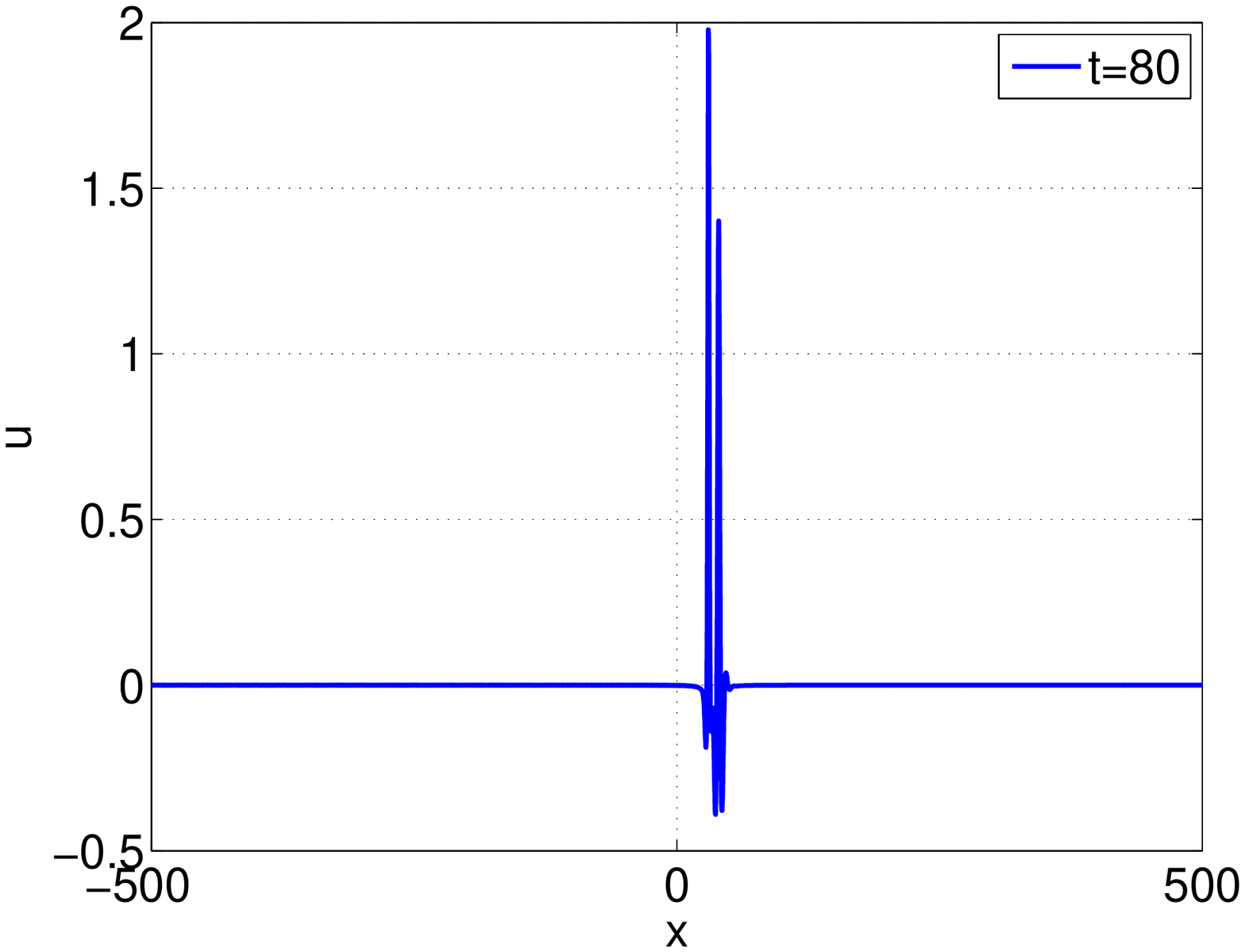}}
\subfigure[]
{\includegraphics[width=6cm]{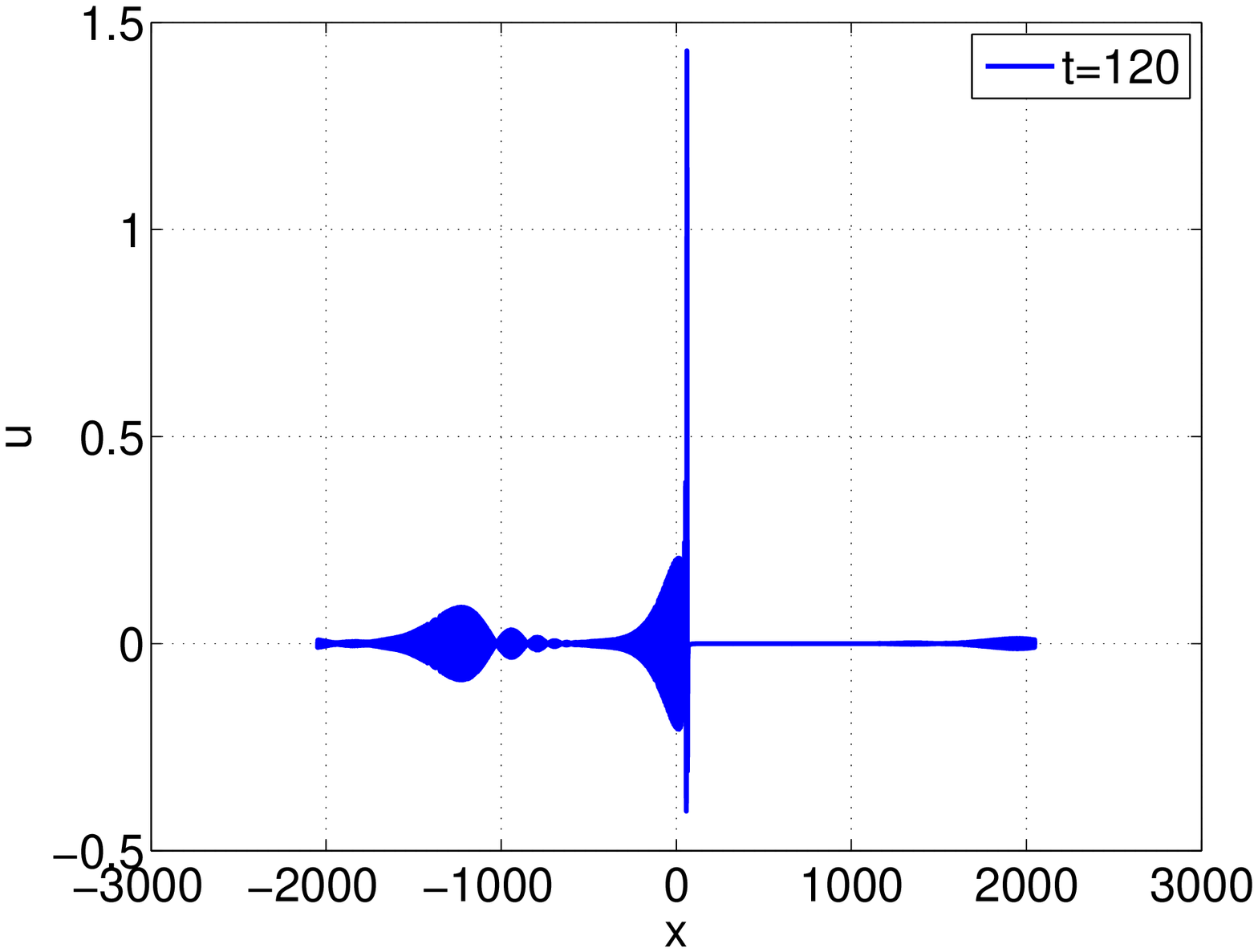}}
\subfigure[]
{\includegraphics[width=6cm]{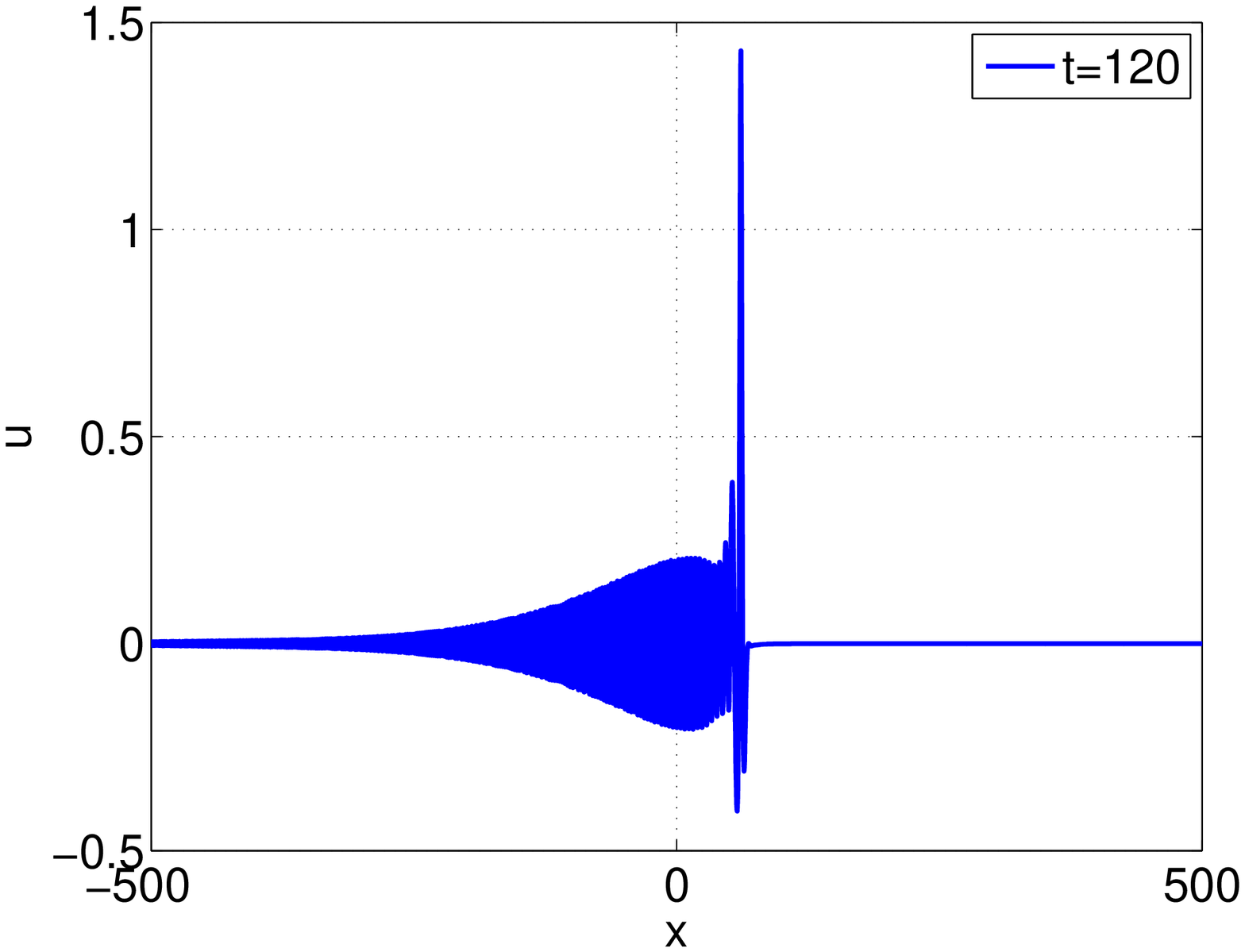}}
\subfigure[]
{\includegraphics[width=6cm]{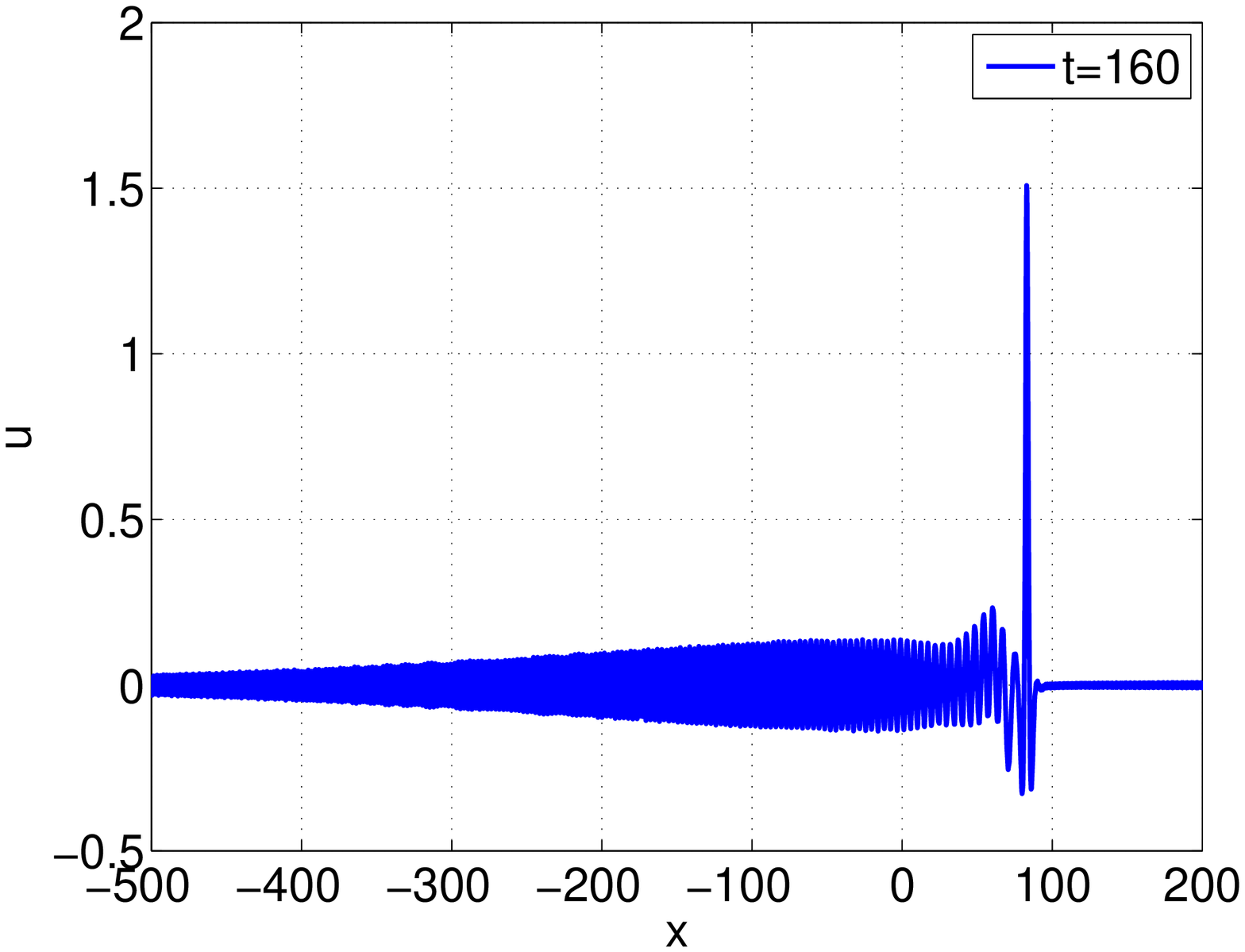}}
\subfigure[]
{\includegraphics[width=6cm]{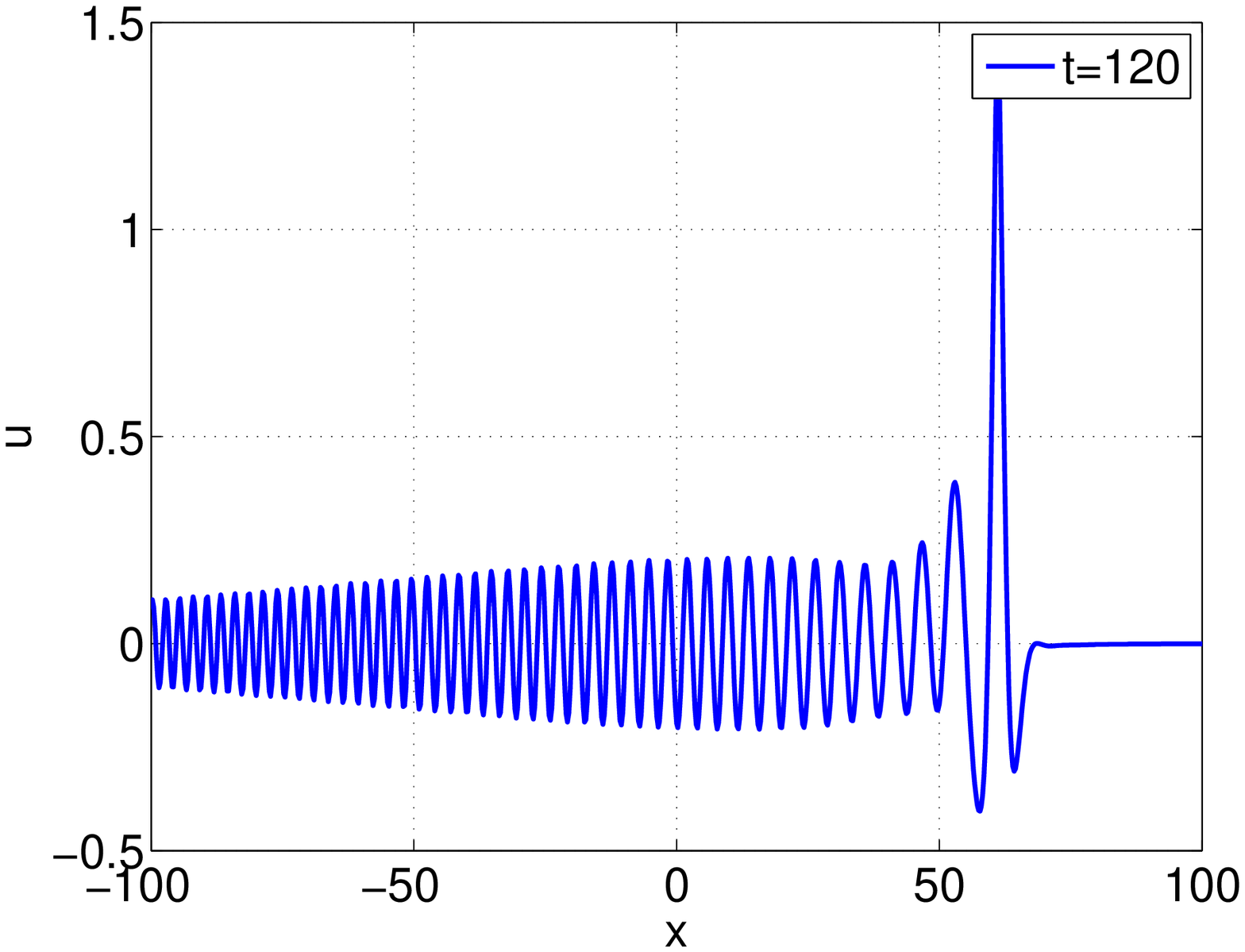}}
\caption{Interactions of solitary waves. Case (G1) $q=4$ of gBenjamin equation. (a)-(f) Numerical approximation at $t=0,40,80,120,160$; (g) is a magnification of (e).}
\label{gbenfig54G13}
\end{figure}
\begin{figure}[!htbp]
\centering
\subfigure[]
{\includegraphics[width=6cm]{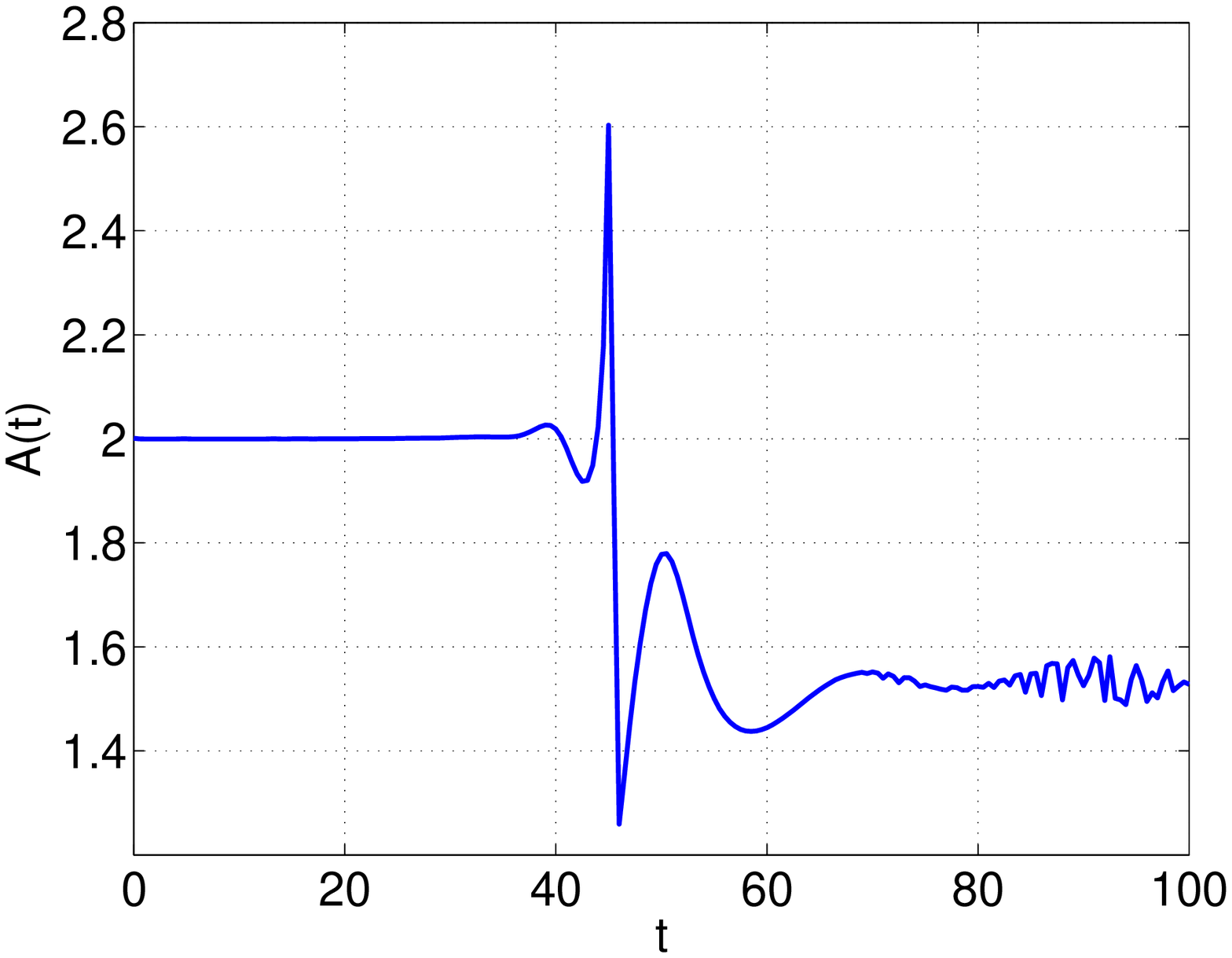}}
\subfigure[]
{\includegraphics[width=6cm]{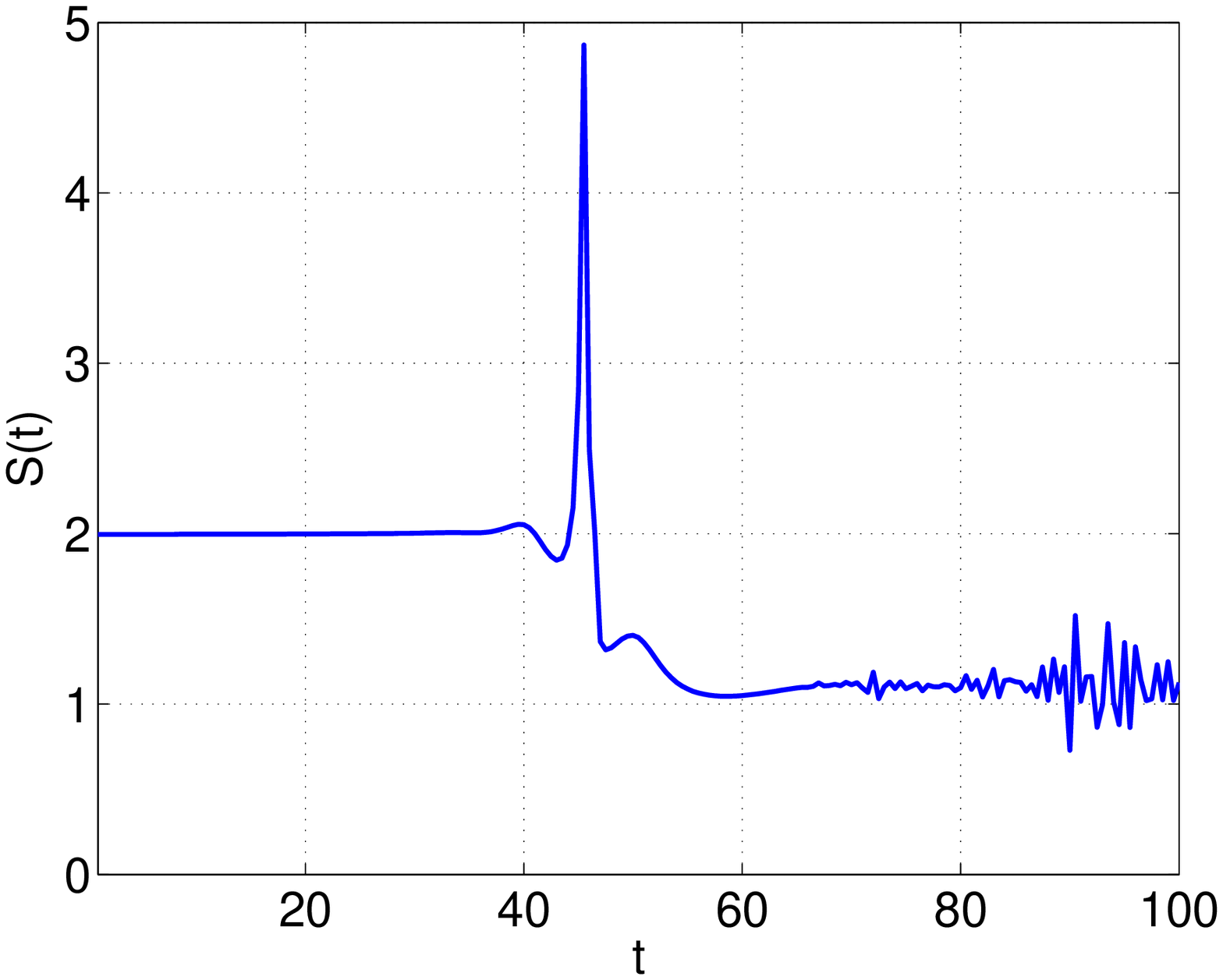}}
\caption{Interactions of solitary waves. Case (G1) $q=4$ of gBenjamin equation. (a) Evolution of the amplitude of the main numerical pulse; (b) Evolution of speed.}
\label{gbenfig54G13b}
\end{figure}

The results of (G2)  for $q=2,3$ give similar conclusions. The interactions look more \lq inelastic\rq\ when $q=3$, in the sense of a larger change of the parameters. On the other hand, when $q=4$ and in both (G1) and (G2) runs the experiments suggest instability (Figures \ref{gbenfig54G13}, \ref{gbenfig54G13b} correspond to (G1).) The dynamics observed in the experiments reminds that of Figure \ref{gbenfig516_1} (obtained from a small perturbation of a single-pulse profile), with the formation of ripples behind the emerging solitary waves that do not seem to disperse as time evolves.

\begin{figure}[!htbp]
\centering
\subfigure[]
{\includegraphics[width=6cm]{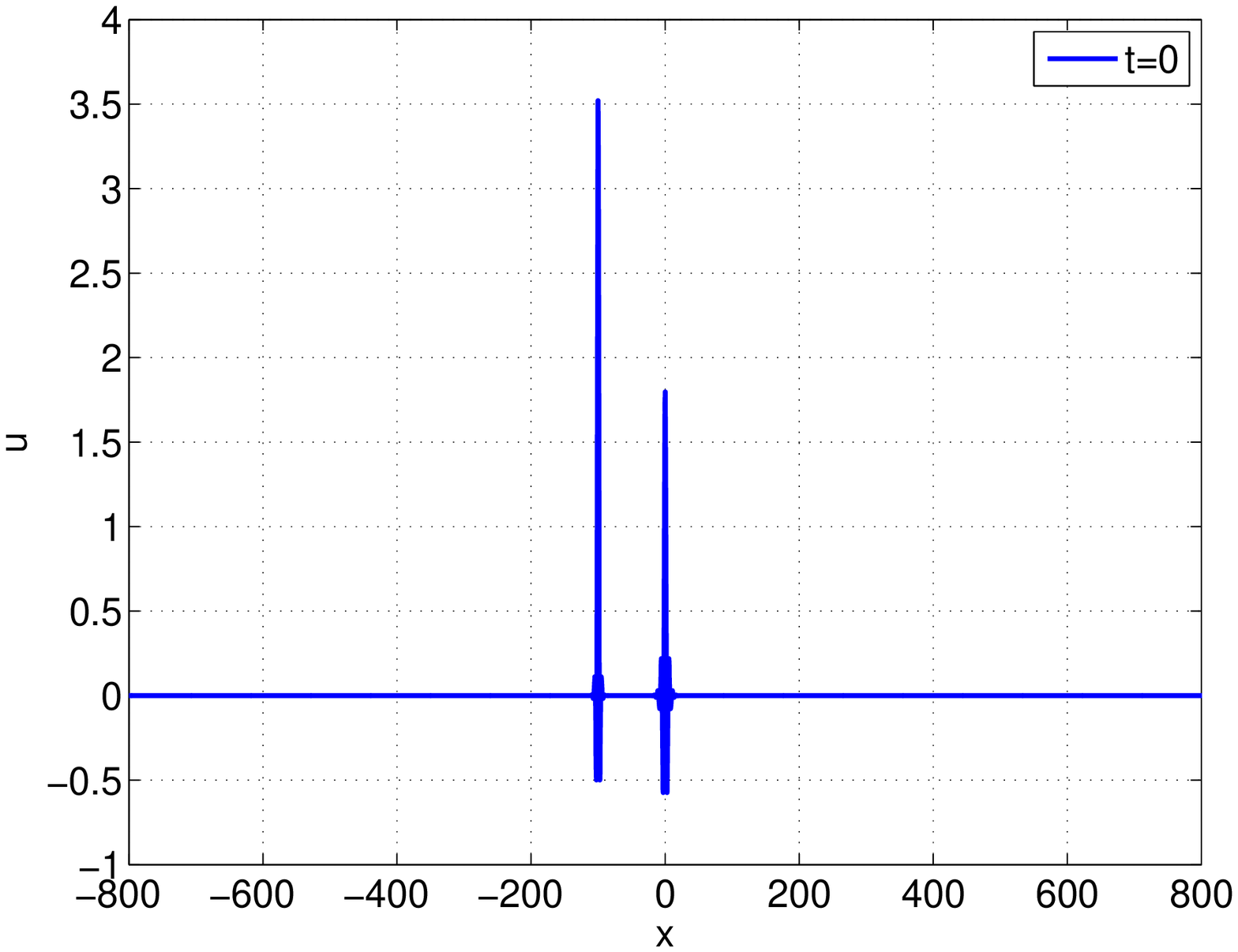}}
\subfigure[]
{\includegraphics[width=6cm]{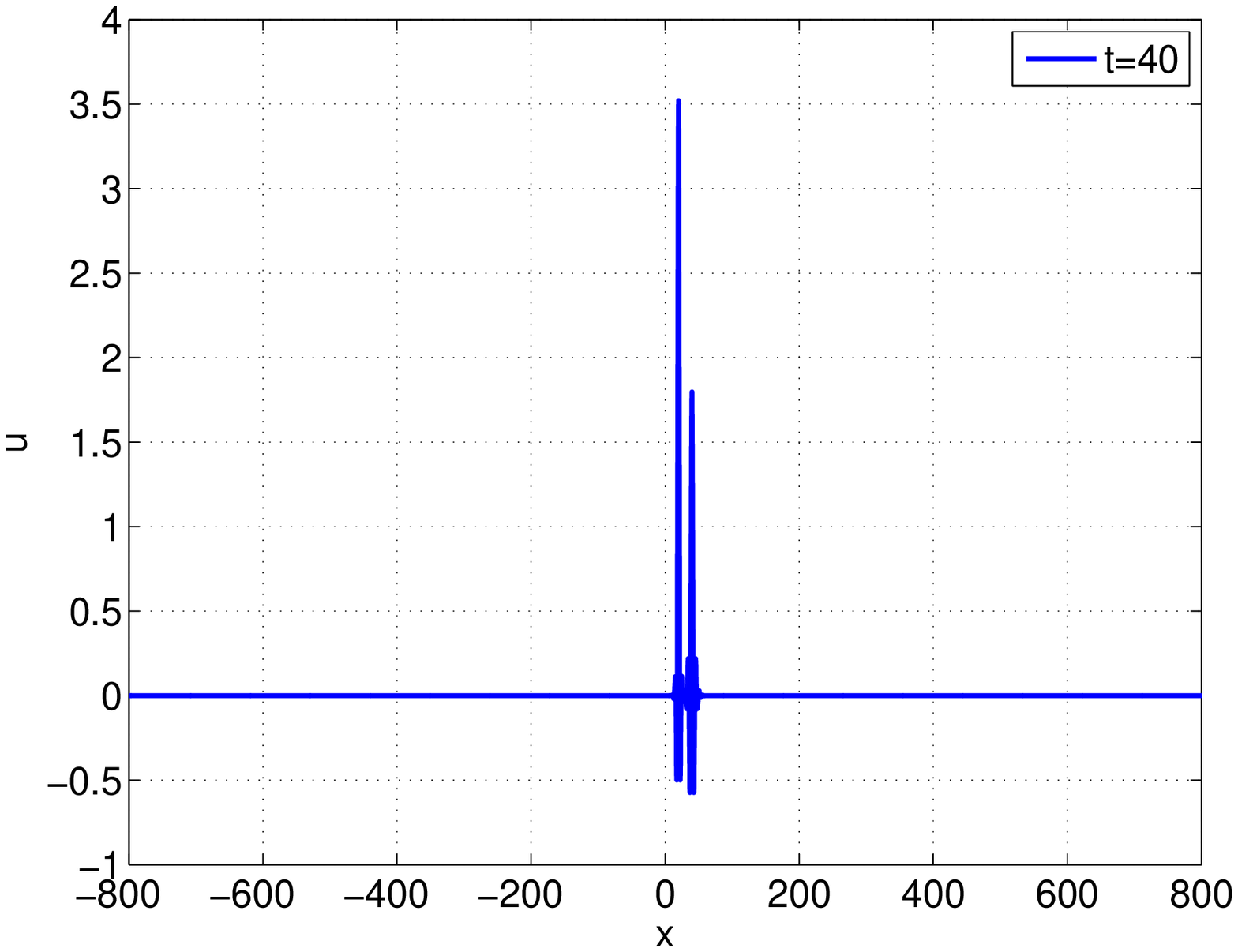}}
\subfigure[]
{\includegraphics[width=6cm]{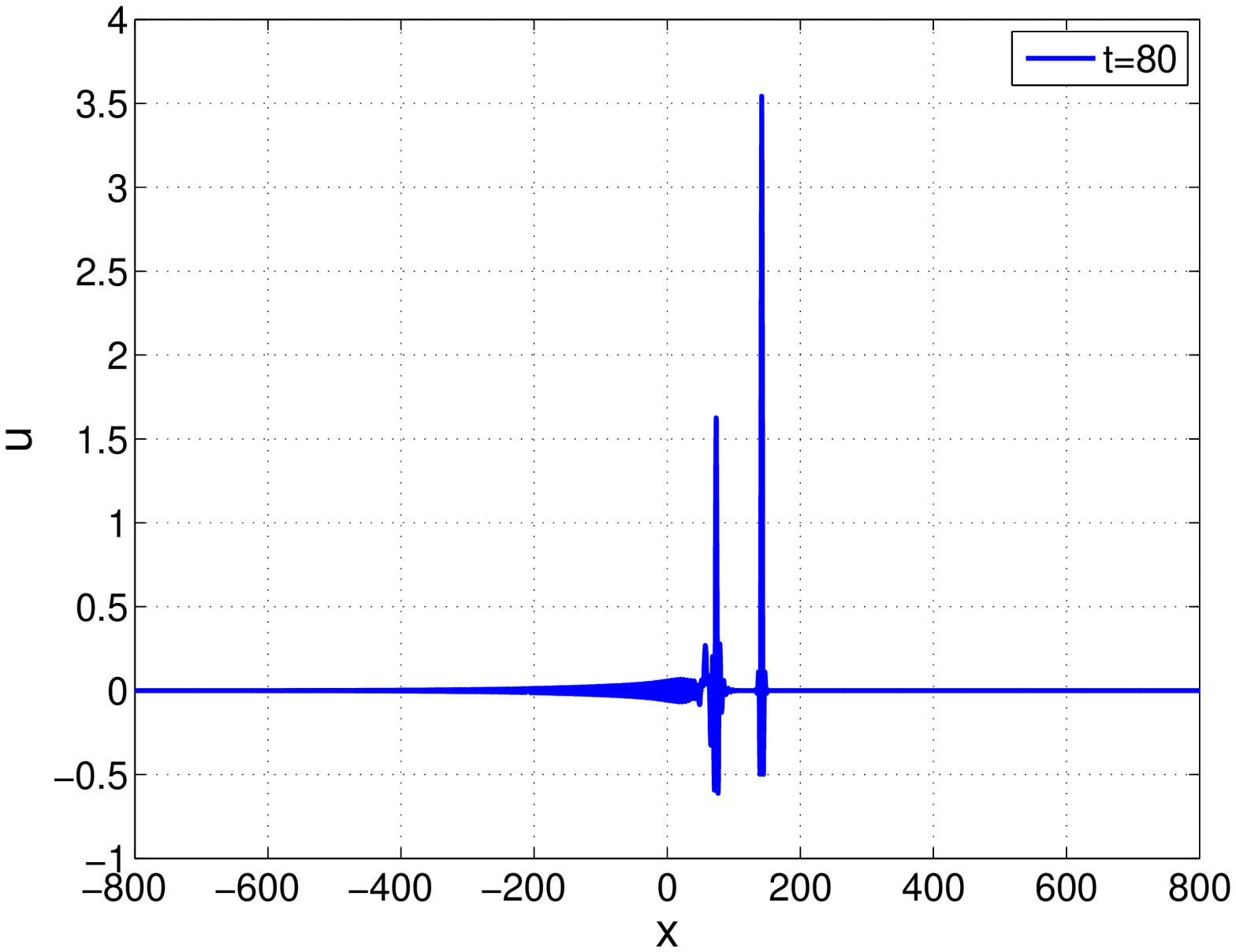}}
\subfigure[]
{\includegraphics[width=6cm]{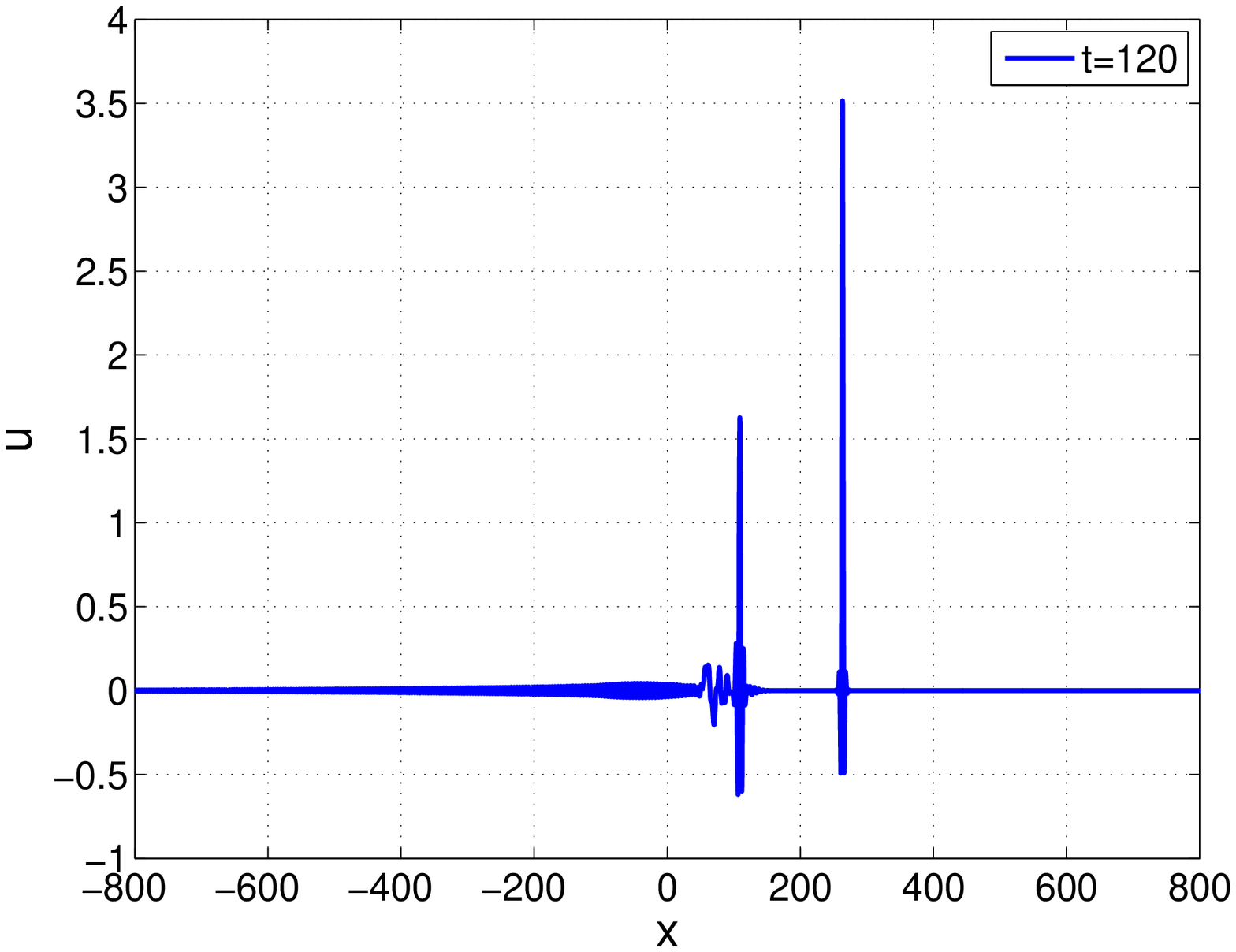}}
\subfigure[]
{\includegraphics[width=6cm]{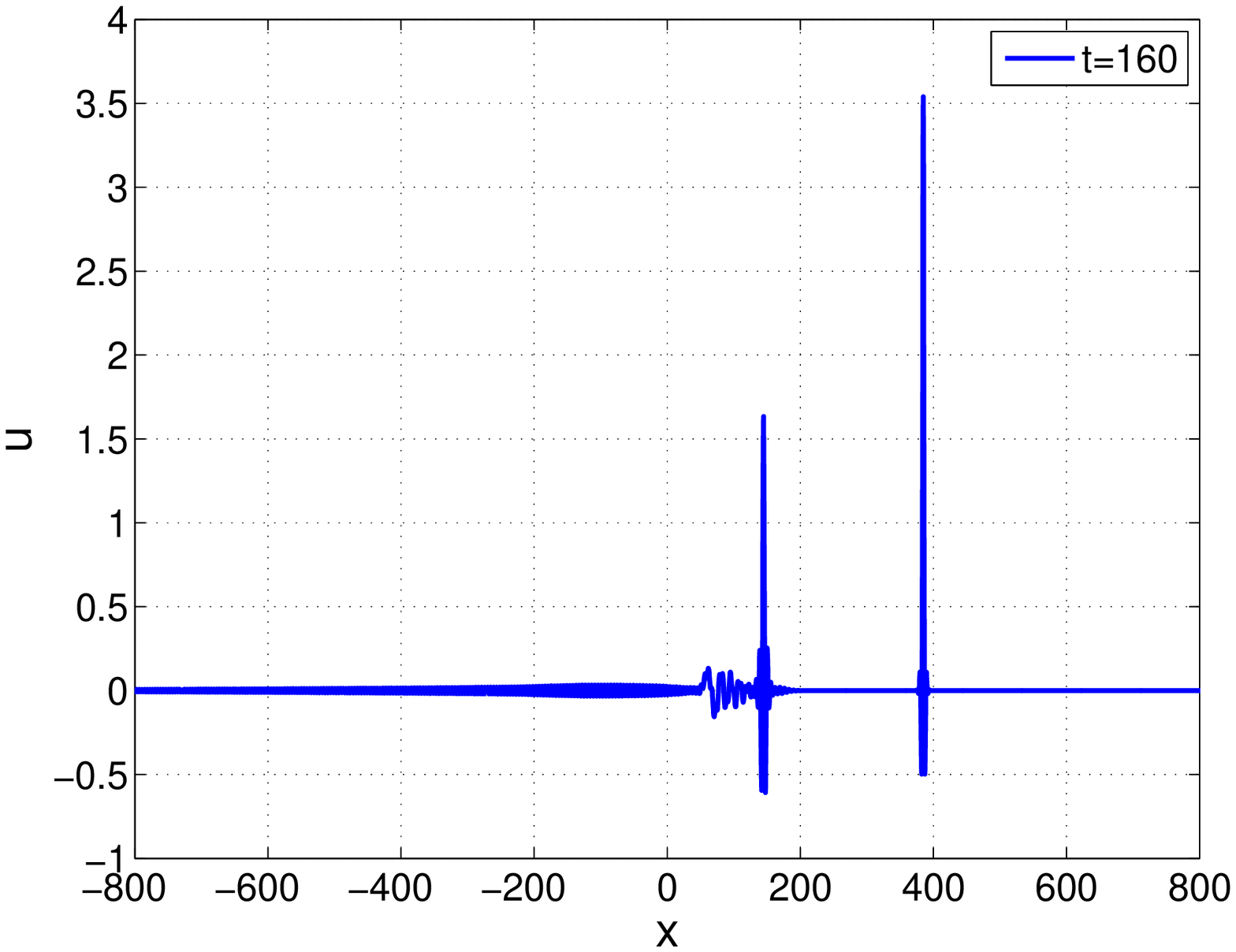}}
\subfigure[]
{\includegraphics[width=6cm]{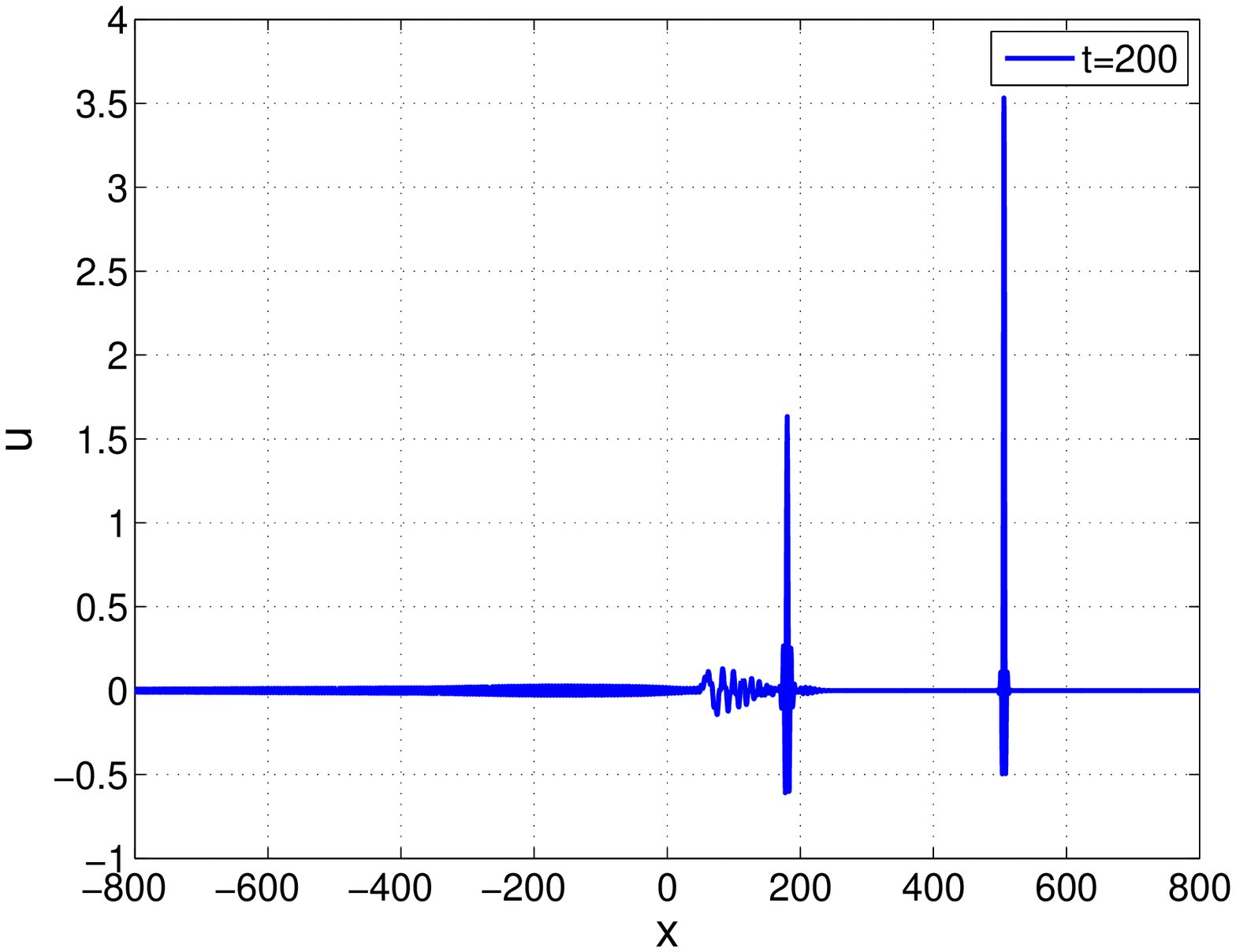}}
\caption{Interactions of solitary waves. General case: $r=3/2, m=2, \delta=1, \gamma=1.5, q=2$. (a)-(f) Numerical approximation at $t=0,40,80,120,160,200$.}
\label{gbenfig54G14}
\end{figure}

\begin{figure}[htbp]
\centering
\subfigure[]
{\includegraphics[width=6cm]{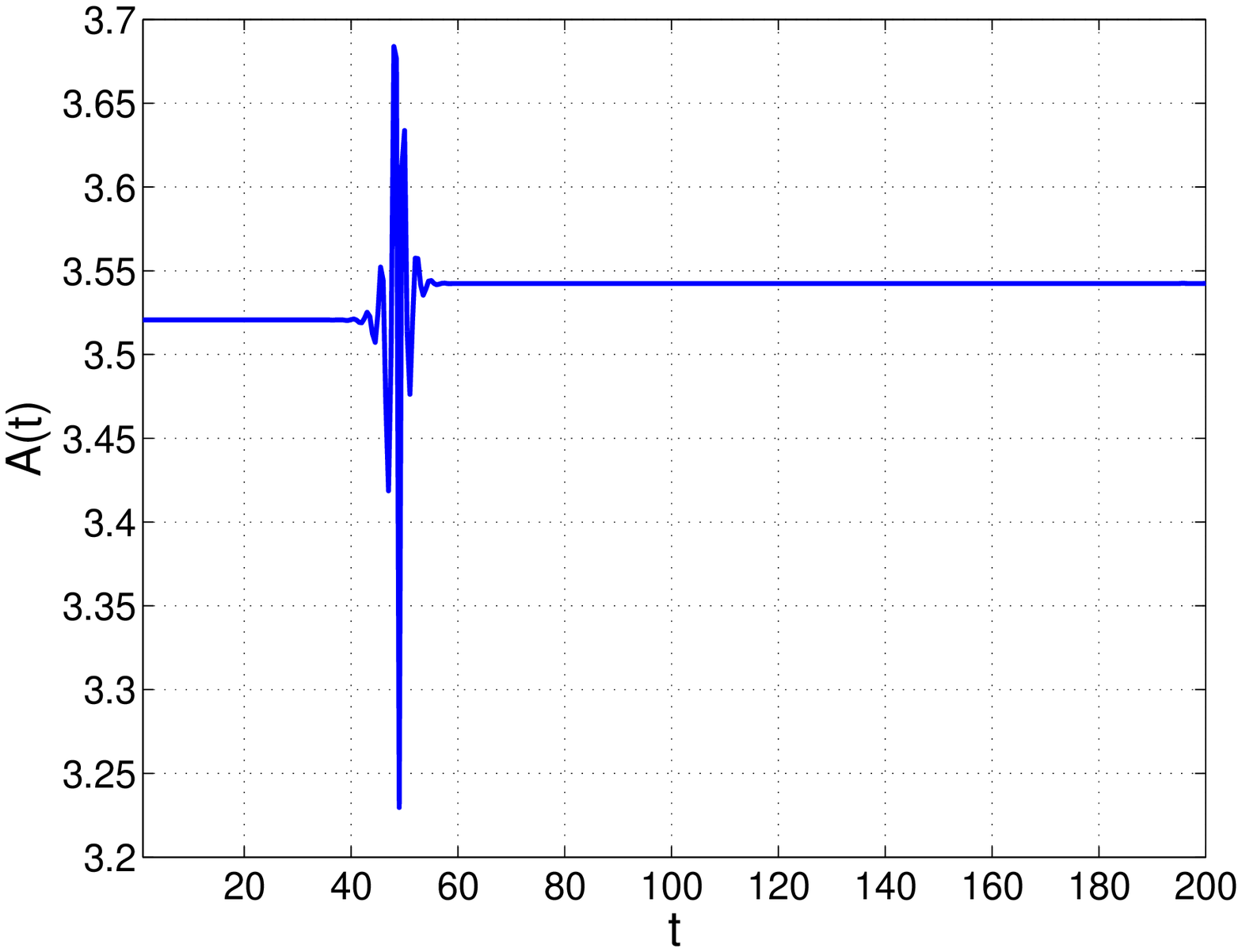}}
\subfigure[]
{\includegraphics[width=6cm]{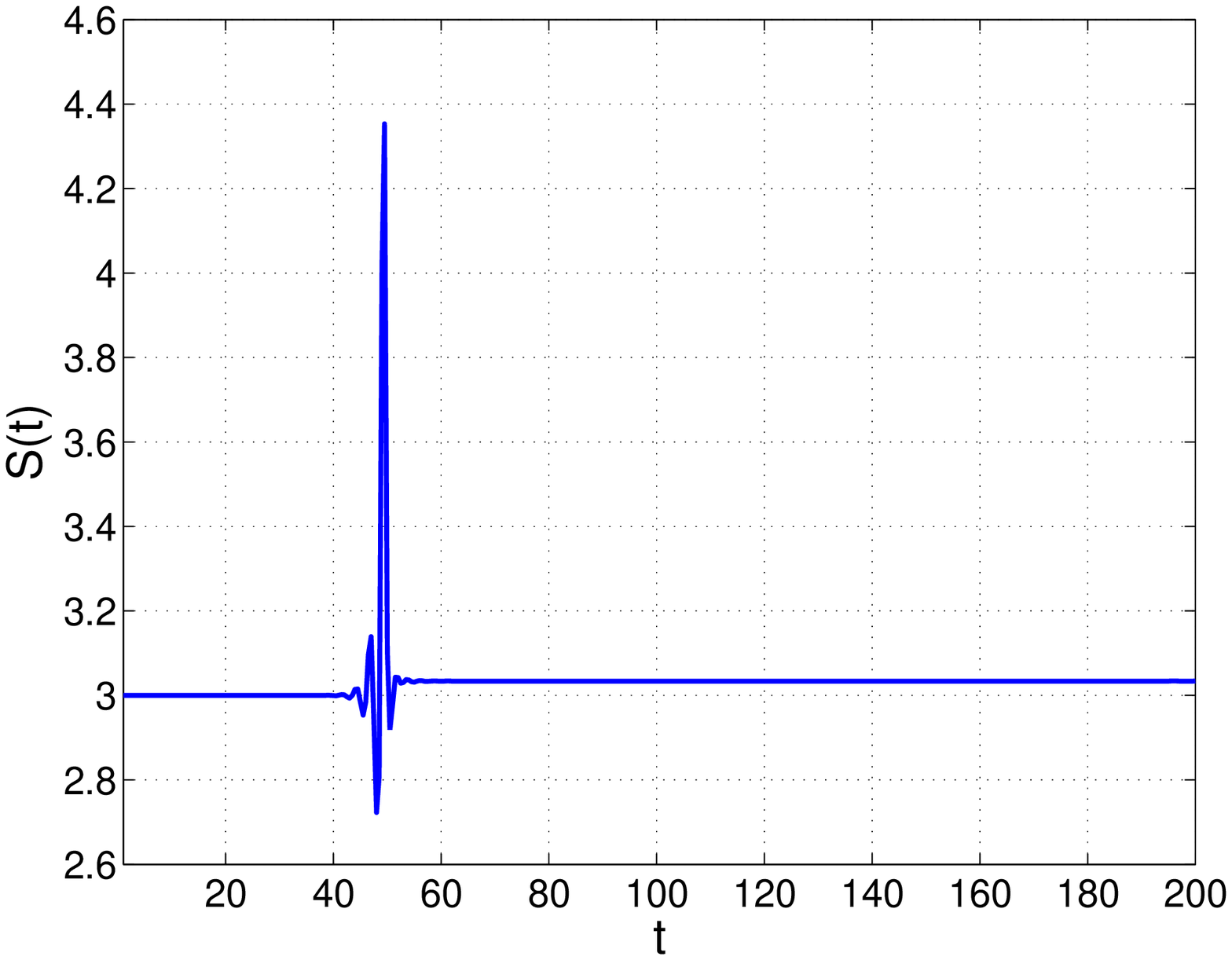}}
\caption{Interactions of solitary waves. General case: $r=3/2, m=2, \delta=1, \gamma=1.5, q=2$. (a) Evolution of the amplitude of the main numerical pulse; (b) Evolution of speed.}
\label{gbenfig54G14b}
\end{figure}

The general case is illustrated with $r=3/2, m=2, \delta=1, \gamma=1.5$.
Two initial profiles of speeds $c_{s}^{(1)}=1, c_{s}^{(2)}=3$ centered at $x_{0}^{(1)}=0, x_{0}^{(2)}=-100$ are taken. The final time is $t=200$ and two values $q=2,3$ have been considered. In this case, $\gamma_{max}(c_{s}^{(1)})\approx 1.754765$ and $\gamma_{max}(c_{s}^{(2)})=2.309401$. The corresponding amplitudes are displayed in Table \ref{gbentab542}; see Figures \ref{gbenfig54G14}, \ref{gbenfig54G14b}, corresponding to the case $q=2$. 
The case $q=4$ leads to similarly unstable results.

\begin{table}
\smallskip
\begin{center}
\begin{tabular}{|c|c|}
\hline\hline
$q=2$&$u_{b}^{(1)}\approx 1.797379$,$u_{a}^{(1)}\approx 1.633485$\\
&$c_{sb}^{(2)}=3, c_{sa}^{(2)}\approx 3.033020$\\
&$u_{b}^{(2)}\approx 3.520598$,$u_{a}^{(2)}\approx 3.542342$\\
\hline
$q=3$&$u_{b}^{(1)}\approx 1.698971$,$u_{a}^{(1)}\approx 1.483980$\\
&$c_{sb}^{(2)}=3, c_{s}^{(2)}\approx 3.098329$\\
&$u_{b}^{(2)}\approx 2.656828$,$u_{a}^{(2)}\approx 2.685788$\\
\hline
\end{tabular}
\end{center}
\caption{Interactions of solitary waves, General case: $r=3/2, m=2, \delta=1, \gamma=1.5$.
Two initial profiles of speeds $c_{s}^{(1)}=1, c_{s}^{(2)}=3$. Amplitudes before ($u_{b}$) and after ($u_{a}$) the collision; speeds ($c_{sb}, c_{sa}$) of the tallest wave.}\label{gbentab542}
\end{table}

\section{Concluding remarks}
\label{sec4}
In this paper a computational study of solitary-wave solutions of generalized versions of the Benjamin equation is carried out. The generalized character of the equations under study, with respect to the the Benjamin equation, involves the dispersion, the surface tension and the nonlinear effects. The first two are modeled by incorporating a linear operator with high-order derivatives in dispersive and nonlocal terms, while the nonlinearity is a homogeneous function of degree $q\geq 2$.

The paper is focused on two aspects concerning the solitary-wave solutions of the equations. First, the numerical generation of the solitary profiles (which is necessary because of the lack of explicit formulas) is performed by the Petviashvili iterative method, combined with the minimal polynomial extrapolation method to accelerate the convergence and to improve, in comparison with previous numerical procedures,  the generation of highly oscillatory solitary pulses. The accuracy of the proposed technique allows us to study some properties of the solitary waves, such as the speed-amplitude relation, as well as to verify the known results of the asymptotic decay depending on the parameters of the equations, \cite{BonaC1998}. 

Once the numerical solitary profiles are accurately computed, they are used to investigate the stability of the solitary waves. The experiments of the second part of the numerical study aim at analyzing the approximate evolution of the solutions from small and large perturbations of solitary-wave profiles and from superpositions of solitary-wave pulses that will interact. The code to simulate the corresponding periodic initial-value problem is based on a Fourier pseudospectral collocation method as spatial discretization and a time integration with a fourth-order, diagonally implicit Runge-Kutta composition method. The computations suggest that the single-pulse solitary waves are stable under small perturbations for $q<4$ while some sort of instability is observed for $q\geq 4$. Larger perturbations of the initial solitary-wave profile seem to evolve into solutions with more solitary pulses, like a resolution property, for $q<4$. We finally examine the interactions of solitary waves. The inelastic character of the collisions seems to grow with $q<4$, in the sense that a larger change of the parameters of the emerging waves is observed. For $q\geq 4$, the interactions of solitary waves apparently lead to similar instabilities.

These preliminary results are being subjected to further computational investigation and, along with the analysis of convergence of the codes, will be the main point in further research.

\section*{Acknowledgements}
This work was supported by Spanish Ministerio de Econom\'{\i}a y Competitividad under the Research Grant MTM2014-54710-P. A. Duran was also supported by Junta de Castilla y Le\'on under the Grant VA041P17.

%

%

%

%
%




\newpage




%


\end{document}